\documentclass[12pt]{ociamthesis}  % default square logo 
\usepackage[utf8]{inputenc}
\usepackage{amssymb}
\usepackage{amsmath,esint}
\usepackage{booktabs}
\usepackage{mathrsfs}
\usepackage{mathtools}
\usepackage{pifont}% http://ctan.org/pkg/pifont
\newcommand{\cmark}{\ding{51}}%
\newcommand{\xmark}{\ding{55}}%
\usepackage{color}
\usepackage{overpic}
\usepackage{url}
\usepackage{upquote,enumitem}
\usepackage[hidelinks]{hyperref}
\usepackage[noabbrev,capitalize]{cleveref}
\usepackage{url}
\usepackage{subfiles}
\usepackage[nottoc]{tocbibind} % Correct numbering of biblio

\usepackage{algorithm}
\usepackage[noend]{algpseudocode}
\algnewcommand\algorithmicinput{\textbf{Input:}}
\algnewcommand\INPUT{\item[\algorithmicinput]}

\usepackage{amsthm}
\newtheorem{theorem}{Theorem}[chapter]
\newtheorem{corollary}{Corollary}[chapter]
\newtheorem{definition}{Definition}[chapter]

\newtheorem{lemma}{Lemma}[chapter]

\newtheorem{proposition}{Proposition}[chapter]
\newtheorem{remark}{Remark}[chapter]

\DeclareMathOperator{\polylog}{{\rm polylog}}
\DeclareMathOperator*{\esssup}{ess\,sup}
\newcommand{\relu}{\text{ReLU}}
\newcommand{\sign}{\text{sign}}
\newcommand{\abs}{\text{abs}}
\newcommand{\m}{\mathbf{m}}
\newcommand{\n}{\mathbf{n}}

\newcommand{\R}{\mathbb{R}}
\renewcommand{\d}{\,\textup{d}}
\newcommand{\Frob}{\textup{F}}
\newcommand{\x}{\mathbf{x}}

\newcommand{\E}{\mathbb{E}}

\renewcommand{\L}{\mathcal{L}}
\newcommand{\F}{\mathscr{F}}
\newcommand{\mtx}[1]{\mathbf{#1}}
\renewcommand{\d}{\,\textup{d}}
\newcommand{\D}{\mathcal{D}}
\newcommand{\N}{\mathcal{N}}
\newcommand{\nrelu}{\mathcal{N}_{\text{ReLU}}}
\newcommand{\nrat}{\mathcal{N}_{\text{Rational}}}
\DeclareMathOperator{\Tr}{Tr}

\DeclareMathOperator{\diam}{diam}
\DeclareMathOperator{\dist}{dist}
\DeclareMathOperator{\cov}{cov}

\DeclareMathOperator{\HS}{HS}

\newcommand{\dobib}{ % Define the command
        \bibliography{references}  %use a bibtex bibliography file refs.bib
		\bibliographystyle{siam}  %use the plain bibliography style
    }

\title{Data-driven discovery of\\[1ex] Green's functions}   %note \\[1ex] is a line break in the title

\author{Nicolas Boull\'e}         %your name
\college{University College}      %your college

\degree{Doctor of Philosophy}     %the degree
\degreedate{Trinity 2022}         %the degree date

\begin{document}

\renewcommand{\dobib}{} % Un-define the command

%this baselineskip gives sufficient line spacing for an examiner to easily
%markup the thesis with comments
\baselineskip=18pt plus1pt

%set the number of sectioning levels that get number and appear in the contents
\setcounter{secnumdepth}{4}
\setcounter{tocdepth}{2}

\maketitle                  % create a title page from the preamble info
%include{dedication}        % include a dedication.tex file

\begin{abstract}
Discovering hidden partial differential equations (PDEs) and operators from data is an important topic at the frontier between machine learning and numerical analysis. Theoretical results and deep learning algorithms are introduced to learn Green's functions associated with linear partial differential equations and rigorously justify PDE learning techniques.

A theoretically rigorous algorithm is derived to obtain a learning rate, which characterizes the amount of training data needed to approximately learn Green's functions associated with elliptic PDEs. The construction connects the fields of PDE learning and numerical linear algebra by extending the randomized singular value decomposition to non-standard Gaussian vectors and Hilbert--Schmidt operators, and exploiting the low-rank hierarchical structure of Green's functions using hierarchical matrices.

Rational neural networks (NNs) are introduced and consist of neural networks with trainable rational activation functions. The highly compositional structure of these networks, combined with rational approximation theory, implies that rational functions have higher approximation power than standard activation functions. In addition, rational NNs may have poles and take arbitrarily large values, which is ideal for approximating functions with singularities such as Green's functions.

Finally, theoretical results on Green's functions and rational NNs are combined to design a human-understandable deep learning method for discovering Green's functions from data. This approach complements state-of-the-art PDE learning techniques, as a wide range of physics can be captured from the learned Green's functions such as dominant modes, symmetries, and singularity locations.

\end{abstract}
          % include the abstract

\begin{acknowledgements}

I would first like to thank my supervisors Patrick Farrell, Marie Rognes, and Alex Townsend for their guidance and suggestions. Their passion and excitement for the field of numerical analysis, as well as their high academic standards, have been a constant source of inspiration and motivation during my DPhil.

I would also like to thank my confirmation examiners, Christoph Reisinger and Justin Sirignano, for their comments and suggestions, as well as Andrew Stuart and Jared Tanner for accepting to be my thesis examiners.

This thesis benefited from discussions with great collaborators, including Efstathios Charalampidis, Vassilios Dallas, Christopher Earls, Ada Ellingsrud, Panayotis Kevrekidis, Seick Kim, Yuji Nakatsukasa, Alberto Paganini, Debasmita Samaddar, Tianyi Shi, and Jonasz S{\l}omka. 

I am grateful to Simula Research Laboratory for co-funding my DPhil along with University College, the Oxford-Radcliffe scholarship, and the InFoMM CDT.

I thank my friends and colleagues from Oxford and Cornell, Boris Andrews,  Francis Aznaran, Pablo Brubeck, Dan Fortunato, Marc Gilles, Gonzalo Gonzalez de Diego, Andrew Horning, Fabian Laakmann, Maike Meier, John Papadopoulos, Alex Puiu, Tianyi Shi, and Heather Wilber, who made this DPhil always enjoyable and fun with great academic and non-academic conversations.

Finally, I am grateful to my family and Tina for their continuous support and encouragements throughout the years.

\end{acknowledgements}
  % include an acknowledgements.tex file

\begin{romanpages}          % start roman page numbering
\tableofcontents            % generate and include a table of contents
%\listoffigures              % generate and include a list of figures
\end{romanpages}            % end roman page numbering

%now include the files of latex for each of the chapters etc

\chapter{Introduction}\label{sec_intro}

This thesis aims at understanding whether partial differential equations (PDEs) can be discovered from data by connecting standard mathematical fields, such as numerical linear algebra, probability, and PDE analysis, with modern deep learning techniques. We focus on learning Green's functions associated with linear PDEs from pairs of forcing functions and solutions. Theoretical bounds exploiting the regularity of the problem are derived and a practical deep learning algorithm is proposed.

\cref{chapt_PDE_learning} derives a theoretically-rigorous scheme for learning Green's functions associated with elliptic PDEs in three dimensions, given input-output pairs. A learning rate is obtained, giving a bound on the number of training pairs needed to learn a Green's function to within a prescribed accuracy with high probability. Along the way, the randomized singular value decomposition (SVD) is extended from matrices to Hilbert--Schmidt (HS) operators, and a quantity is introduced to measure the quality of the training forcing terms to learn Green's functions. The randomized SVD is a popular and effective algorithm for computing a near-best rank $k$ approximation of a matrix using matrix-vector products with standard Gaussian vectors.

\cref{chap_svd} extends the randomized SVD to multivariate Gaussian vectors, allowing one to incorporate prior knowledge of the matrix into the algorithm. This enables us to explore the continuous analogue of the randomized SVD for HS operators using operator-function products with functions drawn from a Gaussian process (GP). A new covariance kernel for GPs, based on weighted Jacobi polynomials, is constructed to rapidly sample the GP and control the smoothness of the randomly generated functions. Numerical examples on matrices and HS operators demonstrate the applicability of the algorithm.

\cref{chapt_rational} considers neural networks with rational activation functions. The choice of the nonlinear activation function in deep learning architectures is crucial and heavily impacts the performance of a neural network. We establish optimal bounds in terms of network complexity and prove that rational neural networks approximate smooth functions more efficiently than networks with Rectified Linear Unit (ReLU) activation functions with exponentially smaller depth. The flexibility and smoothness of rational activation functions make them an attractive alternative to ReLU, as demonstrated by numerical experiments.

\cref{chap_data_green} develops a data-driven approach for learning Green's functions using deep learning. By collecting physical system responses under excitations drawn from a Gaussian process, we train rational neural networks to learn Green's functions of hidden linear PDEs. These functions reveal human-understandable properties and features, such as linear conservation laws and symmetries, along with shock and singularity locations, boundary effects, and dominant modes. The technique is illustrated on several examples and allows us to capture a range of physics, including advection-diffusion, viscous shocks, and Stokes flow in a lid-driven cavity.

\section{Deep learning}

Deep learning has become an important topic across many domains of science due to its recent successes in image recognition, speech recognition, and drug discovery~\cite{hinton2012deep,krizhevsky2012imagenet,lecun2015deep,ma2015deep}. Deep learning techniques are based on objects called artificial neural networks (NNs), which apply a succession of mathematical transformations on an input variable $x$ to output a variable $y$, where $\mathcal{N}(x) = y$ and $\mathcal{N}$ denotes the neural network. An example of simple data fitting task is to assign labels $0$ or $1$ to points in $\R^2$, where, in this case, $x\in\R^2$ and $y\in\{0,1\}$~\cite{higham2019deep}. A large number of NN architectures, characterized by the type of mathematical operations used, have been proposed over the past decades for performing different tasks, such as convolutional NNs for classifying images~\cite{krizhevsky2012imagenet,lecun1989handwritten}, recurrent and long short-term memory neural networks for speech recognition~\cite{graves2014towards,hochreiter1997long,rumelhart1986learning}, and generative adversarial network to generate realistic images~\cite{goodfellow2014generative,karras2019style}. 

We consider one of the most standard types of deep learning model called feedforward neural networks or multilayer perceptrons~\cite[Chapt.~6]{goodfellow2016deep}. Let $L\geq 1$ be an integer and $n_1,n_L$ be the respective dimension of the input and output data. A feedforward network $\mathcal{N}:\R^{n_1}\to\R^{n_L}$, mapping from $\R^{n_1}$ to $\R^{n_L}$, with $L$ layers consists of a composition of $L-1$ functions $f_1,\ldots, f_{L-1}$ of the form
\[\mathcal{N}(x) = f_{L-1}\circ \cdots\circ f_1(x),\quad x\in \R^{n_1}.\]
At a given layer $1\leq i\leq L-1$, the nonlinear transformation $f_i:\R^{n_i}\to \R^{n_{i+1}}$ determines the output of the neural network. The layers $2\leq i\leq L-1$ are called the hidden layers of the network and their dimensionality determines the width of the network, while the number of layers is referred to as the depth~\cite[Chapt.~6]{goodfellow2016deep}. For $1\leq i\leq L-1$, we choose the function $f_i$ to be of the form 
\[f_i:x\mapsto \sigma(W_i x+b_i),\quad x\in \R^{n_i},\]
where $W_i\in \R^{n_i\times n_{i+1}}$ is a matrix called the weight matrix, $b_i\in \R^{n_{i+1}}$ is a bias vector, and $\sigma$ is a nonlinear function called the activation function (also called activation unit). The weight matrices and bias vectors are trainable parameters of the network, and their coefficients are usually obtained using a gradient-based optimization algorithm, such as stochastic gradient descent, applied to a training dataset containing examples of inputs and expected associated outputs of the network~\cite[Chapt.~6.2]{goodfellow2016deep}. In this thesis, we will measure the network complexity using its total number of parameters (\emph{i.e.,}~size) and number of layers (depth), which are standard measures in theoretical deep learning~\cite{anthony1999neural}.

In \cref{chapt_rational}, we will consider NNs with rational activation functions and derive theoretical results that quantify the size needed to approximate smooth functions within a prescribed accuracy. We will establish a connection between standard approximation theory for rational functions and the highly compositional structure of neural networks to show that rational neural networks require fewer parameters than ReLU networks to approximate smooth functions. We expect that the smoothness of rational neural networks away from their poles, together with their potential singularities, make them an interesting alternative to standard activation functions for physics-informed machine learning applications.

\section{Physics-informed machine learning}

Over the past decades, there has been spectacular progress in numerical techniques for solving PDEs, such as finite element methods, finite differences, and spectral methods~\cite{karniadakis2021physics}. However, solving inverse problems to identify parameters of a model or learn a physical model from real-world data remains highly challenging due to missing and noisy data~\cite{arridge2019solving,stuart2010inverse}. Hence, such problems are often ill-posed and require a data-driven approach. Recently, the fields of numerical analysis and machine learning have successfully converged towards physics-informed machine learning, which integrates partial data and prior knowledge on governing physical laws to solve inverse problems using neural networks~\cite{karniadakis2021physics}. The flexibility of the networks, due to the large potential choices of architectures, along with their generalization ability in the presence of big data, either generated by numerical simulations or acquired via experiments, makes them ideal for such tasks. On the other hand, the selection of a specific architecture is a challenging task and is difficult to justify mathematically due to the complexity of the models.

One example of problems that can be tackled by deep learning is to solve a PDE by training a NN on initial and boundary training data. Two popular approaches are physics-informed neural networks (PINNs)~\cite{raissi2019physics} and the deep Galerkin method~\cite{sirignano2018dgm}, which, in their original formulation, aim to solve PDEs of the form
\begin{equation} \label{eq_time_PDE}
\frac{\partial u}{\partial t} + \L(u) = 0,\quad x\in D\subset\R^d,\quad t\in [0,T],
\end{equation}
where the partial differential operator $\L$ is potentially nonlinear. The left-hand side of \cref{eq_time_PDE} is denoted by $f(x,t)$, \emph{i.e.},~$f\coloneqq u_t+\L(u)$. These techniques are attractive because they are mesh-free as they do not require a spatial discretization of the domain and can be applied in high dimensions. The PINN approach consists of approximating the solution $u$ to \cref{eq_time_PDE} by a neural network. This results in a physics-informed neural network $f$, which can be evaluated using chain rule and automatic differentiation~\cite{baydin2018automatic,raissi2019physics}. The loss function is expressed as a sum of a supervised loss of data measurements at the boundary and an unsupervised loss of PDE~\cite{karniadakis2021physics,raissi2019physics}:
\[\textup{Loss} = w_{\textup{data}}\L_{\textup{data}} + w_{\textup{PDE}}\L_{\textup{PDE}},\]
where $w_{\textup{data}}$ and $w_{\textup{PDE}}$ are weights balancing the two terms and $\L_{\textup{data}}$, $\L_{\textup{PDE}}$ are defined as
\[\L_{\textup{data}} = \frac{1}{N_{\textup{data}}}\sum_{i=1}^{N_{\textup{data}}}|u(x_i^{bdr},t_i^{bdr})-u_i^{bdr}|^2,\quad \L_{\textup{PDE}} = \frac{1}{N_{\textup{PDE}}}\sum_{j=1}^{N_{\textup{PDE}}}|f(x_j^{dom},t_j^{dom})|^2.\]
Here, $\{(x_i^{bdr},t_i^{bdr})\}$ are points sampled at the initial and boundary locations, while the points $\{(x_j^{dom},t_j^{dom})\}$ are sampled on the entire domain, and $\mathcal{L}_{\textrm{PDE}}$ is the average of the squared residual of the PDE evaluated at $\{(x_j^{dom},t_j^{dom})\}$. These methods have been generalized since their introductions to tackle a wide range of PDEs such as integro-differential equations~\cite{lu2021deepxde}, fractional PDEs~\cite{pang2019fpinns}, and stochastic PDEs~\cite{zhang2020learning}, and have been applied to problems in fluid mechanics~\cite{raissi2020hidden}, geophysics~\cite{li2020coupled}, and materials science~\cite{shukla2020physics}.

This thesis focuses on another aspect of physics-informed machine learning called PDE learning, whose aim is to discover, or learn, a mathematical model from data. We consider stationary PDEs of the form:
\[\L(u) = f,\]
where $\L$ is a partial differential operator, $f$ is called the forcing term, and $u$ the associated solution of the PDE. The approaches that dominate the PDE learning literature focus on the ``forward'' problem and aim to discover properties of the differential operator $\L$. As an example, sparsity-promoting techniques~\cite{Brunton,Rudy,zhang2018robust} consist of building a library of states $u$ and its spatio-temporal derivatives $u_t, u_x, u_{xx}, u_y,\ldots$ to identify parameters (or coefficients) and discover the main contributing terms in $\L$. Another method aims to find a symbolic expression for $\L$ and identify its dominant coefficients by solving a regression problem~\cite{Udrescu2020,udrescu2020ai2}. Finally, one can also project the operator $\L$ onto a low-dimensional subspace to build a reduced-order model and to significantly speed up standard numerical solvers~\cite{qian2021reduced,qian2020lift}.

An alternative approach, which we will consider, is to study the ``inverse'' problem and directly approximate the PDE solution operator, $\L^{-1}:f\mapsto u$, by an artificial neural network $\mathcal{N}$ from training pairs of forcing terms and solutions $\{f_j,u_j\}_{j=1}^N$~\cite{gin2020deepgreen,kovachki2021neural,li2020neural,li2020fourier,li2020multipole,lu2021learning,wang2021learning}. 
The network $\N$ takes a forcing term $f$ evaluated at a finite number of sensors $\{y_i\}_{i=1}^{N_f}$ and a point $x$ in the domain of $\L^{-1}(f)$ and outputs a real number approximating the solution $u$ to the PDE $\L(u) = f$ evaluated at $x$:
\[\N\left(\begin{bmatrix}
f(y_1) & \cdots & f(y_{N_f})
\end{bmatrix}^\top,x\right)\approx u(x).\]
The NN is then trained by minimizing the following loss function using stochastic gradient descent algorithms:
\[\text{Loss } = \frac{1}{N N_u N_f}\sum_{k=1}^N\sum_{i=1}^{N_u}\sum_{j=1}^{N_f}\left|\N\left(\begin{bmatrix}
f_k(y_1) & \cdots & f_k(y_{N_f})
\end{bmatrix}^\top,x_i\right)-u_k(x_i)\right|^2,\]
where $\{x_i\}_{i=1}^{N_u}$ are spatial points at which the solutions are measured. Unlike coefficient discovery techniques, this approach provides a fast solver for PDEs, which may outperform state-of-the-art numerical solvers~\cite{li2020fourier}. However, the physical interpretation of the learned solution operator remains highly challenging due to the mathematical complexity of the neural network that approximates it. Several black-box deep learning techniques are proposed to approximate the solution operator, which maps forcing terms $f$ to observations of the associated system's responses $u$ such that $\L(u) = f$. These methods are based on the concept of neural operators~\cite{kovachki2021neural}, which generalize neural networks to learn maps between infinite-dimensional function spaces, and mainly differ in their choice of the neural network architecture that is used to approximate the solution map. For example, Fourier neural operator~\cite{li2020fourier} uses a Fourier transform at each layer, while DeepONet~\cite{lu2021learning} contains a concatenation of `trunk' and `branch' networks to enforce additional structure.

On the theoretical side, most of the research has focused on the approximation theory of infinite-dimensional operators by NNs, such as the generalization of the universal approximation theorem~\cite{cybenko1989approximation} to shallow and deep NNs~\cite{chen1995universal,lu2021learning} as well as error estimates for Fourier neural operators and DeepONets with respect to the network width and depth~\cite{kovachki2021universal,kovachki2021neural,lanthaler2021error}. Other approaches aim to approximate the matrix of the discretized Green's functions associated with elliptic PDEs from matrix-vector multiplications by exploiting sparsity patterns or hierarchical structure of the matrix~\cite{lin2011fast,schafer2021sparse}. In addition,~\cite{de2021convergence} derived convergence rates for learning linear self-adjoint operators based on the assumption that the target operator is diagonal in the basis of the Gaussian prior.

In this thesis, we focus on learning linear partial differential operators $\L$ for which the solution operator can be written as an integral operator,
\[\L^{-1}(f)(x) = \int_{D}G(x,y)f(y)\d y = u(x),\]
whose kernel $G$ is known as the Green's function. Our approach contrasts with prior works because we aim to approximate the Green's function instead of the integral operator. As we will see in \cref{chapt_PDE_learning,chap_data_green}, imposing a prior structure on the solution operator offers theoretical and practical advantages over recent PDE learning techniques. First, standard mathematical techniques from elliptic PDE theory and numerical analysis can be exploited to derive rigorous results that quantify the amount of training data needed to learn the solution operator to within a prescribed accuracy. These types of results are notoriously challenging to obtain for deep learning algorithms due to the high nonlinearity of neural network architectures and the complexity of the optimization procedure. Secondly, unlike black-box deep learning techniques, it is possible to extract physical features of the original PDE from the associated Green's function, which is a well-understood mathematical object.

\section{Green's functions}

Throughout this thesis, we consider linear boundary value problems defined on a bounded domain $D\subset\R^d$, with $d\geq 1$, of the form:
\begin{subequations}
\begin{alignat*}{2}
\L u &= f, \quad &&\text{in }D,\\
u &= 0, \quad &&\text{on }\partial D,
\end{alignat*}
\end{subequations}
where $\L$ is a linear partial differential operator and $f:D\to \R$ is a given forcing function. A typical example of such problems is the Poisson equation in one dimension:
\begin{equation} \label{eq_Laplace}
-\frac{d^2 u}{dx^2} = f, \quad x \in (0, 1), \quad u(0)=u(1)=0.
\end{equation}
\cref{eq_Laplace} can be solved for any forcing function $f$ by introducing a kernel $G:[0,1]\times[0,1]\to\R$ so that the solution $u$ can be expressed as the following integral~\cite{evans10,green1854essay,roach1982green},
\begin{equation} \label{eq_int_Green}
u(x) = \int_0^1 G(x,y) f(y) \d y,\quad x\in[0,1].
\end{equation}
The function $G$ is called the Green's function and is a solution to the equation $\L G(x,y) = \delta(x-y)$, where $\delta$ is the Dirac delta function and $x,y\in[0,1]$. 

Green's functions are useful because they are independent of the forcing terms and only characterize the partial differential operators and boundary conditions. Once the Green's function has been determined, then the solution to \cref{eq_Laplace} with any forcing term can be obtained by computing the integral in \cref{eq_int_Green}, which is numerically easier than solving the original PDE and imposing the appropriate boundary conditions~\cite{roach1982green}. Additionally, several properties of the PDE can be recovered from the Green's function, such as symmetries or eigenvalues.

Traditional methods for finding Green's functions can be summarized as deriving analytical formulas, computing eigenvalue expansions, or numerically solving a singular PDE~\cite{evans10,roach1982green}. This is difficult when the geometry of the domain is complex or when the PDE has variable coefficients. Moreover, it requires knowledge of the partial differential operator, which may not be accessible in real applications~\cite{karniadakis2021physics}. Other works study properties of Green's functions and provide theoretical results such as decay bounds along the diagonal of the domain~\cite{cho2012global,gruter1982green,hofmann2007green,kang2010global} or low-rank structure on separable domains~\cite{bebendorf2003existence,boulle2022parabolic,engquist2018approximate}.

In this thesis, we aim to approximate Green's functions from pairs of forcing terms and system's responses $\{(f_j,u_j)\}_{j=1}^N$ by exploiting their low-rank structure on well-separated domains~\cite{bebendorf2003existence}, and combining it with randomized numerical linear algebra~\cite{halko2011finding}.

\section{Low-rank approximation} \label{sec_low_rank}

Let $\mtx{A}$ be an $m\times n$ real matrix with $m\geq n$ and $k\leq n$ be an integer. The best rank $k$ approximation to $\mtx{A}$ in the Frobenius norm is the $m\times n$ real matrix $\mtx{A}_k$, which is solution to the following minimization problem:
\begin{equation} \label{eq_low_rank_pb}
\min_{\mtx{A}_k\in \R^{m\times n}}\|\mtx{A}-\mtx{A}_k\|_\Frob\quad \text{subject to}\quad \textup{rank}(\mtx{A}_k)\leq k,
\end{equation}
where $\|\cdot\|_\Frob$ denotes the Frobenius norm defined as $\|\mtx{A}\|_\Frob=\sqrt{\Tr(\mtx{A}\mtx{A}^*)}$. The Eckart--Young theorem~\cite{eckart1936approximation} states that \eqref{eq_low_rank_pb} has a unique solution given by the truncation of the singular value decomposition of $\mtx{A}$ to the $k$th term. The SVD of an $m\times n$ real matrix $\mtx{A}$, with $m\geq n$, is a factorization of the form $\mtx{A} = \mtx{U} \mtx{\Sigma} \mtx{V}^*$, where $\mtx{U}$ is an $m \times m$ orthogonal matrix of left singular vectors, $\mtx{\Sigma}$ is an $m \times n$ diagonal matrix with entries $\sigma_1(\mtx{A})\geq \cdots \geq \sigma_{n}(\mtx{A})\geq 0$, and $\mtx{V}$ is an $n\times n$ orthogonal matrix of right singular vectors~\cite{golub2013matrix}. Then,
\[\min_{\substack{\mtx{A}_k\in \R^{m\times n}\\ \textup{rank}(\mtx{A_k})\leq k}}\|\mtx{A}-\mtx{A}_k\|_\Frob = \left(\sum_{j=k+1}^{n}\sigma_{j}(\mtx{A})^2\right)^{1/2},\]
where
\[\mtx{A}_k=\sum_{j=1}^k \sigma_j(\mtx{A}) u_j v_j^*.\]
Here, $u_j$ and $v_j$ denote the $j$th column of $\mtx{U}$ and $\mtx{V}$, respectively.

This result can be generalized to functions~\cite{schmidt1907theorie,townsend2014computing} and, in particular, Green's functions of the form $G:D_1\times D_2\to \R$, where $D_1,D_2\subset\R^d$. As an example, if $D_1=[a,b]$ and $D_2=[c,d]$ are two real intervals, and $G$ is square-integrable, then it can be written as the following infinite series, which converges in the $L^2(D_1\times D_2)$ sense to $G$,
\[G(x,y) = \sum_{\substack{j=1\\\sigma_j>0}}^\infty\sigma_j u_j(x) v_j(y),\quad x\in D_1,\,y\in D_2,\]
where $\{u_j\}_{j\geq 1}$ and $\{v_j\}_{j\geq 1}$ form an orthonormal basis of $L^2(D_1)$ and $L^2(D_2)$, and $\sigma_1\geq \sigma_2\geq \cdots\geq 0$ are called the singular values of $G$. This series is referred to as the SVD of $G$. Similar to matrices, the best rank $k$ approximant to $G$ is obtained by truncating its SVD after $k$ terms to obtain a separable approximation
\[G_k(x,y) = \sum_{j=1}^k\sigma_j u_j(x) v_j(y),\quad x\in D_1,\, y\in D_2.\]
By the Eckart--Young theorem, $G_k$ is solution to the following minimization problem:
\[\min_{\substack{f_j\in L^2(D_1)\\ g_j\in L^2(D_2)}}\|G-\sum_{j=1}^k f_j g_j\|_{L^2(D_1\times D_2)} = \|G-G_k\|_{L^2(D_1\times D_2)}=\left(\sum_{j=k+1}^\infty\sigma_j^2\right)^{1/2}.\]

Let $0<\epsilon<1$. If there exists an integer $k>0$ and a separable expression satisfying
\[\|G-\sum_{j=1}^k f_j g_j\|_{L^2(D_1\times D_2)}\leq \epsilon \|G\|_{L^2(D_1\times D_2)},\quad f_j\in L^2(D_1),\,g_j\in L^2(D_2),\]
then we say that $G$ has numerical rank smaller than $k$. We remark that one can easily obtain a bound on the tail of the singular values of $G$ by applying the Eckart--Young theorem as follows,
\[\left(\sum_{j=k+1}^\infty\sigma_j^2\right)^{1/2}=\min_{\substack{f_j\in L^2(D_1)\\ g_j\in L^2(D_2)}}\|G-\sum_{j=1}^k f_j g_j\|_{L^2(D_1\times D_2)}\leq \epsilon \|G\|_{L^2(D_1\times D_2)}.\]
When $k = \mathcal{O}(\log^\delta(1/\epsilon))$ for some small $\delta\in \mathbb{N}$ as $\epsilon \to 0$, then we say that $G$ has exponentially decaying singular values on $D_1\times D_2$.

\section{Randomized singular value decomposition} 

Computing the SVD of a matrix is a fundamental linear algebra task in machine learning~\cite{paterek2007improving}, statistics~\cite{wold1987principal}, and signal processing~\cite{alter2000singular,van1993subspace}. As we saw in \cref{sec_low_rank}, the SVD plays a central role in numerical linear algebra because truncating it after $k$ terms provides the best rank $k$ approximation to $\mtx{A}$ in the spectral and Frobenius norms~\cite{eckart1936approximation,mirsky1960symmetric}. Since computing the SVD of a large matrix can be computationally infeasible, there are various principal component analysis (PCA)~\cite{abdi2010principal,hotelling1933analysis,pearson1901liii} algorithms that perform dimensionality reduction by computing near-best rank $k$ matrix approximations from matrix-vector products~\cite{halko2011finding,martinsson2020randomized,nakatsukasa2020fast,nystrom1930praktische,williams2001using}. The randomized SVD uses matrix-vector products with random test vectors and is one of the most popular algorithms for constructing a low-rank approximation to $\mtx{A}$~\cite{halko2011finding,martinsson2020randomized}. While the error analysis performed in~\cite{halko2011finding} for the randomized SVD uses standard Gaussian random vectors, other random embedding techniques have been considered such as random permutations~\cite{ailon2009fast}, sparse sign matrices~\cite{clarkson2017low,meng2013low,Nelson2013OSNAP,urano2013fast}, and subsampled randomized trigonometric transforms (SRTTs)~\cite{ailon2006approximate,ailon2009fast,Parker95randombutterfly,woolfe2008fast} to mitigate the computational cost of Gaussian vectors in practical applications. Throughout this thesis, we will focus on Gaussian vectors because they yield a more precise error analysis (cf.~\cite[Sec.~8.3]{martinsson2020randomized}).

First, one performs the matrix-vector products $y_1 = \mtx{A}x_1,\,\ldots, \,y_{k+p} = \mtx{A}x_{k+p}$, where $x_1,\ldots,x_{k+p}$ are standard Gaussian random vectors with identically and independently distributed entries and $p\geq 1$ is an oversampling parameter. Then, one computes the economized QR factorization $\begin{bmatrix} y_1 & \cdots & y_{k+p} \end{bmatrix} = \mtx{Q}\mtx{R}$, before forming the rank $\leq k+p$ approximant $\mtx{Q}\mtx{Q}^*\mtx{A}$. Note that if $\mtx{A}$ is symmetric, one can form $\mtx{Q}\mtx{Q}^*\mtx{A}$ by computing $\mtx{Q}(\mtx{A}\mtx{Q})^*$ via matrix-vector products involving $\mtx{A}$; otherwise it requires the adjoint $\mtx{A}^*$. The quality of the rank $\leq k+p$ approximant $\mtx{Q}\mtx{Q}^*\mtx{A}$ is characterized by the following bound for $u,t\geq 1$~\cite[Thm.~10.7]{halko2011finding},
\begin{equation} \label{eq:RandomizedSVDBound}
\| \mtx{A} - \mtx{Q}\mtx{Q}^*\mtx{A} \|_\Frob \leq \left(1 + t\sqrt{\frac{3k}{p+1}}\,\right) \sqrt{\sum_{j=k+1}^n \sigma_j^2(\mtx{A})}+ut\frac{\sqrt{k+p}}{p+1}\sigma_{k+1}(\mtx{A}),
\end{equation}
with failure probability at most $2t^{-p}+e^{-u^2}$. The squared tail of the singular values of $\mtx{A}$, \emph{i.e.},~$\smash{\sqrt{\sum_{j=k+1}^n\sigma_j^2(\mtx{A})}}$, gives the best rank $k$ approximation error to $\mtx{A}$ in the Frobenius norm. This result shows that the randomized SVD can compute a near-best low-rank approximation to $\mtx{A}$ with high probability. In \cref{chapt_PDE_learning,chap_svd}, we will generalize this result to random vectors sampled from a multivariate normal distribution with any covariance matrix, and Hilbert--Schmidt operators.

\section{Hilbert--Schmidt operators} \label{sec_HS}

Hilbert--Schmidt operators generalize the notion of matrices acting on vectors to infinite dimensions with linear operators acting on functions~\cite[Ch.~4]{hsing2015theoretical}. First, let $D_1,D_2\subset\R^d$ be two domains with $d\geq 1$. For $1\leq p\leq \infty$, we denote by $L^p(D_1)$ the space of measurable functions defined on the domain $D_1$ with finite $L^p$ norm, where 
\begin{align*}
\|f\|_{L^p(D_1)} &= \left(\int_{D_1} |f(x)|^p\d x\right)^{1/p} \,\,\, \text{if } p<\infty,\\
\|f\|_{L^\infty(D_1)} &= \inf\left\{C>0,\,|f(x)|\leq C \text{ for almost every }x\in D_1\right\}.
\end{align*}
Since the space of square-integrable functions, $L^2(D_{1})$, is a separable Hilbert space, it admits a complete orthonormal basis $\{e_j\}_{j=1}^\infty$. 

A linear operator $\F: L^2(D_1)\to L^2(D_2)$ is an HS operator~\cite[Def.~4.4.2]{hsing2015theoretical} if it has finite HS norm, $\|\F\|_{\HS}$, defined as
\[
\|\F\|_{\HS} \coloneqq \left(\sum_{j=1}^\infty \|\F e_j\|_{L^2(D_2)}^2\right)^{1/2}<\infty.
\]
This norm does not depend on the choice of the basis~\cite[Thm.~4.4.1]{hsing2015theoretical}. The archetypical example of an HS operator is an integral operator $\F:L^2(D_1)\to L^2(D_2)$ defined as
\[
(\F f)(x)=\int_{D_1} G(x,y)f(y)\d y, \quad f\in L^2(D_1),\, x\in D_2,
\]
where $G\in L^2(D_2\times D_1)$ is the kernel of $\F$ and $\|\F\|_{\HS}=\|G\|_{L^2(D_2\times D_1)}$. The adjoint operator $\F^*:L^2(D_2)\to L^2(D_1)$ is defined as
\[
(\F^* g)(y)=\int_{D_2} G(x,y)g(x)\d x, \quad g\in L^2(D_2),\, y\in D_1.
\]

Since HS operators are compact operators, they have an SVD~\cite[Thm.~4.3.1]{hsing2015theoretical}. That is, that for any $f\in L^2(D_1)$ we have
\begin{equation} 
\mathscr{F}f = \sum_{j=1}^{\infty}\sigma_j\langle q_{1j}, f\rangle q_{2j},
\label{eq:HS_SVD} 
\end{equation} 
where the equality holds in the $L^2(D_2)$ sense. Here, $\sigma_1\geq \sigma_2\geq \cdots\geq 0$ denote the square roots of the eigenvalues of the self-adjoint operator $\F^*\F$, $\{q_{1j}\}$ are the orthonormal eigenvectors of $\F^*\F$, and $\{q_{2j}\}$ are the orthonormal eigenvectors of $\F\F^*$. We refer to $\{(\sigma_j,q_{1j},q_{2j})\}_{j=1}^\infty$ as the singular system of $\F$. When the HS operator is an integral operator, we refer to its singular values as the singular values of the underlying kernel.

Moreover, one finds that $\|\F\|_{\HS}^2 = \sum_{j=1}^\infty \sigma_j^2$, which shows that the HS norm is an infinite dimensional analogue of the Frobenius matrix norm $\|\cdot \|_{\textup{F}}$. In the same way that truncating the SVD after $k$ terms gives the best rank $k$ matrix approximation, truncating~\cref{eq:HS_SVD} gives the best rank $k$ approximation in the HS norm. That is,~\cite[Thm.~4.4.7]{hsing2015theoretical}
\[
\min_{u_j\in L^2(D_1), v_j\in L^2(D_2)}\|\F-\sum_{j=1}^k \langle u_j,\cdot\rangle v_j\|_{\HS} = \|\F-\F_k\|_{\HS}=\left(\sum_{j=k+1}^\infty \sigma_j^2\right)^{1/2},
\]
where the operator $\F_k$ is defined as
\[\F_kf = \sum_{j=1}^{k}\sigma_j\langle q_{1j},f\rangle q_{2j}, \quad f\in L^2(D_1).\]
This result is known as the Eckart--Young--Mirsky theorem~\cite{eckart1936approximation,mirsky1960symmetric}. We will exploit this theorem in \cref{chapt_PDE_learning} to extend the randomized SVD to HS operators and learn Green's functions.

\section{Quasimatrices} \label{sec_Quasimatrices} 
Quasimatrices are an infinite dimensional analogue of tall-skinny matrices~\cite{townsend2015continuous}.  Let $D_1,D_2\subseteq\R^d$ be two domains with $d\geq 1$, we say that $\mtx{\Omega}$ is a $D_1\times k$ quasimatrix, if $\mtx{\Omega}$ is a matrix with $k$ columns where each column is a function in $L^2(D_1)$. That is,
\[
\mtx{\Omega} = \begin{bmatrix} \omega_1 \, | & \! \cdots \! & | \, \omega_k \end{bmatrix}, \quad \omega_j\in L^2(D_1).
\]
Quasimatrices are useful to define analogues of matrix operations for HS operators~\cite{de1991alternative,stewart1998matrix,townsend2015continuous,trefethen1997numerical}. For example, if $\F:L^2(D_1)\to L^2(D_2)$ is an HS operator, then we write $\F\mtx{\Omega}$ to denote the quasimatrix obtained by applying $\F$ to each column of $\mtx{\Omega}$. Moreover, we write $\mtx{\Omega}^*\mtx{\Omega}$ and $\mtx{\Omega}\mtx{\Omega}^*$ to mean the following:
\[
\mtx{\Omega}^*\mtx{\Omega} = \begin{bmatrix}\langle \omega_1,\omega_1 \rangle & \cdots & \langle \omega_1,\omega_k \rangle\\ \vdots & \ddots &  \vdots\\
\langle \omega_k,\omega_1 \rangle & \cdots & \langle \omega_k,\omega_k \rangle \end{bmatrix}, \quad \mtx{\Omega}\mtx{\Omega}^* = \sum_{j=1}^k \omega_j(x)\omega_j(y),
\]
where $\langle \cdot, \cdot \rangle$ is the $L^2(D_1)$ inner-product.  Many operations for rectangular matrices in linear algebra can be generalized to quasimatrices such as the SVD, QR, LU, and Cholesky factorizations~\cite{townsend2015continuous}.

Throughout this thesis, the HS operator denoted by $\mtx{\Omega}\mtx{\Omega}^*\F : L^2(D_1)\to L^2(D_2)$ is given by $\mtx{\Omega}\mtx{\Omega}^*\F f = \sum_{j=1}^k \langle \omega_j,\F f\rangle \omega_j$. Moreover, if $\mtx{\Omega}$ has full column rank then $\mtx{P}_{\mtx{\Omega}}\F \coloneqq \mtx{\Omega}(\mtx{\Omega}^*\mtx{\Omega})^{\dagger}\mtx{\Omega}^*\F$ is the orthogonal projection of the range of $\F$ onto the column space of $\mtx{\Omega}$. Here, $(\mtx{\Omega}^*\mtx{\Omega})^{\dagger}$ is the pseudo-inverse of $\mtx{\Omega}^*\mtx{\Omega}$. This notation is convenient to state the generalization of the randomized SVD in infinite dimensions.

\section{Gaussian processes} \label{sec_GP}
A Gaussian process is an infinite dimensional analogue of a multivariate Gaussian distribution and a function drawn from a GP is analogous to a randomly generated vector. If $K: D\times D\to\R$ is a continuous symmetric positive semi-definite kernel, where $D\subseteq\mathbb{R}^d$ is a domain, then a GP is a stochastic process $\{X_t,\,t\in D\}$ such that for every finite set of indices $t_1,\ldots,t_n\in D$ the vector of random variables $(X_{t_1},\ldots,X_{t_n})$ is a multivariate Gaussian distribution with mean $(0,\ldots,0)$ and covariance $K_{ij} = K(t_i,t_j)$ for $1\leq i,j\leq n$. We denote a GP with mean $(0,\ldots,0)$ and covariance kernel $K$ by $\mathcal{GP}(0,K)$.

Since $K$ is a continuous symmetric positive semi-definite kernel, it has nonnegative eigenvalues $\lambda_1\geq \lambda_2\geq \cdots \geq 0$ and there is an orthonormal basis of eigenfunctions $\{\psi_j\}_{j=1}^\infty$ of $L^2(D)$ such that~\cite[Thm.~4.6.5]{hsing2015theoretical}:
\begin{equation} \label{eq_covariance_kernel}
K(x,y)=\sum_{j=1}^\infty\lambda_j\psi_j(x)\psi_j(y),\quad \int_{D} K(x,y)\psi_j(y)\d y= \lambda_j \psi_j(x),\quad x,y\in D,
\end{equation}
where the infinite sum is absolutely and uniformly convergent~\cite{mercer1909fun}. Note that the eigenvalues of $K$ are the ones of the integral operator with kernel $K$. In addition, we define the trace of the covariance kernel $K$ by $\smash{\Tr(K)\coloneqq \sum_{j=1}^\infty\lambda_j}<\infty$. The eigendecomposition of $K$ gives an algorithm for sampling functions from $\smash{\mathcal{GP}(0,K)}$.  In particular, if 
\[\omega = \sum_{j=1}^{\infty} \sqrt{\lambda_j} c_j\psi_j,\] 
where the coefficients $\{c_j\}_{j=1}^\infty$ are independent and identically distributed (i.i.d.) standard Gaussian random variables and the series converges in mean-square and uniformly, then $\omega\sim\mathcal{GP}(0,K)$. This is known as the Karhunen--Lo\`eve theorem~\cite{karhunen1946lineare,loeve1946functions}. We also have~\cite[Thm.~7.2.5]{hsing2015theoretical} % See https://arxiv.org/pdf/1807.02582.pdf sec. 4.1.3
\[
\E\!\left[\|\omega\|_{L^2(D)}^2\right]=\sum_{j=1}^\infty\lambda_j\E\!\left[c_j^2\right]\|\psi_j\|_{L^2(D)}^2=\sum_{j=1}^\infty \lambda_j=\int_{D}K(y,y)\,\d y<\infty,
\]
where the last equality is analogous to the fact that the trace of a matrix is equal to the sum of its eigenvalues. In this thesis, we restrict our attention to GPs with positive definite covariance kernels so that the eigenvalues of $K$ are strictly positive.

\section{Contribution}

The material of \cref{chapt_PDE_learning} to \cref{chap_data_green} is based on the following four papers with collaborators:
\begin{itemize}
\item \textbf{Learning elliptic PDEs with randomized linear algebra}\\
Nicolas Boull\'e and Alex Townsend\\
\textit{Foundations of Computational Mathematics}, 2022
\item \textbf{A generalization of the randomized singular value decomposition}\\
Nicolas Boull\'e and Alex Townsend\\
\textit{International Conference on Learning Representations}, 2022
\item \textbf{Rational neural networks}\\
Nicolas Boull\'e, Yuji Nakatsukasa, and Alex Townsend\\
\textit{Neural Information Processing Systems}, 2020
\item \textbf{Data-driven discovery of Green's functions with human-understandable deep learning}\\
Nicolas Boull\'e, Christopher J. Earls, and Alex Townsend\\
\textit{Scientific Reports}, 2022
\end{itemize}
My co-authors had advisory roles; I proved the main theoretical results, performed the numerical experiments, and was the lead author in writing the papers.

\dobib

\newcommand{\linfigwidth}{0.2\textwidth}
\setcounter{chapter}{1}
\renewcommand{\thefootnote}{\fnsymbol{footnote}}
\chapter[Learning elliptic PDEs with randomized linear algebra]{Learning elliptic PDEs with randomized linear algebra\footnotemark} \label{chapt_PDE_learning}

\footnotetext{This chapter is based on a paper with Alex Townsend~\cite{boulle2021learning}, published in Foundations of Computational Mathematics. Townsend had an advisory role; I proved the theoretical results and was the lead author in writing the paper.}

\renewcommand*{\thefootnote}{\arabic{footnote}}
\setcounter{footnote}{0}

Can one learn a differential operator from pairs of solutions and righthand sides? If so, how many pairs are required? These two questions have received significant research attention~\cite{feliu2020meta,li2020fourier,long2018pde,pang2019neural}. From data, one hopes to eventually learn physical laws of nature or conservation laws that elude scientists in the biological sciences~\cite{yazdani2020systems}, computational fluid dynamics~\cite{raissi2020hidden}, and computational physics~\cite{raissi2018deep}. The literature contains many highly successful practical schemes based on deep learning techniques~\cite{meng2020ppinn,raissi2019physics}. However, the challenge remains to understand when and why deep learning is effective theoretically. This chapter describes the first theoretically-justified scheme for discovering scalar-valued elliptic partial differential equations (PDEs) in three variables from input-output data and provides a rigorous learning rate. While our novelties are mainly theoretical, we hope to motivate future practical choices in PDE learning.

Let $D\subset\R^3$ be a bounded domain with Lipschitz smooth boundary, $L^2(D)$ be the space of square-integrable functions defined on $D$, $\mathcal{H}^k(D)$ be the space of $k$ times weakly differentiable functions in the $L^2$-sense, and $\mathcal{H}^1_0(D)$ be the closure of $\mathcal{C}_c^\infty(D)$ in $\mathcal{H}^1(D)$. Here, $\mathcal{C}_c^\infty(D)$ is the space of infinitely differentiable compactly supported functions on $D$. Roughly speaking, $\mathcal{H}^1_0(D)$ are the functions in $\mathcal{H}^1(D)$ that are zero on the boundary of $D$. We suppose that there is an unknown second-order uniformly elliptic linear PDE operator $\mathcal{L}:\mathcal{H}^2(D)\cap\mathcal{H}_0^1(D) \to L^2(D)$~\cite{evans10}, which takes the form
\begin{equation} 
\mathcal{L}u(x) = -\nabla \cdot \left(A(x) \nabla u\right)+c(x)\cdot\nabla u+d(x)u, \quad x\in D, \quad u|_{\partial D} = 0.
\label{eq:PDEForm}
\end{equation}
Here, for every $x\in D$, we have that $A(x)\in\mathbb{R}^{3\times 3}$ is a symmetric positive definite matrix with bounded coefficient functions so that $A_{ij}\in L^{\infty}(D)$, $c\in L^r(D)$ with $r\geq 3$, $d\in L^s(D)$ for $s\geq 3/2$, and $d(x)\geq 0$~\cite{kim2019green}. We emphasize that the regularity requirements on the variable coefficients are quite weak. 

The goal of PDE learning is to discover the operator $\mathcal{L}$ from $N\geq 1$ input-output pairs, \emph{i.e.}, $\{(f_j,u_j)\}_{j=1}^N$, where $\mathcal{L}u_j = f_j$ and $u_j|_{\partial D} = 0$ for $1\leq j\leq N$. There are two main types of PDE learning tasks: (1) Experimentally-determined input-output pairs, where one must do the best one can with the predetermined information and (2) Algorithmically-determined input-output pairs, where the data-driven learning algorithm can select $f_1,\ldots,f_N$ for itself. In this chapter, we focus on the PDE learning task where we have algorithmically-determined input-output pairs and aim to provide an upper bound on the sample complexity of the Green's function $G$ associated with $\L$, \emph{i.e.} characterize the number of pairs needed to learn $G$ within a prescribed accuracy. In particular, we suppose that the functions $f_1,\ldots,f_N$ are generated at random and are drawn from a Gaussian process (GP) (see~\cref{sec_GP}). Note that alternative strategies analogue to a power scheme in randomized numerical linear algebra~\cite{golub2013matrix,halko2011finding,rokhlin2010randomized,roweis1997algorithms} to generate forcing terms iteratively might lead to better approximation errors. To keep our theoretical statements manageable, we restrict our attention to PDEs of the form: 
\begin{equation} 
\mathcal{L}u = -\nabla \cdot \left(A(x) \nabla u\right), \quad x\in D, \quad u|_{\partial D} = 0.
\label{eq:PDEsimple}
\end{equation}
Lower-order terms in~\cref{eq:PDEForm} should cause few theoretical problems~\cite{bebendorf2003existence}, though our algorithm and our bounds get far more complicated. 

The approach that dominates the PDE learning literature is to directly learn $\mathcal{L}$ by either (1) learning parameters in the PDE~\cite{bonito2017diffusion,zhao2020learning}, (2) using neural networks (NNs) to approximate the action of the PDE on functions~\cite{raissi2018deep,raissi2018hidden,Karniadakis3,raissi2019physics,raissi2020hidden}, or (3) deriving a model from a library of operators via sparsity considerations~~\cite{Brunton,maddu2019stability,Rudy,schaeffer2017learning,voss2004nonlinear,wang2019variational}.  Instead of trying to learn the unbounded, closed operator $\mathcal{L}$ directly, we follow~\cite{boulle2021data,feliu2020meta,gin2020deepgreen} and discover the Green's function associated with $\mathcal{L}$. That is, we attempt to learn the function $G:D\times D\rightarrow\mathbb{R}^+\cup\{\infty\}$ such that~\cite{evans10}
\begin{equation} 
u_j(x) = \int_{D} G(x,y)f_j(y) \d y,\quad x\in D, \quad 1\leq j\leq N.
\label{eq:IntegralOperator}
\end{equation} 
Seeking $G$, as opposed to $\mathcal{L}$, has several theoretical benefits:
\begin{enumerate}[leftmargin=*,noitemsep]

\item The integral operator in~\cref{eq:IntegralOperator} is compact~\cite{edmunds2013bounded}, while $\mathcal{L}$ is only closed~\cite{edmunds2018spectral}. This allows $G$ to be rigorously learned by input-output pairs $\{(f_j,u_j)\}_{j=1}^N$, as its range can be approximated by finite-dimensional spaces (see~\cref{th_Green}). 

\item It is known that $G$ has a hierarchical low-rank structure~\cite[Thm.~2.8]{bebendorf2003existence}: for $0<\epsilon<1$, there exists a function $G_k(x,y) = \sum_{j=1}^k g_j(x)h_j(y)$ with $k = \mathcal{O}(\log^4(1/\epsilon))$ such that~\cite[Thm.~2.8]{bebendorf2003existence}
\[
\left\|G - G_k\right\|_{L^2(X\times Y)}\leq \epsilon\left\|G\right\|_{L^2(X\times \hat{Y})},
\]
where $X,Y\subseteq D$ are sufficiently separated domains, and $Y\subseteq\hat{Y}\subseteq D$ denotes a larger domain than $Y$ (see \cref{theo_bebendorf} for the definition). The further apart $X$ and $Y$, the faster the singular values of $G$ decay. Moreover, $G$ also has an off-diagonal decay property~\cite{gruter1982green,kang2010global}:
\[
G(x,y) \leq \frac{c}{\| x - y\|_2}\|G\|_{L^2(D\times D)}, \quad x\neq y,\, x \in D,\, y \in D,
\]
where $c$ is a constant independent of $x$ and $y$. Exploiting these structures of $G$ leads to a rigorous algorithm for constructing a global approximant to $G$ (see~\cref{sec_approx_Green}). 

\item The function $G$ is smooth away from its diagonal, allowing one to efficiently approximate it~\cite{gruter1982green}.

\end{enumerate} 
Once a global approximation $\tilde{G}$ has been constructed for $G$ using input-output pairs, given a new righthand side $f$ one can directly compute the integral in~\cref{eq:IntegralOperator} to obtain the corresponding solution $u$  to~\cref{eq:PDEForm}. Usually, numerically computing the integral in~\cref{eq:IntegralOperator} must be done with sufficient care as $G$ possesses a singularity when $x = y$. However, our global approximation $\tilde{G}$ has a hierarchical structure and is constructed as $0$ near the diagonal. Therefore, for each fixed $x\in D$, we simply recommend that $\int_{D} \tilde{G}(x,y) f_j(y)\d y$ is partitioned into the panels that corresponds to the hierarchical decomposition, and then discretized each panel with a quadrature rule.

There are two main contributions in this chapter: (1) the generalization of the randomized singular value decomposition (SVD) algorithm for learning matrices from matrix-vector products to Hilbert--Schmidt (HS) operators and (2) a theoretical learning rate for discovering Green's functions associated with PDEs of the form~\cref{eq:PDEsimple}. These contributions are summarized in \cref{th_tropp_random_svd_Frob,th_Green}.

\cref{th_tropp_random_svd_Frob} says that, with high probability, one can recover a near-best rank $k$ HS operator using $k+p$ operator-function products, for a small integer $p$. In the bound of the theorem, a quantity, denoted by $0< \gamma_k \leq 1$, measures the quality of the input-output training pairs (see~\cref{sec:TwoWarnings,sec_quality_kernel}). We then combine~\cref{th_tropp_random_svd_Frob} with the theory of Green's functions for elliptic PDEs to derive a theoretical learning rate for PDEs.

In \cref{th_Green}, we show that Green's functions associated with uniformly elliptic PDEs in three dimensions can be recovered using $N=\mathcal{O}(\epsilon^{-6}\log^4(1/\epsilon))$ input-output pairs $(f_j,u_j)_{j=1}^N$  to within an accuracy of $\mathcal{O}(\Gamma_\epsilon^{-1/2}\log^3(1/\epsilon)\epsilon)$ with high probability, for $0<\epsilon<1$. Our learning rate associated with uniformly elliptic PDEs in three variables is therefore $\mathcal{O}(\epsilon^{-6}\log^4(1/\epsilon))$. The quantity $0< \Gamma_\epsilon\leq 1$ (defined in~\cref{eq_define_gamma_eps}) measures the quality of the GP used to generate the random functions $\{f_j\}_{j=1}^N$ for learning $G$. We emphasize that the number of training pairs is small only if the GP's quality is high. The probability bound in \cref{th_Green} implies that the constructed approximation is close to $G$ with high probability and converges almost surely to the Green's function as $\epsilon\to 0$.

\section{Low-rank approximation of Hilbert--Schmidt operators} \label{sec_random_SVD}
In a landmark paper, Halko, Martinsson, and Tropp proved that one could learn the column space of a finite matrix---to high accuracy and with a high probability of success---by using matrix-vector products with standard Gaussian random vectors~\cite{halko2011finding}. We now set out to generalize this from matrices to HS operators. Alternative randomized low-rank approximation techniques such as the generalized Nystr\"om method~\cite{nakatsukasa2020fast} might also be generalized in a similar manner.
Since the proof is relatively long, we state our final generalization now.

\begin{theorem} \label{th_tropp_random_svd_Frob}
Let $D_1,D_2\subseteq \mathbb{R}^d$ be domains with $d\geq 1$ and $\F:L^2(D_1)\to L^2(D_2)$ be an HS operator. Select a target rank $k\geq 1$, an oversampling parameter $p\geq 2$, and a $D_1\times (k+p)$ quasimatrix $\mtx{\Omega}$ such that each column is i.i.d. and drawn from $\mathcal{GP}(0,K)$, where $K:D_1\times D_1\to\R$ is a continuous symmetric positive definite kernel with eigenvalues $\lambda_1\geq \lambda_2\geq \cdots>0$. If $\mtx{Y}=\F\mtx{\Omega}$, then
\begin{equation}
\label{eq:MainExpectationBound}
\mathbb{E}\!\left[\|\F - \mtx{P}_\mtx{Y}\F\|_{\HS}\right]\leq \left(1+\sqrt{\frac{1}{\gamma_k}\frac{k(k+p)}{p-1}}\,\right)\left(\sum_{j=k+1}^{\infty}\sigma_j^2\right)^{1/2},
\end{equation}
where $\gamma_k = k/(\lambda_1\Tr(\mtx{C}^{-1}))$ with $\mtx{C}_{ij}=\int_{D_1\times D_1}v_i(x)K(x,y)v_j(y)\d x\d y$ for $1\leq i,j\leq k$. Here, $\mtx{P}_\mtx{Y}$ is the orthogonal projection onto the vector space spanned by the columns of $\mtx{Y}$, $\sigma_j$ is the $j$th singular value of $\F$, and $v_j$ is the $j$th right singular vector of $\F$. 

Assume further that $p\geq 4$, then for any $s,t\geq 1$, we have
\begin{equation}
\label{eq:MainProbabilityBound}
\|\F-\mtx{P}_{\mtx{Y}}\F\|_{\HS}\leq \sqrt{1+ t^2s^2 \frac{3}{\gamma_k}\frac{k(k+p)}{p+1}\sum_{j=1}^\infty \frac{\lambda_j}{\lambda_1}}\,\left(\sum_{j=k+1}^\infty\sigma_j^2\right)^{1/2},
\end{equation}
with probability $\geq 1 - t^{-p}-[s e^{-(s^2-1)/2}]^{k+p}$.
\end{theorem}

We remark that the term $[s e^{-(s^2-1)/2}]^{k+p}$ in the statement of \cref{th_tropp_random_svd_Frob} is bounded by $e^{-s^2}$ for $s\geq 2$ and $k+p\geq 5$. The term $0\leq \gamma_k\leq 1$ is discussed in \cref{sec_quality_kernel} and is bounded by the inverse of the harmonic mean of $k$ eigenvalues of the covariance kernel under some conditions on the kernel eigenvectors and the right singular vectors of the Hilbert--Schmidt operator $\F$ (see~\cref{lem_bound_sigma}). In the rest of the section, we prove this theorem.

\subsection{Three caveats that make the generalization non-trivial}\label{sec:TwoWarnings}
One might imagine that the generalization of the randomized SVD algorithm from matrices to HS operators is trivial, but this is not the case due to three caveats.

First, the randomized SVD on finite matrices always uses matrix-vector products with standard Gaussian random vectors~\cite{halko2011finding}. However, for GPs, one must always have a continuous kernel $K$ in $\mathcal{GP}(0,K)$, which discretizes to a non-standard multivariate Gaussian distribution. Therefore, we must extend~\cite[Thm.~10.5]{halko2011finding} to allow for non-standard multivariate Gaussian distributions. The discrete version of our extension is the following:

\begin{corollary} \label{cor_finite_dim}
Let $\mtx{A}$ be a real $n_2\times n_1$ matrix with singular values $\sigma_1\geq \cdots \geq \sigma_{\min\{n_1,n_2\}}$. Choose a target rank $k\geq 1$ and an oversampling parameter $p\geq 2$. Draw an $n_1\times (k+p)$ Gaussian matrix, $\mtx{\Omega}$, with independent columns where each column is i.i.d. from a multivariate Gaussian distribution with mean $(0,\ldots,0)^\top$ and positive definite covariance matrix $\mtx{K}$. If $\mtx{Y}=\mtx{A}\mtx{\Omega}$, then the expected approximation error is bounded by
\begin{equation} \label{eq_expect_cor}
\mathbb{E}\left[\|\mtx{A}-\mtx{P}_\mtx{Y}\mtx{A}\|_{\textup{F}}\right]\leq \left(1+\sqrt{\frac{k+p}{p-1}\sum_{j=n_1-k+1}^{n_1}\frac{\lambda_1}{\lambda_j}}\,\right)\left(\sum_{j=k+1}^\infty \sigma_j^2\right)^{1/2},
\end{equation}
where $\lambda_1 \geq \cdots \geq \lambda_{n_1} > 0$ are the eigenvalues of $\mtx{K}$ and $\mtx{P}_\mtx{Y}$ is the orthogonal projection onto the vector space spanned by the columns of $\mtx{Y}$. Assume further that $p\geq 4$, then for any $s,t\geq 1$, we have
\[\|\mtx{A}-\mtx{P}_\mtx{Y}\mtx{A}\|_{\textup{F}}\leq \left(\!1+ ts\cdot \sqrt{\frac{3(k+p)}{p+1}\left(\sum_{j=1}^{n_1}\lambda_j\right)\sum_{j=n_1-k+1}^{n_1}\frac{1}{\lambda_j}}\,\right)\!\left(\sum_{j=k+1}^\infty\sigma_j^2\right)^{1/2},\]
with probability $\geq 1 - t^{-p}-[s e^{-(s^2-1)/2}]^{k+p}$.
\end{corollary}

Choosing a covariance matrix $\mtx{K}$ with eigenvalue decay so that $\lim_{n_1\rightarrow\infty}\sum_{j=1}^{n_1} \lambda_j <\infty$ allows $\E[\|\mtx{\Omega}\|_{\textup{F}}^2]$ to remain bounded as $n_1\to\infty$. This is of interest when applying the randomized SVD algorithm to extremely large matrices and is critical for HS operators. A stronger statement of this result (see \cref{th_svd_main}) shows that prior information on $\mtx{A}$ can be incorporated into the covariance matrix to achieve lower approximation error than the randomized SVD with standard Gaussian vectors.

Secondly, we need an additional essential assumption. The kernel in $\mathcal{GP}(0,K)$ is ``reasonable" for learning $\F$, where reasonableness is measured by the quantity $\gamma_k$ in \cref{th_tropp_random_svd_Frob}. If the first $k$ right singular functions of the HS operator $v_1,\ldots,v_k$ are spanned by the first $k+m$ eigenfunctions of $K$ $\psi_1,\ldots,\psi_{k+m}$, for some $m\in \mathbb{N}$, then (see~\cref{eq_lower_bound_gamma} and~\cref{lem_bound_sigma})
\[
\frac{1}{k}\sum_{j=1}^k\frac{\lambda_1}{\lambda_j}\leq \frac{1}{\gamma_k} \leq \frac{1}{k}\sum_{j=m+1}^{k+m}\frac{\lambda_1}{\lambda_j}.
\] 
In the matrix setting, this assumption always holds with $m=n_1-k$ (see~\cref{cor_finite_dim}) and one can have $\gamma_k = 1$ when $\lambda_1=\cdots=\lambda_{n_1}$~\cite[Thm.~10.5]{halko2011finding}.

Finally, probabilistic error bounds for the randomized SVD in~\cite{halko2011finding} are derived using tail bounds for functions of standard Gaussian matrices~\cite[Sec.~5.1]{ledoux2001concentration}. Unfortunately, we are not aware of tail bounds for non-standard Gaussian quasimatrices. This results in a weaker bounds by a factor of $\sqrt{k+p}$ in \cref{cor_finite_dim} compared to~\cite[Thm.~10.7]{halko2011finding}.

\subsection{Deterministic error bound}
Apart from the three caveats, the proof of \cref{th_tropp_random_svd_Frob} follows the outline of the argument in~\cite[Thm.~10.5]{halko2011finding}. We define two quasimatrices $\mtx{U}$ and $\mtx{V}$ containing the left and right singular functions of $\F$ so that the $j$th column of $\mtx{V}$ is $v_j$. We also denote by $\mtx{\Sigma}$ the infinite diagonal matrix with the singular values of $\F$, \emph{i.e.}, $\sigma_1\geq\sigma_2\geq\cdots\geq0$, on the diagonal. Finally, for a fixed $k\geq 1$, we define the $D_1\times k$ quasimatrix as the truncation of $\mtx{V}$ after the first $k$ columns and $\mtx{V}_2$ as the remainder. Similarly, we split $\mtx{\Sigma}$ into two parts:
\[\begin{array}{@{}c@{}c@{}c@{}c@{}c@{}c}
        && k & \infty &\\
    \mtx{\Sigma} =  &\left. \begin{array}{c} \\ \\ \end{array} \!\!\! \right( &
    \begin{array}{c} \mtx{\Sigma}_1 \\ 0 \end{array} &
    \begin{array}{c} 0 \\ \mtx{\Sigma}_2 \end{array} &
    \left. \!\!\!\begin{array}{c} \\ \\ \end{array}  \right) &
    \begin{array}{c} k \\ \infty \\ \end{array}
\end{array}.\]
We are ready to prove an infinite dimensional analogue of~\cite[Thm.~9.1]{halko2011finding} for HS operators.

\begin{theorem}[Deterministic error bound] \label{th_tropp_deter_svd}
Let $\F:L^2(D_1)\to L^2(D_2)$ be an HS operator with SVD given in~\cref{eq:HS_SVD}. Let $\mtx{\Omega}$ be a $D_1\times k$ quasimatrix and $\mtx{Y}=\F\mtx{\Omega}$. If $\mtx{\Omega}_1=\mtx{V}_1^*\mtx{\Omega}$ and $\mtx{\Omega}_2=\mtx{V}_2^*\mtx{\Omega}$, then assuming $\mtx{\Omega}_1$ has full rank, we have
\[
\|\F-\mtx{P}_\mtx{Y}\F\|_{\HS}^2\leq \|\mtx{\Sigma}_2\|_{\HS}^2+\|\mtx{\Sigma}_2\mtx{\Omega}_2\mtx{\Omega}_1^{\dagger}\|_{\HS}^2,
\]
where $\mtx{P}_\mtx{Y}=\mtx{Y}(\mtx{Y}^*\mtx{Y})^\dagger\mtx{Y}^*$ is the orthogonal projection onto the space spanned by the columns of $\mtx{Y}$ and $\smash{\mtx{\Omega}_1^{\dagger} = (\mtx{\Omega}_1^*\mtx{\Omega}_1)^{-1}\mtx{\Omega}_1^*}$. 
\end{theorem}
\begin{proof}
First, note that because $\mtx{U}\mtx{U}^*$ is the orthonormal projection onto the range of $\F$ and $\mtx{U}$ is a basis for the range, we have
\[
\|\F-\mtx{P}_\mtx{Y}\F\|_{\HS} = \|\mtx{U}\mtx{U}^*\F-\mtx{P}_{\mtx{Y}}\mtx{U}\mtx{U}^*\F\|_{\HS}.
\]
By Parseval's theorem~\cite[Thm.~4.18]{rudin1987real}, we have
\[ 
\|\mtx{U}\mtx{U}^* \F-\mtx{P}_{\mtx{Y}}\mtx{U}\mtx{U}^* \F\|_{\HS}= \|\mtx{U}^*\mtx{U}\mtx{U}^* \F-\mtx{U}^*\mtx{P}_{\mtx{Y}}\mtx{U}\mtx{U}^* \F \mtx{V}\|_{\HS}.
\]
Moreover, we have the equality $\|\F-\mtx{P}_{\mtx{Y}}\F\|_{\HS}=\|(\mtx{I}-\mtx{P}_{\mtx{U}^*\mtx{Y}})\mtx{U}^*\F\mtx{V}\|_{\HS}$ 
because the inner product $\langle\sum_{j=1}^\infty \alpha_j u_j,\sum_{j=1}^\infty \beta u_j\rangle=0$ if and only if $\sum_{j=1}^\infty \alpha_j \beta_j=0$. We now take $\mtx{A}=\mtx{U}^*\F\mtx{V}$, which is a bounded infinite matrix such that $\|\mtx{A}\|_{\textup{F}}=\|\F\|_{\HS}<\infty$. The statement of the theorem immediately follows from the proof of~\cite[Thm.~9.1]{halko2011finding}.
\end{proof}

This theorem shows that the bound on the approximation error $\|\F-\mtx{P}_\mtx{Y}\F\|_{\HS}$ depends on the singular values of the HS operator and the test matrix $\mtx{\Omega}$.

\subsection{Probability distribution of \texorpdfstring{$\mtx{\Omega}_1$}{Omega\_1}} \label{sec_proba_omega_1}

If the columns of $\mtx{\Omega}$ are independent and identically distributed as $\mathcal{GP}(0,K)$, then the matrix $\mtx{\Omega}_1$ in~\cref{th_tropp_deter_svd} is of size $k\times \ell$ with entries that follow a Gaussian distribution.  To see this, note that 
\[
\mtx{\Omega}_1 = \mtx{V}_1^*\mtx{\Omega} = 
\left(
\begin{array}{ccc}
\langle v_1,\omega_1\rangle & \cdots & \langle v_1, \omega_{\ell}\rangle\\
\vdots & \ddots & \vdots \\
\langle v_k,\omega_1\rangle & \cdots & \langle v_k,\omega_{\ell}\rangle
\end{array}
\right), \quad \omega_j\sim\mathcal{GP}(0,K). 
\]
If $\omega \sim \mathcal{GP}(0,K)$ with $K$ given in \cref{eq_covariance_kernel}, then we find that 
\[\langle v,\omega\rangle \sim \mathcal{N}\left(0,\sum_{j=1}^\infty \lambda_j \langle v,\psi_j\rangle^2\right)\]
so we conclude that $\mtx{\Omega}_1$ has Gaussian entries with zero mean. Finding the covariances between the entries is more involved. 
\begin{lemma} \label{prop_cov_matrix}
With the same setup as~\cref{th_tropp_deter_svd}, suppose that the columns of $\mtx{\Omega}$ are independent and identically distributed as $\mathcal{GP}(0,K)$. Then, the matrix $\mtx{\Omega}_1=\mtx{V}_1^*\mtx{\Omega}$ in~\cref{th_tropp_deter_svd} has independent columns and each column is identically distributed as a multivariate Gaussian with positive definite covariance matrix $\mtx{C}$ given by
\begin{equation} \label{eq_def_sigma_k}
\mtx{C}_{ij} = \int_{D_1\times D_1} v_i(x)K(x,y)v_{j}(y)\d x\d y, \quad 1\leq i,j\leq k,
\end{equation}
where $v_i$ is the $i$th column of $\mtx{V}_1$.
\end{lemma}
\begin{proof}
We already know that the entries are Gaussian with mean $0$. Moreover, the columns are independent because $\omega_1,\ldots,\omega_\ell$ are independent. Therefore, we focus on the covariance matrix. Let $1\leq i,i'\leq k$, $1\leq j,j'\leq \ell$, then since $\mathbb{E}\!\left[\langle v_i,\omega_j\rangle\right] = 0$ we have
\[
\cov(\langle v_i,\omega_j\rangle, \langle v_{i'},\omega_{j'}\rangle) 
= \E\left[ \langle v_i,\omega_j\rangle\, \langle v_{i'},\omega_{j'}\rangle \right] = \mathbb{E}\left[X_{ij}X_{i'j'}\right],
\]
where $X_{ij} = \langle v_i,\omega_j\rangle$. Since $\langle v_i,\omega_j\rangle\sim\sum_{n=1}^\infty \sqrt{\lambda_n}c_n^{(j)} \langle v_i,\psi_n\rangle$, where $c_n^{(j)}\sim\mathcal{N}(0,1)$, we have
\[
\cov(\langle v_i,\omega_j\rangle, \langle v_{i'},\omega_{j'}\rangle) = \E\left[\lim_{m_1,m_2\to\infty}X_{ij}^{m_1} X_{i'j'}^{m_2}\right],\quad X_{ij}^{m_1}\coloneqq\sum_{n=1}^{m_1} \sqrt{\lambda_n}c_n^{(j)}  \langle v_i,\psi_n\rangle.
\]
We first show that $\lim_{m_1,m_2\to\infty}\left|\E\!\left[\!X_{ij}^{m_1}X_{i'j'}^{m_2}\right]-\E\!\left[X_{ij}X_{i'j'}\right]\right|=0$. For any $m_1,m_2\geq 1$, we have by the triangle inequality,
\begin{align*}
\left|\E\!\left[\!X_{ij}^{m_1}X_{i'j'}^{m_2}\!\right]-\E\!\left[X_{ij}X_{i'j'}\right]\right| \!\!&\leq \E\!\left[\left|X_{ij}^{m_1}X_{i'j'}^{m_2}-X_{ij}X_{i'j'}\right|\right]\\
&\leq \E\!\left[\left|(X_{ij}^{m_1}-X_{ij})X_{i'j'}^{m_2}\right|\right] \!\!+\E\!\left[\left|X_{ij}(X_{i'j'}^{m_2}-X_{i'j'})\right|\right]\\
&\!\!\!\!\!\!\!\!\!\!\!\!\!\!\!\!\!\!\!\!\!\!\!\!\!\!\!\!\!\!\!\!\!\!\!\!
\leq 
\E\!\left[\left|X_{ij}^{m_1}-X_{ij}\right|^2\right]^{\tfrac{1}{2}}\!\E\!\left[\left|X_{i'j'}^{m_2}\right|^2\right]^{\tfrac{1}{2}}\!\!
+
\E\!\left[\left|X_{i'j'}-X_{i'j'}^{m_2}\right|^2\right]^{\tfrac{1}{2}}\!\E\!\left[\left|X_{ij}\right|^2\right]^{\tfrac{1}{2}}\!,
\end{align*}
where the last inequality follows from the Cauchy--Schwarz inequality. We now set out to show that both terms in the last inequality converge to zero as $m_1,m_2\to\infty$. The terms $\smash{\E[|X_{i'j'}^{m_2}|^2]}$ and $\smash{\E[|X_{ij}|^2]}$ are bounded by $\sum_{n=1}^\infty \lambda_n<\infty$, using the Cauchy--Schwarz inequality. Moreover, we have
\[
\E\left[\left|X_{ij}^{m_1}-X_{ij}\right|^2\right] = \E\left[\left|\sum_{n=m_1+1}^\infty \sqrt{\lambda_n}c_n^{(j)} \langle v_i,\psi_n\rangle\right|^2\right]\leq \sum_{n=m_1+1}^\infty \lambda_n \xrightarrow[m_1 \to \infty]{} 0,
\]
because $X_{ij} - X_{ij}^{m_1}\sim \mathcal{N}(0,\sum_{n=m_1+1}^\infty\lambda_n\langle v_i,\psi_n\rangle^2)$. Then, we find that $\cov(X_{ij}, X_{i'j'}) = \lim_{m_1,m_2\to\infty}\E[X_{ij}^{m_1}X_{i'j'}^{m_2}]$ and we obtain
\begin{align*}
\cov(X_{ij},X_{i'j'})
& =\lim_{m_1,m_2\to\infty}\E\left[\sum_{n=1}^{m_1}\sum_{n'=1}^{m_2} \sqrt{\lambda_n \lambda_{n'}}c_n^{(j)}c_{n'}^{(j')} \langle v_i,\psi_n\rangle\langle v_{i'},\psi_{n'}\rangle\right] \\
&=\lim_{m_1,m_2\to\infty}\sum_{n=1}^{m_1}\sum_{n'=1}^{m_2} \sqrt{\lambda_n \lambda_{n'}}\E[c_n^{(j)}c_{n'}^{(j')}] \langle v_i,\psi_n\rangle\langle v_{i'},\psi_{n'}\rangle.
\end{align*}
The latter expression is zero if $n\neq n'$ or $j\neq j'$ because then $c_n^{(j)}$ and $c_{n'}^{(j')}$ are independent random variables with mean $0$. Since $\E[(c_n^{(j)})^2] = 1$, we have
\[
\cov(X_{ij},X_{i'j'})= 
\begin{cases}
\sum_{n=1}^\infty \lambda_n \langle v_i,\psi_n\rangle\langle v_{i'},\psi_{n}\rangle,  & j=j',\\
0, & \text{otherwise}.
\end{cases}
\]
The result follows as the infinite sum is equal to the integral in~\cref{eq_def_sigma_k}.  To see that $\mtx{C}$ is positive definite, let $a\in \R^k$, then $a^* \mtx{C} a=\E[Z_a^2]\geq 0$, where $Z_a\sim\mathcal{N}(0,\sum_{n=1}^\infty \lambda_n \langle a_1 v_1+\cdots+a_k v_k,\psi_n\rangle^2)$. Moreover, $a^* \mtx{C} a = 0$ implies that $a=0$ because $v_1,\ldots,v_k$ are orthonormal and $\{\psi_n\}$ is an orthonormal basis of $L^2(D_1)$. 
\end{proof}

\cref{prop_cov_matrix} gives the distribution of the matrix $\mtx{\Omega}_1$, which is essential to prove \cref{th_tropp_random_svd_Frob} in \cref{sec_proof_thm_1}. In particular, $\mtx{\Omega}_1$ has independent columns that are each distributed as a multivariate Gaussian with covariance matrix given in \cref{eq_def_sigma_k}.

\subsection{Quality of the covariance kernel} \label{sec_quality_kernel}
To investigate the quality of the kernel, we introduce the Wishart distribution, which is a family of probability distributions over symmetric and nonnegative-definite matrices that often appear in the context of covariance matrices~\cite{wishart1928generalised}. If $\mtx{\Omega}_1$ is a $k\times \ell$ random matrix with independent columns, where each column is a multivariate Gaussian distribution with mean $(0,\ldots,0)^\top$ and covariance $\mtx{C}$, then $\mtx{A}=\mtx{\Omega}_1 \mtx{\Omega}_1^*$ has a Wishart distribution~\cite{wishart1928generalised}. We write $\mtx{A}\sim W_k(\ell,\mtx{C})$. We note that $\|\mtx{\Omega}_1^{\dagger}\|_{\textup{F}}^2=\Tr[(\mtx{\Omega}_1^{\dagger})^*\mtx{\Omega}_1^{\dagger}]=\Tr(\mtx{A}^{-1})$, where the second equality holds with probability one because the matrix $\mtx{A}=\mtx{\Omega}_1\mtx{\Omega}_1^*$ is invertible with probability one (see~\cite[Thm.~3.1.4]{muirhead2009aspects}). By~\cite[Thm.~3.2.12]{muirhead2009aspects} for $\ell-k\geq 2$, we have $\E[\mtx{A}^{-1}]=\frac{1}{\ell-k-1}\mtx{C}^{-1}$, $\E[\Tr(\mtx{A}^{-1})]=\Tr(\mtx{C}^{-1})/(\ell-k-1)$, and conclude that 
\begin{equation} \label{prop_rand_2}
\mathbb{E}\left[\|\mtx{\Omega}_1^{\dagger}\|_{\textup{F}}^2\right]=\frac{1}{\gamma_k\lambda_1}\frac{k}{\ell-k-1},\quad \gamma_k \coloneqq \frac{k}{\lambda_1 \Tr(\mtx{C}^{-1})}. 
\end{equation} 

The quantity $\gamma_k$ can be viewed as measuring the quality of the covariance kernel $K$ for learning the HS operator $\F$ (see~\cref{th_tropp_random_svd_Frob}). First, $1\leq \gamma_k<\infty$ as $\mtx{C}$ is symmetric positive definite. Moreover, for $1\leq j\leq k$, the $j$th largest eigenvalue of $\mtx{C}$ is bounded by the $j$th largest eigenvalue of $K$ as $\mtx{C}$ is a principal submatrix of $\mtx{V}^*K\mtx{V}$~\cite[Sec.~III.5]{kato2013perturbation}. Therefore, the following inequality holds,
\begin{equation} \label{eq_lower_bound_gamma}
\frac{1}{k}\sum_{j=1}^k\frac{\lambda_1}{\lambda_j}\leq \frac{1}{\gamma_k}<\infty,
\end{equation}
and the harmonic mean of the first $k$ scaled eigenvalues of $K$ is a lower bound for $1/\gamma_k$. In the ideal situation, the eigenfunctions of $K$ are the right singular functions of $\F$, \emph{i.e.}, $\psi_n=v_n$, $\mtx{C}$ is a diagonal matrix with entries $\lambda_1,\ldots,\lambda_k$, and $\gamma_k=k/(\sum_{j=1}^k \lambda_1/\lambda_j)$ is as small as possible.

We now provide a useful upper bound on $\gamma_k$ in a more general setting.

\begin{lemma} \label{lem_bound_sigma}
Let $\mtx{V}_1$ be a $D_1\times k$ quasimatrix with orthonormal columns and assume that there exists $m\in\mathbb{N}$ such that the columns of $\mtx{V}_1$ are spanned by the first $k+m$ eigenvectors of the continuous positive definite kernel $K:D_1\times D_1\to \R$. Then 
\[
\frac{1}{\gamma_k} \leq \frac{1}{k}\sum_{j=m+1}^{k+m}\frac{\lambda_1}{\lambda_j},
\]
where $\lambda_1\geq \lambda_2\geq\cdots> 0$ are the eigenvalues of $K$. This bound is tight in the sense that the inequality can be attained as an equality. 
\end{lemma}
\begin{proof}
Let $\mtx{Q} = \left[v_1\,|\,\cdots\,|\,v_k\,|\,q_{k+1}\,|\cdots\,|\,q_{k+m}\right]$ be a quasimatrix with orthonormal columns whose columns form an orthonormal basis for ${\rm Span}(\psi_1,\ldots,\psi_{k+m})$.  Then, $\mtx{Q}$ is an invariant space of $K$ and $\mtx{C}$ is a principal submatrix of $\mtx{Q}^*K\mtx{Q}$, which has eigenvalues $\lambda_1\geq \cdots\geq \lambda_{k+m}$. By~\cite[Thm.~6.46]{kato2013perturbation} the $k$ eigenvalues of $\mtx{C}$, denoted by $\mu_1,\ldots,\mu_{k}$, are greater than the first $k+m$ eigenvalues of $K$: $\mu_{j}\geq \lambda_{m+j}$ for $1\leq j\leq k$, and the result follows as the trace of a matrix is the sum of its eigenvalues. 
\end{proof}

\subsection{Probabilistic error bounds}
As discussed in~\cref{sec:TwoWarnings}, we need to extend the probability bounds of the randomized SVD to allow for non-standard Gaussian random vectors. The following lemma is a generalization of~\cite[Thm.~A.7]{halko2011finding}.

\begin{lemma} \label{thm_propa_bound}
Let $k,\ell\geq 1$ such that $\ell-k\geq 4$ and $\mtx{\Omega}_1$ be a $k\times \ell$ random matrix with independent columns such that each column has mean $(0,\ldots,0)^\top$ and positive definite covariance $\mtx{C}$. For all $t\geq 1$, we have
\[\mathbb{P}\left\{\|\mtx{\Omega}_1^\dagger\|_{\textup{F}}^2>\frac{3\Tr(\mtx{C}^{-1})}{\ell-k+1}\cdot t^2\right\}\leq t^{-(\ell-k)}.\]
\end{lemma}
\begin{proof}
Since $\mtx{\Omega}_1\mtx{\Omega}_1^*\sim W_{k}(\ell,\mtx{C})$, the reciprocals of its diagonal elements follow a scaled chi-square distribution~\cite[Thm.~3.2.12]{muirhead2009aspects}, \emph{i.e.},
\[
\frac{\left((\mtx{\Omega}_1\mtx{\Omega}_1^*)^{-1}\right)_{jj}}{\left(\mtx{C}^{-1}\right)_{jj}} \sim X_j^{-1},\quad X_j\sim \chi_{\ell-k+1}^2,\quad 1\leq j\leq k.
\]
Let $Z=\|\mtx{\Omega}_1^\dagger\|_{\textup{F}}^2=\Tr[(\mtx{\Omega}_1\mtx{\Omega}_1^*)^{-1}]$ and $q=(\ell-k)/2$. Following the proof of~\cite[Thm.~A.7]{halko2011finding}, we have the inequality
\[
\mathbb{P}\left\{|Z|\geq \frac{3\Tr(\mtx{C}^{-1})}{\ell-k+1}\cdot t^2\right\}\leq \left[\frac{3\Tr(\mtx{C}^{-1})}{\ell-k+1}\cdot t^2\right]^{-q}\E\left[|Z|^q\right],\quad t\geq 1.
\]
Moreover, by the Minkowski inequality, we have
\[
\left(\E\left[|Z^q|\right]\right)^{1/q}=\left(\E\left[\left|\sum_{j=1}^k[\mtx{C}^{-1}]_{jj}X_j^{-1}\right|^q\right]\right)^{1/q}\!\!\leq 
\sum_{j=1}^k [\mtx{C}^{-1}]_{jj}\E\left[|X_j^{-1}|^q\right]^{1/q}\leq \frac{3\Tr(\mtx{C}^{-1})}{\ell-k+1},
\]
where the last inequality is from~\cite[Lem.~A.9]{halko2011finding}. The result follows from the argument in the proof of~\cite[Thm.~A.7]{halko2011finding}.
\end{proof}

Under the assumption of \cref{lem_bound_sigma}, we find that~\cref{thm_propa_bound} gives the following bound: 
\[
\mathbb{P}\left\{\|\mtx{\Omega}_1^\dagger\|_{\textup{F}}>t\cdot \sqrt{\frac{3}{\ell-k+1}\sum_{j=m+1}^{k+m}\lambda_j^{-1}}\right\}\leq t^{-(\ell-k)}.
\]
In particular, in the finite dimensional case when $\lambda_1=\cdots=\lambda_n=1$, we recover the probabilistic bound found in~\cite[Thm.~A.7]{halko2011finding}.

To obtain the probability statement found in \cref{eq_proba_Green_bound} we require control of the tail of the distribution of a Gaussian quasimatrix with non-standard covariance kernel (see~\cref{sec_proof_thm_1}). In the theory of the randomized SVD, one relies on the concentration of measure results~\cite[Prop.~10.3]{halko2011finding}. However, we need to employ a different strategy and instead directly bound the HS norm of $\mtx{\Omega}_2$. One difficulty is that the norm of this matrix must be controlled for large dimensions $n$, which leads to a weaker probability bound than~\cite{halko2011finding}. While it is possible to apply Markov's inequality to obtain deviation bounds, we highlight that \cref{prop_control_matrices} provides a Chernoff-type bound, \emph{i.e.}, exponential decay of the tail distribution of $\|\mtx{\Omega}_2\|_{\HS}$, which is crucial to approximate Green's functions (see \cref{sec_proba_fail_Green}).

\begin{lemma} \label{prop_control_matrices}
With the same notation as in~\cref{th_tropp_deter_svd}, let $\ell\geq k \geq 1$. For all $s\geq 1$ we have
\[\mathbb{P}\left\{\|\mtx{\Omega}_2\|_{\HS}^2>\ell s^2\Tr(K)\right\}\leq  \left[s e^{-(s^2-1)/2}\right]^{\ell}.\]
\end{lemma}

\begin{proof}
We first remark that
\begin{equation} \label{eq_distrib_omega}
\|\mtx{\Omega}_2\|_{\HS}^2 \leq \|\mtx{\Omega}\|_{\HS}^2 = \sum_{j=1}^\ell Z_j,\quad Z_j \coloneqq \|\omega_j\|_{L^2(D_1)}^2,
\end{equation}
where the $Z_j$ are i.i.d.~because $\omega_j\sim \mathcal{GP}(0,K)$ are i.i.d. For $1\leq j\leq \ell$, we have (c.f.~\cref{sec_GP}),
\[\omega_j = \sum_{m=1}^\infty c_m^{(j)}\sqrt{\lambda_m}\psi_m,\]
where $c_m^{(j)}\sim \mathcal{N}(0,1)$ are i.i.d.~for $m\geq 1$ and $1\leq j\leq \ell$. First, since the series in \cref{eq_distrib_omega} converges absolutely, we have
\[Z_j = \sum_{m=1}^\infty (c_m^{(j)})^2\lambda_m= \lim_{N\to\infty}\sum_{m=1}^N X_m, \quad X_m = (c_m^{(j)})^2\lambda_m,\]
where the $X_m$ are independent random variables and $X_m \sim \lambda_m\chi^2$ for $1\leq m\leq N$. Here, $\chi^2$ denotes the chi-squared distribution~\cite[Chapt.~4.3]{mood1950introduction}.

Let $N\geq 1$ and $0<\theta<1/(2\Tr(K))$, we can bound the moment generating function of $\sum_{m=1}^N X_m$ as
\begin{align*}
\E\left[e^{\theta \sum_{m=1}^N X_m}\right] &= \prod_{m=1}^N\E\left[e^{\theta X_m}\right] = \prod_{m=1}^N (1-2\theta\lambda_m)^{-1/2}\leq \left(1-2\theta \sum_{m=1}^N \lambda_m\right)^{-1/2}\\
&\leq \left(1-2\theta \Tr(K)\right)^{-1/2},
\end{align*}
because $X_m/\lambda_m$ are independent random variables that follow a chi-squared distribution. Using the monotone convergence theorem, we have
\[\E\left[e^{\theta Z_j}\right]\leq (1-2\theta \Tr(K))^{-1/2}.\]

Let $\tilde{s}\geq 0$ and $0<\theta<1/(2\Tr(K))$. By the Chernoff bound~\cite[Thm.~1]{chernoff1952measure}, we obtain
\begin{align*}
\mathbb{P} \left\{\|\mtx{\Omega}_2\|_{\HS}^2 > \ell (1+\tilde{s})\Tr(K)\right\} &\leq e^{-(1+\tilde{s})\Tr(K)\ell \theta}\E\left[e^{\theta Z_j}\right]^\ell \\
&=e^{-(1+\tilde{s})\Tr(K)\ell \theta}(1-2\theta \Tr(K))^{-\ell/2}.
\end{align*}
We can minimize this upper bound over $0<\theta<1/(2\Tr(K))$ by choosing $\theta = \tilde{s}/(2(1+\tilde{s})\Tr(K))$, which gives
\[\mathbb{P} \left\{\|\mtx{\Omega}_2\|_{\HS}^2 > \ell (1+\tilde{s})\Tr(K)\right\}\leq (1+\tilde{s})^{\ell/2} e^{-\ell \tilde{s}/2}.\]
Choosing $s=\sqrt{1+\tilde{s}}\geq 1$ concludes the proof.
\end{proof}

\cref{prop_control_matrices} can be refined further to take into account the interaction between the Hilbert--Schmidt operator $\F$ and the covariance kernel $K$ (see~\cref{prop_control_matrices_refined}).

\subsection{Randomized SVD algorithm for HS operators}\label{sec_proof_thm_1}

We first prove an intermediary result, which generalizes \cite[Prop.~10.1]{halko2011finding} to HS operators. Note that one may obtain sharper bounds using a suitably chosen covariance kernels that yields a lower approximation error (see~\cref{chap_svd}).

\begin{lemma} \label{prop_rand_1}
Let $\mtx{\Sigma}_2$, $\mtx{V}_2$, and $\mtx{\Omega}$ be defined as in \cref{th_tropp_deter_svd}, and $\mtx{T}$ be an $\ell\times k$ matrix, where $\ell\geq k\geq 1$. Then,
\[
\E\left[\|\mtx{\Sigma}_2 \mtx{V}_2^* \mtx{\Omega} \mtx{T}\|_{\HS}^2\right] \leq \lambda_1\|\mtx{\Sigma}_2\|_{\HS}^2\|\mtx{T}\|_{\textup{F}}^2,
\]
where $\lambda_1$ is the first eigenvalue of $K$.
\end{lemma}
\begin{proof}
Let $\mtx{T} = \mtx{U}_\mtx{T} \mtx{D}_\mtx{T} \mtx{V}_\mtx{T}^*$ be the SVD of $\mtx{T}$. If $\{v_{\mtx{T},i}\}_{i=1}^k$ are the columns of $\mtx{V}_\mtx{T}$, then 
\[
\E\left[\|\mtx{\Sigma}_2 \mtx{V}_2^* \mtx{\Omega}\mtx{T}\|_{\HS}^2\right]= \sum_{i=1}^k\E\left[\|\mtx{\Sigma}_2 \mtx{\Omega}_2\mtx{U}_{\mtx{T}} \mtx{D}_\mtx{T} \mtx{V}_\mtx{T}^* v_{\mtx{T},i}\|_2^2\right],
\]
where $\mtx{\Omega}_2=\mtx{V}_2^*\mtx{\Omega}$. Therefore, we have
\[
\E\left[\|\mtx{\Sigma}_2 \mtx{\Omega}_2\mtx{T}\|_{\HS}^2\right] = \sum_{i=1}^k ((\mtx{D}_\mtx{T})_{ii})^2 \E\left[\|\mtx{\Sigma}_2 \mtx{\Omega}_2\mtx{U}_{\mtx{T}}(:,i)\|_{2}^2\right].
\]
Moreover, using the monotone convergence theorem for non-negative random variables, we have 
\begin{align*}
\E\left[\|\mtx{\Sigma}_2 \mtx{\Omega}_2\mtx{U}_{\mtx{T}}(:,i)\|_{2}^2\right]
&= \E\left[\sum_{n=1}^\infty\sum_{j=1}^{\ell}\sigma_{k+n}^2 \left|\mtx{\Omega}_2(n,j)\right|^2\mtx{U}_{\mtx{T}}(j,i)^2\right]\\
&= \sum_{n=1}^\infty\sum_{j=1}^{\ell}\sigma_{k+n}^2 \mtx{U}_{\mtx{T}}(j,i)^2\E\left[\left|\mtx{\Omega}_2(n,j)\right|^2\right],
\end{align*}
where $\sigma_{k+1},\sigma_{k+2},\ldots$ are the diagonal elements of $\mtx{\Sigma}_2$. Then, the quasimatrix $\mtx{\Omega}_2$ has independent columns and, using \cref{prop_cov_matrix}, we have 
\[\E\left[|\mtx{\Omega}_2(n,j)|^2\right] = \int_{D_1\times D_1}v_{k+n}(x)K(x,y)v_{k+n}(y)\d x\d y,\] where $v_{k+n}$ is the $n$th column of $\mtx{V}_2$. Then, $\E\left[|\mtx{\Omega}_2(n,j)|^2\right]\leq \lambda_1$, as $\E\left[|\mtx{\Omega}_2(n,j)|^2\right]$ is written as a Rayleigh quotient. Finally, we have 
\[\E\left[\|\mtx{\Sigma}_2 \mtx{V}_2^* \mtx{\Omega} \mtx{T}\|_{\HS}^2\right]\leq \lambda_1\sum_{i=1}^k ((\mtx{D}_{\mtx{T}})_{ii})^2\sum_{j=1}^{\ell}\mtx{U}_{\mtx{T}}(j,i)^2\sum_{n=1}^\infty\sigma_{k+n}^2 = \lambda_1 \|\mtx{T}\|_{\textup{F}}^2\|\mtx{\Sigma}_2\|_{\HS}^2,\] 
by orthonormality of the columns on $\mtx{U}_{\mtx{T}}$.
\end{proof}

We are now ready to prove~\cref{th_tropp_random_svd_Frob}, which shows that the randomized SVD can be generalized to HS operators. 
\begin{proof}[Proof of Theorem~\ref{th_tropp_random_svd_Frob}]
Let $\mtx{\Omega}_1, \mtx{\Omega}_2$ be the quasimatrices defined in~\cref{th_tropp_deter_svd}. The $k\times (k+p)$ matrix $\mtx{\Omega}_1$ has full rank with probability one and by \cref{th_tropp_deter_svd}, we have
\begin{align*}
\E\left[\|(\mtx{I}-\mtx{P}_\mtx{Y})\F\|_{\HS}\right] &\leq \E\left[\left(\|\mtx{\Sigma}_2\|_{\HS}^2+ \|\mtx{\Sigma}_2\mtx{\Omega}_2\mtx{\Omega}_1^{\dagger}\|_{\HS}^2\right)^{1/2}\right]\leq \|\mtx{\Sigma}_2\|_{\HS} + \E \|\mtx{\Sigma}_2\mtx{\Omega}_2\mtx{\Omega}_1^{\dagger}\|_{\HS}\\
&\leq \|\mtx{\Sigma}_2\|_{\HS} + \E \left[\|\mtx{\Sigma}_2\mtx{\Omega}_2\|_{\HS}^2\right]^{1/2}\E \left[\|\mtx{\Omega}_1^\dagger\|_{\textup{F}}^2\right]^{1/2},
\end{align*}
where the last inequality follows from Cauchy--Schwarz inequality. Then, combining \cref{prop_rand_1} and \cref{prop_rand_2}, we have
\[\E \left[\|\mtx{\Sigma}_2\mtx{\Omega}_2\|_{\HS}^2\right]\leq \lambda_1(k+p)\|\mtx{\Sigma}_2\|_{\HS}^2\quad\text{and} \quad \E\left[\|\mtx{\Omega}_1\|^2_{\textup{F}}\right]\leq \frac{1}{\gamma_k\lambda_1}\frac{k}{p-1},\]
where $\gamma_k$ is defined in \cref{sec_quality_kernel}.
The observation that $\|\mtx{\Sigma}_2\|_{\HS}^2=\sum_{j=k+1}^\infty \sigma_j^2$ concludes the proof of \cref{eq:MainExpectationBound}.

For the probabilistic bound in \cref{eq:MainProbabilityBound}, we note that by \cref{th_tropp_deter_svd} we have,
\[
\|\F-\mtx{P}_\mtx{Y}\F\|_{\HS}^2\leq \|\mtx{\Sigma}_2\|_{\HS}^2+\|\mtx{\Sigma}_2\mtx{\Omega}_2\mtx{\Omega}_1^{\dagger}\|_{\HS}^2\leq (1+\|\mtx{\Omega}_2\|_{\HS}^2\|\mtx{\Omega}_1^{\dagger}\|_{\textup{F}}^2)\|\mtx{\Sigma}_2\|_{\HS}^2,
\]
where the second inequality uses the submultiplicativity of the HS norm. The bound follows from bounding $\|\mtx{\Omega}_1^\dagger\|_{\textup{F}}^2$ and $\|\mtx{\Omega}_2\|_{\HS}^2$ using \cref{thm_propa_bound,prop_control_matrices}, respectively.
\end{proof}

\begin{remark}
The expectation bound~\eqref{eq:MainExpectationBound} in \cref{th_tropp_random_svd_Frob} does not control the square of the HS norm and therefore cannot be used to obtain an expectation bound for the randomized scheme for learning Green's functions described in \cref{sec_approx_Green}.
\end{remark}

The following proposition provides an expectation bound for the randomized SVD of the HS norm squared.

\begin{proposition}
With the notations of \cref{th_tropp_random_svd_Frob}, we have
\[\mathbb{E}\!\left[\|\F - \mtx{P}_\mtx{Y}\F\|_{\HS}^2\right]\leq \left(1+\frac{3\sqrt{2}}{\gamma_k}\,\frac{k(k+p)}{p+1}\sum_{j=1}^\infty\frac{\lambda_j}{\lambda_1}\right) \sum_{j=k+1}^{\infty}\sigma_j^2.\]
\end{proposition}

\begin{proof}
Let $\mtx{\Omega}_1, \mtx{\Omega}_2$ be the quasimatrices defined in~\cref{th_tropp_deter_svd}. We combine \cref{th_tropp_deter_svd} with the submultiplicativity of the HS norm and Cauchy--Schwarz inequality to obtain
\[\E[\|\F-\mtx{P}_\mtx{Y}\F\|_{\HS}^2]\!\leq \!\|\mtx{\Sigma}_2\|_{\HS}^2+\E[\|\mtx{\Sigma}_2\mtx{\Omega}_2\mtx{\Omega}_1^{\dagger}\|_{\HS}^2]\!\leq \!(1+\E[\|\mtx{\Omega}_2\|_{\HS}^4]^{\frac{1}{2}}\E[\|\mtx{\Omega}_1^{\dagger}\|_{\textup{F}}^4]^{\frac{1}{2}})\|\mtx{\Sigma}_2\|_{\HS}^2.\]
We can then control both terms $\E[\|\mtx{\Omega}_2\|_{\HS}^4]^{1/2}$ and $\E[\|\mtx{\Omega}_1^{\dagger}\|_{\textup{F}}^4]^{1/2}$ independently.

First, since $\mtx{\Omega}_1\mtx{\Omega}_1^{*}\sim W_{k}(k+p,\mtx{C})$, where $\mtx{C}$ is defined in \cref{prop_cov_matrix}, we have (cf.~the proof of \cref{thm_propa_bound})
\[
\frac{\left((\mtx{\Omega}_1\mtx{\Omega}_1^*)^{-1}\right)_{jj}}{\left(\mtx{C}^{-1}\right)_{jj}} \sim X_j^{-1},\quad X_j\sim \chi_{p+1}^2,\quad 1\leq j\leq k.
\]
Therefore,
\begin{align*}
\E[\|\mtx{\Omega}_1^\dagger\|^4]^{1/2} &= \E[\Tr((\mtx{\Omega}_1\mtx{\Omega}_1^*)^{-1})^2]^{1/2}= \E\left[\left(\sum_{j=1}^k (\mtx{C}^{-1})_{jj}X_j^{-1}\right)^2\right]^{1/2}\\
&\leq \Tr(\mtx{C}^{-1}) \E^2[X_1^{-1}],
\end{align*}
by the triangle inequality for the norm defined as $\E^2(Z) \coloneqq \E[|Z|^2]^{1/2}$ (see~\cite[Sec.~A.3.1]{halko2011finding}). Finally, using~\cite[Lem.~A.9]{halko2011finding}, we have $\E^2[X_1^{-1}] = 3/(p+1)$, which gives
\[\E[\|\mtx{\Omega}_1^\dagger\|^4]^{1/2} \leq \frac{3}{p+1}\Tr(\mtx{C}^{-1}).\]
The second term can be bounded as
\[\|\mtx{\Omega}_2\|_{\HS}^2 \leq \|\mtx{\Omega}\|_{\HS}^2 = \sum_{j=1}^{k+p} Z_j,\quad Z_j = \|\omega_j\|_{L^2(D_1)}^2,\]
where the $Z_j$ are i.i.d.. Therefore, by the triangle inequality for the $\E^2$-norm applied to the random variable $\|\mtx{\Omega}\|_{\HS}^2$, we have
\[\E(\|\mtx{\Omega}_2\|_{\HS}^4)^{1/2}\leq (k+p)\E[Z_1^2]^{1/2}.\]
We then characterize the distribution of $Z_1$ following the proof of \cref{prop_control_matrices} as
\[Z_1 = \sum_{m=1}^\infty \lambda_m Y_m,\quad Y_m\sim \chi^2.\]
Applying triangle inequality to $\E^2(Z_1)$ yields
\[\E[Z_1^2]^{1/2}\leq \sum_{m=1}^\infty \lambda_m \E[Y_m^2]^{1/2}=\Tr(K)\E[Y_1^2]^{1/2}=\sqrt{2}\Tr(K),\]
which concludes the proof.
\end{proof}

\section{Recovering the Green's function from input-output pairs} \label{sec_approx_Green}

It is known that the Green's function associated with~\cref{eq:PDEsimple} always exists, is unique, and is a nonnegative function $G : D \times D \to \R^+ \cup \{\infty\}$ such that 
\[
u(x) =\int_D G(x,y)f(y)\d y, \quad  f \in \mathcal{C}_c^\infty(D).
\]
For each $y\in\Omega$ and any $r>0$, we have $G(\cdot, y) \in \mathcal{H}^1(D \setminus B_r(y))\cap \mathcal{W}_0^{1,1}(D)$~\cite{gruter1982green}. Here, $B_r(y) = \{z\in\mathbb{R}^3 : \|z - y\|_2 < r\}$, $\mathcal{W}^{1,1}(D)$ is the space of weakly differentiable functions in the $L^1$-sense, and $\mathcal{W}^{1,1}_0(D)$ is the closure of $\mathcal{C}_c^\infty(D)$ in $\mathcal{W}^{1,1}(D)$. Since the PDE in~\cref{eq:PDEsimple} is self-adjoint, we also know that for almost every $x,y\in D$, we have $G(x,y) = G(y,x)$~\cite{gruter1982green}. 

We now state \cref{th_Green}, which shows that if $N = \mathcal{O}(\epsilon^{-6}\log^4(1/\epsilon))$ and one has $N$ input-output pairs $\{(f_j,u_j)\}_{j=1}^N$ with algorithmically-selected $f_j$, then the Green's function associated with $\mathcal{L}$ in \cref{eq:PDEsimple} can be recovered to within an accuracy of $\mathcal{O}(\Gamma_\epsilon^{-1/2}\log^3(1/\epsilon)\epsilon)$ with high probability. Here, the quantity $0<\Gamma_\epsilon\leq 1$ measures the quality of the random input functions $\{f_j\}_{j=1}^N$ (see~\cref{sec_green_reconst_adm}).

\begin{theorem} \label{th_Green}
Let $0<\epsilon<1$, $D\subset \R^3$ be a bounded Lipschitz domain, and $\L$ given in~\cref{eq:PDEsimple}. If $G$ is the Green's function associated with $\L$, then there is a randomized algorithm that constructs an approximation $\tilde{G}$ of $G$ using $\mathcal{O}(\epsilon^{-6}\log^4(1/\epsilon))$ input-output pairs such that, as $\epsilon\rightarrow 0$, we have
\begin{equation} \label{eq_proba_Green_bound}
\|G-\tilde{G}\|_{L^2(D\times D)} = \mathcal{O} \left(\Gamma_\epsilon^{-1/2}\log^{3}(1/\epsilon)\epsilon\right)\|G\|_{L^2(D\times D)},
\end{equation}
with probability $\geq 1 - \mathcal{O}(\epsilon^{\log(1/\epsilon)-6})$. The term $\Gamma_\epsilon$ is defined by~\cref{eq_define_gamma_eps}.
\end{theorem}

For simplicity, we have not reported the dependence of the bound in \cref{eq_proba_Green_bound} with respect to the spectral condition number, $\kappa_C = \lambda_{\max}/\lambda_{\min}$\footnote{Here, $\lambda_{\max}$ is defined as $\sup_{x\in D}\lambda_{\max}(A(x))$ and $\lambda_{\min} = \inf_{x\in D}\lambda_{\min}(A(x))>0$.}, of the coefficient matrix $A(x)$ in~\cref{eq:PDEsimple}.

Our algorithm that leads to the proof of \cref{th_Green} relies on the extension of the randomized SVD to HS operators (see \cref{sec_random_SVD}) and a hierarchical partition of the domain of $G$ into ``well-separated" domains. The scheme described in this section is summarized by \cref{alg_green}.

\renewcommand{\algorithmicrequire}{\textbf{Input:}}
\renewcommand{\algorithmicensure}{\textbf{Output:}}
\begin{algorithm}
\caption{Approximation of the Green's function}\label{alg_green}
\begin{algorithmic}[1]
\Require Action of the integral operator with kernel $G$
\Ensure Approximation $\tilde{G}$ of $G$
\State Construct a hierarchical partition of the domain $D\times D$
\State Approximate the Green's function on the admissible domains with the randomized SVD
\State Neglect $G$ on the non-admissible domains using a decay bound for the Green's function near the diagonal
\end{algorithmic}
\end{algorithm}

\subsection{Recovering the Green's function on admissible domains} \label{sec_exist_Green}

Roughly speaking, as $\|x-y\|_2$ increases $G$ becomes smoother about $(x,y)$, which can be made precise using so-called admissible domains~\cite{ballani2016matrices,bebendorf2008hierarchical,hackbusch2015hierarchical}. For $X,Y\subset\R^3$, let $\diam X\coloneqq\sup_{x,y\in X}\|x-y\|_2$ be the diameter of $X$, and $\dist(X,Y)\coloneqq\inf_{x\in X,y\in Y}\|x-y\|_2$ be the shortest distance between $X$ and $Y$. Admissible domains are defined as follows.
\begin{definition}
For a fixed parameter $\rho>0$, we say that two bounded and non-empty domains $X,Y\subset \R^3$ are admissible if
\[\dist(X,Y)\geq \rho \max \{\diam X, \diam Y\}.\]
Otherwise, we say that $X\times Y$ is non-admissible.
\end{definition}
\noindent There exists a weaker definition of admissible domains, which only requires that $\dist(X,Y)\geq \rho \min \{\diam X, \diam Y\}$~\cite[p.~59]{hackbusch2015hierarchical}, but we do not consider it.

\subsubsection{Approximation theory on admissible domains}\label{sec:ApproxTheoryAdmissibleDomain}
It turns out that the Green's function associated with~\cref{eq:PDEsimple} has exponentially decaying singular values when restricted to admissible domains. Roughly speaking, if $X,Y\subset D$ are such that $X\times Y$ is an admissible domain, then $G$ is well-approximated by a function of the form~\cite{bebendorf2003existence}
\begin{equation} \label{eq_approx_G}
G_k(x,y) = \sum_{j=1}^k g_j(x)h_j(y), \quad (x,y)\in X\times Y,
\end{equation}
for some functions $g_1,\ldots,g_k\in L^2(X)$ and $h_1,\ldots,h_k\in L^2(Y)$. This is summarized in \cref{theo_bebendorf}, which is a corollary of~\cite[Thm.~2.8]{bebendorf2003existence}. 
\begin{theorem} \label{theo_bebendorf}
Let $G$ be the Green's function associated with~\cref{eq:PDEsimple} and $\rho>0$. Let $X,Y\subset D$ such that $\dist(X,Y)\geq \rho\max\{\diam X, \diam Y\}$.  Then, for any $0<\epsilon<1$, there exists $k\leq k_\epsilon\coloneqq\lceil c(\rho,\diam D,\kappa_C)\rceil \lceil\log (1/\epsilon)\rceil^{4}+\lceil\log (1/\epsilon)\rceil$ and an approximant, $G_k$, of $G$ in the form given in~\cref{eq_approx_G} such that
\[
\|G-G_k\|_{L^2(X\times Y)}\leq \epsilon \|G\|_{L^2(X\times \hat{Y})}, \quad \hat{Y}\coloneqq\{y\in D,\, \dist(y,Y)\leq \frac{\rho}{2}\diam Y\},
\]
where $\kappa_C = \lambda_{\max}/\lambda_{\min}$ is the spectral condition number of the coefficient matrix $A(x)$ in~\cref{eq:PDEsimple} and $c$ is a constant that only depends on $\rho$, $\diam D$, $\kappa_C$. 
\end{theorem}
\begin{proof} 
In~\cite[Thm.~2.8]{bebendorf2003existence}, it is shown that if $Y=\tilde{Y}\cap D$ and $\tilde{Y}$ is convex, then there exists $k\leq c_{\rho/2}^3\lceil\log (1/\epsilon)\rceil^{4}+\lceil\log (1/\epsilon)\rceil$ and an approximant, $G_k$, of $G$ such that 
\begin{equation} \label{sec_sep_approx}
\|G(x,\cdot)-G_k(x,\cdot)\|_{L^2(Y)}\leq \epsilon \|G(x,\cdot)\|_{L^2(\hat{Y})}, \quad x\in X,
\end{equation}
where $\hat{Y}\coloneqq\{y\in D,\, \dist(y,Y)\leq \frac{\rho}{2}\diam Y\}$ and $c_{\rho/2}$ is a constant that only depends on $\rho$, $\diam Y$, and $\kappa_C$. As remarked by~\cite{bebendorf2003existence}, $\tilde{Y}$ can be included in a convex of diameter $\diam D$ that includes $D$ to obtain the constant $c(\rho,\diam D,\kappa_C)$. The statement follows by integrating the error bound in~\cref{sec_sep_approx} over $X$.
\end{proof} 

Since the truncated SVD of $G$ on $X\times Y$ gives the best rank $k_\epsilon\geq k$ approximation to $G$,~\cref{theo_bebendorf} also gives bounds on singular values:
\begin{equation} \label{eq_tail_SVD_error}
\left(\sum\nolimits_{j=k_\epsilon+1}^\infty\sigma_{j,X\times Y}^2\right)^{1/2} \leq \|G-G_{k}\|_{L^2(X\times Y)} \leq \epsilon \|G\|_{L^2(X\times \hat{Y})},
\end{equation}
where $\sigma_{j,X\times Y}$ is the $j$th singular value of $G$ restricted to $X\times Y$. Since $k_\epsilon = \mathcal{O}(\log^4(1/\epsilon))$, we conclude that the singular values of $G$ restricted to admissible domains $X\times Y$ rapidly decay to zero. 

\subsubsection{Randomized SVD for admissible domains} \label{sec_adm_domain}

Since $G$ has exponentially decaying singular values on admissible domains $X\times Y$, we use the randomized SVD for HS operators to learn $G$ on $X\times Y$ with high probability (see~\cref{sec_random_SVD}).

We start by defining a GP on the domain $Y$. Let $\mathcal{R}_{Y\times Y}K$ be the restriction\footnote{We denote the restriction operator by $\mathcal{R}_{Y\times Y}:L^2(D\times D)\to L^2(Y\times Y)$.} of the covariance kernel $K$ to the domain $Y\times Y$, which is a continuous symmetric positive definite kernel so that $\mathcal{GP}(0,\mathcal{R}_{Y\times Y}K)$ defines a GP on $Y$. We choose a target rank $k\geq 1$, an oversampling parameter $p\geq 2$, and form a quasimatrix $\mtx{\Omega} = \begin{bmatrix}f_1\,|\,\cdots \,|\, f_{k+p}\end{bmatrix}$ such that $f_j\in L^2(Y)$ and $f_j\sim \mathcal{GP}(0,\mathcal{R}_{Y\times Y} K)$ are identically distributed and independent. We then extend by zero each column of $\mtx{\Omega}$ from $L^2(Y)$ to $L^2(D)$ by $\mathcal{R}_Y^*\mtx{\Omega}=\begin{bmatrix}\mathcal{R}_Y^* f_1\,|\,\cdots \,|\, \mathcal{R}_Y^* f_{k+p}\end{bmatrix}$, where $\mathcal{R}_Y^* f_j\sim \mathcal{GP}(0,\mathcal{R}_{Y\times Y}^*\mathcal{R}_{Y\times Y}K)$. The zero extension operator $\mathcal{R}_Y^*:L^2(Y)\to L^2(D)$ is the adjoint of $\mathcal{R}_Y:L^2(D)\to L^2(Y)$.

Given the training data, $\mtx{Y} = \begin{bmatrix}u_1\,|\,\cdots \,|\, u_{k+p} \end{bmatrix}$ such that $\mathcal{L}u_j = \mathcal{R}_Y^* f_j$ and $u_j|_{\partial D} = 0$, we now construct an approximation to $G$ on $X\times Y$ using the randomized SVD (see~\cref{sec_random_SVD}). Following \cref{th_tropp_random_svd_Frob}, we have the following approximation error for $t\geq 1$ and $s\geq 2$:
\begin{equation} \label{eq_approx_error_loc_Green}
\|G-\tilde{G}_{X\times Y}\|_{L^2(X\times Y)}^2 \leq \left(1+t^2s^2\frac{3}{\gamma_{k,X\times Y}}\frac{k(k+p)}{p+1}\sum_{j=1}^\infty\frac{\lambda_j}{\lambda_1}\,\right)\left(\sum\nolimits_{j=k+1}^\infty\sigma_{j,X\times Y}^2\right)^{1/2},
\end{equation}
with probability greater than $1-t^{-p}-e^{-s^2(k+p)}$. Here, $\lambda_1\geq \lambda_2\geq\cdots>0$ are the eigenvalues of $K$, $\tilde{G}_{X\times Y} = \mtx{P}_{\mathcal{R}_{X}\mtx{Y}}\mathcal{R}_{X}\F\mathcal{R}_{Y}^*$ and $\mtx{P}_{\mathcal{R}_{X}\mtx{Y}} = \mathcal{R}_{X}\mtx{Y}((\mathcal{R}_{X}\mtx{Y})^*\mathcal{R}_{X}\mtx{Y})^{\dagger}(\mathcal{R}_{X}\mtx{Y})^*$ is the orthogonal projection onto the space spanned by the columns of $\mathcal{R}_{X}\mtx{Y}$. Moreover, $\gamma_{k,X\times Y}$ is a measure of the quality of the covariance kernel of $\mathcal{GP}(0,\mathcal{R}_{Y\times Y}^*\mathcal{R}_{Y\times Y}K)$ (see \cref{sec_quality_kernel}) and, for $1\leq i,j\leq k$, defined as $
\gamma_{k,X\times Y} = k/(\lambda_1\Tr(\mtx{C}_{X\times Y}^{-1}))$, where
\[[\mtx{C}_{X\times Y}]_{ij} = \int_{D\times D} \mathcal{R}_Y^*v_{i,X\times Y}(x)K(x,y) \mathcal{R}_Y^*v_{j,X\times Y}(y)\d x\d y,\]
and $v_{1,X\times Y},\ldots, v_{k,X\times Y}\in L^2(Y)$ are the first $k$ right singular functions of $G$ restricted to $X\times Y$. 

Unfortunately, there is a big problem with the formula $\tilde{G}_{X\times Y} = \mtx{P}_{\mathcal{R}_X\mtx{Y}}\mathcal{R}_{X}\F\mathcal{R}_{Y}^*$. It cannot be formed because we only have access to input-output data, so we have no mechanism for composing $\mtx{P}_{\mathcal{R}_X\mtx{Y}}$ on the left of $\mathcal{R}_{X}\F\mathcal{R}_{Y}^*$. Instead, we note that since the partial differential operator in~\cref{eq:PDEsimple} is self-adjoint, $\F$ is self-adjoint, and $G$ is itself symmetric. That means we can use this to write down a formula for $\tilde{G}_{Y\times X}$ instead. That is, 
\[
\tilde{G}_{Y\times X} = \tilde{G}_{X\times Y}^* = \mathcal{R}_{Y}\F\mathcal{R}_{X}^*\mtx{P}_{\mathcal{R}_X\mtx{Y}}, 
\]
where we used the fact that $\mtx{P}_{\mathcal{R}_X\mtx{Y}}$ is also self-adjoint. This means we can construct $\tilde{G}_{Y\times X}$ by asking for more input-output data to assess the quasimatrix $\F(\mathcal{R}_X^*\mathcal{R}_X\mtx{Y})$. Of course, to compute $\tilde{G}_{X\times Y}$, we can swap the roles of $X$ and $Y$ in the above argument. 

With a target rank of $k=k_\epsilon = \lceil c(\rho,\diam D,\kappa_C)\rceil\lceil\log (1/\epsilon)\rceil^{4}+\lceil\log (1/\epsilon)\rceil$ and an oversampling parameter of $p = k_\epsilon$, we can combine~\cref{theo_bebendorf} and \cref{eq_tail_SVD_error,eq_approx_error_loc_Green} to obtain the bound 
\[\|G-\tilde{G}_{X\times Y}\|_{L^2(X\times Y)}^2 \leq \left(1+t^2s^2\frac{6k_\epsilon}{\gamma_{k_{\epsilon},X\times Y}}\sum_{j=1}^\infty\frac{\lambda_j}{\lambda_1}\,\right)\epsilon^2\|G\|_{L^2(X\times \hat{Y})}^2 ,\]
with probability greater than $1-t^{-k_{\epsilon}}-e^{-2s^2 k_\epsilon}$. A similar approximation error holds for $\tilde{G}_{Y\times X}$ without additional evaluations of $\F$. We conclude that our algorithm requires $N_{\epsilon, X\times Y} \!= 2(k_\epsilon+p) =  \mathcal{O}\!\left(\log^4(1/\epsilon)\right)$ input-output pairs to learn an approximant to $G$ on $X\times Y$ and $Y\times X$.

\subsection{Ignoring the Green's function on non-admissible domains} \label{sec_lp_estimates}
When the Green's function is restricted to non-admissible domains, its singular values may not decay. Instead, to learn $G$ we take advantage of the off-diagonal decay property of $G$.
It is known that for almost every $x\neq y\in D$ then
\begin{equation} \label{eq_widman_green}
G(x,y)\leq \frac{c_{\kappa_C}}{\|x-y\|_2}\|G\|_{L^2(D\times D)},
\end{equation}
where $c_{\kappa_C}$ is an implicit constant that only depends on $\kappa_C$ (see~\cite[Thm.~1.1]{gruter1982green}). Note that we have normalized~\cite[Eq.~1.8]{gruter1982green} to highlight the dependence on $\|G\|_{L^2(D\times D)}$.

If $X\times Y$ is a non-admissible domain, then for any $(x,y)\in X\times Y$, we find that 
\[
\|x-y\|_2\leq \dist(X,Y)+\diam(X)+\diam(Y)< (2+\rho)\max\{\diam X,\diam Y\},
\]
because $\dist(X,Y)<\rho\max\{\diam X,\diam Y\}$. This means that $x\in B_r(y)\cap D$, where $r=(2+\rho)\max\{\diam X,\diam Y\}$. Using~\cref{eq_widman_green}, we have
\begin{align*}
\int_X G(x,y)^2 dx &\leq \int_{B_r(y)\cap D}G(x,y)^2\d x \leq c_{\kappa_C}^2\|G\|_{L^2(D\times D)}^2 \int_{B_r(y)}\|x-y\|_2^{-2}\d x\\
&\leq 4\pi c_{\kappa_C}^2 r\|G\|_{L^2(D\times D)}^2.
\end{align*}
Noting that $\diam(Y)\leq r/(2+\rho)$ and $\int_Y 1 \d y\leq 4\pi ({\rm diam}(Y)/2)^3/3$, we have the following inequality for non-admissible domains $X\times Y$:
\begin{equation} \label{eq_norm_green_non_adm}
\|G\|_{L^2(X\times Y)}^2\leq \frac{2\pi^2}{3(2+\rho)^3} c_{\kappa_C}^2 r^4 \|G\|_{L^2(D\times D)}^2,
\end{equation}
where $r=(2+\rho)\max\{\diam X,\diam Y\}$.
We conclude that the Green's function restricted to a non-admissible domain has a relatively small norm when the domain itself is small. Therefore, in our approximant $\tilde{G}$ for $G$, we ignore $G$ on non-admissible domains by setting $\tilde{G}$ to be zero.

\subsection{Hierarchical admissible partition of domain}\label{sec_hierar_domain}
We now describe a hierarchical partitioning of $D\times D$ so that many subdomains are admissible domains, and the non-admissible domains are all small. For ease of notation, we may assume---without loss of generality---that $\diam D = 1$ and $D\subset [0,1]^3$; otherwise, one should shift and scale $D$. Moreover, partitioning $[0,1]^3$ and restricting the partition to $D$ is easier than partitioning $D$ directly. For the definition of admissible domains, we find it convenient to select $\rho = 1/\sqrt{3}$. 

\begin{figure}[htbp]
\centering
\vspace{0.2cm}
\begin{overpic}[width=\textwidth]{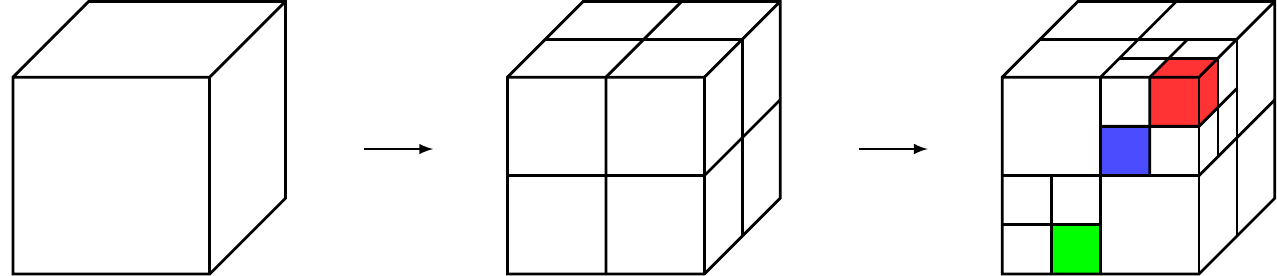}
\end{overpic}
\caption{Two levels of hierarchical partitioning of $[0,1]^3$. The blue and green domains are admissible, while the blue and red domains are non-admissible.}
\label{fig_hierar_tree}
\end{figure}

Let $I = [0,1]^3$. The hierarchical partitioning for $n$ levels is defined recursively as:
\begin{itemize}[leftmargin=*]
\item $I_{1\times 1\times 1}\coloneqq I_1\times I_1\times I_1=[0,1]^3$ is the root for level $L = 0$.
\item At a given level $0\leq L\leq n-1$, if $I_{j_1\times j_2\times j_3}\coloneqq I_{j_1}\times I_{j_2}\times I_{j_3}$ is a node of the tree, then it has $8$ children defined as
\[
\{I_{2 j_1+n_j(1)}\times I_{2 j_2+n_j(2)} \times I_{2j_3+n_j(3)}\mid n_j\in\{0,1\}^3\}.
\]
Here, if $I_j=[a,b]$, $0\leq a<b\leq 1$, then $I_{2j} = \left[a,\frac{a+b}{2}\right]$ and $I_{2j+1} = \left[\frac{a+b}{2},b\right]$.
\end{itemize}

\begin{figure}[htbp]
\centering
\vspace{0.5cm}
\begin{overpic}[width=0.4\textwidth, trim=100 50 100 10, clip]{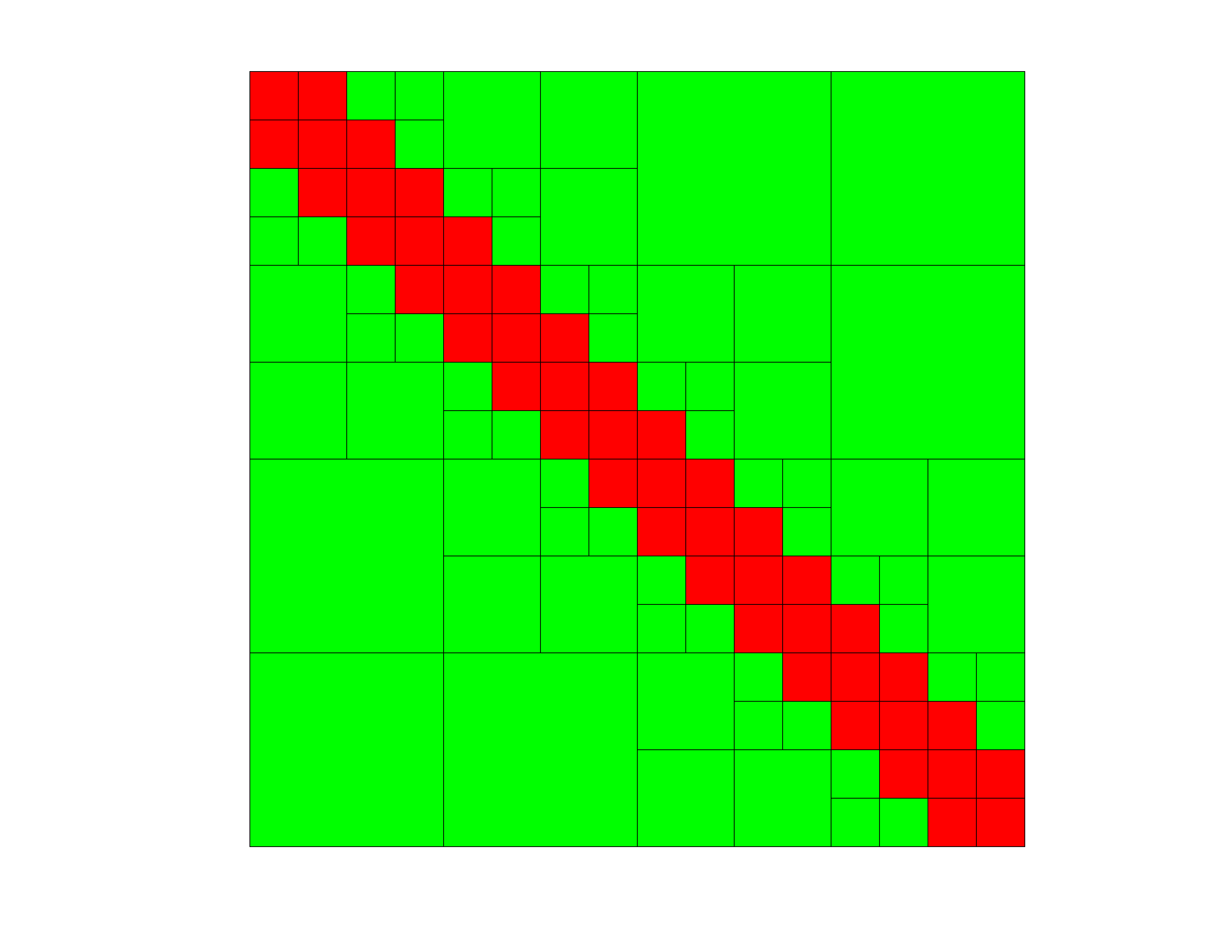}
\put(49,88){1D}
\end{overpic}
\hspace{1cm}
\begin{overpic}[width=0.4\textwidth, trim=100 50 100 10, clip]{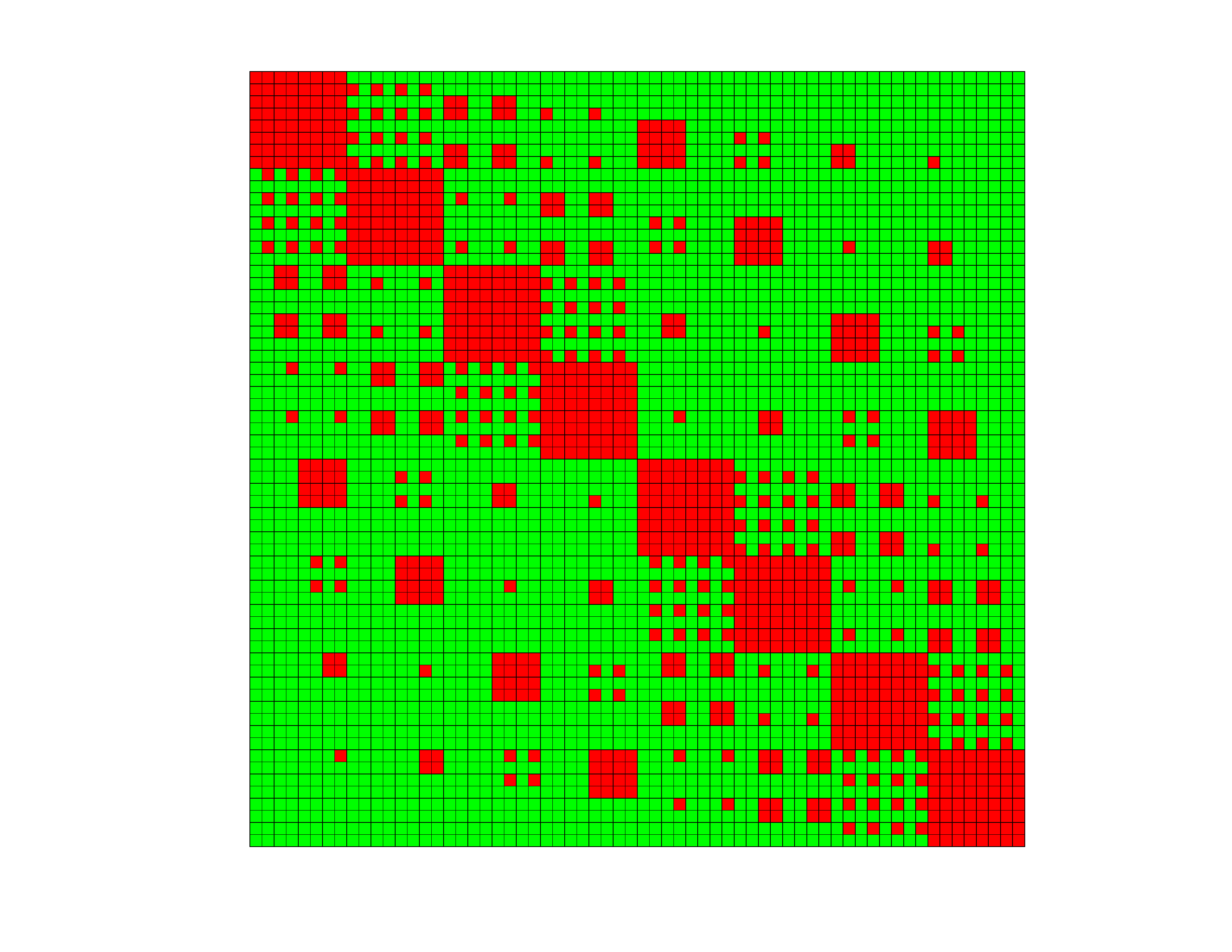}
\put(49,88){3D}
\end{overpic}
\caption{For illustration purposes, we include the hierarchical structure of the Green's functions in 1D after $4$ levels (left) and in 3D after $2$ levels (right). The hierarchical structure in 3D is complicated as this is physically a $6$-dimensional tensor that has been rearranged so it can be visualized.}
\label{fig_hierar_mat}
\end{figure}

The set of non-admissible domains can be given by an unwieldy expression 
\begin{equation} \label{eq_nonadm_d}
P_{\text{non-adm}} = \bigcup_{\substack{\bigwedge_{i=1}^3|j_i-\tilde{j}_i|\leq 1\\ 2^n\leq j_1,j_2,j_3\leq 2^{n+1}-1\\ 2^n\leq \tilde{j}_1,\tilde{j}_2,\tilde{j}_3\leq 2^{n+1}-1}} I_{j_1\times j_2\times j_3}\times I_{\tilde{j}_1\times \tilde{j}_2\times \tilde{j}_3},
\end{equation}
where $\land$ is the logical ``and'' operator. The set of admissible domains is given by 
\begin{equation} \label{eq_adm_d}
P_{\text{adm}} = \bigcup_{L=1}^n \Lambda(P_{\text{non-adm}}(L-1))\backslash P_{\text{non-adm}}(L)),
\end{equation}
where $P_{\text{non-adm}}(L)$ is the set of non-admissible domain for a hierarchical level of $L$ and
\[\Lambda(P_{\text{non-adm}}(L-1))=\bigcup_{\substack{I_{j_1\times j_2 \times j_3}\times I_{\tilde{j}_1\times\tilde{j}_2\times \tilde{j}_3}\\ \in P_{\text{non-adm}}(L-1)}}\, \bigcup_{n_j,n_{\tilde{j}}\in\{0,1\}^3}I_{\bigtimes_{i=1}^3 2j_i+n_j(i)}\times I_{\bigtimes_{i=1}^3 2 \tilde{j}_i+n_{\tilde{j}}(i)}.\]
Using \cref{eq_nonadm_d}-\cref{eq_adm_d}, the number of admissible and non-admissible domains are precisely $|P_{\text{non-adm}}| = (3\times 2^n-2)^3$ and $|P_{\text{adm}}| =  \sum_{\ell=1}^n 2^{6}(3\times 2^{L-1}-2)^3-(3\times 2^{L}-2)^3$.  In particular, the size of the partition at the hierarchical level $0\leq L\leq n$ is equal to $8^L$ and the tree has a total of $(8^{n+1}-1)/7$ nodes (see~\cref{fig_hierar_mat}).

Finally, the hierarchical partition of $D\times D$ can be defined via the partition $P=P_{\text{adm}}\cup P_{\text{non-adm}}$ of $[0,1]^3$ by doing the following: 
\[
D\times D =\bigcup\limits_{\tau\times\sigma\in P} (\tau\cap D)\times(\sigma\cap D).
\]
The sets of admissible and non-admissible domains of $D\times D$ are denoted by $P_{\text{adm}}$ and $P_{\text{non-adm}}$ in the next sections.

\subsection{Recovering the Green's function on the entire domain} \label{sec_proof_approx_green}
We now show that we can recover $G$ on the entire domain $D\times D$.

\subsubsection{Global approximation on the non-admissible set} \label{sec_approx_non_adm}
\begin{figure}[htbp]
\centering
\vspace{0.5cm}
\begin{overpic}[width=0.3\textwidth, trim=100 50 100 50, clip]{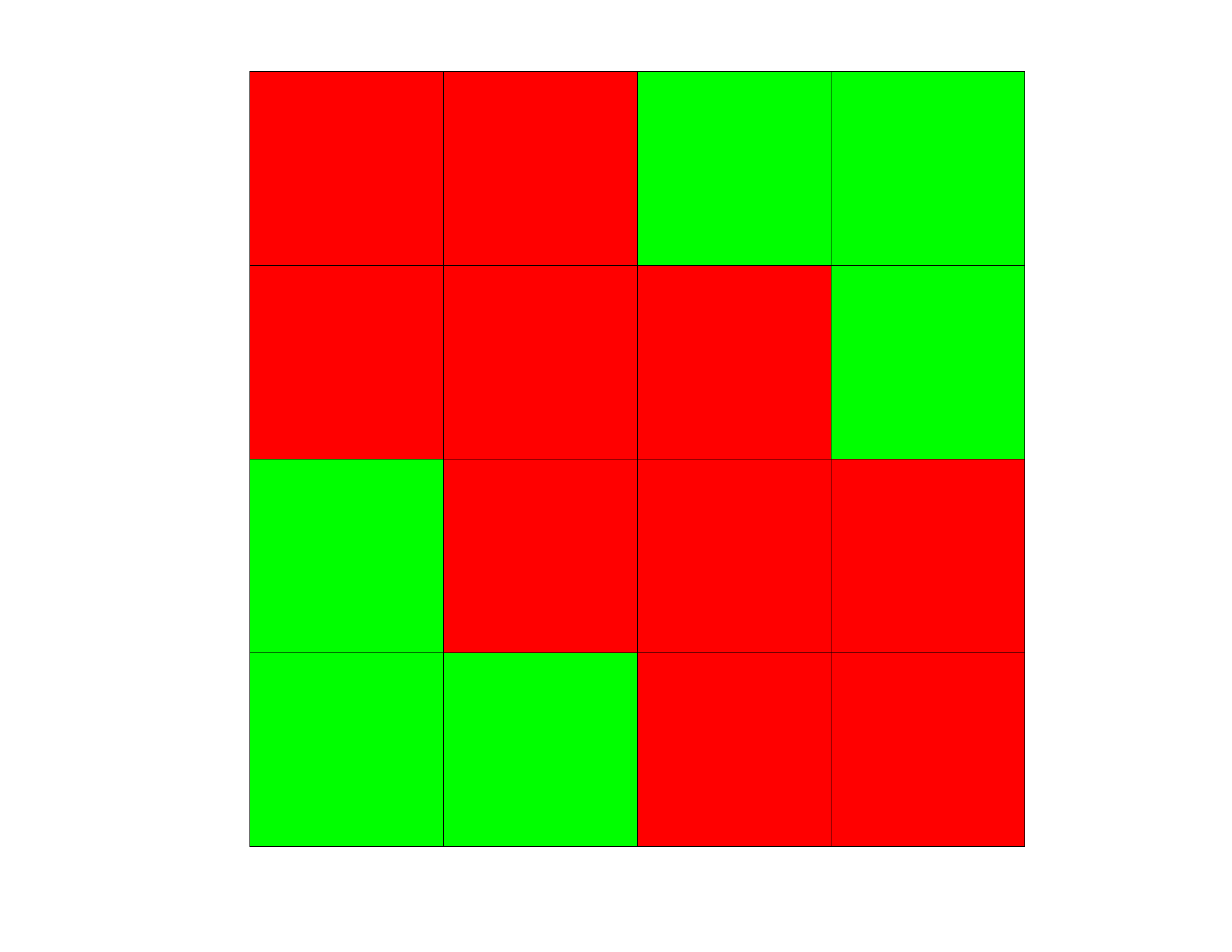}
\put(41,89){Level 2}
\end{overpic}
\begin{overpic}[width=0.3\textwidth, trim=100 50 100 50, clip]{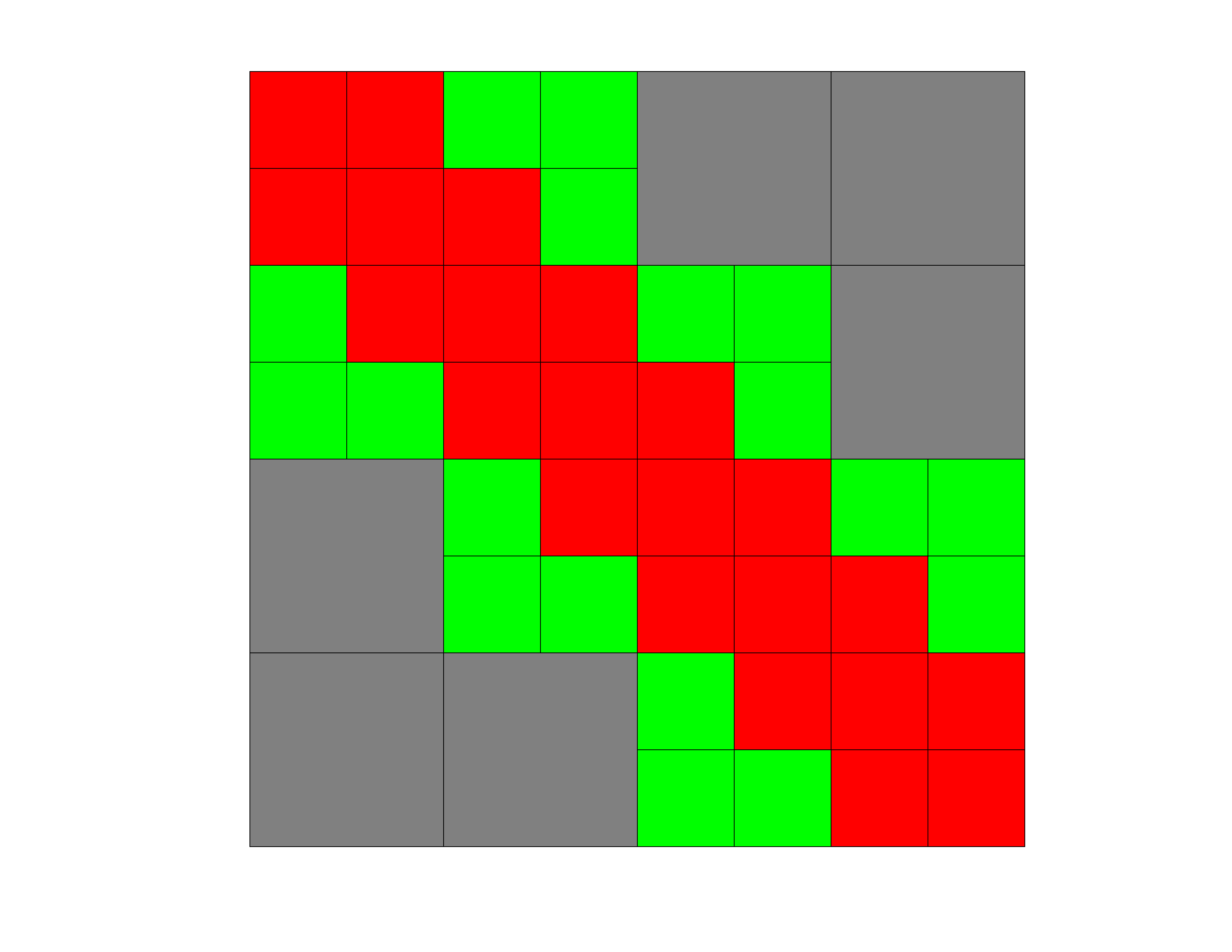}
\put(41,89){Level 3}
\end{overpic}
\begin{overpic}[width=0.3\textwidth, trim=100 50 100 50, clip]{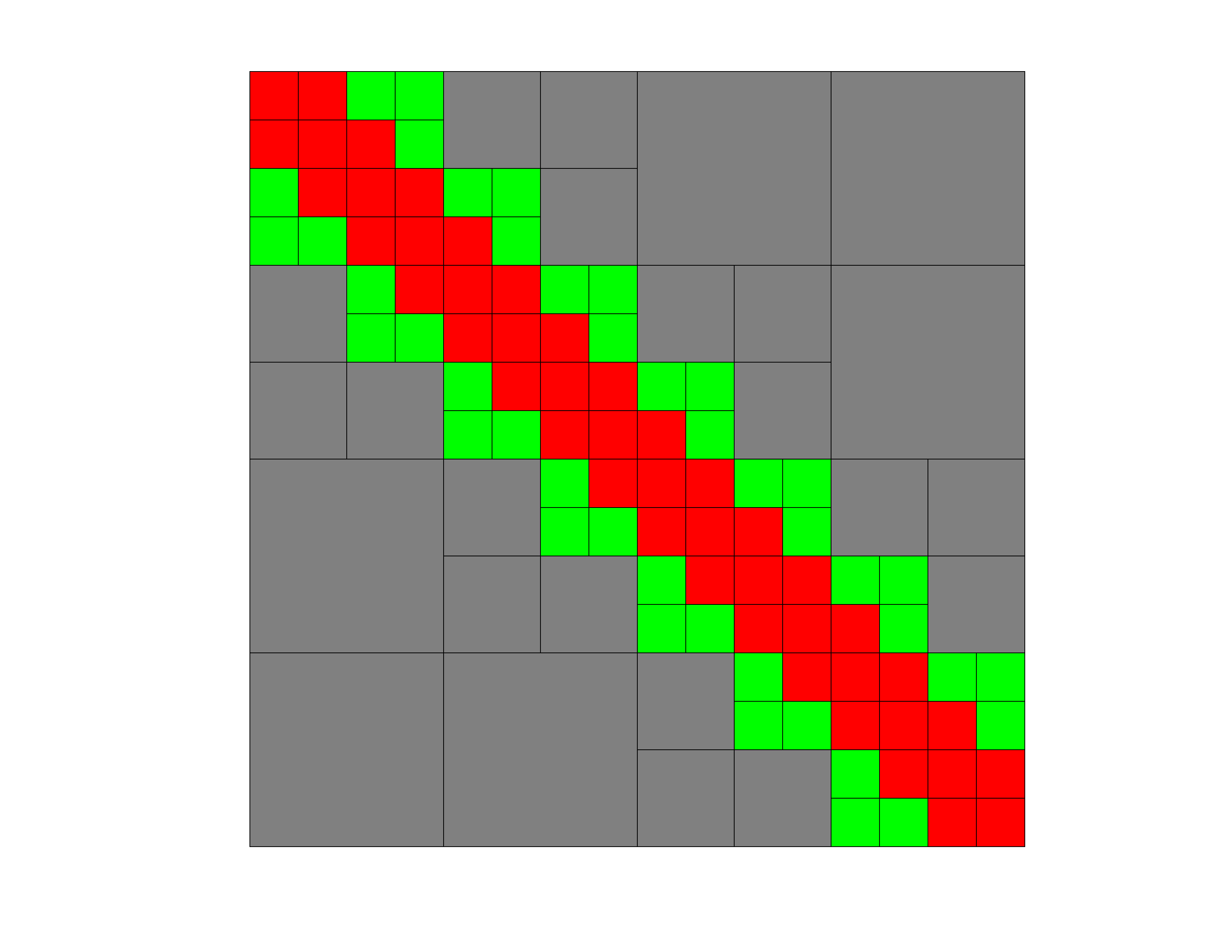}
\put(41,89){Level 4}
\end{overpic}
\caption{For illustration purposes, we include the hierarchical structure of the Green function in 1D. The green blocks are admissible domains at that level, the gray blocks are admissible at a higher level, and the red blocks are the non-admissible domains at that level. The area of the non-admissible domains decreases at deeper levels.}
\label{fig_hierar_mat_1D}
\end{figure}

Let $n_\epsilon$ be the number of levels in the hierarchical partition $D\times D$ (see~\cref{sec_hierar_domain}). We want to make sure that the norm of the Green's function on all non-admissible domains is small so that we can safely ignore that part of $G$ (see~\cref{sec_lp_estimates}). As one increases the hierarchical partitioning levels, the volume of the non-admissible domains get smaller (see~\cref{fig_hierar_mat_1D}). 

Let $X\times Y\in P_{\text{non-adm}}$ be a non-admissible domain, the two domains $X$ and $Y$ have diameter bounded by $\sqrt{3}/2^{n_\epsilon}$ because they are included in cubes of side length $1/2^{n_\epsilon}$ (see \cref{sec_hierar_domain}). Combining this with~\cref{eq_norm_green_non_adm} yields
\[
\|G\|_{L^2(X\times Y)}^2\leq 2\pi^2(6+\sqrt{3})c_{\kappa_C}^2 2^{-4n_\epsilon}\|G\|_{L^2(D\times D)}^2.\]
Therefore, the $L^2$-norm of $G$ on the non-admissible domain $P_{\text{non-adm}}$ satisfies
\[\|G\|_{L^2(P_{\text{non-adm}})}^2 = \sum_{X\times Y\in P_{\text{non-adm}}}\|G\|_{L^2(X\times Y)}^2\leq 54\pi^2(6+\sqrt{3})c_{\kappa_C}^2 2^{-n_\epsilon}\|G\|_{L^2(D\times D)}^2,\]
where we used $|P_{\text{non-adm}}| = (3\times 2^{n_\epsilon}-2)^3\leq 27(2^{3n_\epsilon})$. 
This means that if we select $n_\epsilon$ to be
\begin{equation} \label{eq_level_epsilon}
n_\epsilon = \left\lceil \log_2(54\pi^2(6+\sqrt{3})c_{\kappa_C}^2)+2\log_2(1/\epsilon)\right\rceil \sim 2\log_2(1/\epsilon),
\end{equation}
then we guarantee that $\|G\|_{L^2(P_{\text{non-adm}})}\leq\epsilon\|G\|_{L^2(D\times D)}$.  We can safely ignore $G$ on non-admissible domains---by taking the zero approximant---while approximating $G$ to within $\epsilon$.

\subsubsection{Learning rate of the Green's function}
 \label{sec_green_reconst_adm}

Following \cref{sec_adm_domain}, we can construct an approximant $\tilde{G}_{X\times Y}$ to the Green's function on an admissible domain $X\times Y$ of the hierarchical partitioning using the HS randomized SVD algorithm, which requires $N_{\epsilon,X\times Y}=\smash{\mathcal{O}(\log^4(1/\epsilon))}$ input-output training pairs (see \cref{sec_adm_domain}). Therefore, the number of training input-output pairs needed to construct an approximant to $G$ on all admissible domains is given by 
\[
N_\epsilon = \sum_{X\times Y\in P_{\text{adm}}} N_{\epsilon,X\times Y} = \mathcal{O}\left(|P_{\text{adm}}|\log^4(1/\epsilon)\right),\]
where $|P_{\text{adm}}|$ denotes the total number of admissible domains at the hierarchical level $n_\epsilon$, which is given by \cref{eq_level_epsilon}. Then, we have (see~\cref{sec_hierar_domain}):
\begin{equation} \label{eq_number_admissible}
|P_{\text{adm}}| = \sum_{\ell=1}^{n_\epsilon} 2^{6}(3\times 2^{\ell-1}-2)^3-(3\times 2^{\ell}-2)^3 \leq 6^3 2^{3n_\epsilon},
\end{equation}
and, using~\cref{eq_level_epsilon}, we obtain $|P_{\text{adm}}| = \mathcal{O}(1/\epsilon^6)$. This means that the total number of required input-output training pairs to learn $G$ with high probability is bounded by $N_\epsilon = \mathcal{O}\left(\epsilon^{-6}\log^4(1/\epsilon)\right)$.

\subsubsection{Global approximation error} \label{sec_proba_fail_Green}

We know that with $N_\epsilon = \mathcal{O}(\epsilon^{-6}\log^4(1/\epsilon))$ input-output training pairs, we can construct an accurate approximant to $G$ on each admissible and non-admissible domain. Since the number of admissible and non-admissible domains depends on $\epsilon$, we now check that this implies a globally accurate approximant that we denote by $\tilde{G}$. 

Since $\tilde{G}$ is zero on non-admissible domains and $P_{\text{adm}}\cap P_{\text{non-adm}}$ has measure zero, we have
\begin{equation} \label{eq_ineq_G}
 \|G-\tilde{G}\|_{L^2(D\times D)}^2
\leq \epsilon^2\|G\|_{L^2(D\times D)}^2+\sum_{X\times Y\in P_{\text{adm}}}\|G-\tilde{G}\|_{L^2(X\times Y)}^2.
\end{equation}
Following \cref{sec_green_reconst_adm}, if $X\times Y$ is admissible then the approximation error satisfies 
\[\|G-\tilde{G}_{X\times Y}\|_{L^2(X\times Y)}^2 \leq 12 t^2s^2\frac{k_\epsilon}{\gamma_{k_{\epsilon},X\times Y}}\sum_{j=1}^\infty\frac{\lambda_j}{\lambda_1}\epsilon^2\|G\|_{L^2(X\times \hat{Y})}^2 ,\]
with probability greater than $1-t^{-k_{\epsilon}}-e^{-2s^2 k_\epsilon}$. Here, $\hat{Y}=\{y\in D,\, \dist(y,Y)\leq \diam Y/2\sqrt{3}\}$ (see~\cref{theo_bebendorf} with $\rho=1/\sqrt{3}$). To measure the worst $\gamma_{k_\epsilon, X\times Y}$, we define
\begin{equation} \label{eq_define_gamma_eps}
\Gamma_\epsilon = \min\{\gamma_{k_\epsilon,X\times Y}:  X\times Y\in P_{\text{adm}}\}.
\end{equation}
From \cref{eq_lower_bound_gamma}, we know that $0<\Gamma_\epsilon\leq 1$ and that $1/\Gamma_\epsilon$ is greater than the harmonic mean of the first $k_\epsilon$ scaled eigenvalues of the covariance kernel $K$, \emph{i.e.},
\begin{equation} \label{eq_gamma_eps_bound}
\frac{1}{\Gamma_\epsilon} \geq \frac{1}{k_\epsilon}\sum_{j=1}^{k_\epsilon}\frac{\lambda_1}{\lambda_j},
\end{equation}

Now, one can see that $X\times \hat{Y}$ is included in at most $5^3=125$ neighbours including itself. 
Assuming that all the probability bounds hold on the admissible domains, this implies that
\begin{align*}
\sum_{X\times Y\in P_{\text{adm}}}\|G-\tilde{G}\|_{L^2(X\times Y)}^2 &\leq 
\sum_{X\times Y\in P_{\text{adm}}}\|G-\tilde{G}\|_{L^2(X\times Y)}^2\\
&\leq 12 t^2s^2\frac{k_\epsilon}{\lambda_1\Gamma_{\epsilon}}\Tr(K)\epsilon^2\sum_{X\times Y\in P_{\text{adm}}}\|G\|_{L^2(X\times\hat{Y})}^2\\
&\leq 1500 t^2s^2\frac{k_\epsilon}{\lambda_1\Gamma_{\epsilon}}\Tr(K)\epsilon^2\|G\|^2_{L^2(D\times D)}.
\end{align*}
We then choose $t=e$ and $s=k_\epsilon^{1/4}$ so that the approximation bound on each admissible domain holds with probability of failure less than $2e^{-\sqrt{k_\epsilon}}$. Finally, using \cref{eq_ineq_G} we conclude that as $\epsilon \to 0$, the approximation error on $D\times D$ satisfies
\[\|G-\tilde{G}\|_{L^2(D\times D)} = \mathcal{O}\left(\Gamma_\epsilon^{-1/2}\log^3(1/\epsilon)\epsilon\right)\|G\|_{L^2(D\times D)},\]
with probability $\geq (1-2e^{-\sqrt{k_\epsilon}})^{6^3 2^{3n_\epsilon}}=1-\mathcal{O}(\epsilon^{\log(1/\epsilon)-6})$, where $n_\epsilon$ is given by \cref{eq_level_epsilon}. We conclude that the approximant $\tilde{G}$ is a good approximation to $G$ with very high probability.

\section{Discussion} \label{sec_further_work}

There are several possible extensions of the results presented in this chapter related to the recovery of hierarchical matrices, the study of other partial differential operators, and practical deep learning applications, which we discuss further in this section.

\subsection{Fast and stable reconstruction of hierarchical matrices} \label{sec_stable_H_reconst}

We described an algorithm for reconstructing Green's function on admissible domains of a hierarchical partition of $D\times D$ that requires performing the HS randomized SVD $\mathcal{O}(\epsilon^{-6})$ times. We want to reduce it to a factor that is $\mathcal{O}(\text{polylog}(1/\epsilon))$. A polylogarithmic function in $x$ is any polynomial in $\log(x)$ and is denoted by $\polylog(x)$.

For $n\times n$ hierarchical matrices, there are several existing algorithms for recovering the matrix based on matrix-vector products~\cite{boukaram2019randomized,lin2011fast,martinsson2011fast,martinsson2016compressing}. There are two main approaches: (1) the ``bottom-up'' approach: one begins at the lowest level of the hierarchy and moves up and (2) the ``top-down'' approach: one updates the approximant by peeling off the off-diagonal blocks and going down the hierarchy. The bottom-up approach requires $\mathcal{O}(n)$ applications of the randomized SVD algorithm~\cite{martinsson2011fast}. There are lower complexity alternatives that only require $\mathcal{O}(\log(n))$ matrix-vector products with random vectors~\cite{lin2011fast}. However, the algorithm in~\cite{lin2011fast} is not yet proven to be theoretically stable as errors from low-rank approximations potentially accumulate exponentially, though this is not observed in practice. For symmetric positive semi-definite matrices, it may be possible to employ a sparse Cholesky factorization~\cite{schafer2021sparse,schafer2017compression}. This leads us to formulate the following challenge: 

\fbox{\begin{minipage}[t][1.2\height][c]{\dimexpr\textwidth-15\fboxsep-2\fboxrule\relax}
\centering
{\bf Algorithmic challenge:} Design a provably stable algorithm that can recover an $n\times n$ hierarchical matrix using $\mathcal{O}(\log(n))$ matrix-vector products with high probability. 
\end{minipage}}

\medskip

If one can design such an algorithm and it can be extended to HS operators, then the $\mathcal{O}(\epsilon^{-6}\log^4(1/\epsilon))$ term in~\cref{th_Green} may improve to~$\mathcal{O}(\text{polylog}(1/\epsilon))$. This means that the learning rate of partial differential operators of the form of \cref{eq:PDEsimple} will be a polynomial in $\log(1/\epsilon)$ and grow sublinearly with respect to $1/\epsilon$.

\subsection{Extension to other partial differential operators}
Our learning rate for elliptic partial differential operators (PDOs) in three variables (see~\cref{sec_approx_Green}) depends on the decay of the singular values of the Green's function on admissible domains~\cite{bebendorf2003existence}. We expect that one can also find the learning rate for other PDOs. 

It is known that the Green's functions associated to elliptic PDOs in two dimensions exist and satisfy the following pointwise estimate~\cite{dong2009green}:
\begin{equation} \label{eq_estimate_2d}
|G(x,y)|\leq C\left(\frac{1}{\gamma R^2}+\log\left(\frac{R}{\|x-y\|_2}\right)\right), \quad \|x-y\|_2\leq R\coloneqq\frac{1}{2}\max(d_x,d_y),
\end{equation}
where $d_x=\dist(x,\partial D)$, $\gamma$ is a constant depending on the size of the domain $D$, and $C$ is an implicit constant. One can conclude that $G(x,\cdot)$ is locally integrable for all $x\in D$ with $\|G(x,\cdot)\|_{L^p(B_r(x)\cap D)}<\infty$ for $r>0$ and $1\leq p<\infty$.  We believe that the pointwise estimate in~\cref{eq_estimate_2d} implies the off-diagonal low-rank structure of $G$ here, as suggested in~\cite{bebendorf2003existence}. Therefore, we expect that the results in this chapter can be extended to elliptic PDOs in two variables. It should also be possible to characterize the learning rate for elliptic PDOs with lower order terms (under reasonable conditions)~\cite{dong2020green,hwang2020green,kim2019green} as the associated Green's functions have similar regularity and pointwise estimates. The main task is to extend~\cite[Thm.~2.8]{bebendorf2003existence} to construct separable approximations of the Green's functions on admissible domains.

PDOs in four or more variables are far more challenging since we rely on the following bound on the Green's function on non-admissible domains~\cite{gruter1982green}:
\[
G(x,y)\leq \frac{c(d,\kappa_C)}{\lambda_{\min}}\|x-y\|_2^{2-d}, \quad x\neq y\in D,
\]
where $D\subset\R^d$, $d\geq 3$ is the dimension, and $c$ is a constant depending only on $d$ and $\kappa_C$. This inequality implies that the $L^p$-norm of $G$ on non-admissible domains is finite when $0\leq p < d/(d-2)$.  However, for a dimension $d\geq 4$, we have $p<2$ and one cannot ensure that the $L^2$ norm of $G$ is finite. Therefore, the Green's function may not be compatible with the HS randomized SVD.

The low-rank theory of Bebendorf and Hackbush has been recently extended from elliptic to parabolic operators~\cite{bebendorf2003existence} and combined with pointwise estimates for Green's functions~\cite{kim2020green} to obtain a learning rate parabolic PDEs, expressed in the $L^1$-norm~\cite{boulle2022parabolic}. In contrast, we believe that deriving a theoretical learning rate for hyperbolic PDOs remains a significant research challenge for many reasons. The first roadblock is that the Green's function associated with hyperbolic PDOs do not necessarily lie in $L^2(D\times D)$. For example, the Green's function associated with the wave equation in three variables, \emph{i.e.}, $\L=\partial_t^2-\nabla^2$, is not square-integrable as
\[
G(x,t,y,s) = \frac{\delta(t-s-\|x-y\|_2)}{4\pi \|x-y\|_2},\quad (x,t),(y,s)\in \R^3\times [0,\infty),
\]
where $\delta(\cdot)$ is the Dirac delta function.

Finally, while the extension to nonlinear dynamical systems seem out of reach of the technique presented in this chapter, characterizing the sample complexity of such systems and learn finite-dimensional approximations of the dynamics using Koopman operator theory~\cite{alexander2020operator,brunton2021modern,koopman1931hamiltonian} would be an interesting future research direction.

\subsection{Connection with neural networks}

As a concluding remark, we emphasize that the algorithm described in this chapter to learn Green's functions is not meant to be applied in practice. The proof of~\cref{th_Green} relies on the construction of a hierarchical partition of the domain $D\times D$ and the HS randomized SVD algorithm applied on each admissible domain. While this gives an algorithm for approximating Green's functions with high probability, it would be prohibitively computationally expensive to employ the hierarchical scheme and the generalization of the randomized SVD to large-scale three-dimensional problems. 

However, there are more practical approaches based on deep learning that currently do not yet have theoretical guarantees~\cite{feliu2020meta,gin2020deepgreen}. As we will see in \cref{chap_data_green}, deep learning techniques may be more competitive due to their ability to learn non self-adjoint problems and the fast optimization algorithms, based on stochastic gradient descent, for training neural networks. There are many possible connections between the work presented in this chapter and neural networks from practical and theoretical viewpoints.

A promising opportunity that we will explore in \cref{chap_data_green} is to design a NN that can learn and approximate Green's functions using input-output training pairs $\{(f_j,u_j)\}_{j=1}^N$. Once a neural network $\mathcal{N}$ has been trained such that $\|\mathcal{N}-G\|_{L^2}\leq \epsilon \|G\|_{L^2}$, the solution to $\L u = f$ can be obtained by computing the integral $u(x) = \int_D \mathcal{N}(x,y)f(y)\d y$. Therefore, this may give an efficient computational approach for discovering operators since a NN is only trained once. Incorporating a priori knowledge of the Green's function into the network architecture design could be particularly beneficial. As an example, one might exploit the low-rank structure by performing dimensionality reduction with an autoencoder~\cite{gin2020deepgreen,goodfellow2016deep}. One could also wrap the selection of the kernel in the GP for generating random functions and training data into a Bayesian framework. We expect that the theory developed in this chapter could guide deep learning experiments regarding the choice of training data and the type neural network architectures used to take advantage of the hierarchical structure of Green's functions and their singularity along the diagonal.

Finally, we wonder how many parameters in a NN are needed to approximate a Green's function associated with elliptic PDOs within a tolerance of $0<\epsilon<1$. Can one exploit the off-diagonal low-rank structure of Green's functions to reduce the number of parameters? We expect the recent work on the characterization of ReLU NNs' approximation power is useful~\cite{guhring2019error,petersen2018optimal,yarotsky2017error}. The use of NNs with high approximation power such as rational NNs might also be of interest to approximate the singularities of the Green's function near the diagonal, as we shall see in \cref{chapt_rational,chap_data_green}.

\dobib

 % Randomized SVD

\setcounter{chapter}{2}
\renewcommand{\thefootnote}{\fnsymbol{footnote}}
\chapter[A generalization of the randomized singular value decomposition]{A generalization of the randomized singular value decomposition\footnotemark} \label{chap_svd}

\footnotetext{This chapter is based on a paper with Alex Townsend~\cite{boulle2021generalization}, published in ICLR 2022. Townsend had an advisory role; I proved the theoretical results, performed the numerical experiments, and was the lead author in writing the paper.}
\renewcommand*{\thefootnote}{\arabic{footnote}}
\setcounter{footnote}{0}

The theory behind the randomized SVD has been extended in \cref{chapt_PDE_learning} to nonstandard covariance matrices and HS operators. However, the probability bounds, generalizing~\cite[Thm.~10.7]{halko2011finding}, are not sharp enough to emphasize the improved performance of covariance matrices with prior information over the standard randomized SVD. In this chapter, we improve the bounds obtained in \cref{chapt_PDE_learning} when the matrix-vector products are with multivariate Gaussian random vectors. Our theory allows for multivariate Gaussian random input vectors that have a general symmetric positive semi-definite covariance matrix. A key novelty of this work is that prior knowledge of the matrix $\mtx{A}$ can be exploited to design covariance matrices that achieve lower approximation errors than the randomized SVD with standard Gaussian vectors. We then design a practical algorithm for learning Hilbert--Schmidt (HS) operators using random input functions, sampled from a Gaussian process (GP). Examples of applications include learning integral kernels such as Green's functions associated with linear partial differential equations, as discussed in the previous chapter.

The choice of the covariance kernel in the GP is crucial and impacts both the theoretical bounds and numerical results of the randomized SVD. This leads us to introduce a new covariance kernel based on weighted Jacobi polynomials for learning HS operators. One of the main advantages of this kernel is that it is directly expressed as a Karhunen--Lo\`eve expansion~\cite{karhunen1946lineare,loeve1946functions} so that it is faster to sample functions from the associated GP than using a standard squared-exponential kernel. In addition, we show that the smoothness of the functions sampled from a GP with the Jacobi kernel can be controlled as it is related to the decay rate of the kernel's eigenvalues.

\section{Theoretical bounds for non-standard covariance matrices} \label{sec_random_svd}

In this section we provide new probability bounds for GPs with nonstandard covariance matrices. Let $m\geq n\geq 1$ and $\mtx{A}$ be an $m\times n$ real matrix with singular value decomposition $\mtx{A}=\mtx{U}\mtx{\Sigma}\mtx{V}^*$, where $\mtx{U}$ and $\mtx{V}$ are orthonormal matrices, and $\mtx{\Sigma}$ be an $m\times n$ diagonal matrix with entries $\sigma_1(\mtx{A})\geq \cdots\geq \sigma_n(\mtx{A})\geq 0$. For a fixed target rank $k\geq 1$, we define $\mtx{\Sigma}_1$ and $\mtx{\Sigma}_2$ to be the $k\times k$ and $(n-k)\times (n-k)$ diagonal matrices, which respectively contain the first $k$ singular values of $\mtx{A}$: $\sigma_1(\mtx{A})\geq \cdots \geq \sigma_k(\mtx{A})$, and the remaining singular values. Let $\mtx{V}_1$ be the $n\times k$ matrix obtained by truncating $\mtx{V}_1$ after $k$ columns and $\mtx{V}_2$ the remainder. In this section, $\mtx{K}$ denotes a symmetric positive semi-definite $n\times n$ matrix with $k$th largest eigenvalue $\lambda_k>0$ and $\mtx{\Omega}\in \R^{n\times \ell}$ a Gaussian random matrix with $\ell\geq k$ independent columns sampled from a multivariate normal distribution with covariance matrix $\mtx{K}$. Finally, we define $\mtx{\Omega}_1\coloneqq \mtx{V}_1^*\mtx{\Omega}$ and $\mtx{\Omega}_2\coloneqq \mtx{V}_2^*\mtx{\Omega}$. The following theorem is a refinement of \cref{th_tropp_random_svd_Frob}. While it is formulated with matrices, the same result holds for HS operators in infinite dimensions.

\begin{theorem} \label{th_svd_main}
Let $\mtx{A}$ be an $m\times n$ matrix, $k\geq 1$ an integer, and choose an oversampling parameter $p\geq 4$. If $\mtx{\Omega}\in\mathbb{R}^{n\times (k+p)}$ is a Gaussian random matrix, where each column is i.i.d. from a multivariate Gaussian distribution with covariance matrix $\mtx{K}\in\mathbb{R}^{n\times n}$, and $\mtx{Q}\mtx{R} = \mtx{A}\mtx{\Omega}$ is the economized QR decomposition of $\mtx{A}\mtx{\Omega}$, then for all $u,t\geq 1$,
\begin{equation} \label{eq_svd_ProbabilityBound}
\|\mtx{A}-\mtx{Q}\mtx{Q}^*\mtx{A}\|_{\textup{F}} \leq \left(1+ut\sqrt{(k+p)\frac{3k}{p+1}\frac{\beta_k}{\gamma_k}}\,\right)\sqrt{\sum_{j=k+1}^n\sigma_j^2(\mtx{A})},
\end{equation}
with failure probability at most $t^{-p}+[ue^{-(u^2-1)/2}]^{k+p}$. 
Here, the covariance quality factors are denoted by $\smash{\gamma_k \!=\! k/(\lambda_1 \Tr((\mtx{V}_1^*\mtx{K}\mtx{V}_1)^{-1})))}$ and $\smash{\beta_k \!=\! \Tr(\mtx{\Sigma}_2^2\mtx{V}_2^*\mtx{K}\mtx{V}_2)/(\lambda_1\|\mtx{\Sigma}_2\|_\Frob^2)}$, where $\lambda_1$ is the largest eigenvalue of $\mtx{K}$.
\end{theorem}

This result differs from \cref{th_tropp_random_svd_Frob} and~\cite[Thm.~10.5]{halko2011finding} due to the additional factors $\gamma_k$ and $\beta_k$, which measure the quality of the covariance matrix to learn $\mtx{A}$ in \cref{th_svd_main}. They can be respectively bounded (\cref{lem_bound_sigma,lemm_beta_bound}) using the eigenvalues $\lambda_1\geq \cdots\geq \lambda_n$ of the covariance matrix $\mtx{K}$ and the singular values of $\mtx{A}$ as:
\begin{equation} \label{eq_bound_gamma_beta}
\frac{1}{\gamma_k}\leq \frac{1}{k}\sum_{j=n-k+1}^{n}\frac{\lambda_1}{\lambda_j},\qquad \beta_k\leq \sum_{j=k+1}^{n}\frac{\lambda_{j-k}}{\lambda_1}\sigma_j^2(\mtx{A})\bigg/ \sum_{j=k+1}^n\sigma_j^2(\mtx{A}).
\end{equation}
This shows that the performance of the generalized randomized SVD depends on the decay rate of the sequence $\{\lambda_j\}$. The quantities $\gamma_k$ and $\beta_k$ depend on how much prior information of the $k+1,\ldots,n$ right singular vectors of $\mtx{A}$ is encoded in $\mtx{K}$. In the ideal situation where these singular vectors are known, then one can define $\mtx{K}$ such that $\beta_k=0$ for $\lambda_{k+1}=\cdots=\lambda_n=0$. Unlike the weaker but more explicit bound proven in \cref{sec_random_SVD}, this highlights that a suitably chosen covariance matrix can outperform the randomized SVD with standard Gaussian vectors (see \cref{sec_exp_matrix} for a numerical example).

The proof of \cref{th_svd_main} will require bounding $\|\mtx{\Omega}_1^\dagger\|_{\textup{F}}^2$, which we achieve using \cref{thm_propa_bound}, as well as the term $\|\mtx{\Sigma}_2\mtx{\Omega}_2\|_{\textup{F}}^2$, which is done in the following lemma.

\begin{lemma} \label{prop_control_matrices_refined}
With the notations introduced at the beginning of the section, for all $s\geq 0$, we have
\[\mathbb{P}\left\{\|\mtx{\Sigma}_2\mtx{\Omega}_2\|_{\Frob}^2>\ell (1+s) \Tr(\mtx{\Sigma}_2^2\mtx{V}_2^*\mtx{K}\mtx{V}_2)\right\}\leq  (1+s)^{\ell/2} e^{-s\ell/2}.\]
\end{lemma}

\begin{proof}
Let $\omega_j$ be the $j$th column of $\mtx{\Omega}$ for $1\leq j\leq \ell$ and $v_1,\ldots,v_n$ be the $n$ columns of the orthonormal matrix $\mtx{V}$. We first remark that 
\[\|\mtx{\Omega}_2\|_{\Frob}^2 = \sum_{j=1}^\ell Z_j,\quad Z_j \coloneqq \sum_{n_1=1}^{n-k}\sigma_{k+n_1}^2(\mtx{A}) (v_{k+n_1}^*\omega_j)^2,\]
where the $Z_j$ are i.i.d.~because $\omega_j\sim \mathcal{N}(0,\mtx{K})$ are i.i.d. Let $\lambda_1\geq \cdots\geq  \lambda_n\geq 0$ be the eigenvalues of $\mtx{K}$ with eigenvectors $\psi_1,\ldots,\psi_n\in\R^n$. For $1\leq j\leq \ell$, we have,
\[\omega_j = \sum_{i=1}^n (c_i^{(j)})^2\sqrt{\lambda_i}\psi_i,\]
where $c_i^{(j)}\sim \mathcal{N}(0,1)$ are i.i.d.~for $1\leq i\leq n$ and $1\leq j\leq \ell$. Then,
\[Z_j = \sum_{i=1}^n (c_i^{(j)})^2 \lambda_i\sum_{n_1=1}^{n-k}\sigma_{k+n_1}^2(\mtx{A}) (v_{k+n_1}^*\psi_i)^2
= \sum_{i=1}^n X_i\]
where the $X_i$ are independent. Let $\gamma_i = \lambda_i\sum_{n_1=1}^{n-k} \sigma_{k+n_1}^2(\mtx{A})(v_{k+n_1}^*\psi_i)^2$, then $X_i \sim \gamma_i\chi^2$ for $1\leq i\leq n$. 

Let $0<\theta<1/(2\sum_{i=1}^n\gamma_i)$. We can bound the moment generating function of $\sum_{i=1}^n X_i$ as
\[\E\left[e^{\theta \sum_{i=1}^n X_i}\right] = \prod_{i=1}^n\E\left[e^{\theta X_i}\right] = \prod_{i=1}^n (1-2\theta\gamma_i)^{-1/2}\leq \left(1-2\theta \sum_{i=1}^n \gamma_i\right)^{-1/2}\]
because the $X_i/\gamma_i$ are independent and follow a chi-squared distribution. The right inequality is obtained by showing by recurrence that, if $a_1,\ldots,a_n\geq 0$ are such that $\sum_{i=1}^n a_i\leq 1$, then $\prod_{i=1}^n (1-a_i)\geq 1-\sum_{i=1}^n a_i$. For convenience, we define $C_1 \coloneqq \sum_{i=1}^n \gamma_n$, we have shown that 
\[\E\left[e^{\theta Z_j}\right]\leq (1-2\theta C_1)^{-1/2}.\]
Moreover, we find that
\begin{align*}
C_1 &= \sum_{n_1=1}^{n-k} \sigma_{k+n_1}^2 v_{k+n_1}^*\left(\sum_{i=1}^n \psi_i^*\lambda_i\psi_i\right)v_{k+n_1}=\sum_{n_1=1}^{n-k} \sigma_{k+n_1}^2(\mtx{A}) v_{k+n_1}^*\mtx{K} v_{k+n_1}\\
&=\Tr(\mtx{\Sigma}_2^2\mtx{V}_2^*\mtx{K}\mtx{V}_2).
\end{align*}

Let $s\geq 0$ and $0<\theta<1/(2C_1)$. By the Chernoff bound~\cite[Thm.~1]{chernoff1952measure}, we obtain
\begin{align*}
\mathbb{P} \left\{\|\mtx{\Sigma}_2\mtx{\Omega}_2\|_{\Frob}^2 > \ell (1+s)C_1\right\} &\leq e^{-(1+s)C_1\ell \theta}\E\left[e^{\theta Z_j}\right]^\ell \\
&=e^{-(1+s)C_1\ell \theta}(1-2\theta C_1)^{-\ell/2}.
\end{align*}
We minimize the bound over $0<\theta<1/(2\Tr(K))$ by choosing $\theta = s/(2(1+s)C_1)$, which gives
\[\mathbb{P} \left\{\|\mtx{\Sigma}_2\mtx{\Omega}_2\|_{\Frob}^2 > \ell (1+s)C_1\right\}\leq (1+s)^{\ell/2} e^{-\ell s/2}.\]
\end{proof}

We now prove \cref{th_svd_main}, which provides a refined probability bound for the performance of the generalized randomized SVD on matrices.

\begin{proof}[Proof of \cref{th_svd_main}]
Using \cref{th_tropp_deter_svd} and the submultiplicativity of the Frobenius norm, we have
\begin{equation} \label{eq_ineq_frob}
\|\mtx{A}-\mtx{Q}\mtx{Q}^*\mtx{A}\|_{\Frob}^2\leq \|\mtx{\Sigma}_2\|_{\Frob}^2+\|\mtx{\Sigma}_2\mtx{\Omega}_2\|_{\Frob}^2\|\mtx{\Omega}_1^\dagger\|_{\Frob}^2.
\end{equation}
Let $\ell = k+p$ with $p\geq 4$. Combining \cref{thm_propa_bound,prop_control_matrices_refined} to bound the terms $\|\mtx{\Sigma}_2\mtx{\Omega}_2\|_{\Frob}^2$ and $\|\mtx{\Omega}_1^\dagger\|_{\Frob}^2$ in \cref{eq_ineq_frob} yields the following probability estimate:
\begin{align*}
\|\mtx{A}-\mtx{Q}\mtx{Q}^*\mtx{A}\|_{\textup{F}}^2 &\leq \|\mtx{\Sigma}_2\|_{\Frob}^2+3t^2(1+s)\frac{k+p}{p+1}\Tr((\mtx{V}_1^*\mtx{K}\mtx{V}_1)^{-1})\Tr(\mtx{\Sigma}_2^2\mtx{V}_2^*\mtx{K}\mtx{V}_2)\\
&\leq \left(1+3t^2(1+s)\frac{(k+p)k}{p+1}\frac{\beta_k}{\gamma_k}\right)\sum_{j=k+1}^n\sigma_j^2(\mtx{A}),
\end{align*}
with failure probability at most $t^{-p}+(1+s)^{(k+p)/2}e^{-s(k+p)/2}$. Note that we introduced $\gamma_k\coloneqq k/(\lambda_1 \Tr((\mtx{V}_1^*\mtx{K}\mtx{V}_1)^{-1})))$ and $\beta_k \coloneqq \Tr(\mtx{\Sigma}_2^2\mtx{V}_2^*\mtx{K}\mtx{V}_2)/(\lambda_1\|\mtx{\Sigma}_2\|_\Frob^2)$. We conclude the proof by defining $u=\sqrt{1+s}\geq 1$.
\end{proof}

The following Lemma provides an estimate of the quantity $\beta_k$ introduced in the statement of \cref{th_svd_main}.

\begin{lemma} \label{lemm_beta_bound}
Let $\beta_k = \Tr(\mtx{\Sigma}_2^2\mtx{V}_2^*\mtx{K}\mtx{V}_2)/(\lambda_1\|\mtx{\Sigma}_2\|_\Frob^2)$, then the following inequality holds
\[\beta_k\leq \sum_{j=k+1}^{n}\frac{\lambda_{j-k}}{\lambda_1}\sigma_j^2(\mtx{A})\bigg/ \sum_{j=k+1}^n\sigma_j^2(\mtx{A}).\]
\end{lemma}

\begin{proof}
Let $\mu_1\geq \cdots\geq \mu_{n-k}$ be the eigenvalues of the matrix $\mtx{V}_2^*\mtx{K}\mtx{V}_2$. Using von Neumann's trace inequality~\cite{mirsky1975trace,von1937some}, we have
\[\Tr(\mtx{\Sigma}_2^2\mtx{V}_2^*\mtx{K}\mtx{V}_2)\leq \sum_{j=k+1}^n\mu_{j-k}\sigma^2_j(\mtx{A}).\]
Then, the matrix $\mtx{V}_2^*\mtx{K}\mtx{V}_2$ is a principal submatrix of $\mtx{V}^*\mtx{K}\mtx{V}$, which has the same eigenvalues of $K$. Therefore, by~\cite[Thm.~6.46]{kato2013perturbation}, the eigenvalues of $\mtx{V}_2^*\mtx{K}\mtx{V}_2$ are individually bounded by the eigenvalues of $\mtx{K}$, \emph{i.e.}, $\mu_j \leq \lambda_{j}$ for $1\leq j\leq n-k$, which concludes the proof.
\end{proof}

Finally, we highlight that the statement of~\cref{th_svd_main} can be simplified by choosing $p=5$, $t = 4$, and $u=3$ to highlight the difference with the standard bounds for the randomized SVD.

\begin{corollary}[Generalized randomized SVD] \label{cor_gen_svd}
Let $\mtx{A}$ be an $m\times n$ matrix and $k\geq 1$ an integer. If $\mtx{\Omega}\in\mathbb{R}^{n\times (k+5)}$ is a Gaussian random matrix, where each column is i.i.d. from a multivariate Gaussian distribution with symmetric positive semi-definite covariance matrix $\mtx{K}\in\mathbb{R}^{n\times n}$, and $\mtx{Q}\mtx{R} = \mtx{A}\mtx{\Omega}$ is the economized QR decomposition of $\mtx{A}\mtx{\Omega}$, then 
\[\mathbb{P}\left[ \| \mtx{A} - \mtx{Q}\mtx{Q}^*\mtx{A} \|_\Frob \leq \left(1+ 9 \sqrt{k(k+5)\frac{\beta_k}{\gamma_k}} \right) \sqrt{\sum_{j=k+1}^n \sigma_j^2(\mtx{A})} \,\right] \geq 0.999.\]
\end{corollary}

\noindent In contrast, a simplification of the theorem for the randomized SVD~\cite[Thm.~10.7]{halko2011finding} by choosing $t=6$ and $u=4$ gives the following result.

\begin{corollary}[Randomized SVD]
Let $\mtx{A}$ be an $m\times n$ matrix and $k\geq 1$ an integer. If $\mtx{\Omega}\in\mathbb{R}^{n\times (k+5)}$ is a standard Gaussian random matrix and $\mtx{Q}\mtx{R} = \mtx{A}\mtx{\Omega}$ is the economized QR decomposition of $\mtx{A}\mtx{\Omega}$, then 
\[\mathbb{P}\left[ \| \mtx{A} - \mtx{Q}\mtx{Q}^*\mtx{A} \|_\Frob \leq \left(1 + 16\sqrt{k+5} \right) \sqrt{\sum_{j=k+1}^n \sigma_j^2(\mtx{A})} \,\right] \geq 0.999.\]
\end{corollary}

The following proposition bounds the expected approximation error of the randomized SVD with multivariate Gaussian inputs.

\begin{proposition}\label{prop_exp_random_svd}
Let $\mtx{A}$ be an $m\times n$ matrix, $k\geq 1$ an integer, and choose an oversampling parameter $p\geq 2$. If $\mtx{\Omega}\in\mathbb{R}^{n\times (k+p)}$ is a Gaussian random matrix, where each column is sampled from a multivariate Gaussian distribution with covariance matrix $\mtx{K}\in\mathbb{R}^{n\times n}$, and $\mtx{Q}\mtx{R} = \mtx{A}\mtx{\Omega}$ is the economized QR decomposition of $\mtx{A}\mtx{\Omega}$, then,
\[\E\left[\|\mtx{A} - \mtx{Q}\mtx{Q}^*\mtx{A}\|_{\Frob}\right]\leq \left(1+\sqrt{\frac{\beta_k}{\gamma_k}\frac{k(k+p)}{p-1}}\,\right)\sqrt{\sum_{j=k+1}^n\sigma^2_j(\mtx{A})},\]
where $\gamma_k = k/(\lambda_1 \Tr((\mtx{V}_1^*\mtx{K}\mtx{V}_1)^{-1})))$ and $\beta_k = \Tr(\mtx{\Sigma}_2^2\mtx{V}_2^*\mtx{K}\mtx{V}_2)/(\lambda_1\|\mtx{\Sigma}_2\|_\Frob^2)$.
\end{proposition}

We remark that for standard Gaussian inputs, we have $\gamma_k=\beta_k=1$ in \cref{prop_exp_random_svd}, and we recover the average Frobenius error of the randomized SVD~\cite[Thm.~10.5]{halko2011finding} up to a factor of $(k+p)$ due to the non-independence of $\mtx{\Omega}_1$ and $\mtx{\Omega}_2$ in general. The proof of \cref{prop_exp_random_svd} consists of combining the proof of \cref{th_tropp_random_svd_Frob} with the following lemma, which is a refinement of~\cref{prop_rand_1}.

\begin{lemma}\label{lem_1}
Let $\ell\geq 1$, $\mtx{\Omega}\in\mathbb{R}^{n\times \ell}$ be a Gaussian random matrix, where each column is sampled from a multivariate Gaussian distribution with covariance matrix $\mtx{K}$, and $\mtx{T}$ be an $\ell\times k$ matrix. Then,
\begin{equation}
\E[\|\mtx{\Sigma}_2\mtx{V}^*_2\mtx{\Omega}\mtx{T}\|_{\textup{F}}^2]= \Tr(\mtx{\Sigma}_2^2\mtx{V}_2^*\mtx{K}\mtx{V}_2)\|\mtx{T}\|_{\textup{F}}^2.
\end{equation}
\end{lemma}

\begin{proof}
Let $\mtx{K}=\mtx{Q}_\mtx{K}\mtx{\Lambda}\mtx{Q}_\mtx{K}^*$ be the eigenvalue decomposition of $\mtx{K}$, where $\mtx{Q}_\mtx{K}$ is orthonormal and $\mtx{\Lambda}$ is a diagonal matrix containing the eigenvalues of $\mtx{K}$ in decreasing order: $\lambda_1\geq \cdots\geq \lambda_n\geq 0$. We note that $\mtx{\Omega}$ can be expressed as $\mtx{\Omega} = \mtx{Q}_\mtx{K}\mtx{\Lambda}^{1/2}\mtx{G}$, where $\mtx{G}$ is a standard Gaussian matrix. Let $\mtx{S}=\mtx{\Sigma}_2 \mtx{V}_2^* \mtx{Q}_\mtx{K}\mtx{\Lambda}^{1/2}$, the proof follows from~\cite[Prop.~A.1]{halko2011finding}, which shows that $\E\|\mtx{S}\mtx{G}\mtx{T}\|_\Frob^2=\|\mtx{S}\|_\Frob^2\|\mtx{T}\|_\Frob^2$.
\end{proof}

Note that one can bound the term $\Tr(\mtx{\Sigma}_2^2\mtx{V}_2^*\mtx{K}\mtx{V}_2)$ by $\lambda_1\|\mtx{\Sigma}_2\|_{\textup{F}}^2$, where $\lambda_1$ is the largest eigenvalue of $\mtx{K}$ (cf.~\cref{prop_rand_1}). While this provides a simple upper bound, it does not demonstrate that the use of a covariance matrix containing prior information on the singular vectors of $\mtx{A}$ can outperform the randomized SVD with standard Gaussian inputs. 

\section{Randomized SVD for Hilbert--Schmidt operators} \label{sec_approx_HS}

We now describe the randomized SVD for learning HS operators (see \cref{alg_SVD}). The algorithm is implemented in the Chebfun software system~\cite{driscoll2014chebfun}, which is a MATLAB package for computing with functions. The Chebfun implementation of the randomized SVD for HS operators uses Chebfun's capabilities, which offer continuous analogues of several matrix operations like the QR decomposition and numerical integration. Indeed, the continuous analogue of a matrix-vector multiplication $\mtx{A}\mtx{\Omega}$ for an HS integral operator $\F$ (see~\cref{sec_HS} for definitions and properties of HS operators), with kernel $G:D\times D\to \R$, is
\[(\F f)(x) = \int_D G(x,y)f(y)\d y, \quad x\in D,\, f\in L^2(D),\]
where $D\subset\R^d$ with $d\geq 1$.

\renewcommand{\algorithmicrequire}{\textbf{Input:}}
\renewcommand{\algorithmicensure}{\textbf{Output:}}
\begin{algorithm}
\caption{Randomized SVD for HS operators}\label{alg_SVD}
\begin{algorithmic}[1]
\Require HS integral operator $\F$ with kernel $G(x,y)$, number of samples $k>0$
\Ensure Approximation $G_k$ of $G$
\State Define a GP covariance kernel $K$
\State Sample the GP $k$ times to generate a quasimatrix of random functions $\Omega=[f_1 \ldots  f_k]$
\State Evaluate the integral operator at $\Omega$, $Y = [\F(f_1)\ldots\F(f_k)]$
\State Orthonormalize the columns of $Y$, $Q=\text{orth}(Y)=[q_1\ldots q_k]$
\State Compute an approximation to $G$ by evaluating the adjoint of $\F$
\State Initialize $G_k(x,y)$ to $0$
\For{$i=1:k$}
\State $G_k(x,y) \gets G_k(x,y) + q_i(x)\int_D G(z,y)q_i(z)\d z$
\EndFor
\end{algorithmic}
\end{algorithm}

The algorithm takes as input an integral operator that we aim to approximate. Note that we focus here on learning an integral operator, but other HS operators would work similarly.
The first step of the randomized SVD for HS operators consists of generating a $D\times k$ quasimatrix $\Omega$ by sampling a GP $k$ times, where $k$ is the target rank (see \cref{sec_cov_kernel}). Therefore, each column of $\Omega$ is an object, consisting of a polynomial approximation of a smooth random function sampled from the GP in the Chebyshev basis. After evaluating the HS operator at $\Omega$ to obtain a quasimatrix $Y$, we use the QR algorithm~\cite{townsend2015continuous} to obtain an orthonormal basis $Q$ for the range of the columns of $Y$. Then, the randomized SVD for HS operators requires the left-evaluation of the operator $\F$ or, equivalently, the evaluation of its adjoint $\F_t$ satisfying:
\[(\F_t f)(x) = \int_D G(y,x)f(y)\d y, \quad x\in D.\]
We evaluate the adjoint of $\F$ at each column vector of $Q$ to construct an approximation $G_k$ of $G$. 
Finally, the approximation error between the operator kernel $G$ and the learned kernel $G_k$ can be computed in the $L^2$-norm, corresponding to the HS norm of the integral operator.

\section{Covariance kernels}  \label{sec_cov_kernel}

\begin{figure}[htbp]
\vspace{0.5cm}
\centering
\begin{overpic}[width=0.95\textwidth]{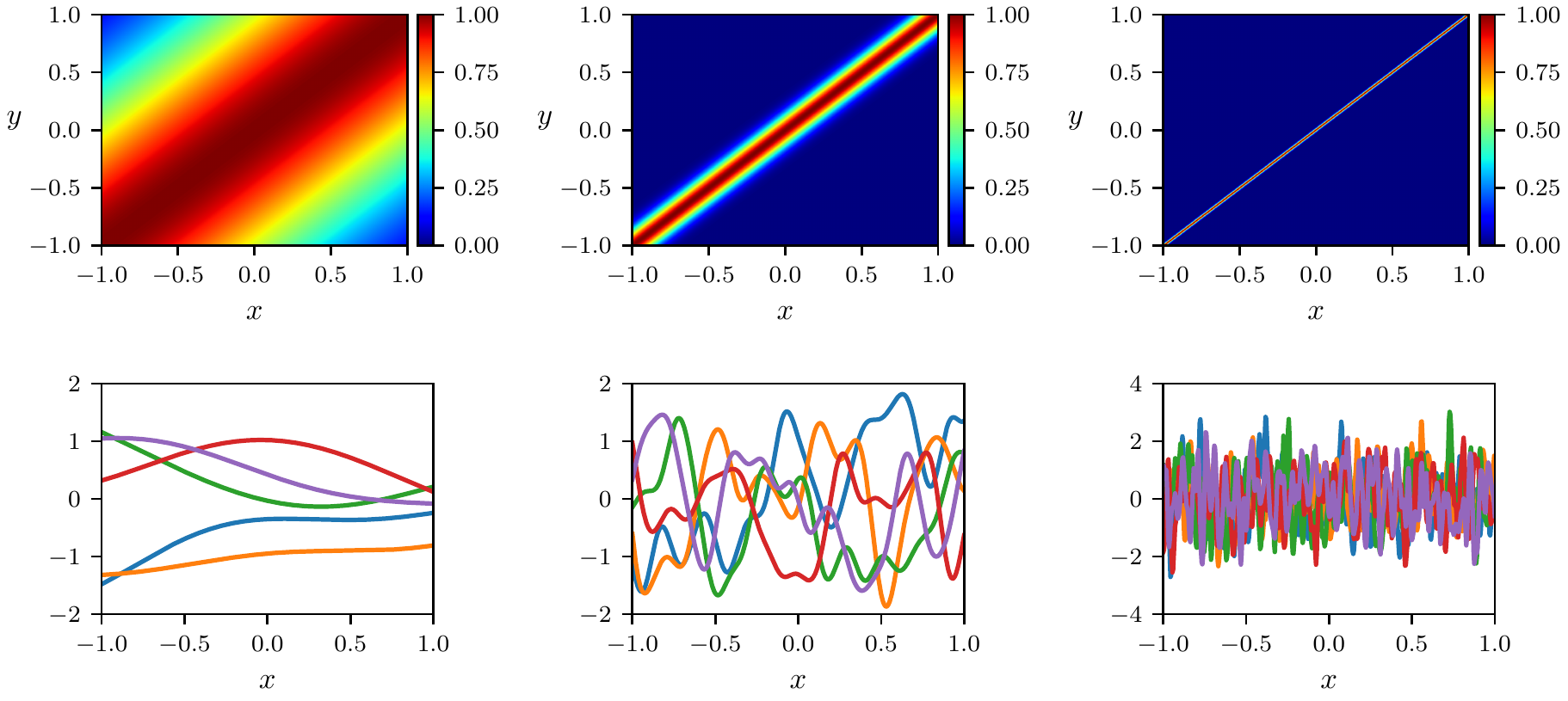}
\put(14,46){$\ell=1$}
\put(47,46){$\ell=0.1$}
\put(80,46){$\ell=0.01$}
\end{overpic}
\caption{Squared-exponential covariance kernel $K_{\text{SE}}$ with parameter $\ell=1,0.1,0.01$ (top row) and five functions sampled from $\mathcal{GP}(0, K_{\text{SE}})$ (bottom row).}
\label{fig_examples_GP}
\end{figure}

To generate the random input functions $f_1,\ldots,f_k$ for the randomized SVD for HS operators, we draw them from a GP, denoted by $\mathcal{GP}(0,K)$, for a certain covariance kernel $K$. A widely employed covariance kernel is the squared-exponential function $K_{\text{SE}}$~\cite{rasmussen2006gaussian} given by
\begin{equation} \label{eq_SE_kernel}
K_{\text{SE}}(x,y) = \exp\left(-|x-y|^2/(2\ell^2)\right),\quad x,y\in D,
\end{equation}
where $\ell>0$ is a parameter controlling the length-scale of the GP. This kernel is isotropic as it only depends on $|x-y|$, is infinitely differentiable, and its eigenvalues decay supergeometrically to $0$. Since the bound in~\cref{th_svd_main} degrades as the ratio $\lambda_1/\lambda_j$ increases for $j\geq k+1$ (cf.~\cref{eq_bound_gamma_beta}), the randomized SVD for learning HS operators prefers covariance kernels with slowly decaying eigenvalues. 
Our randomized SVD cannot hope to learn HS operators where the range of the operator has a rank greater than $\tilde{k}$, where $\tilde{k}$ is such that the $\tilde{k}$th eigenvalue of $K_{\text{SE}}$ reaches machine precision. In \cref{fig_examples_GP}, we display the squared-exponential kernel with length-scale parameters $\ell = 1,0.1,0.01$ together with sampled functions from $\mathcal{GP}(0,K_{\text{SE}})$. We observe that the functions become more oscillatory as the length-scale parameter $\ell$ decreases and hence the numerical rank of the kernel increases or, equivalently, the associated eigenvalues $\{\lambda_j\}$ decay more slowly to zero.

Other popular kernels for GPs include the Mat\'ern kernel~\cite{lindgren2011explicit,rasmussen2006gaussian} and Brownian bridge~\cite{nelsen2020random}. Prior information on the HS operator can also be enforced through the choice of the covariance kernel. For instance, one can impose the periodicity of the samples by using the following squared-exponential periodic kernel:
\[K_\text{Per}(x,y) = \exp\left(-\frac{2}{\ell^2}\sin^2\left(\frac{x-y}{2}\right)\right),\quad x,y\in D,\]
where $\ell>0$ is the length-scale parameter.

\subsection{Sample random functions from a Gaussian process} \label{sec_numer_GP}

In finite dimensions, a random vector $u\sim \mathcal{N}(0,\mtx{K})$, where $\mtx{K}\in \R^{n\times n}$ is a covariance matrix with Cholesky factorization $\mtx{K}=\mtx{L}\mtx{L}^*$, can be generated from the matrix-vector product $u = \mtx{L} c$. Here, $c\in \R^n$ is a vector whose entries follow the standard Gaussian distribution. We now detail how this process extends to infinite dimensions with a continuous covariance kernel. Let $K$ be a continuous symmetric positive-definite covariance function defined on the domain $[a,b]\times [a,b]\subset \R^2$ with $-\infty < a < b < \infty$. We consider the continuous analogue of the Cholesky factorization to write $K$ as~\cite{townsend2015continuous}
\[K(x,y) = \sum_{j=1}^\infty r_j(x) r_j(y) = L_c(x)L_c^*(y), \quad x,y\in [a,b],\]
where $r_j$ is the $j$th row of $L_c$, which ---in Chebfun's terminology--- is a lower-triangular quasimatrix. In practice, we truncate the series after $n$ terms, either arbitrarily or when the $n$th largest kernel eigenvalue, $\lambda_n$, falls below machine precision. Then, if $c\in\R^n$ follows the standard Gaussian distribution, a function $u$ can be sampled from $\mathcal{GP}(0,K)$ as $u=L_c c$. That is,
\[u(x) = \sum_{j=1}^n c_j r_j(x), \quad x\in [a,b].\]
The continuous Cholesky factorization is implemented in Chebfun2~\cite{townsend2013extension}, which is the extension of Chebfun for computing with two-dimensional functions. As an example, the polynomial approximation, which is accurate up to essentially machine precision, of the squared-exponential covariance kernel $K_{\text{SE}}$ with parameter $\ell = 0.01$ on $[-1,1]^2$ yields a numerical rank of $n = 503$. The functions sampled from $\mathcal{GP}(0,K_{\text{SE}})$ become more oscillatory as the length-scale parameter $\ell$ decreases and hence the numerical rank of the kernel increases or, equivalently, the associated eigenvalues sequence $\{\lambda_j\}$ decays more slowly to zero.

\subsection{Influence of the kernel's eigenvalues and Mercer's representation} \label{sec_choice_lambda}

The covariance kernel can also be defined from its Mercer's representation as
\begin{equation} \label{eq:KLexpansion}
K(x,y)=\sum_{j=1}^\infty\lambda_j \psi_j(x)\psi_j(y),\quad x,y\in D,
\end{equation} 
where $\{\psi_j\}$ is an orthonormal basis of $L^2(D)$ and $\lambda_1\geq \lambda_2\geq \cdots> 0$~\cite[Thm.~4.6.5]{hsing2015theoretical}. We prefer to construct $K$ directly from Mercer's representations for several reasons. First, one can impose prior knowledge of the kernel of the HS operator on the eigenfunctions of $K$ (such as periodicity or smoothness). Then, one can often generate samples from $\mathcal{GP}(0,K)$ efficiently using~\cref{eq:KLexpansion}. Finally, one can control the decay rate of the eigenvalues of $K$. 

Hence, the quantity $\gamma_k$ in the probability bound of~\cref{th_svd_main} measures the quality of the covariance kernel $K$ in $\smash{\mathcal{GP}(0,K)}$ to generate random functions that can learn the HS operator $\mathcal{F}$. To minimize $1/\gamma_k$ we would like to select the eigenvalues $\lambda_1\geq \lambda_2\geq \cdots >0$ of $K$ so that they have the slowest possible decay rate while maintaining $\smash{\sum_{j=1}^\infty \lambda_j <\infty}$. One needs $\smash{\{\lambda_j\}\in \ell^1}$ to guarantee that $\smash{\omega \sim\mathcal{GP}(0,K)}$ has finite expected squared $L^2$-norm, \emph{i.e.},~$\smash{\E[\|\omega\|_{L^2(D)}^2]=\sum_{j=1}^\infty\lambda_j<\infty}$. The best sequence of eigenvalues we know that satisfies this property is called the Rissanen sequence~\cite{rissanen1983universal} and is given by $\smash{\lambda_j = R_j \coloneqq 2^{-L(j)}}$, where  
\[
L(j) = \log_2(c_0)+\log_2^*(j), \quad \log_2^*(j) = \sum_{i=2}^\infty \max(\log_2^{(i)}(j),0), \quad c_0 = \sum_{i=2}^\infty 2^{-\log_2^*(i)},
\]
and $\log_2^{(i)}(j)=\log_2\circ\cdots\circ\log_2(j)$ is the composition of $\log_2(\cdot)$ $i$ times. Other options for the choice of eigenvalues include any sequence of the form $\lambda_j  = j^{-\nu}$ for $\nu>1$.

\subsection{Jacobi covariance kernel} \label{sec_Jacobi_kernel}

If $D=[-1,1]$, then a natural choice of orthonormal basis of $L^2(D)$ to define the Mercer's representation of the kernel are weighted Jacobi polynomials~\cite{deheuvels2008karhunen,olver2010nist}. That is, for a weight function $w_{\alpha,\beta}(x) = (1-x)^\alpha(1+x)^\beta$ with $\alpha,\beta>-1$, and any positive eigenvalue sequence $\{\lambda_j\}$, we consider the Jacobi kernel
\begin{equation} \label{eq:JacobiKernel}
K_{\text{Jac}}^{(\alpha,\beta)}(x,y) = \sum_{j=0}^\infty \lambda_{j+1}w_{\alpha,\beta}^{1/2}(x)\tilde{P}^{(\alpha,\beta)}_j(x)w_{\alpha,\beta}^{1/2}(y)\tilde{P}^{(\alpha,\beta)}_j(y), \quad x,y\in[-1,1],
\end{equation} 
where $\smash{\tilde{P}^{(\alpha,\beta)}_j}$ is the scaled Jacobi polynomial of degree $j$ and parameters $(\alpha,\beta)$ where $P_j^{(\alpha,\beta)}$ is defined by Rodrigues' formula~\cite[Eq.~4.3.1]{Szego1939orthogonal} as
\[w_{\alpha,\beta}(x) P_j^{(\alpha,\beta)}(x)=\frac{(-1)^j}{2^j n!}\frac{d^j}{dx^j}\left\{w_{\alpha,\beta}(x)(1-x^2)^j\right\}.\] 
The polynomials $\smash{\tilde{P}^{(\alpha,\beta)}_j}$ are normalized such that $\smash{\|w_{\alpha,\beta}^{1/2}\tilde{P}^{(\alpha,\beta)}_j\|_{L^2([-1,1])}=1}$ and $\{\lambda_j\}$ is chosen such that $K_{\text{Jac}}^{(\alpha,\beta)}\in L^2([-1,1]^2)$. In this case, a random function can be sampled as
\[u(x) = \sum_{j=0}^\infty \sqrt{\lambda_{j+1}} c_j w_{\alpha,\beta}^{1/2} \tilde{P}_j^{(\alpha,\beta)}(x), \quad x\in [-1,1],\]
where $c_j\sim \mathcal{N}(0,1)$ for $0\leq j\leq \infty$.

A desirable property of a covariance kernel is to be unbiased towards one spatial direction, \emph{i.e.}, $K(x,y)=K(-y,-x)$ for $x,y\in [-1,1]$, which motivates us to always select $\alpha=\beta$. Moreover, it is desirable to have the eigenfunctions of $\smash{K_{\text{Jac}}^{(\alpha,\beta)}}$ to be polynomial so that one can generate samples from $\mathcal{GP}(0,K)$ efficiently. This leads us to choose $\alpha$ and $\beta$ to be even integers. The choice of $\alpha=\beta=0$ gives the Legendre kernel~\cite{foster2020optimal,habermann2019semicircle}. In the rest of this chapter, we will use~\cref{eq:JacobiKernel} with $\alpha = \beta = 2$ to ensure that functions sampled from the associated GP satisfy homogeneous Dirichlet boundary conditions (see~\cref{fig_examples_legendre_GP}). We emphasize that covariance kernels on higher dimensional domains of the form $D=[-1,1]^d$, for $d\geq 2$, can be defined using tensor products of weighted Jacobi polynomials.

\subsection{Smoothness of functions sampled from a GP with Jacobi kernel} \label{sec_smooth}

We now connect the decay rate of the eigenvalues of the Jacobi covariance kernel $K_{\text{Jac}}^{(2,2)}$ to the smoothness of the samples from $\mathcal{GP}(0,K_{\text{Jac}}^{(2,2)})$. Hence, the Jacobi covariance function allows the control of the decay rate of the eigenvalues $\{\lambda_j\}$ as well as the smoothness of the resulting randomly generated functions. First, \cref{th_regularity_series} asserts that if the coefficients of an infinite polynomial series have sufficient decay, then the resulting series is smooth with regularity depending on the decay rate. This result can be seen as a converse to~\cite[Thm.~7.1]{trefethen2019approximation}.

\begin{lemma} \label{th_regularity_series}
Let $\{p_j\}$ be a family of polynomials such that $\max_{x\in[-1,1]} |p_j(x)| = 1$ and $\textup{deg}(p_j)\leq j$. If $f_n(x)=\sum_{j=0}^n a_j p_j(x)$ with $|a_j|\leq j^{-\nu}$ for $\nu>1$, then $f_n$ converges uniformly to $f(x)=\sum_{j=0}^\infty a_j p_j(x)$ and $f$ is $\mu$ times continuously differentiable for any integer $\mu$ such that $\mu<(\nu-1)/2$. 
\end{lemma}

\begin{proof}
By Markov brothers' inequality~\cite{markov1890question}, for all $j\geq 0$ and $0\leq \mu\leq j$, we have $\max_{x\in[-1,1]} |p_j^{(\mu)}(x)| \leq j^{2\mu}$. 
Therefore, $|f_n^{(\mu)}(x)|\leq \sum_{j=0}^n |a_j| \|p_j^{(\mu)}\|_{\infty}\leq \sum_{j=0}^n j^{2\mu-\nu}$ so $|f_n^{(\mu)}(x)|<\infty$ if $\mu<(\nu-1)/2$. The result follows from a standard result on uniform convergence and differentiation~\cite[Thm.~7.17]{rudin1976principles}.
\end{proof}

Note that the main application of this lemma occurs when $\textup{deg}(p_j)=j$ for all $j\geq 0$. We then prove a bound on ultraspherical polynomials in order to apply \cref{th_regularity_series} to functions sampled from the GP with the Jacobi covariance kernel $K_{\text{Jac}}^{(2,2)}$. First, note that $\tilde{P}_j^{(2,2)}$ is a scaled ultraspherical polynomial $\tilde{C}^{(5/2)}_j$ with parameter $5/2$ and degree $j\geq 0$ so it can be bounded by the following proposition.
 
\begin{proposition} \label{prop_bound_ultra}
Let $\tilde{C}^{(5/2)}_j$ be the ultraspherical polynomial of degree $j$ with parameter $5/2$, normalized such that $\int_{-1}^1 (1-x^2)^2\tilde{C}^{(5/2)}_j(x)^2\d x=1$. Then,
\begin{equation} \label{eq_bound_ultra}
\max_{x\in [-1,1]}|(1-x^2)\tilde{C}^{(5/2)}_j(x)| \leq 2\sqrt{j+5/12}, \quad j\geq 0.
\end{equation}
\end{proposition}
\begin{proof}
Let $j\geq 0$ and $x\in[-1,1]$, according to~\cite[Table~18.3.1]{olver2010nist},
\begin{equation} \label{eq_norm_ultra}
\tilde{C}^{(5/2)}_j(x) = 3\sqrt{\frac{j+5/2}{(j+1)(j+2)(j+3)(j+4)}}C^{(5/2)}_j(x),
\end{equation}
where $C^{(5/2)}_j(x)$ is the standard ultraspherical polynomial. Using~\cite[(18.9.8)]{olver2010nist}, we have
\[(1-x^2)C^{(5/2)}_j(x) = \frac{(j+3)(j+4)C_j^{(3/2)}(x)-(j+1)(j+2)C_{j+2}^{(3/2)}(x)}{6(j+5/2)}.\]
By using~\cite[(18.9.7)]{olver2010nist}, we have $(C_{j+2}^{(3/2)}(x)-C_{j}^{(3/2)}(x))/2 = (j+5/2)C_{j+2}^{(1/2)}(x)$ and hence, 
\[(1-x^2)C^{(5/2)}_j(x) = \frac{2}{3}C_j^{(3/2)}(x)-\frac{(j+1)(j+2)}{3}C_{j+2}^{(1/2)}(x).\]
We bound the two terms with~\cite[(18.14.4)]{olver2010nist} to obtain the following inequalities:
\[|C_{j}^{(3/2)}(x)|\leq \frac{(j+1)(j+2)}{2},\quad |C_{j+2}^{(1/2)}(x)|\leq 1.\]
Hence, $|(1-x^2)C^{(5/2)}_j(x)| \leq 2(j+1)(j+2)/3$ and following \cref{eq_norm_ultra} we obtain
\[|(1-x^2)\tilde{C}^{(5/2)}_j(x)|\leq 2\sqrt{\frac{(j+1)(j+2)(j+5/2)}{(j+3)(j+4)}}\leq 2\sqrt{j+5/12},\]
which concludes the proof.
\end{proof}

\begin{figure}[htbp]
\centering
\vspace{0.2cm}
\begin{overpic}[width=0.8\textwidth]{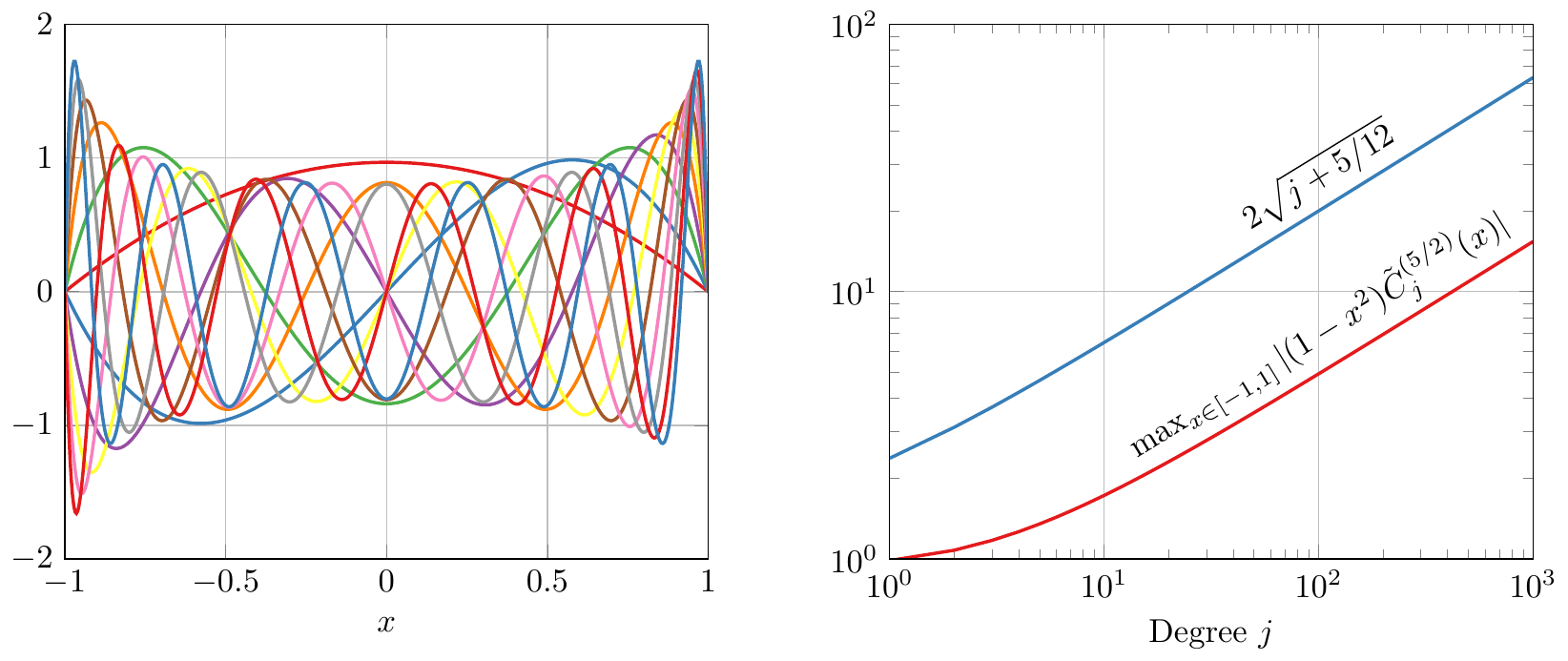}
\put(-3,39){(a)}
\put(47.5,39){(b)}
\end{overpic}
\caption{(a) Normalized ultraspherical polynomials $\tilde{C}^{(5/2)}_j$ up to degree $10$. (b) Theoretical bound (blue line) for the maximum of weighted ultraspherical polynomials on $[-1,1]$, obtained in \cref{prop_bound_ultra}, against the one observed numerically (red line).}
\label{fig_bound_ultra}
\end{figure}

The bound given in \cref{prop_bound_ultra} differs initially by a factor of $4/3$ from the numerically observed upper bound ($1.5\sqrt{j+5/12}$\,) as shown by \cref{fig_bound_ultra}. We now state the following theorem about the regularity of functions sampled from $\smash{\mathcal{GP}(0,K_{\text{Jac}}^{(2,2)})}$, which guarantees that if the eigenvalues are chosen such that $\lambda_j=\mathcal{O}(1/j^{\nu})$ with $\nu>3$, then $\smash{f\sim \mathcal{GP}(0,K_{\text{Jac}}^{(2,2)})}$ is almost surely continuous. Moreover, a faster decay of the eigenvalues of $\smash{K_{\text{Jac}}^{(2,2)}}$ implies higher regularity of the sampled functions, in an almost sure sense.

\begin{theorem} \label{th_regularity_GP}
Let $\{\lambda_j\}\in \ell^1(\R^+)$ be a positive sequence such that $\lambda_j = \mathcal{O}(j^{-\nu})$ for $\nu>3$. If $f$ is sampled from $\mathcal{GP}(0,K_{\text{Jac}}^{(2,2)})$, then $f\in \mathcal{C}^\mu([-1,1])$ almost surely for any integer $\mu < (\nu-3)/2$.
\end{theorem}

\begin{proof}
Since $f\sim \mathcal{GP}(0,K_{\text{Jac}}^{(2,2)})$, $f \sim \sum_{j=0}^\infty c_j\sqrt{\lambda_{j+1}}(1-x^2)\tilde{P}^{(2,2)}_j(x)$, where $c_j\sim\mathcal{N}(0,1)$ for $j\geq 0$. Let $f_n$ denote the truncation of $f$ after $n$ terms. By letting $M>0$ be the constant such that $\lambda_{j+1}\leq M (j+1)^{-\nu}$, we find that 
\[
\|f-f_n\|_\infty\leq S_n, \quad S_n \coloneqq 2\sqrt{M} \sum_{j=n+2}^\infty |c_{j-1}| j^{(1-\nu)/2},
\]
where we used $\max_{x\in[-1,1]} |(1-x^2)\tilde{P}_j^{(2,2)}(x)| \leq 2\sqrt{j+1}$ (cf.~\cref{prop_bound_ultra}). Thus, we have 
\[\mathbb{P}\left(\lim_{n\to\infty}\|f-f_n\|_\infty=0\right)\geq \mathbb{P}\left(\lim_{n\to\infty}S_n=0\right).\] Here, $S_n\sim X_n=\sum_{j=n+2}^\infty Y_j j^{(1-\nu)/2}$, where $Y_j$ follows a half-normal distribution~\cite{leone1961folded} with parameter $\sigma=1$ and the $(Y_j)_j$ are independent. We want to show that $X_n\xrightarrow{a.s.}0$. For $\epsilon>0$, using Chebyshev's inequality, we have:
\[\sum_{n=0}^\infty \mathbb{P}(|X_n|\geq \epsilon) \leq \frac{1}{\epsilon^2}\sum_{n=0}^\infty \left(1-\frac{2}{\pi}\right)\sum_{j=n+2}^\infty \frac{1}{j^{\nu-1}} \leq \frac{1}{\epsilon^2}\left(1-\frac{2}{\pi}\right)\frac{1}{\nu-2}\sum_{n=1}^\infty \frac{1}{n^{\nu-2}},\]
which is finite if $\nu>3$. Therefore, using the Borel--Cantelli Lemma~\cite[Chapt.~2.3]{durrett2019probability}, $X_n$ converges to $0$ almost surely and
$\mathbb{P}(\lim_{n\to\infty}X_n=0)=1$. Finally,
\[\mathbb{P}\left(\lim_{n\to\infty}\|f-f_n\|_\infty=0\right)\geq \mathbb{P}\left(\lim_{n\to\infty}X_n=0\right)=1,\]
which proves that $\{f_n\}$ converges uniformly and hence $f$ is continuous with probability one. The statement for higher order derivatives follows the proof of \cref{th_regularity_series}.
\end{proof}

\begin{figure}[htbp]
\vspace{0.5cm}
\centering
\begin{overpic}[width=0.95\textwidth]{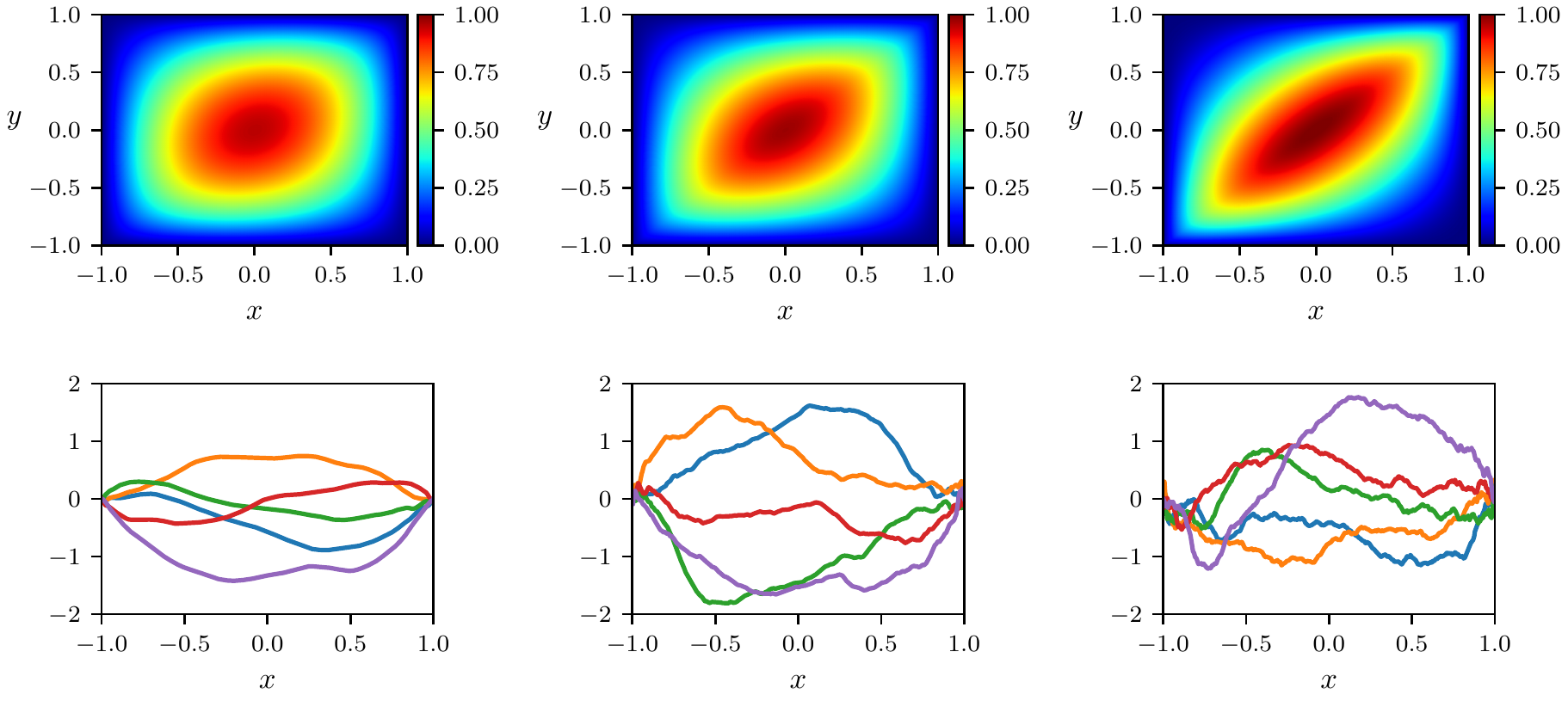}
\put(11.5,46){$\lambda_j=1/j^4$}
\put(45,46){$\lambda_j=1/j^3$}
\put(78,46){$\lambda_j=R_j/j$}
\end{overpic}
\caption{Covariance kernel $K_{\text{Jac}}^{(2,2)}$ constructed using Jacobi polynomials of type $(2,2)$ with $\lambda_j=1/j^4$, $1/j^3$, and $R_j/j$, where $R_j$ is the Rissanen sequence (top). The bottom panels illustrate functions sampled from $\smash{\mathcal{GP}(0, K_{\text{Jac}}^{(2,2)})}$ with the different eigenvalue sequences. The series for generating the random functions are truncated to $n=500$.}
\label{fig_examples_legendre_GP}
\end{figure}

This theorem can be seen as a particular case of Driscoll's zero-one law~\cite{driscoll1973reproducing}, which characterizes the regularity of functions samples from GPs (see also~\cite{kanagawa2018gaussian}). Hence, one must have $\sum_{j=1}^\infty j\lambda_{j} < \infty$ so that the series of functions in \cref{eq:JacobiKernel} converges uniformly and $\smash{K_{\text{Jac}}^{(2,2)}}$ is a continuous kernel. Under this additional constraint, the best choice of eigenvalues is given by a scaled Rissanen sequence: $\lambda_j = R_j/j$, for $j\geq 1$ (cf.~\cref{sec_choice_lambda}). In \cref{fig_examples_legendre_GP}, we display the Jacobi kernel of type $(2,2)$ with functions sampled from the corresponding GP. We selected eigenvalue sequences of different decay rates: from the faster $1/j^4$ to the slower Rissanen sequence $R_j/j$ (\cref{sec_choice_lambda}). For $\lambda_j=1/j^3$ and $\lambda_j=R_j/j$, we observe a large variation of the randomly generated functions near $x=\pm 1$, indicating a potential discontinuity of the samples at these two points as $n\to\infty$. This is in agreement with \cref{th_regularity_GP}, which only guarantees continuity (with probability one) of the randomly generated functions if $\lambda_j \sim 1/j^\nu$ with $\nu>3$. 

\section{Numerical experiments} \label{HS_random_SVD}

We now perform several numerical experiments with the randomized SVD to learn matrices using random vectors sampled for a multivariate Gaussian distribution and HS operators.

\subsection{Covariance matrix with prior knowledge} \label{sec_exp_matrix}

The approximation error bound in~\cref{th_svd_main} depends on the eigenvalues of the covariance matrix, which dictates the distribution of the column vectors of the input matrix $\mtx{\Omega}$. Roughly speaking, the more prior knowledge of the matrix $\mtx{A}$ that can be incorporated into the covariance matrix, the better. In this numerical example, we investigate whether the standard randomized SVD, which uses the identity as its covariance matrix, can be improved by using a different covariance matrix. We then attempt to learn the discretized $2000\times 2000$ matrix, \emph{i.e.},~the discrete Green's function, of the inverse of the following differential operator:
\[\mathcal{L} u = d^2u/dx^2-100\sin(5\pi x)u,\quad x\in[0,1].\]
We vary the number of columns (\emph{i.e.} samples from the GP) in the input matrix $\mtx{\Omega}$ from $1$ to $2000$. 

\begin{figure}[htbp]
\centering
\vspace{0.5cm}
\begin{overpic}[width=\textwidth]{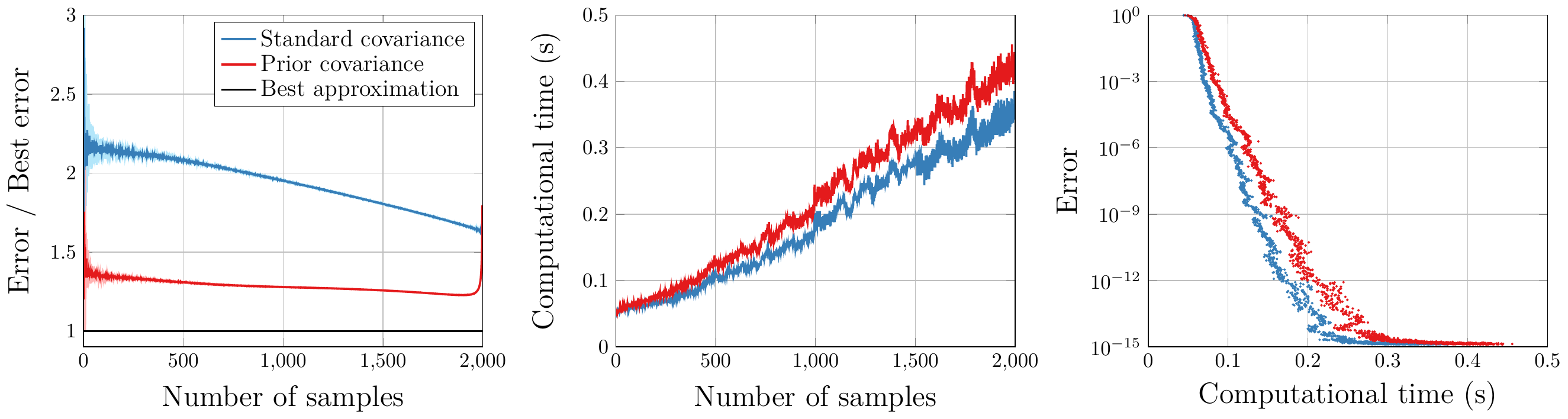}
\put(0,26){(a)}
\put(33,26){(b)}
\put(66.5,26){(c)}
\end{overpic}
\caption{(a) Ratio between the average randomized SVD approximation error (over 10 runs) of the $2000\times 2000$ matrix of the inverse of the differential operator $\mathcal{L} u = d^2u/dx^2-100\sin(5\pi x)u$ on $[0,1]$, and the best approximation error. The error bars in light colour (blue and red) illustrate one standard deviation. (b) Average computational time of the algorithm (over 10 runs). The eigenvalue decomposition of the covariance matrix has been precomputed offline. (c) Randomized SVD approximation error with standard and prior covariance matrices with respect to the computational time.}
\label{fig_SVD_matrix}
\end{figure}

In \cref{fig_SVD_matrix}(a), we compare the ratios between the relative error in the Frobenius norm given by the randomized SVD and the best approximation error, obtained by truncating the SVD of $\mtx{A}$. The prior covariance matrix $\mtx{K}$ consists of the discretized $2000\times 2000$ matrix of the Green's function of the negative Laplace operator $\mathcal{L}u = -d^2u/dx^2$ on $[0,1]$ to incorporate knowledge of the diffusion term in the matrix $\mtx{A}$. We see that a nonstandard covariance matrix leads to a higher approximation accuracy, with a reduction of the error by a factor of $1.3$-$1.6$ compared to the standard randomized SVD.

At the same time, the procedure is only $20\%$ slower\footnote{Timings were performed on an Intel Xeon CPU E5-2667 v2 @ 3.30GHz using MATLAB R2020b without explicit parallelization.} on average (\cref{fig_SVD_matrix}(b)) as one can precompute the eigenvalue decomposition of the covariance matrix. Hence, sampling a random vector from a multivariate normal distribution with an arbitrary covariance matrix $\mtx{K}$ can be computationally expensive when the dimension, $n$, of the matrix is large as it requires the computation of a Cholesky factorization, which can be done in $\mathcal{O}(n^3)$ operations. We highlight that this step can be precomputed once, such that the overhead of the generalized SVD can be essentially expressed as the cost of an extra matrix-vector multiplication. Then, the difference in timings between standard and prior covariance matrices is marginal as shown by \cref{fig_SVD_matrix}(b).

We observe in \cref{fig_SVD_matrix}(c) that using a standard covariance matrix offers a better trade-off between error and computational time. However, choosing a prior covariance matrix is of interest in applications where the sampling time is much higher than the numerical linear algebra costs to maximize the accuracy of the approximation matrix from a limited number of samples.

Additionally, we would like to highlight that prior covariance matrices can be designed and derived using physical knowledge of the problem, such as its diffusive nature, which can also significantly decrease the precomputation cost. In this example, we employ the discretized Green's function of the negative Laplacian operator with homogeneous Dirichlet boundary conditions, given by $\mathcal{L}u=-d^2u/dx^2$ on $[0,1]$, for which we know the eigenvalue decomposition. Hence, the eigenvalues and normalized eigenfunctions are respectively given by
\[\lambda_n=\frac{1}{\pi^2n^2},\quad \psi_n(x)=\sqrt{2}\sin(n\pi x),\quad x\in[0,1],\,n\geq 1.\]
Therefore, one can employ Mercer's representation (see~\cref{eq:KLexpansion}) to sample the random vectors and precompute the covariance matrix in $\mathcal{O}(n^2)$ operations. For a problem of size $n=2000$, it takes $0.16$s to precompute the matrix.

\subsection{Randomized SVD for Hilbert--Schmidt operators}

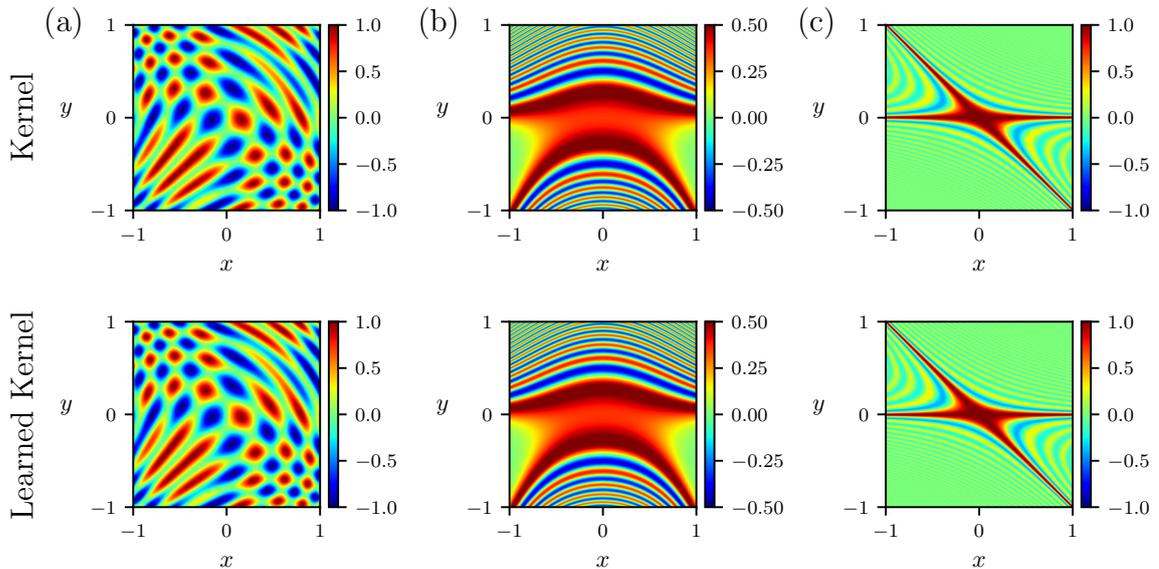
\begin{figure}[htbp]
\centering
\vspace{0.1cm}
\hspace{0.5cm}
\begin{overpic}[width=0.95\textwidth,trim = 0 0 0 0, clip]{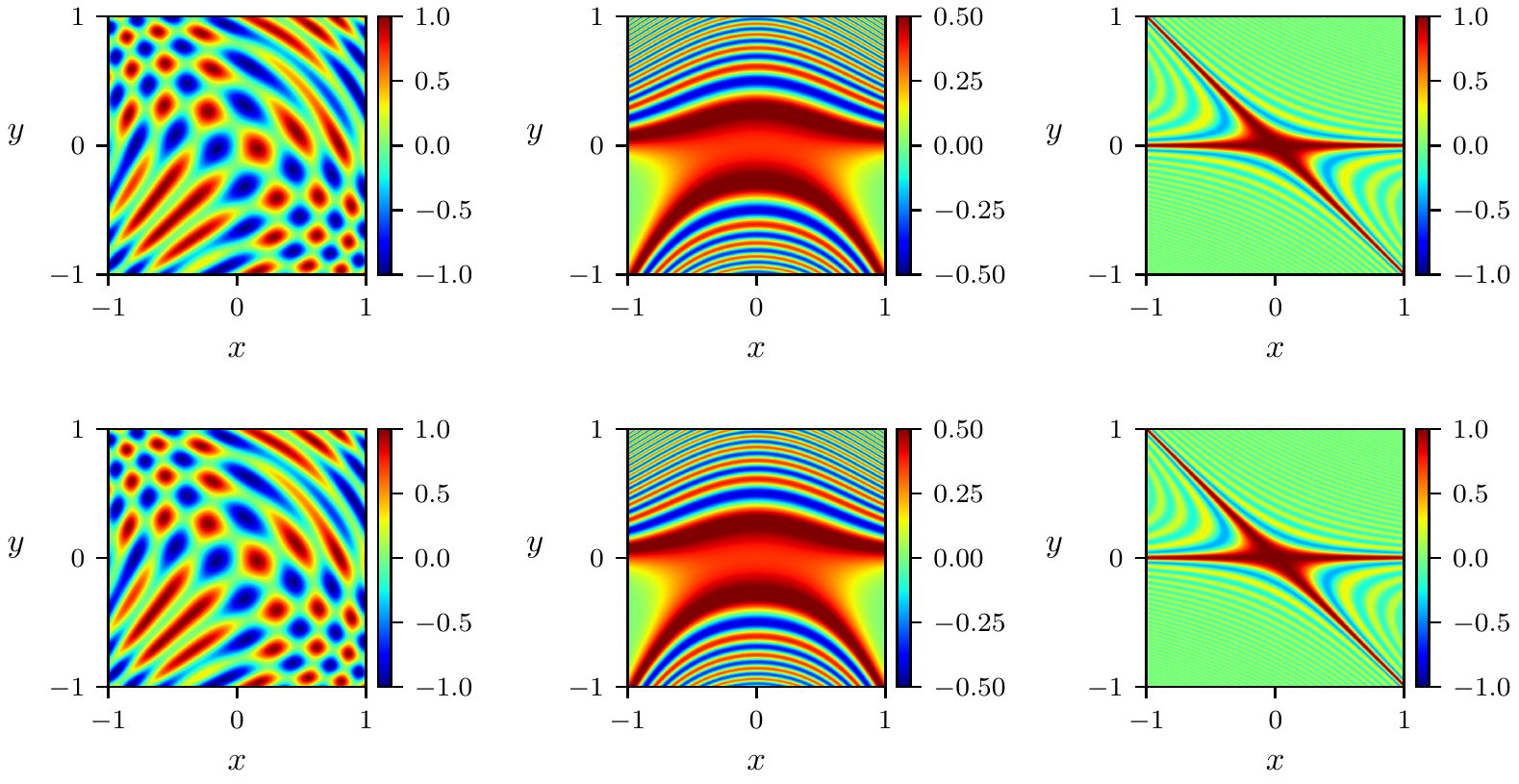}
\put(-1,50){(a)}
\put(33,50){(b)}
\put(67.5,50){(c)}
\put(-4,38){\rotatebox{90}{Kernel}}
\put(-4,5.2){\rotatebox{90}{Learned Kernel}}
\end{overpic}
\caption{Kernels of three HS operators (top) together with the kernels learned by the randomized SVD for HS operators (bottom), using the squared-exponential covariance kernel $K_{\text{SE}}$ with parameter $\ell=0.01$ and one hundred functions sampled from $\mathcal{GP}(0,K_{\text{SE}})$.}
\label{fig_examples_SVD}
\end{figure}

We now apply the randomized SVD for HS operators to learn kernels of integral operators. In this first example, the kernel is defined as~\cite{townsend2013example} 
\[G(x,y) = \cos(10(x^2+y))\sin(10(x+y^2)), \quad x,y\in [-1,1],\]
and is displayed in \cref{fig_examples_SVD}(a). We employ the squared-exponential covariance kernel $K_{\text{SE}}$ with parameter $\ell=0.01$ and $k=100$ samples (see \cref{eq_SE_kernel}) to sample random functions from the associated GP. The learned kernel $G_k$ is represented on the bottom panel of \cref{fig_examples_SVD}(a) and has an approximation error around machine precision.

As a second application of the randomized SVD for HS operators, we learn the kernel $G(x,y) = \text{Ai}(-13(x^2y+y^2))$ for $x,y\in[-1,1]$, where $\text{Ai}$ is the Airy function~\cite[Chapt.~9]{olver2010nist} defined by
\[\text{Ai}(x) = \frac{1}{\pi}\int_0^\infty\cos\left(\frac{t^3}{3}+xt\right)\d t, \quad x\in\R.\]
We plot the kernel and its low-rank approximant given by the randomized SVD for HS operators in \cref{fig_examples_SVD}(b) and obtain an approximation error (measured in the $L^2$-norm) of $5.04\times 10^{-14}$. The two kernels have a numerical rank equal to $42$.

The last example consists of learning the HS operator associated with the kernel $G(x,y)=J_0(100(xy+y^2))$ for $x,y\in [-1,1]$, where $J_0$ is the Bessel function of the first kind~\cite[Chapt.~10]{olver2010nist} defined as
\[J_0(x) = \frac{1}{\pi}\int_0^{\pi}\cos(x\sin t)\d t,\quad x\in\R,\]
and plotted in \cref{fig_examples_SVD}(c). The rank of this kernel is equal to $91$ while its approximation is of rank $89$ and the approximation error is equal to $4.88\times 10^{-13}$. We observe that in the three numerical examples displayed in \cref{fig_examples_SVD}, the differences between the learned and the original kernels are not visually perceptible.

\begin{figure}[htbp]
\centering
\vspace{0.1cm}
\begin{overpic}[width=0.9\textwidth]{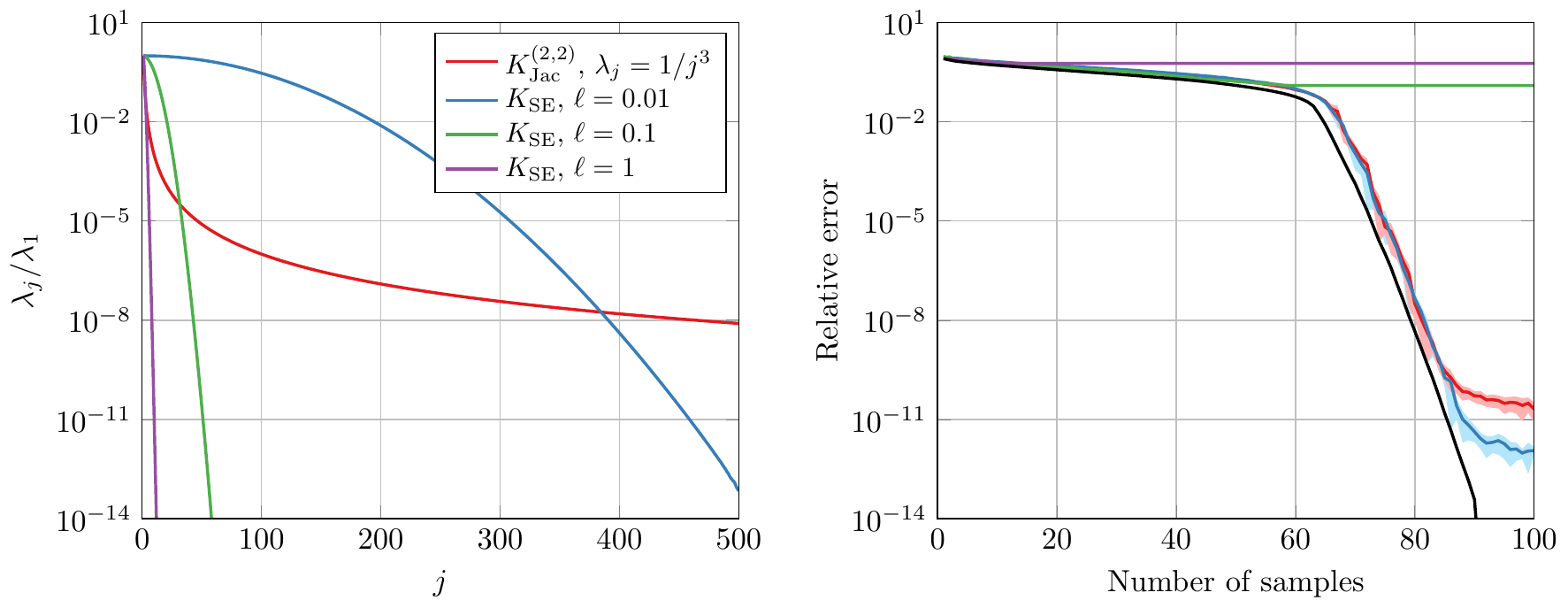}
\end{overpic}
\caption{Left: Scaled eigenvalues of the Jacobi covariance kernel $\smash{K_{\text{Jac}}^{(2,2)}}$ with sequence $\lambda_j=1/j^3$ and squared-exponential kernels $K_{\text{SE}}$ with parameters $\ell=0.01,0.1,1$, respectively. Right: Average (over 10 runs) relative approximation error in the $L^2$-norm between the Bessel kernel $G(x,y)=J_0(100(xy+y^2))$ and its low-rank approximation $G_k(x,y)$, obtained from the randomized SVD by sampling the GPs $k$ times. The error bars in light colour (blue and red) illustrate one standard deviation and the black line indicates the best approximation error given by the tail of the singular values of $G$.}
\label{fig_convergence_SVD}
\end{figure}

Finally, we evaluate the influence of the choice of covariance kernel and number of samples in \cref{fig_convergence_SVD}. Here, we vary the number of samples from $k=1$ to $k=100$ and use the randomized SVD for HS operators with four different covariance kernels: the squared-exponential $\smash{K_{\text{SE}}}$ with parameters $\ell = 0.01,0.1,1$, and the Jacobi kernel $\smash{K_{\text{Jac}}^{(2,2)}}$ with eigenvalues $\lambda_j = 1/j^3$, for $j\geq 1$. In the left panel of \cref{fig_convergence_SVD}, we represent the eigenvalue ratio $\lambda_j/\lambda_1$ of the four kernels and observe that this quantity falls below machine precision for the squared-exponential kernel with $\ell=1$ and $\ell=0.1$ at $j=13$ and $j=59$, respectively. In \cref{fig_convergence_SVD}(right), we observe that these two kernels fail to approximate kernels of high numerical rank. The other two kernels have a much slower decay of eigenvalues and can capture (or learn) more complicated kernels. We then see in the right panel of \cref{fig_convergence_SVD} that the relative approximation errors obtained using $\smash{K_{\text{Jac}}^{(2,2)}}$ and $\smash{K_{\text{SE}}}$ are close to the best approximation error given by the squared tail of the singular values of the integral kernel $G(x,y)$, \emph{i.e.},~$(\sum_{j\geq k+1}\sigma_j^2)^{1/2}$. The overshoot in the error at $k=100$ compared to the machine precision is due to the decay of the eigenvalues of the covariance kernels. Hence, spatial directions associated with small eigenvalues are harder to learn accurately. This issue does not arise in finite dimensions with the standard randomized SVD because the covariance kernel used there is isotropic, \emph{i.e.},~all its eigenvalues are equal to one. However, this choice is no longer possible for learning HS integral operators as the covariance kernel $K$ must be squared-integrable. The relative approximation errors at $k=100$ (averaged over 10 runs) using $\smash{K_{\text{Jac}}^{(2,2)}}$ and $K_{\text{SE}}$ with $\ell = 0.01$ are $\smash{\text{Error}(K_{\text{Jac}}^{(2,2)}) \approx 2.6\times 10^{-11}}$, and $\smash{\text{Error}(K_{\text{SE}})\approx 5.7\times 10^{-13}}$, which gives a ratio of
\begin{equation} \label{eq_ratio_rel}
\text{Error}(K_{\text{Jac}}^{(2,2)})/\text{Error}(K_{\text{SE}}) \approx 45.6.
\end{equation}
However, the square-root of the ratio of the quality of the two kernels for $k=91$ is equal to 
\begin{equation} \label{eq_ratio_gamma}
\sqrt{\gamma_{91}(K_{\text{SE}})/\gamma_{91}(K_{\text{Jac}}^{(2,2)})}\approx 117.8,
\end{equation}
which is of the same order of magnitude of \cref{eq_ratio_rel} as predicted by~\cref{th_svd_main}. In \cref{eq_ratio_gamma}, $\smash{\gamma_{91}(K_{\text{SE}})\approx 5.88\times 10^{-2}}$ and $\smash{\gamma_{91}(K_{\text{Jac}}^{(2,2)})\approx 4.24\times 10^{-6}}$ are both computed using Chebfun.

In conclusion, this section provides numerical insights to motivate the choice of the covariance kernel to learn HS operators. Following \cref{fig_convergence_SVD}, a kernel with slowly decaying eigenvalues is preferable and yields better approximation errors or higher learning rate with respect to the number of samples, especially when learning a kernel with a large numerical rank. The optimal choice from a theoretical viewpoint is to select a covariance kernel whose eigenvalues have a decay rate similar to the Rissanen sequence~\cite{rissanen1983universal}, but other choices may be preferable in practice to ensure smoothness of the sample functions (cf.~\cref{sec_smooth}).

\dobib
 % PDE learning

\setcounter{chapter}{3}
\renewcommand{\thefootnote}{\fnsymbol{footnote}}
\chapter[Rational neural networks]{Rational neural networks\footnotemark} \label{chapt_rational}

\footnotetext{This chapter is based on a paper with Yuji Nakatsukasa and Alex Townsend~\cite{boulle2020rational}, published in NeurIPS 2020. Nakatsukasa and Townsend had advisory roles and proved \cref{th_abs_approx}. I proved the other theoretical results, performed the numerical experiments, and was the lead author in writing the paper.}

\renewcommand*{\thefootnote}{\arabic{footnote}}
\setcounter{footnote}{0}

A key question in designing deep learning architectures is the choice of the activation function to reduce the number of trainable parameters of the network while keeping the same approximation power~\cite{goodfellow2016deep}. While smooth activation functions such as sigmoid, logistic, or hyperbolic tangent are widely used, they suffer from the ``vanishing gradient problem''~\cite{bengio1994learning} because their derivatives are zero for large inputs. Neural networks (NNs) based on polynomial activation functions are an alternative~\cite{cheng2018polynomial,daws2019polynomial, goyal2019learning, guarnieri1999multilayer,ma2005constructive,vecci1998learning}, but can be numerically unstable due to large gradients for large inputs~\cite{bengio1994learning}. Moreover, polynomials do not approximate non-smooth functions efficiently~\cite{trefethen2019approximation}, which can lead to optimization issues in classification problems. A popular choice of activation function is the Rectified Linear Unit (ReLU) defined as $\relu(x)=\max(x,0)$~\cite{jarrett2009best,nair2010rectified}. It has numerous advantages, such as being fast to evaluate and zero for many inputs~\cite{glorot2011deep}. Many theoretical studies characterize and understand the expressiveness of shallow and deep ReLU neural networks from the perspective of approximation theory~\cite{devore1989optimal,liang2016deep,mhaskar1996neural,telgarsky2016benefits,yarotsky2017error}.

ReLU networks also suffer from drawbacks, which are most evident during training. The main disadvantage is that the gradient of ReLU is zero for negative real numbers. Therefore, its derivative is zero if the activation function is saturated~\cite{maas2013rectifier}. Several adaptations to ReLU have been proposed over the past few years, such as Leaky ReLU~\cite{maas2013rectifier}, Exponential Linear Unit (ELU)~\cite{clevert2015fast}, Parametric Linear Unit (PReLU)~\cite{he2015delving}, and Scaled Exponential Linear Unit (SELU)~\cite{klambauer2017self}, to improve the initialization and optimization of neural networks and avoid the use of batch normalization layers~\cite{ioffe2015batch}. These modifications outperform ReLU in image classification applications, and some of these activation functions have trainable parameters, which are learned by gradient descent at the same time as the hyperparameters of the network. To obtain significant benefits for image classification and partial differential equation (PDE) solvers, one can also perform an exhaustive search over trainable activation functions constructed from standard units~\cite{jagtap2020adaptive,ramachandran2017searching}. 

In this chapter, we are motivated by designing activation functions with greater approximation power than ReLU. We study rational neural networks, which are neural networks with activation functions that are trainable rational functions. These will have improved theoretical guarantees on expressivity compared to ReLU, as we shall see in \cref{sec_th_result}.

\section{Definitions}

We consider neural networks whose activation functions consist of rational functions with trainable coefficients $a_i$ and $b_j$, \emph{i.e.}, functions of the form:
\begin{equation} \label{eq_rational}
F(x) = \frac{P(x)}{Q(x)}=\frac{\sum_{i=0}^{r_P} a_ix^i}{\sum_{j=0}^{r_Q} b_jx^j}, \quad a_P\neq 0,\,b_Q\neq 0,
\end{equation}
where $r_P$ and $r_Q$ are the polynomial degrees of the numerator and denominator, respectively.  We say that $F(x)$ is of type $(r_P,r_Q)$ and degree $\max(r_P,r_Q)$.

The use of rational functions in deep learning is motivated by the theoretical work of Telgarsky, who proved error bounds on the approximation of ReLU neural networks by high-degree rational functions and vice versa~\cite{telgarsky2017neural}. On the practical side, neural networks based on rational activation functions are considered by Molina et al.~\cite{molina2019pad}, who defined a safe Pad\'e Activation Unit (PAU) as
\[
F(x) = \frac{\sum_{i=0}^{r_P} a_ix^i}{1+|\sum_{j=1}^{r_Q} b_jx^j|}.
\]
The denominator is selected so that $F(x)$ does not have poles located on the real axis. PAU networks can learn new activation functions and are competitive with state-of-the-art neural networks for image classification. However, this choice results in a non-smooth activation function and makes the gradient expensive to evaluate during training. In a closely related work, Chen et al.~\cite{chen2018rational} propose high-degree rational activation functions in a neural network, which have benefits in terms of approximation power.  However, this choice can significantly increase the number of parameters in the network, causing the training stage to be computationally expensive. 

In this chapter, we use low-degree rational functions as activation functions, which are then composed together by the neural network to build high-degree rational functions. In this way, we can leverage the approximation power of high-degree rational functions without making training expensive. We highlight the approximation power of rational networks and provide optimal error bounds to demonstrate that rational neural networks theoretically outperform ReLU networks. Motivated by our theoretical results, we consider rational activation functions of type $(3,2)$, \emph{i.e.}, $r_P=3$ and $r_Q=2$. This type appears naturally in the theoretical analysis due to the composition property of Zolotarev sign functions (see \cref{sec_approx_relu_rat}): the degree of the overall rational function represented by the rational neural network is an enormous $3^{\#\textup{layers}}$, while the number of trainable parameters only grows linearly with respect to the depth of the network. A low-degree activation function keeps the number of trainable parameters small, while the implicit composition in a neural network gives us the approximation power of high-degree rationals. This choice is also motivated empirically, and we do not claim that the type $(3,2)$ is the best choice for all situations as the configurations may depend on the application as shown later by \cref{fig_rational_loss_2d}.
Our experiments\footnote{All code and hyperparameters are publicly available at~\cite{boulleGitrational}.} on the approximation of smooth functions and generative adversarial networks (GANs) suggest that rational neural networks are an attractive alternative to ReLU networks (see \cref{sec_experiments}).

\section{Theoretical results on rational neural networks} \label{sec_th_result}
Here, we demonstrate the theoretical benefit of using neural networks based on rational activation functions due to their superiority over ReLU in approximating functions. We derive optimal bounds in terms of the total number of trainable parameters (also called size) needed by rational networks to approximate ReLU networks as well as functions in the Sobolev space $\mathcal{W}^{n,\infty}([0,1]^d)$, where $n,d\geq 1$ are integers. Throughout this chapter, we take $\epsilon$ to be a small parameter with $0<\epsilon<1$. We first show that an $\epsilon$-approximation on the domain $[-1,1]^d$ of a ReLU network ($\nrelu$) by a rational neural network ($\nrat$) must have the following size (indicated in brackets):
\begin{equation} \label{eq_rat_approx_relu}
\nrat [\Omega(\log(\log(1/\epsilon)))] \leq \nrelu \leq \nrat [\mathcal{O}(\log(\log(1/\epsilon)))],
\end{equation}
where the constants only depend on the size and depth of the ReLU network. Here, the upper bound means that all ReLU networks can be approximated to within $\epsilon$ by a rational network of size $\mathcal{O}(\log(\log(1/\epsilon)))$. The lower bound means that there is a ReLU network that cannot be $\epsilon$-approximated by a rational network of size less than $C\log(\log(1/\epsilon))$, for some constant $C>0$. In comparison, the size needed by a ReLU network to approximate a rational neural network within the tolerance of $\epsilon$ is given by the following inequalities:
\begin{equation} \label{eq_relu_approx_rat}
\nrelu [\Omega(\log(1/\epsilon))] \leq \nrat \leq \nrelu [\mathcal{O}(\log(1/\epsilon))^3],
\end{equation}
where the constants only depend on the size and depth of the rational neural network. This means that all rational networks can be approximated to within $\epsilon$ by a ReLU network of size $\mathcal{O}(\log(1/\epsilon))^3$, while there is a rational network that cannot be $\epsilon$-approximated by a ReLU network of size less than $\Omega(\log(1/\epsilon))$. A comparison between~\eqref{eq_rat_approx_relu} and~\eqref{eq_relu_approx_rat} suggests that rational networks could be more expressive than ReLU.

\subsection{Approximation of ReLU networks by rational neural networks} \label{sec_approx_relu_rat}

Telgarsky showed that neural networks and rational functions can approximate each other in the sense that there exists a rational function of degree $\mathcal{O}(\polylog(1/\epsilon))$ that is $\epsilon$-close to a ReLU network~\cite[Thm.~1.1]{telgarsky2017neural}, where $\epsilon>0$ is a small number.

\begin{theorem}[Telgarsky] \label{thm_telgarsky}
Let $0<\epsilon <1$ and let $\|\cdot\|_1$ denote the vector 1-norm. The following two statements hold: 
\begin{enumerate}[leftmargin=0cm,itemindent=.5cm,noitemsep]
\item Let $k$ be a nonnegative integer and $p:[0,1]^d\to [-1,1]$, $q:[0,1]^d\to [2^{-k},1]$ be polynomials of degree $\leq r$, each with $\leq s$ monomials. Then, there exists a ReLU network $\nrelu:[0,1]^d\to\mathbb{R}$ of size
\[
\mathcal{O}\left(k^7\log(1/\epsilon)^3+\min\left\{srk\log(sr/\epsilon),sdk^2\log(dsr/\epsilon)^2\right\}\right),
\]
such that 
\[\sup_{x\in[0,1]^d} \left|\nrelu(x)-\frac{p(x)}{q(x)}\right| \leq \epsilon.\]
\item  Let $\nrelu:[-1,1]^d \to \mathbb{R}$ be a ReLU network with $M$ layers and at most $k$ nodes per layer, where each node computes $x\mapsto\relu(a^\top x+b)$ and the pair $(a,b)$ (possibly distinct across nodes) satisfies $\|a\|_1+|b| \leq 1$. Then, there exists a rational function $R:[-1,1]^d \to \mathbb{R}$ with degree (maximum of numerator and denominator)
\[
\mathcal{O}\left(k^M\log(M/\epsilon)^M\right),
\]
such that 
\[\sup_{x\in[-1,1]^d} |\nrelu(x)-R(x)| \leq \epsilon.\]
\end{enumerate}
\end{theorem}

To prove this statement, Telgarsky used a rational function constructed with Newman polynomials~\cite{newman1964rational} to obtain a rational approximation to the ReLU function that converges with square-root exponential accuracy. That is, Telgarsky needed a rational function of degree $\Omega(\log(1/\epsilon)^2)$ to achieve a tolerance of $\epsilon$. A degree $r$ rational function can be represented with $2(r+1)$ coefficients, \emph{i.e.}, $a_0,\ldots,a_r$ and $b_0,\ldots, b_r$ in~\cref{eq_rational}. Therefore, the rational approximation to a ReLU network constructed by Telgarsky requires at least $\Omega(\polylog(1/\epsilon))$ parameters. In contrast, for any rational function, Telgarsky showed that there exists a ReLU network of size $\mathcal{O}(\polylog(1/\epsilon))$ that is an $\epsilon$-approximation on $[0,1]^d$.

Our key observation is that by composing low-degree rational functions together, we can approximate a ReLU network much more efficiently in terms of the size (rather than the degree) of the rational network. Our theoretical work is based on a family of rationals called Zolotarev sign functions, which are the best rational approximation in the infinity norm on $[-1,-\ell]\cup[\ell,1]$, with $0<\ell<1$, to the $\sign$ function~\cite{achieser2013theory,petrushev2011rational}, defined as
\[\sign(x) = 
\begin{cases}
-1, & x<0,\\
0, & x=0,\\
1, & x>0.
\end{cases}\]
We first show that a rational function can approximate the absolute value function $|x|$ on $[-1,1]$ with square-root exponential convergence using a composition of Zolotarev functions.

\begin{lemma}  \label{th_abs_approx}
For any integer $k\geq 0$, we have
\[
\min_{r\in\mathcal{R}_{k,k}} \max_{x\in[-1,1]} \left| |x| - xr(x)\right| \leq 4 e^{-\pi\sqrt{k/2}},
\]
where $\mathcal{R}_{k,k}$ is the space of rational functions of type at most $(k,k)$. Thus, $xr(x)$ is a rational approximant to $|x|$ of type at most $(k+1,k)$. Moreover, if $k=\prod_{i=1}^pk_i$ for some $p\geq 1$ and integers $k_1,\ldots,k_p\geq 2$, then $r$ can be written as $r=R_p\circ\cdots\circ R_1$, where $R_i\in\mathcal{R}_{k_i,k_i}$.
\end{lemma}

\begin{proof}
Let $0<\ell<1$ be a real number and consider the sign function on the domain $[-1,-\ell]\cup[\ell,1]$, \emph{i.e.}, 
\[
{\rm \sign}(x) = 
\begin{cases}
-1, & x\in[-1,-\ell], \\
+1, & x\in [\ell,1].
\end{cases}
\]
By~\cite[Equation~(33)]{beckermann2017singular}, we find that for any $k\geq 0$,
\[
\min_{r\in\mathcal{R}_{k,k}} \max_{x\in[-1,-\ell]\cup[\ell,1]} |\sign(x) - r(x)| \leq 4\left[\exp\left(\frac{\pi^2}{2\log(4/\ell)}\right)\right]^{-k}.
\]
Let $r(x)$ be the rational function of type $(k,k)$ that attains the minimum~\cite[Equation~(12)]{beckermann2017singular}. We refer to such $r(x)$ as the Zolotarev sign function. It is given by
\[
r(x) = Mx \frac{\prod_{j=1}^{\lfloor (k-1)/2\rfloor} x^2+c_{2j}}{\prod_{j=1}^{\lfloor k/2\rfloor} x^2+c_{2j-1}}, \quad c_j = \ell^2\frac{{\rm sn^2}(jK(\kappa)/k ; \kappa)}{1-{\rm sn^2}(jK(\kappa)/k ; \kappa)}.
\]
Here, $M$ is a real constant selected so that ${\rm sign}(x)-r(x)$ equioscillates on $[-1,-\ell]\cup[\ell,1]$, $\kappa = \sqrt{1-\ell^{2}}$, ${\rm sn}(\cdot)$ is the first Jacobian elliptic function, and $K$ is the complete elliptic integral of the first kind. Since $|x| = x\cdot{\rm sign}(x)$ we have the following inequality,
\[
\begin{aligned} 
\max_{x\in[-1,-\ell]\cup[\ell,1]}\left| |x| - xr(x)\right| &= \max_{x\in[-1,-\ell]\cup[\ell,1]}\left| x\cdot\sign(x) - xr(x)\right|\\
&\leq \max_{x\in[-1,-\ell]\cup[\ell,1]}\left| \sign(x) - r(x)\right|.
\end{aligned} 
\]
The last inequality follows because $|x|\leq 1$ on $[-1,-\ell]\cup[\ell,1]$. Moreover, since $xr(x)\geq 0$ for $x\in[-1,1]$ (see~\cite[Equation~(12)]{beckermann2017singular}) we have
\[
\max_{x\in[-\ell,\ell]}\left| |x| - xr(x)\right| \leq \max_{x\in[-\ell,\ell]}\left| x\right|\leq \ell. 
\]
Therefore, 
\[
\max_{x\in[-1,1]}\left| |x| - xr(x)\right|\leq \max\left\{\ell, 4\left[\exp\left(\frac{\pi^2}{2\log(4/\ell)}\right)\right]^{-k}\right\}. 
\]
Now, we select $0<\ell<1$ to minimize this upper bound. One finds that $\ell = 4\exp(-\pi\sqrt{k/2})$ and the result follows immediately. 

For the final claim, let $r$ be the Zolotarev sign function $Z_{k}(\cdot\,;\ell)$ of type $(k,k)$ on $[-1,-\ell]\cup[\ell,1]$, with $k=\prod_{i=1}^p k_i$. By definition, $Z_{k}(\cdot;\ell)$ is the best rational approximation of degree $k$ to the $\sign$ function on $[-1,-\ell]\cup[\ell,1]$. We know from~\cite{lebedev1977zolotarev,nakatsukasa2016computing} that there exist $p$ Zolotarev sign functions $R_1,\ldots,R_p$, where each $R_i$ is of type $(k_i,k_i)$, such that
\begin{equation} \label{eq_composition_Zolotarev}
r(x) \coloneqq Z_{k}(x;\ell) = R_p(\cdots(R_2(R_1(x)))\cdots).
\end{equation}
\end{proof}

A composition of $k\geq 1$ Zolotarev sign functions of type $(3,2)$ has type $(3^k,3^k-1)$ but can be represented with $7k$ parameters instead of $2\times 3^k+1$. This property enables the construction of a rational approximation to ReLU using compositions of low-degree Zolotarev sign functions with $\mathcal{O}(\log(\log(1/\epsilon)))$ parameters in \cref{lem_relu_rat}. The proof of \cref{lem_relu_rat} is a direct consequence of the previous lemma and the properties of Zolotarev sign functions.

\begin{lemma} \label{lem_relu_rat}
Let $0<\epsilon <1$. There exists a rational network $\nrat:[-1,1]\to[-1,1]$ of size $\mathcal{O}(\log(\log(1/\epsilon)))$ such that 
\[\|\nrat-\relu\|_{\infty}\coloneqq\max_{x\in[-1,1]}|\nrat(x)-\relu(x)| \leq \epsilon.\] 
Moreover, no rational network of size smaller than $\Omega(\log(\log(1/\epsilon)))$ can achieve this.
\end{lemma}

\begin{proof}
Let $0<\epsilon<1$, $0<\ell<1$, $k\geq 1$, and $r$ be the Zolotarev sign function $Z_{3^k}(\cdot\,;\ell)$ of type $(3^k,3^k-1)$. Again from~\cite{lebedev1977zolotarev,nakatsukasa2016computing}, we see that there exist $k$ Zolotarev sign functions $R_1,\ldots,R_k$ of type $(3,2)$ such that their composition equals $Z_{3^k}(x;\ell)$, \emph{i.e.},
\begin{equation} \label{eq_composition_Zolotarev_2}
r(x) \coloneqq Z_{3^k}(x;\ell) = R_k(\cdots (R_2(R_1(x))\cdots).
\end{equation}
Following the proof of \cref{th_abs_approx}, we have the inequality
\begin{equation} \label{eq_proof_relu_rat}
\max_{x\in[-1,1]} \left| |x| - xr(x)\right| \leq 4 e^{-\pi\sqrt{3^k/2}},
\end{equation}
where we chose $\ell=4\exp(-\pi\sqrt{3^k/2})$. Now, we take
\begin{equation} \label{eq_proof_relu_rat_k}
k=\left\lceil \frac{\ln(2/\pi^2)+2\ln(\ln(4/\epsilon))}{\ln(3)} \right\rceil,
\end{equation}
so that the right-hand side of~\cref{eq_proof_relu_rat} is bounded by $\epsilon$. Finally, we use the identity
\[
\relu(x) = \frac{|x|+x}{2},\quad x\in\mathbb{R},
\]
to define a rational approximation to the ReLU function on the interval $[-1,1]$ as
\[
\tilde{r}(x)=\frac{1}{2}\left(\frac{xr(x)}{1+\epsilon}+x\right).
\]
Therefore, we have the following inequalities for $x\in[-1,1]$,
\begin{align*}
|\relu(x)-\tilde{r}(x)| &= \frac{1}{2}\left||x|-\frac{xr(x)}{1+\epsilon}\right|\leq \frac{1}{2(1+\epsilon)}(||x|-xr(x)|+\epsilon|x|)\\
&\leq \frac{\epsilon}{1+\epsilon}\leq \epsilon.
\end{align*}
Then, $r$ is a composition of $k$ rational functions of type $(3,2)$ and can be represented using at most $7k$ coefficients (see~\cref{eq_composition_Zolotarev}). Moreover, using~\cref{eq_proof_relu_rat_k}, we see that $k=\mathcal{O}(\log(\log(1/\epsilon)))$, which means that $\tilde{r}$ is representable by a rational network of size $\mathcal{O}(\log(\log(1/\epsilon)))$. Finally, $|\tilde{r}(x)|\leq 1$ for $x\in[-1,1]$.

The lower bound on the rational networks size will be proved separately later in \cref{prop_optimal_bound_rat_relu}.
\end{proof}

The upper bound on the complexity of the neural network obtained in \cref{lem_relu_rat} is optimal, as proved by Vyacheslavov~\cite{vyacheslavov1975uniform}.

\begin{theorem}[Vyacheslavov] \label{th_approx_abs}
The following inequalities hold:
\begin{equation}
C_1e^{-\pi\sqrt{k}}\leq \max_{x\in[-1,1]}||x|-r_k(x)| \leq C_2 e^{-\pi\sqrt{k}}, \quad k\geq 0,
\end{equation}
where $r_k$ is the best rational approximation to $|x|$ in $[-1,1]$ from $\mathcal{R}_{k,k}$. Here,  $C_1,C_2>0$ are constants that are independent of $k$.
\end{theorem}

We first deduce the following corollary, giving lower and upper bounds on the optimal rational approximation to the ReLU function.

\begin{corollary} \label{cor_best_relu}
The following inequalities hold: 
\begin{equation}
\frac{C_1}{2}e^{-\pi\sqrt{k}}\leq \|\relu-r_k\|_{\infty} \leq \frac{C_2}{2} e^{-\pi\sqrt{k}}, \quad k\geq 0,
\end{equation}
where $r_k$ is the best rational approximation to ReLU on $[-1,1]$ in $\mathcal{R}_{k,k}$ and $C_1,C_2>0$ are constants given by \cref{th_approx_abs}.
\end{corollary}

\begin{proof}
Let $k$ be an integer and let $r_k\in\mathcal{R}_{k,k}$ be any rational function of degree $\leq k$. Now, define $r_\abs(x)=2r_k(x)-x$. Since $\relu(x)=(|x|+x)/2$, we have
\begin{align*}
\|\relu-r_k\|_{\infty} &= \max_{x\in[-1,1]} \left|\frac{1}{2}(r_\abs(x)+x)-\frac{1}{2}(|x|+x)\right| = \max_{x\in[-1,1]} \frac{1}{2}\left|r_\abs(x)-|x|\right|\\
&\geq \frac{1}{2}C_1e^{-\pi\sqrt{k}},
\end{align*}
where the inequality is from~\cref{th_approx_abs}.  Now, let $r_k\in\mathcal{R}_{k,k}$ be the best rational approximation to $|x|$ on $[-1,1]$. Now, define $r_\relu(x) = (r_k(x)+x)/2$. We find that
\begin{align*}
\|\relu-r_\relu\|_{\infty} &= \max_{x\in[-1,1]}  \left| \frac{1}{2}(|x|+x) - \frac{1}{2}(r_k(x)+x)\right|=  \max_{x\in[-1,1]}\frac{1}{2} \left||x|-r_k(x)\right|\\
&\leq \frac{1}{2}C_2e^{-\pi\sqrt{k}},
\end{align*}
which proves that the best approximation to ReLU satisfies the upper bound.
\end{proof}

We now show that a rational neural network must be at least $\Omega(\log(\log(1/\epsilon)))$ in size (total number of nodes) to approximate the ReLU function to within $\epsilon$. 

\begin{proposition} \label{prop_optimal_bound_rat_relu}
Let $0<\epsilon<1$. A rational neural network that approximates the ReLU function on $[-1,1]$ to within $\epsilon$ has size of at least $\Omega(\log(\log(1/\epsilon)))$. 
\end{proposition}
\begin{proof}
Let $\nrat: [-1,1]\rightarrow \mathbb{R}$ be a rational neural network with $k_1,\ldots,k_M\geq 1$ nodes at each of its $M$ layers, and assume that its activation functions are rational functions of type at most $(r_P,r_Q)$. Let $d_r=\max(r_P,r_Q)$ be the maximum of the degrees of the activation functions of $\nrat$. Such a network has size $\sum_{i=1}^M k_i$. Note that $\nrat$ itself is a rational function of degree $d$, where from additions and compositions of rational functions we have $d \leq d_r^M\prod_{i=1}^M k_i$.  If $\nrat$ is an $\epsilon$-approximation to the ReLU function on $[-1,1]$, we know by~\cref{cor_best_relu} that 
\begin{equation} \label{eq_degree_relu}
\frac{C_1}{2}e^{-\pi\sqrt{d}}\geq\epsilon, \quad d\geq \left(\frac{1}{\pi}\ln\left(\frac{C_1}{2\epsilon}\right)\right)^2.
\end{equation}
The statement follows by minimizing the size of $\nrat$, \emph{i.e.}, $\sum_{i=1}^M k_i$ subject to 
\[d_r^M\prod_{i=1}^M k_i\geq  \left(\frac{1}{\pi}\ln\left(\frac{C_1}{2\epsilon}\right)\right)^2.
\]
That is,
\begin{equation} \label{eq_KKT}
\sum_{i=1}^M\ln(k_i)+M\ln(d_r) \geq 2\ln\left(\ln\left(\frac{C_1}{2\epsilon}\right)\right)-2\ln(\pi).
\end{equation}
We introduce a Lagrange multiplier $\lambda\in\mathbb{R}$ and define the Lagrangian of this optimization problem as
\[
\mathcal{L}(k_1,\ldots,k_M,\lambda)=\sum_{i=1}^M k_i+\lambda\left[2\ln\left(\ln\left(\frac{C_1}{2\epsilon}\right)\right)-2\ln(\pi) - \sum_{i=1}^M\ln(k_i) - M\ln(d_r)\right].
\]
One finds using the Karush--Kuhn--Tucker conditions~\cite{kuhn1951} that $k_1=\cdots=k_M=\lambda$. Then, using \cref{eq_KKT}, we find that $\lambda$ satisfies
\begin{equation} \label{eq_lambda}
\ln(\lambda) \geq \frac{2}{M}\left[\ln\left(\ln\left(\frac{C_1}{2\epsilon}\right)\right)-\ln(\pi)\right] - \ln(d_r) =: \ln(\lambda^*).
\end{equation}
Therefore, the rational network $\nrat$ with $M$ layers that approximates the ReLU function to within $\epsilon$ on $[-1,1]$ has a size of at least $s(M)\coloneqq M\lambda^*$, where $\lambda^*$ is given by~\cref{eq_lambda} and depends on $M$. We now minimize $s(M)$ with respect to the number of layers $M\geq 1$. We remark that minimizing $s$ is equivalent of minimizing $\ln(s)$, where
\[
\ln(s(M)) = \ln(M) + \ln(\lambda^*) = \ln(M)+\frac{2}{M}\left[\ln\left(\ln\left(\frac{C_1}{2\epsilon}\right)\right)-\ln(\pi)\right] - \ln(d_r).
\]
One finds that one should take $k_1 =\cdots=k_M=\lambda^*= \mathcal{O}(1)$ and $M =\Omega(\log(\log(1/\epsilon)))$. The result follows. 
 \end{proof}

\begin{figure}[htbp]
\centering
\vspace{0.2cm}
\begin{overpic}[width=0.9\textwidth]{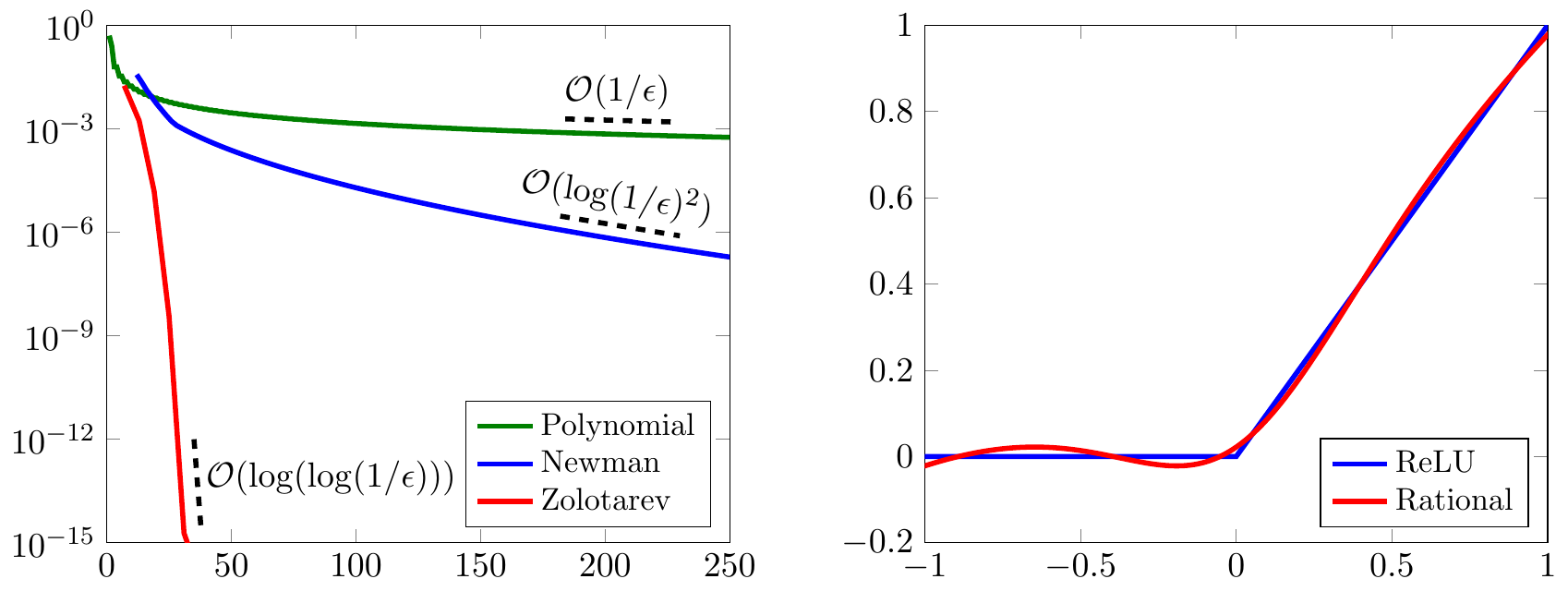}
\put(14,-2.5){Number of parameters}
\put(-3,11){\rotatebox{90}{$\|\relu-r_N\|_{\infty}$}}
\put(78.3,-2.5){$x$}
\end{overpic}
\vspace{0.3cm}
\caption{Left: Approximation error $\|\relu-r_N\|_{\infty}$ of the Newman (blue), Zolotarev sign functions (red), and best polynomial approximation~\cite{pachon2009barycentric} of degree $N-1$ (green) $r_N$ to ReLU with respect to the number of parameters required to represent $r_N$. Right: Best rational function of type $(3,2)$ (red) that approximates the ReLU function (blue). We use this to initialize the rational activation functions when training a rational neural network.}
\label{fig_init_rat}
\end{figure}

The proof of \cref{prop_optimal_bound_rat_relu} shows that the bound obtained in \cref{lem_relu_rat} is optimal in the sense that a rational network requires at least $\Omega(\log(\log(1/\epsilon)))$ parameters to approximate the ReLU function on $[-1,1]$ to within the tolerance $\epsilon>0$. The convergence of the Zolotarev sign functions to the ReLU function is much faster, with respect to the number of parameters, than the rational constructed with Newman polynomials (see~\cref{fig_init_rat}(left)). We also include in this panel the algebraic convergence of $\mathcal{O}(1/\epsilon)$ obtained by polynomials~\cite{trefethen2019approximation} as a comparison.

The converse of \cref{lem_relu_rat}, which is a consequence of a theorem proved by Telgarsky~\cite[Theorem~1.1]{telgarsky2017neural}, shows that any rational function can be approximated by a ReLU network of size at most $\mathcal{O}(\log(1/\epsilon)^3)$.

\begin{lemma} \label{lem_relu_approx}
Let $0<\epsilon<1$. If $R:[-1,1]\to[-1,1]$ is a rational function, then there exists a ReLU network $\nrelu:[-1,1]\to[-1,1]$ of size $\mathcal{O}(\log(1/\epsilon)^3)$ such that $\|R-\nrelu\|_{\infty}\leq \epsilon$.
\end{lemma}

\begin{proof}
Let $0<\epsilon<1$ and $R:[-1,1]\to[-1,1]$ be a rational function. Take $\tilde{R}(x) = R( 2x-1)$, which is still a rational function. Without loss of generality, we can assume that $\tilde{R}$ is an irreducible rational function (otherwise cancel factors till it is irreducible). Since $\tilde{R}$ is a rational, it can be written as $\tilde{R}=p/q$ with $\max_{x\in[0,1]}|q(x)|=1$. Moreover, we know that $\tilde{R}(x)\in[-1,1]$ for $x\in[0,1]$ so we can assume that $q(x)\geq 0$ for $x\in [0,1]$ (it is either positive or negative by continuity). Since $R$ is continuous on $[-1,1]$, there is an integer $n\geq 1$ such that $q(x)\in[2^{-n},1]$ for $x\in[0,1]$. Furthermore, we find that $|p(x)|\leq 1$ for $x\in[0,1]$ because $|R(x)|\leq 1$ and $|q(x)|\leq 1$ for $x\in[0,1]$. By~\cite[Theorem~1.1]{telgarsky2017neural}, there exists a ReLU network $\nrelu:[0,1]\to \mathbb{R}$ of size $\mathcal{O}(n^7\log(1/\epsilon)^3)$ such that
\[
\max_{x\in[0,1]}\left|\nrelu(x)-\frac{p(x)}{q(x)}\right|\leq \frac{\epsilon}{2}.
\]
We now define a scaled ReLU network $\tilde{\N}_{\relu}(x)=\nrelu(x)/(1+\epsilon/2)$ such that $|\tilde{\N}_{\relu}(x)|\leq 1$ for $x\in[0,1]$.
Therefore, for all $x\in[0,1]$,
\[
\left|\tilde{\N}_{\relu}(x)-\tilde{R}(x)\right| = \left|\frac{\nrelu(x)}{1+\epsilon/2}-\frac{p(x)}{q(x)}\right|
\leq \frac{1}{1+\epsilon/2}\left(\left|\nrelu(x)-\frac{p(x)}{q(x)}\right|+\frac{\epsilon}{2}\left|\frac{p(x)}{q(x)}\right|\right)
\leq \epsilon.
\]
Therefore, $x\mapsto \tilde{\N}_{\relu}((x+1)/2)$ is a ReLU neural network of size $\mathcal{O}(\log(1/\epsilon)^3)$ that is an $\epsilon$-approximation to $R$ on $[-1,1]$.
\end{proof}

To demonstrate the improved approximation power of rational neural networks over ReLU networks ($\mathcal{O}(\log(\log(1/\epsilon)))$ versus $\mathcal{O}(\log(1/\epsilon)^3)$), it is known that a ReLU network that approximates $x^2$, which is rational, to within $\epsilon$ on $[-1,1]$ must be of size at least $\Omega(\log(1/\epsilon))$~\cite[Theorem~11]{liang2016deep}. 

We can now state our main theorem based on~\cref{lem_relu_rat,lem_relu_approx}. \cref{th_rat_network} provides bounds on the approximation power of ReLU networks by rational neural networks and vice versa. We regard~\cref{th_rat_network} as an analogue of~\cite[Thm.~1.1]{telgarsky2017neural} for our Zolotarev sign functions, where we are counting the number of training parameters instead of the degree of the rational functions. In particular, our rational networks have high degrees but can be represented with few parameters due to compositions, making training more computationally efficient. While Telgarsky required a rational function with $\mathcal{O}(k^M\log(M/\epsilon)^M)$ parameters to approximate a ReLU network with fewer than $k$ nodes in each of $M$ layers to within a tolerance of $\epsilon$, we construct a rational network that only has size $\mathcal{O}(kM\log(\log(M/\epsilon)))$.

\begin{theorem} \label{th_rat_network}
Let $0<\epsilon <1$ and let $\|\cdot\|_1$ denote the vector 1-norm. The following two statements hold: 
\begin{enumerate}[leftmargin=0cm,itemindent=.5cm,noitemsep]
\item Let $\nrat:[-1,1]^d \to [-1,1]$ be a rational network with $M$ layers and at most $k$ nodes per layer, where each node computes $x\mapsto r(a^\top x+b)$ and $r$ is a rational function with Lipschitz constant $L$ ($a$, $b$, and $r$ are possibly distinct across nodes). Suppose further that $\|a\|_1+|b|\leq 1$ and $r:[-1,1]\rightarrow[-1,1]$. Then, there exists a ReLU network $\nrelu:[-1,1]^d\to[-1,1]$ of size
\[
\mathcal{O}\left(kM\log(ML^M/\epsilon)^3\right)
\]
such that $\max_{x\in[-1,1]^d} |\nrat(x)-\nrelu(x)| \leq \epsilon$.  
\item  Let $\nrelu:[-1,1]^d \to [-1,1]$ be a ReLU network with $M$ layers and at most $k$ nodes per layer, where each node computes $x\mapsto\relu(a^\top x+b)$ and the pair $(a,b)$ (possibly distinct across nodes) satisfies $\|a\|_1+|b| \leq 1$. Then, there exists a rational network $\nrat:[-1,1]^d \to [-1,1]$ of size
\[
\mathcal{O}(kM\log(\log(M/\epsilon)))
\]
such that $\max_{x\in[-1,1]^d} |\nrelu(x)-\nrat(x)| \leq \epsilon$.\end{enumerate}
\end{theorem}

\begin{proof} \leavevmode
The statement of~\cref{th_rat_network} comes in two parts, and we prove them separately.  The structure of the proof closely follows~\cite[Lemma~1.3]{telgarsky2017neural}.

1.~~Consider the subnetwork $H$ of the rational network $\nrat$, consisting of the layers of $\nrat$ up to the $J$th layer for some $1\leq J\leq M-1$. Let $H_\relu$ denote the ReLU network obtained by replacing each rational function $r_{ij}$ in $H$ by a ReLU network approximation $f_{r_{ij}}$ at a given tolerance $\epsilon_j>0$ for $1\leq j\leq J$ and $1\leq i\leq k_j$, such that $|H_\relu(x)|\leq 1$ for $x\in[-1,1]$ (see Lemma~\ref{lem_relu_approx}). Let $x\mapsto r_{i,J+1}(a_{i,J+1}^\top H(x)+b_{i,J+1})$ be the output of the rational network $\nrat$ at layer $J+1$ and node $i$ for $1\leq i\leq k_J$. Now, approximate node $i$ in the $(J+1)$st layer by a ReLU network $f_{r_i,J+1}$ with tolerance $\epsilon_{J+1}>0$ (see Lemma~\ref{lem_relu_approx}). The approximation error $E_{i,J+1}$ between the rational and the approximating ReLU network at layer $J+1$ and node $i$ satisfies
\begin{align*}
E_{i,J+1} &= |f_{r_{i,J+1}}(a_{i,J+1}^\top H_\relu(x)+b_{i,J+1}) - r_{i,J+1}(a_{i,J+1}^\top H(x)+b_{i,J+1})| \\
& \leq \underbrace{|f_{r_{i,J+1}}(a_{i,J+1}^\top H_\relu(x)+b_{i,J+1}) - r_{i,J+1}(a_{i,J+1}^\top H_\relu(x)+b_{i,J+1})|}_{(1)}\\
& + \underbrace{|r_{i,J+1}(a_{i,J+1}^\top H_\relu(x)+b_{i,J+1}) - r_{i,J+1}(a_{i,J+1}^\top H(x)+b_{i,J+1})|}_{(2)}. 
\end{align*}
The first term is bounded by 
\[
(1) \leq \max_{x\in[-1,1]} \left| r_{i,J+1}(x)-f_{r_{i,J+1}}\right| \leq \epsilon_{J+1},
\]
since $\left|a_{i,J+1}^\top H_\relu(x)+b_{i,J+1}\right| \leq \|a_{i,J+1}\|_1+|b_{i,J+1}|\leq 1$ by assumption. The second term is bounded as the Lipschitz constant of $r_{i,J+1}$ is at most $L$. That is,  
\[
(2) \leq L\|a_{i,J+1}\|_1\max_{x\in[-1,1]^d} \left\|H_\relu(x)-H(x)\right\|_\infty \leq  L\max_{x\in[-1,1]^d} \left\|H_\relu(x)-H(x)\right\|_\infty,
\]
where we used the fact that $\|a_{i,J+1}\|_1\leq 1$ and $\|H_\relu(x)\|_\infty \leq 1$ for $x\in[-1,1]^d$. We find that we have the following set of inequalities:  
\[\max_{1\leq i\leq k_{j+1}} E_{i,j+1}\leq L \max_{1\leq i\leq k_{j}} E_{i,j}+\epsilon_{j+1}, \quad 1\leq i\leq k_j, \quad 1\leq j\leq J+1,\] 
with $E_{i,0} = 0$.  If we select $\epsilon_j=\epsilon L^{j-J-1}/(J+1)$, then we find that $\max_{1\leq i\leq k_{J+1}} E_{i,J+1} \!\leq\! \epsilon$. When $J = M-1$, the ReLU network approximates the original rational network, $\nrat$, and the ReLU network has size
\[
\mathcal{O}\left(k \sum_{j=1}^M\log\left(\frac{M}{L^{j-M}\epsilon}\right)^3\right).
\]
where we used the fact that $k_j\leq k$ for $1\leq j\leq M$. This can be simplified a little since 
\[
\sum_{j=1}^M \log\left(\frac{M}{L^{j-M}\epsilon}\right)^3=\sum_{j=1}^M\left(\log(ML^M/\epsilon)+j\log(1/L)\right)^3=\mathcal{O}\!\left(M\log(ML^M/\epsilon)^3\right).
\]

2.~~Telgarsky proved in~\cite[Lemma~1.3]{telgarsky2017neural} that if $H_R$ is a neural network obtained by replacing all the ReLU activation functions in $\nrelu$ by rational functions $R$ for $1\leq j\leq M$, which satisfies $R(x)\in[-1,1]$ and $|R(x)-\relu(x)|\leq \epsilon/M$ for $x\in[-1,1]$, then 
\[
\max_{x\in[-1,1]^d}|\nrelu(x)-H_R(x)|\leq \epsilon.
\] 
Let $\tilde{R}$ be a rational neural network approximating ReLU with a tolerance of $\epsilon/M$, constructed by~\cref{lem_relu_rat}. Then, $\tilde{R}$ is rational network of size $\mathcal{O}(\log(\log(M/\epsilon)))$ and thus, $H_{\tilde{R}}$ is a rational neural network of size $\mathcal{O}(Mk\log(\log(M/\epsilon)))$.
\end{proof}

\Cref{th_rat_network} highlights the improved approximation power of rational neural networks over ReLU networks. ReLU networks of size $\mathcal{O}(\polylog(1/\epsilon))$ are required to approximate rational networks while rational networks of size only $\mathcal{O}(\log(\log(1/\epsilon)))$ are sufficient to approximate ReLU networks.

\subsection{Approximation of functions by rational networks} \label{sec_func_rat}

A important question is the required size and depth of deep neural networks to approximate smooth functions~\cite{liang2016deep,montanelli2019deep,yarotsky2017error}. In this section, we consider the approximation theory of rational networks. In particular, we consider the approximation of functions in the Sobolev space $\mathcal{W}^{n,\infty}([0,1]^d)$, where $n\geq 1$ is the regularity of the functions and $d\geq 1$. The norm of a function $f\in\mathcal{W}^{n,\infty}([0,1]^d)$ is defined as
\[
\|f\|_{\mathcal{W}^{n,\infty}([0,1]^d)}=\max_{|\mathbf{n}|\leq n} \esssup_{\mathbf{x}\in[0,1]^d}|D^{\mathbf{n}}f(\mathbf{x})|,
\]
where $\mathbf{n}$ is the multi-index $\mathbf{n}=(n_1,\ldots,n_d)\in\{0,\ldots,n\}^d$, and $D^{\mathbf{n}}f$ is the corresponding weak derivative of $f$. In this section, we consider the approximation of functions from
\[
F_{d,n} \coloneqq \{f\in \mathcal{W}^{n,\infty}([0,1]^d),\quad \|f\|_{\mathcal{W}^{n,\infty}([0,1]^d)}\leq 1\}.
\]
By the Sobolev embedding theorem~\cite{brezis2010functional}, $F_{d,n}$ contains the functions in $\mathcal{C}^{n-1}([0,1]^d)$, which is the class of functions whose first $n-1$ derivatives are Lipschitz continuous. 
Yarotsky derived upper bounds on the size of neural networks with piecewise linear activation functions needed to approximate functions in $F_{d,n}$~\cite[Thm.~1]{yarotsky2017error}. In particular, Yarotsky constructed an $\epsilon$-approximation to functions in $F_{d,n}$ with a ReLU network of size at most $\mathcal{O}(\epsilon^{-d/n}\log(1/\epsilon))$ and depth smaller than $\mathcal{O}(\log(1/\epsilon))$.

\begin{theorem}[Yarotsky] \label{th_yarotsky}
Let $d\geq 1$, $n\geq 1$, $0<\epsilon<1$, and $f\in F_{d,n}$. There exists a ReLU neural network $\nrelu$ of size
\[\mathcal{O}(\epsilon^{-d/n}\log(1/\epsilon))\]
and maximum depth $\mathcal{O}(\log(1/\epsilon))$ such that $\|f-\nrelu\|_{\infty}\leq \epsilon$.
\end{theorem}

The term $\epsilon^{-d/n}$ in \cref{th_yarotsky} is introduced by a local Taylor approximation, while the $\log(1/\epsilon)$ term is the size of the ReLU network needed to approximate monomials, \emph{i.e.}, $x^j$ for $j\geq 0$, in the Taylor series expansion. We now present an analogue of \cref{th_yarotsky} for a rational neural network.

\begin{theorem} \label{th_approx_smooth_rat}
Let $d\geq 1$, $n\geq 1$, $0<\epsilon<1$, and $f\in F_{d,n}$. There exists a rational neural network $\nrat$ of size
\[\mathcal{O}(\epsilon^{-d/n}\log(\log(1/\epsilon)))\]
and maximum depth $\mathcal{O}(\log(\log(1/\epsilon)))$ such that $\|f-\nrat\|_{\infty}\leq \epsilon$.
\end{theorem}

The proof of \cref{th_approx_smooth_rat} consists of approximating $f$ by a local Taylor expansion. One needs to approximate the piecewise linear functions and monomials arising in the Taylor expansion by rational networks. The main distinction between Yarotsky's argument and the proof of~\cref{th_approx_smooth_rat} is that monomials can be represented by rational neural networks with a size that does not depend on the accuracy of $\epsilon$. In contrast, ReLU networks require $\mathcal{O}(\log(1/\epsilon))$ parameters. Meanwhile, while ReLU neural networks can exactly approximate piecewise linear functions with a constant number of parameters, rational networks can approximate them with a size of a most $\mathcal{O}(\log(\log(1/\epsilon)))$ (see~\cref{lem_relu_rat}). That is, rational neural networks approximate piecewise linear functions much faster than ReLU networks approximate polynomials. This allows the existence of a rational network approximation to $f$ with exponentially smaller depth ($\mathcal{O}(\log(\log(1/\epsilon)))$) than the ReLU networks constructed by Yarotsky.

We first show that the construction in~\cref{lem_relu_rat} can approximate any piecewise linear function on $[-1,1]$.

\begin{proposition} \label{th_approx_piecewise_rat}
Let $0<\epsilon<1$ and let $g:[0,1]\rightarrow\mathbb{R}$ be any continuous piecewise linear function with $m\geq 1$ breakpoints and Lipschitz constant $L>0$. Then, there exists a rational neural network $\nrat:[0,1]\rightarrow \mathbb{R}$ of size at most 
\[
\mathcal{O}(m\log(\log(L/\epsilon)))
\]
such that $\max_{x\in[0,1]} |g(x)-\nrat(x)|\leq \epsilon$.
\end{proposition}

\begin{proof}
Let $0\leq b_1<\cdots<b_M\leq 1$ be the breakpoints of $g$. In a similar way to the proof of~\cite[Proposition 1]{yarotsky2017error}, we first express $\rho$ as the following sum:
\begin{equation} \label{eq_proof_piecewise_rat}
g(x)=c_0\relu(b_1-x)+\sum_{j=1}^mc_j\relu(x-b_j) + c_{m+1},
\end{equation}
for some constants $c_0,\ldots,c_{m+1}\in\mathbb{R}$. Therefore, $g$ can be exactly represented using a ReLU network with $m+1$ nodes and one layer, \emph{i.e.}, 
\[
g(x) = \begin{pmatrix}
c_0 & c_1& \cdots & c_m
\end{pmatrix}
\begin{pmatrix}
\relu(-x+b_1)\\
\relu(x-b_1)\\
\vdots\\
\relu(x-b_m)
\end{pmatrix} + c_{m+1}.
%x\mapsto
%\relu
%\begin{pmatrix}
%a_1-x\\
%x-a_1\\
%\vdots\\
%x-a_M
%\end{pmatrix}
%\to
%\begin{pmatrix}
%c_0 & \cdots & c_M
%\end{pmatrix}
%\begin{pmatrix}
%\sigma(a_1-x)\\
%\sigma(x-a_1)\\
%\vdots\\
%\sigma(x-a_M)
%\end{pmatrix}
%+h
%=\rho(x).
\]
Since $g$ has a Lipschitz constant of $L$, we find that $|c_0|\leq L$ and $\sum_{j=1}^m |c_j| \leq L$. Using~\cref{lem_relu_rat} we can approximate a ReLU function on $[-1,1]$ with tolerance $\epsilon/(2L)$ by a rational network $R_{\relu}$ of size $\mathcal{O}(\log(\log(2L/\epsilon)))$. Now, we construct $\nrat:[0,1]\rightarrow \mathbb{R}$ as a rational network obtained by replacing the ReLU functions in $g$ by $R_{\relu}$. We have the following error estimate: 
\[
\max_{x\in[0,1]} |g(x)-\nrat(x)| \leq |c_0|\|\relu-R_{\relu}\|_{\infty}+\sum_{j=1}^m |c_j| \|\relu-R_{\relu}\|_{\infty} \leq \frac{\epsilon}{2}+\frac{\epsilon}{2}\leq \epsilon.
\]
The result follows as $\nrat$ is of size $\mathcal{O}(m\log(\log(L/\epsilon)))$.
\end{proof}

We remark that the size of the rational network required to approximate a piecewise linear function depends on $\epsilon$. In contrast, ReLU neural networks can represent piecewise linear functions exactly. In the next proposition, we show that a rational neural network can represent $x^n$, for some integer $n$, exactly.

\begin{proposition} \label{prop_rat_x_n}
Let $n\geq 1$, $r_P\geq 2$, and $r_Q\geq 0$.  There exists a rational network $\nrat$, with rational activation functions of type $(r_P,r_Q)$, of size at most $5\lfloor\log_{r_P}(n)\rfloor^2+1$ such that $\nrat(x) = x^n$ for all $x\in\mathbb{R}$.
\end{proposition}
\begin{proof}
We start by expressing $n$ in base $r_P$, \emph{i.e.}, 
\[
n=\sum_{\ell=0}^{\lfloor\log_{r_P}\!(n)\rfloor}c_\ell r_P^\ell, \qquad c_\ell\in\{0,1,\ldots,r_P-1\}.
\] 
This means we can represent $x^n$ as
\begin{equation} \label{eq_epr_x_n}
x^n=\prod_{\ell=0}^{\lfloor\log_{r_P}\!(n)\rfloor} x^{c_\ell r_P^\ell}.
\end{equation}
Note that $x^{c_\ell r_P^\ell}$ is just $x^{r_P}$ composed $\ell$ times as well as composed with $x^{c_\ell}$ so can be represented by a rational neural network with $\ell+1$ layers, each with one node. Therefore, all the $x^{c_\ell r_P^\ell}$ terms can be represented in rational networks that in total have size
\[
\sum_{\ell=0}^{\lfloor\log_{r_P}\!(n)\rfloor} \!\!(\ell+1) = \frac{1}{2}(\lfloor\log_{r_P}\!(n)\rfloor)^2 + \frac{3}{2}\lfloor\log_{r_P}\!(n)\rfloor + 1.
\]
The function $x^n$ can be formed by multiplying all the $x^{c_\ell r_P^\ell}$ terms together. Since $xy = (x^2 + y^2 - (x-y)^2)/2$, there is a rational network with one layer and three nodes that represents the multiplication operation. Therefore, multiplying all the terms together requires a rational network of size at most $3\lfloor\log_{r_P}\!(n)\rfloor$ (see \cref{eq_epr_x_n}). The result follows by noting that $x^2/2+9x/2+1\leq 5x^2+1$ for $x\geq 1$. 
\end{proof}

\begin{figure}[htbp]
\centering
\vspace{0.2cm}
\begin{overpic}[width=0.4\textwidth]{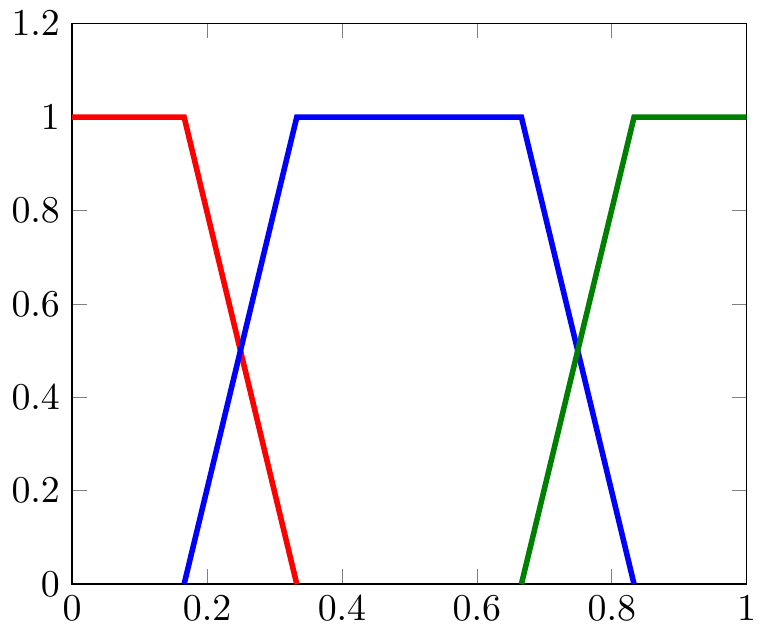}
\put(52,-5){$x$}
\vspace{0.2cm}
\end{overpic}
\caption{Partition of unity: $\psi_0$ (red), $\psi_1$ (blue), and $\psi_2$ (green), for $N=2$.}
\label{fig_partition}
\end{figure}

We can now prove \cref{th_approx_smooth_rat} using the two previous propositions.

\begin{proof}[Proof of \cref{th_approx_smooth_rat}]
The proof is based on the proof of~\cite[Theorem~1]{yarotsky2017error} and consists of replacing the piecewise linear functions and monomials arising in the local Taylor approximation of the function $f$ by rational networks using the previous approximation results.

Let $N\geq 1$ be an integer and consider a partition of unity of $(N+1)^d$ functions $\phi_\m$ on the domain $[0,1]^d$, \emph{i.e.},
\[\sum_{\m\in\{0,\ldots,N\}^d} \phi_\m(\x)=1, \quad \phi_\m(\x)=\prod_{k=1}^d\psi_{m_k}(x_k), \quad \x=(x_1,\ldots,x_d),\]
where $\m=(m_1,\ldots,m_d)$, and $\psi_{m_k}$ is given by
\[\psi_{m_k}(x)=
\begin{cases}
1, \quad & \text{if }\left|x_k-\frac{m_k}{N}\right|<\frac{1}{3N},\\
0, \quad & \text{if }\left|x_k-\frac{m_k}{N}\right|>\frac{2}{3N},\\
2-3N\left|x_k-\frac{m_k}{N}\right|, \quad & \text{otherwise}.
\end{cases}
\]
Examples of the functions $\psi_{m_k}$ are shown in~\cref{fig_partition} when $N=2$. We now define a local Taylor approximation of $f$ by 
\[
f_N(\x) =\sum_{\m\in\{0,\ldots,N\}^d}\phi_\m(\x) P_\m(\x),
\]
where $P_\m$ denotes the degree $n-1$ Taylor polynomial of $f$ at $\x=\m/N$. That is, 
\begin{equation} \label{eq_proof_th_4_expansion}
P_\m(\x)=\sum_{|\n|<n} \frac{D^\n f(\tfrac{\m}{N})}{\n !} \left(\x-\frac{\m}{N}\right)^\n,
\end{equation}
where $|\n|=\sum_{k=1}^d n_k$, $\n !=\prod_{k=1}^d n_k!$, and $(\x-\m/N)^\n=\prod_{k=1}^d(x_k-m_k/N)^{n_k}$.  Let $\x\in[0,1]^d$ and note that
\[
\text{support}(\phi_\m)\subset \left\{\x = (x_1,\ldots,x_d):\left|x_k-\frac{m_k}{N}\right|<\frac{1}{N}\right\}, \quad \m\in\{0,\ldots,N\}^d.
\]
Hence, the approximation error between $f$ and its local Taylor approximation satisfies
\begin{align*}
|f(\x)-f_N(\x)| &= \left|\sum_{\m\in\{0,\ldots,N\}^d} \phi_\m (f(\x)-P_{\m}(\x))\right|\leq \sum_{\m:\left|x_k-\frac{m_k}{N}\right|<\frac{1}{N}}|f(\x)-P_\m(\x)|\\
&\leq \frac{2^d d^n}{n!}\left(\frac{1}{N}\right)^n \max_{|\n|=n}\esssup_{\x\in[0,1]^d}|D^\n f(\x)| \leq \frac{2^d d^n}{n!}\left(\frac{1}{N}\right)^n.
\end{align*}
We now select (see~\cite[Theorem~1]{yarotsky2017error} for a similar idea)
\[
N=\left\lceil\left(\frac{n!}{2^dd^n}\frac{\epsilon}{2}\right)^{-1/n}\right\rceil,
\]
so that 
\begin{equation} \label{eq_taylor_approx}
\max_{\x\in[0,1]^d} \left|f(\x)-f_N(\x)\right| \leq \epsilon/2.
\end{equation}
We now approximate the function $f_n$ by a rational network using~\cref{th_approx_piecewise_rat,prop_rat_x_n}. First, we write $f_N$ as
\begin{equation} \label{eq_proof_th_4_expansion_2}
f_N(\x) = \sum_{\m\in\{0,\ldots,N\}^d}\sum_{|\n|<n}a_{\m,\n}\phi_\m(\x)\left(\x-\frac{\m}{N}\right)^\n,
\end{equation}
where $|a_{\m,\n}|\!\leq\! 1$ and the monomials are uniformly bounded by $1$ (see \cref{eq_proof_th_4_expansion}). \cref{eq_proof_th_4_expansion_2} consists of at most $d^n(N+1)^d$ terms of the form $\phi_\m(\x)(\x-\m/N)^\n$. The monomial part $(\x-\m/N)^\n$ in \cref{eq_proof_th_4_expansion_2} is representable by a rational network of size $\mathcal{O}(d\log(n)^2)$ using~\cref{prop_rat_x_n}, including the fact that the multiplication is a rational network with one layer and three nodes. Let $0<\delta<1$ be a small number, for each $m_k\in\{0,\ldots,N\}$ the piecewise linear function $\psi_{m_k}$ has a Lipschitz constant of $L=3N$. Therefore, it can be approximated with a tolerance $\delta$ by a rational network $\tilde{\psi}_{m_k}$ of size $\mathcal{O}(\log(\log(N/\delta)))$ (see \cref{th_approx_piecewise_rat}). We can assume $\|\tilde{\psi}_{m_k}\|_\infty=1$ by increasing the size of the network by a constant. This yields the following approximation error between a term in~\cref{eq_proof_th_4_expansion_2} and the rational network constructed using $\tilde{\psi}_{m_k}$:
\begin{multline*}
\left|\phi_\m(\x)\left(\x-\frac{\m}{N}\right)^\n-\prod_{k=1}^d\tilde{\psi}_{m_k}(x_k)\left(\x-\frac{\m}{N}\right)^\n\right|
\leq \left|\prod_{k=1}^d \psi_{m_k}(x_k)-\prod_{k=1}^d\tilde{\psi}_{m_k}(x_k)\right|\\
\begin{aligned}
&\leq \left|\psi_{m_1}(x_1)-\tilde{\psi}_{m_1}(x_1)\right| \left|\prod_{k=2}^d \psi_{m_k}(x_k)\right|
+\left|\tilde{\psi}_{m_1}(x_1)\right| \left|\prod_{k=2}^d \psi_{m_k}(x_k)-\prod_{k=2}^d\tilde{\psi}_{m_k}(x_k)\right|\\
&\leq \left|\psi_{m_1}(x_1)-\tilde{\psi}_{m_1}(x_1)\right| 
+\left|\prod_{k=2}^d \psi_{m_k}(x_k)-\prod_{k=2}^d\tilde{\psi}_{m_k}(x_k)\right|\\
&\leq \delta +\left|\prod_{k=2}^d \psi_{m_k}(x_k)-\prod_{k=2}^d\tilde{\psi}_{m_k}(x_k)\right| \leq d\delta.
\end{aligned}
\end{multline*}
Here, the final inequality is derived by repeating the argument of the previous inequalities for $x_2,\ldots,x_d$. If we denote by $\nrat$ the rational network approximation to $f_N$ constructed above, then, for all $\x\in[0,1]^d$, we have
\begin{align*}
|f_N(\x)-\nrat(\x)|&\leq \!\!\!\sum_{\m\in\{0,\ldots,N\}^d}\sum_{|\n|<n}|a_{\m,\n}|\left|\phi_\m(\x)\!\left(\x-\frac{\m}{N}\right)^\n-\prod_{k=1}^d\tilde{\psi}_{m_k}(x_k)\!\left(\x-\frac{\m}{N}\right)^\n\right|\\
&\leq 2^d d^{n+1}\delta.
\end{align*}
Therefore, we select $\delta = \epsilon/(2^{d+1} d^{n+1})$ so that $\max_{\x\in [0,1]^d} |f_N(\x)-\tilde{f}_N(\x)|\leq \epsilon/2$. Then, by \cref{eq_taylor_approx}, we have
\[
\max_{\x\in[0,1]^d} \left|f(\x)-\nrat(\x)\right| \leq \frac{\epsilon}{2}+\frac{\epsilon}{2}\leq \epsilon.
\]
The statement of the theorem follows as the rational network $\nrat$ has size at most
\[
\mathcal{O}(d^n(N+1)^d\log(\log(N/\delta)))\!=\!\mathcal{O}(\epsilon^{-d/n}\log(\log(1/\epsilon^{1+1/n})))\!=\!\mathcal{O}(\epsilon^{-d/n}\log(\log(1/\epsilon))).
\]
\end{proof}

A theorem proved by DeVore et al.~\cite{devore1989optimal} gives a lower bound of $\Omega(\epsilon^{-d/n})$ on the number of parameters needed by a neural network to express any function in $F_{d,n}$ with an error $\epsilon$, under the assumption that the weights are chosen continuously. Comparing $\mathcal{O}(\epsilon^{-d/n}\log(\log(1/\epsilon)))$ and $\mathcal{O}(\epsilon^{-d/n}\log(1/\epsilon))$, we find that rational neural networks require exponentially fewer nodes than ReLU networks with respect to the optimal bound of $\Omega(\epsilon^{-d/n})$ to approximate functions in $F_{d,n}$.

\section{Experiments using rational neural networks} \label{sec_experiments}

In this section, we consider neural networks with trainable rational activation functions of type $(3,2)$. We select the type $(3,2)$ based on empirical performance; roughly, a low-degree (but higher than $1$) rational function is ideal for generating high-degree rational functions by composition, with a small number of parameters. The rational activation units can be easily implemented in the open-source TensorFlow library~\cite{abadi2016tensorflow} by using the \texttt{polyval} and \texttt{divide} commands for function evaluations. The coefficients of the numerators and denominators of the rational activation functions are trainable parameters, determined at the same time as the weights and biases of the neural network by backpropagation and a gradient descent optimization algorithm.

One crucial question is the initialization of the coefficients of the rational activation functions~\cite{chen2018rational,molina2019pad}. A badly initialized rational function might contain poles on the real axis, leading to exploding values, or converge to a local minimum in the optimization process. Our experiments, supported by the empirical results of Molina et al.~\cite{molina2019pad}, show that initializing each rational function with the best rational approximation to the ReLU function (as described in \cref{lem_relu_rat}) produces good performance. The underlying idea is to initialize rational networks near a network with ReLU activation functions, widely used for deep learning. Then, the adaptivity of the rational functions allows for further improvements during the training phase. We represent the initial rational function used in our experiments in \cref{fig_init_rat}(right). The coefficients of this function are obtained by using the \texttt{minimax} command, available in the Chebfun software~\cite{driscoll2014chebfun,filip2018rational} for numerically computing rational approximations, and are given in \cref{table_coeffs}.

\begin{table}[htbp]
  \caption{Initialization coefficients of the rational activation functions.}
  \label{table_coeffs}
  \centering
  \vspace{0.5cm}
  \begin{tabular}{ccccccc}
    \toprule
    $a_0$ & $a_1$ & $a_2$ & $a_3$ & $b_0$ & $b_1$ & $b_2$ \\
    \midrule
    $1.1915$ & $1.5957$ & $0.5000$ & $0.0218$ & $2.3830$ & $0.0000$ & $1.0000$ \\
    \bottomrule
  \end{tabular}
\end{table}

In the following experiments, we use a single rational activation function of type $(3,2)$ at each layer, instead of different functions at each node to reduce the number of trainable parameters and the computational training expense. This adds 7 degrees of freedom per layer.

\subsection{Approximation of functions} \label{sec_approx_fun}

Raissi, Perdikaris, and Karniadakis~\cite{raissi2018deep,raissi2019physics} introduce a framework called \emph{deep hidden physics models} for discovering nonlinear partial differential equations (PDEs) from observations. This technique requires to solving the following interpolation problem: given the observation data $(u_i)_{1\leq i\leq N}$ at the spatio-temporal points $(x_i,t_i)_{1\leq i\leq N}$, find a neural network $\mathcal{N}$ (called the identification network), that minimizes the loss function
\begin{equation} \label{eq_loss_approx}
\mathcal{L}=\frac{1}{N}\sum_{i=1}^N|\mathcal{N}(x_i,t_i)-u_i|^2.
\end{equation}
This technique has successfully discovered hidden models in fluid mechanics~\cite{raissi2020hidden}, solid mechanics~\cite{haghighat2020deep}, and nonlinear PDEs such as the Korteweg--de Vries (KdV) equation~\cite{raissi2019physics}. Raissi et al.~use an identification network, consisting of $4$ layers and $50$ nodes per layer, to interpolate samples from a solution to the KdV equation. Moreover, they observe that networks based on smooth activation functions, such as the hyperbolic tangent ($\tanh(x)$) or the sinusoid ($\sin(x)$), outperform ReLU neural networks~\cite{raissi2018deep,raissi2019physics}. However, the performance of these smooth activation functions highly depends on the application.

Moreover, these functions might not be adapted to approximate non-smooth or highly oscillatory solutions. Recently, Jagtap, Kawaguchi, and Karniadakis~\cite{jagtap2020adaptive} proposed and analyzed different adaptive activation functions to approximate smooth and discontinuous functions with physics-informed neural networks. More specifically, they use an adaptive version of classical activation functions such as sigmoid, hyperbolic tangent, ReLU, and Leaky ReLU. The choice of these trainable activation functions introduces another parameter in the design of the neural network architecture, which may not be ideal for use for a black-box data-driven PDE solver. 

\begin{figure}[htbp]
\centering
\vspace{0.2cm}
\begin{overpic}[width=0.45\textwidth]{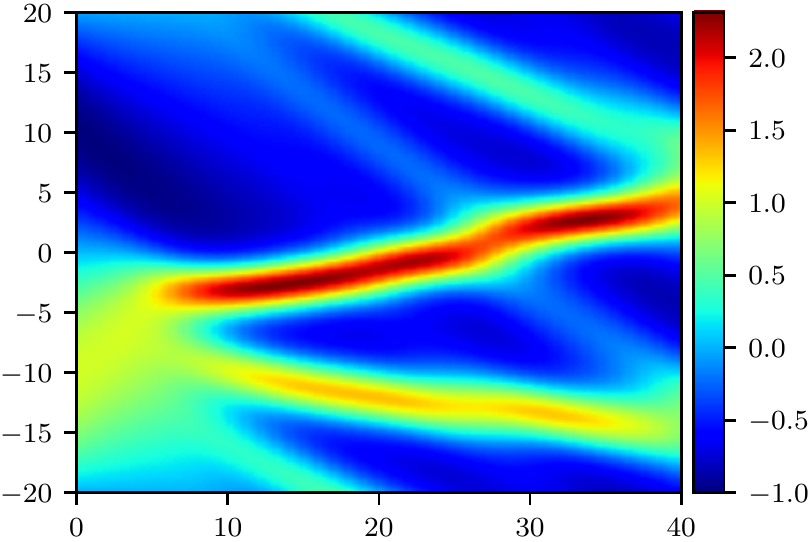}
\put(46,-6){$t$}
\put(-1,34.3){$x$}
\end{overpic}
\hspace{0.7cm}
\begin{overpic}[width=0.44\textwidth, trim=30 15 45 25,clip]{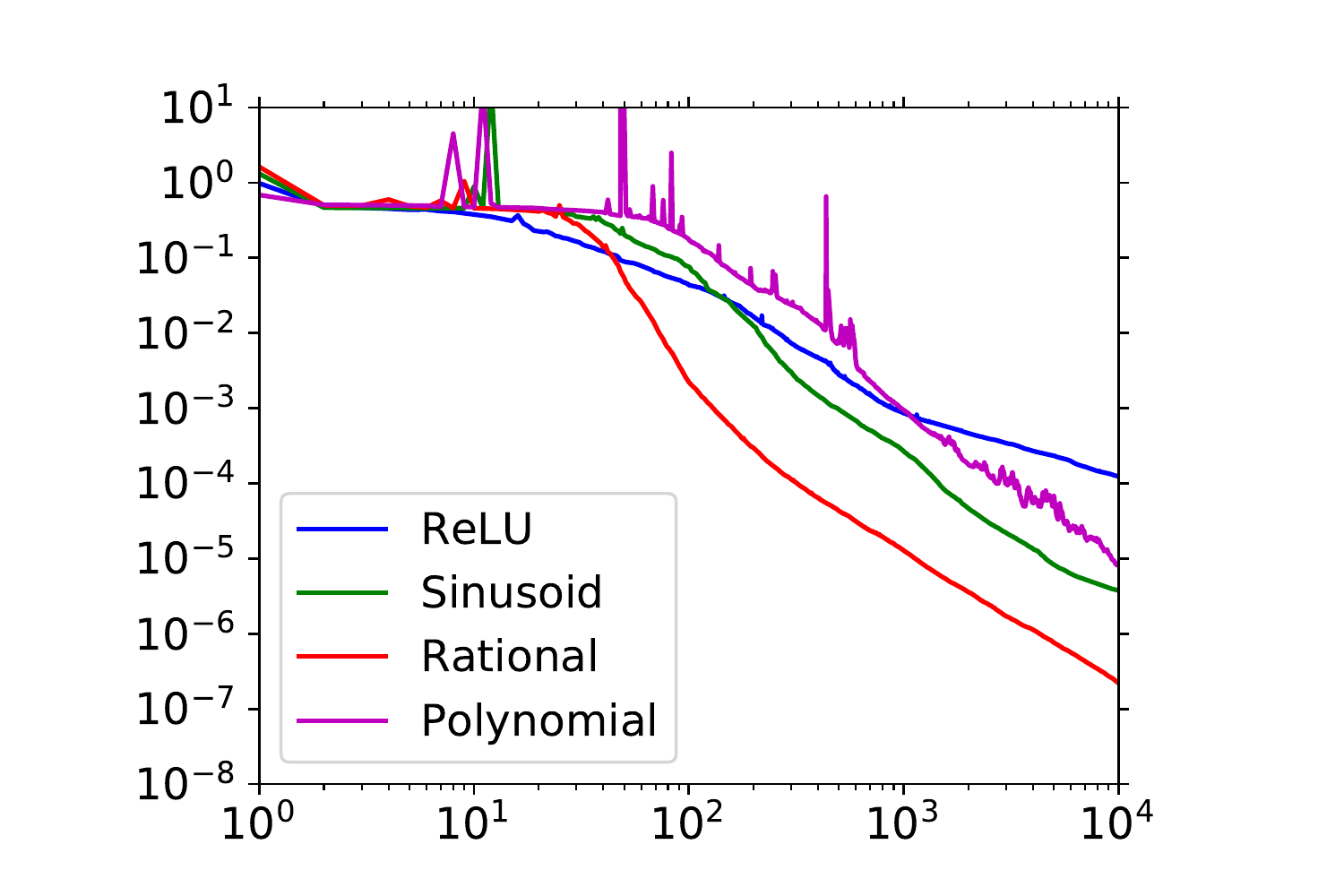}
\put(45,-7){Epochs}
\put(-7,17){\rotatebox{90}{Validation loss}}
\end{overpic}
\vspace{0.4cm}
\caption{Solution to the KdV equation used as training data (left) and validation loss of a ReLU (blue), sinusoid (green), rational (red), and polynomial (purple) neural networks with respect to the number of optimization steps (right).}
\label{fig_loss_2d}
\end{figure}

We illustrate that rational neural networks can address the issues mentioned above due to their adaptivity and approximation power (see \cref{sec_th_result}). Similarly to Raissi~\cite{raissi2018deep}, we use a solution $u$ to the KdV equation:
\[
u_t=-uu_x-u_{xxx}, \quad u(x,0) = -\sin(\pi x/20),
\]
as training data for the identification network (see the left panel of \cref{fig_loss_2d}). We use the TensorFlow implementation\footnote{We adapt the code that is publicly available~\cite{raissiGit}.} of the deep hidden physics model framework to build and train the identifier network $\mathcal{N}$ that approximates a solution $u$ to the KdV equation.  The true solution is computed on the domain $(x,t)\in[-20,20]\times[0,40]$ by Raissi~\cite{raissi2018deep} using the Chebfun package~\cite{driscoll2014chebfun} with a spectral Fourier discretization of $512$ and a time-step of $\Delta t = 10^{-4}$. Moreover, the solution is stored after every $2000$ time steps, giving a testing data set of approximately $10^5$ spatio-temporal points in $[-20,20]\times[0,40]$. We then constituted the training and validation sets (of $10^4$ points each) by randomly subsampling the solution at $2\times 10^4$ points in $[-20,20]\times[0,40]$.

In a similar manner to~\cite{raissi2018deep}, we use a fully connected identification network to approximate $u$ with $4$ hidden layers with $50$ nodes per layer. The network is trained using the L-BFGS optimization algorithm with $10,\!000$ iterations. We train and compare four networks with the following activation functions: ReLU, sinusoid, trainable rational functions of type $(3,2)$, and trainable polynomials of degree $3$. Furthermore, the rational activation functions are initialized to be the best approximation to the ReLU function, using the initial coefficients reported in~\cref{table_coeffs}.

The mean squared error (MSE) of the neural networks on the validation set throughout the training phase is reported in the right panel of \cref{fig_loss_2d}. We observe that the rational neural network outperforms the sinusoid network, despite having the same asymptotic convergence rate. The network with polynomial activation functions (chosen to be of degree 3 in this example) is harder to train than the rational network, as shown by the non-smooth validation loss (see the right panel of \cref{fig_loss_2d}). We highlight that rational neural networks are never much bigger in terms of trainable parameters than ReLU networks since the increase is only linear with respect to the number of layers. Here, the ReLU network has $8000$ parameters (consisting of weights and biases), while the rational network has $8000+7\times\#\textup{layers}=8035$. The  ReLU, sinusoid, rational, and polynomial networks achieve the following mean square errors after $10^4$ epochs:
\begin{align*}
&\text{MSE}(u_{\text{ReLU}}) = 1.9\times 10^{-4}, &&\text{MSE}(u_{\text{Sinusoid}}) = 3.3\times 10^{-6},\\ 
&\text{MSE}(u_{\text{Rational}}) = 1.2\times 10^{-7}, && \text{MSE}(u_{\text{Polynomial}}) = 3.6\times 10^{-5}.
\end{align*}

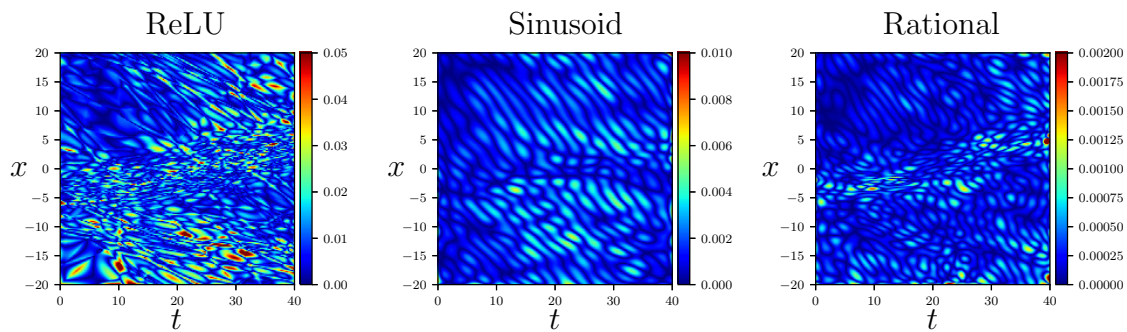
\begin{figure}[htbp]
\centering
\vspace{0.3cm}
\begin{overpic}[width=0.95\textwidth]{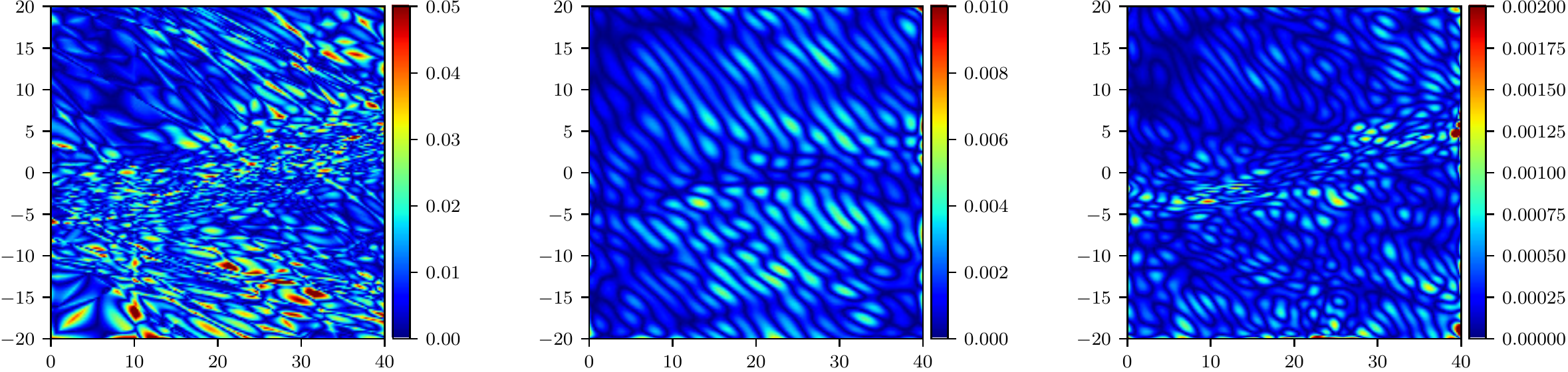}
\put(11,25){ReLU}
\put(13.3,-2){$t$}
\put(-1.3,12){$x$}
\put(44,25){Sinusoid}
\put(47.7,-2){$t$}
\put(33,12){$x$}
\put(78,25){Rational}
\put(82,-2){$t$}
\put(67,12){$x$}
\end{overpic}
\vspace{0.1cm}
\caption{Approximation errors of the neural networks with ReLU, sinusoid, and rational activation layers. Note the different scales of the errors.}
\label{fig_2d_approx}
\end{figure}

The rational neural network is approximately five times more accurate than the sinusoid network used by Raissi and twenty times more accurate than the ReLU network. The absolute approximation errors between the different neural networks and the exact solution to the KdV equation is illustrated in \cref{fig_2d_approx}. We find that the approximation errors made by the ReLU network are not uniformly distributed in space and time and located in specific regions, indicating that a network with non-smooth activation functions is not appropriate to resolve smooth solutions to PDEs.

\begin{figure}[htbp]
\centering
\vspace{0.2cm}
\begin{overpic}[width=0.44\textwidth, trim=30 15 45 25,clip]{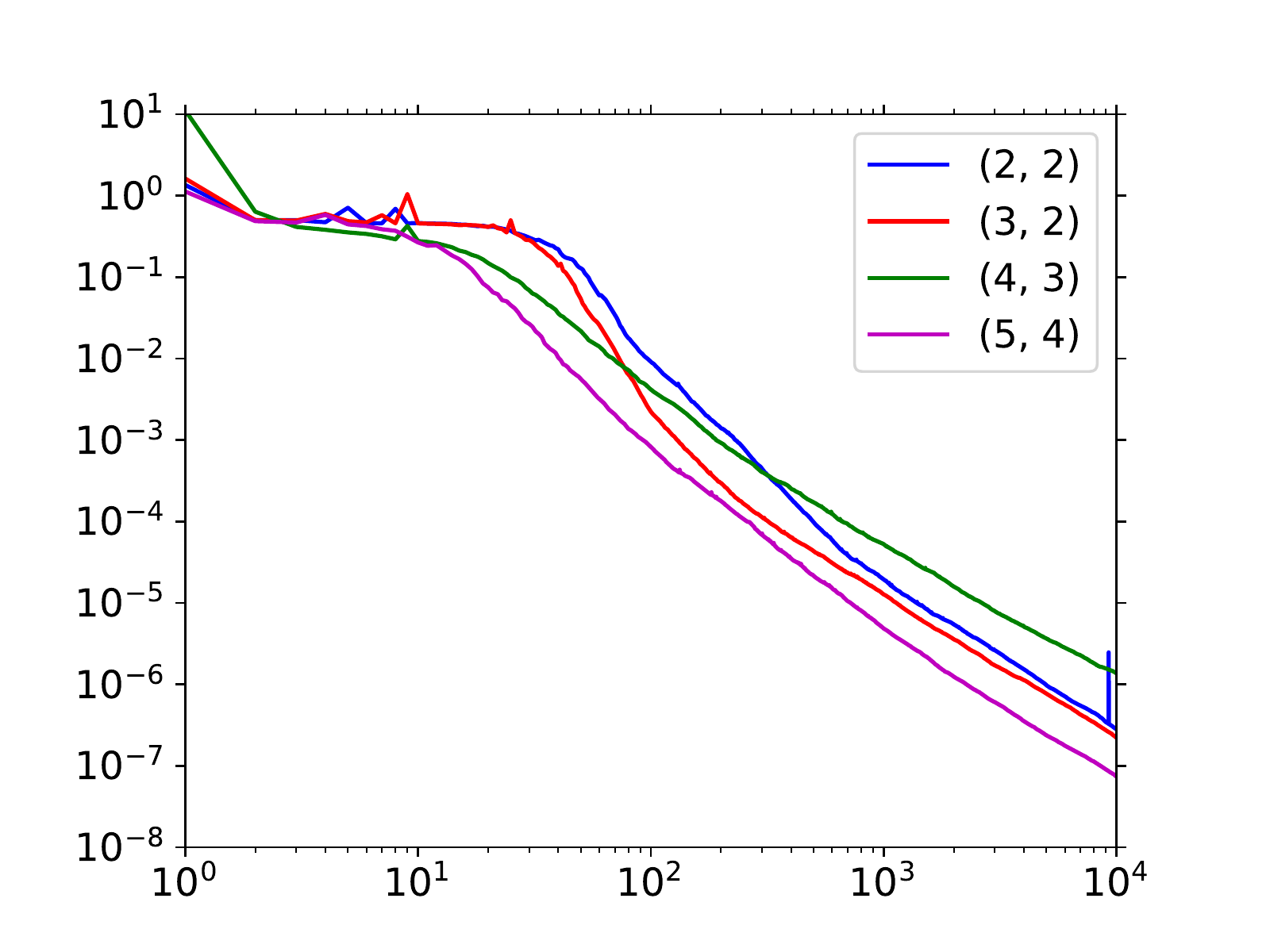}
\put(45,-5){Epochs}
\put(-9,23){\rotatebox{90}{Validation loss}}
\end{overpic}
\vspace{0.2cm}
\caption{Validation loss of rational networks of types $(2,2)$, $(3,2)$, $(4,3)$, and $(5,4)$ with respect to the number of epochs.}
\label{fig_rational_loss_2d}
\end{figure}

Finally, in \cref{fig_rational_loss_2d}, we compare rational neural networks with different degree activation functions (each initialized to approximate the ReLU function using the MATLAB code 
\texttt{initial\_rational\_coeffs.m} available at~\cite{boulleGitrational}) and find that they all performed better than ReLU networks. While a type $(3,2)$ rational offers a good trade-off between the number of parameters and quality of approximation according to the theoretical results presented in \cref{sec_th_result}, the type of rational function might well depend on the application considered.

\subsection{Generative adversarial networks}

Generative adversarial networks are used to generate synthetic examples from an existing dataset~\cite{goodfellow2014generative}. They consist of two networks: a generator to produce synthetic samples and a discriminator to evaluate the samples of the generator with the training dataset. Radford et al.~\cite{radford2015unsupervised} describe deep convolutional generative adversarial networks (DCGANs) to build good image representations using convolutional architectures. They evaluate their model on the MNIST and ImageNet image datasets~\cite{deng2009imagenet,lecun1998gradient}. 

This section highlights the simplicity of using rational activation functions in existing neural network architectures by training an Auxiliary Classifier GAN (ACGAN)~\cite{odena2017conditional} on the MNIST dataset. In particular, the neural network, referred to as the ReLU network in this section, consists of convolutional generator and discriminator networks with ReLU and Leaky ReLU~\cite{maas2013rectifier} activation units (respectively) and is used as a reference GAN. We adapt the Keras example in~\cite{KerasDoc} to train an Auxiliary Classifier GAN with rational activation functions on the MNIST. The hyper-parameters used for the GAN experiment are given in~\cref{table_GAN}. Moreover, the GAN is trained on $20$ epochs with a batch size of $100$ by Adam's optimization algorithm~\cite{kingma2014adam} and the following parameters: $\alpha=0.0002$ and $\beta_1=0.5$, as suggested by~\cite{radford2015unsupervised}.

\begin{table}[htbp]
  \caption{Hyper-parameters of the GAN experiment, BN denotes the presence of a Batch normalization layer. The Generator and Discriminator networks are trained with ReLU and rational activation functions, initialized with the coefficients reported in \cref{table_coeffs}. Transposed convolution layers and rational activation functions are respectively abbreviated as ``Transp. Conv.'' and ``Rat.''.}
  \label{table_GAN}
  \centering
  \vspace{0.5cm}
  \begin{tabular}{lcccccc}
    \toprule
     \multicolumn{1}{c}{Operation} & Kernel & Strides & Features & BN & Dropout & Activation \\
    \midrule
    \multicolumn{1}{c}{Generator} \\    
    Linear & N/A & N/A & 3456 & \xmark & 0.0 & ReLU / Rat.\\
    Transp. Conv. & $5\times 5$ & $1\times 1$ & 192 & \cmark & 0.0 & ReLU / Rat.\\
    Transp. Conv. & $5\times 5$ & $2\times 2$ & 96 & \cmark & 0.0 & ReLU / Rat.\\
    Transp. Conv. & $5\times 5$ & $2\times 2$ & 1 & \xmark & 0.0 & Tanh\\
    \multicolumn{1}{c}{Discriminator} \\    
    Convolution & $3\times 3$ & $2\times 2$ & 32 & \xmark & 0.3 & Leaky ReLU / Rat.\\
    Convolution & $3\times 3$ & $1\times 1$ & 64 & \xmark & 0.3 & Leaky ReLU / Rat.\\
    Convolution & $3\times 3$ & $2\times 2$ & 128 & \xmark & 0.3 & Leaky ReLU / Rat.\\
    Convolution & $3\times 3$ & $1\times 1$ & 256 & \xmark & 0.3 & Leaky ReLU / Rat.\\
    Linear & N/A & N/A & 11 & \xmark & 0.0 & Soft-Sigmoid\\
    \bottomrule
  \end{tabular}
\end{table}

As in the experiment described in \cref{sec_approx_fun}, we replace the activation units of the generative and discriminator networks by a rational function with trainable coefficients (see~\cref{fig_init_rat}). We initialize the activation functions in the training phase with the best rational function that approximates the ReLU function on $[-1,1]$.

\begin{figure}[htbp]
\centering
\begin{minipage}{0.8\textwidth}
\centering
\begin{overpic}[width=0.225\textwidth,trim={0 952px 140px 0},clip]{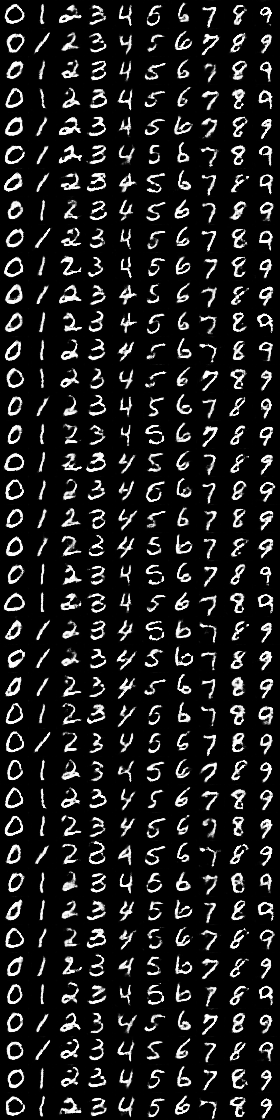}
\put(23,-11){epoch 5}
\put(-15,35){\rotatebox{90}{ReLU}}
\end{overpic}
\begin{overpic}[width=0.225\textwidth,trim={0 952px 140px 0},clip]{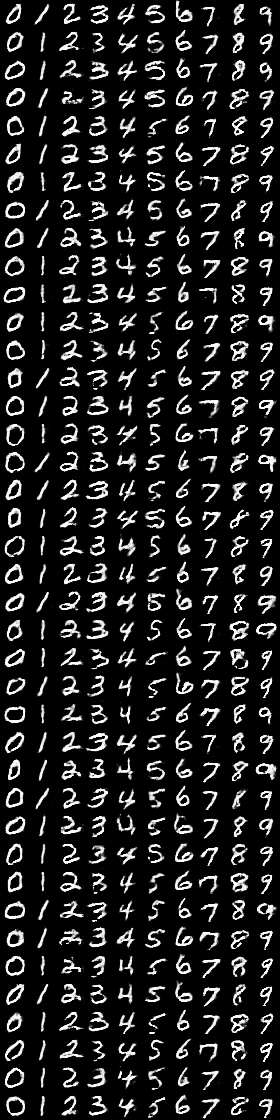}
\put(19,-11){epoch 10}
\end{overpic}
\begin{overpic}[width=0.225\textwidth,trim={0 952px 140px 0},clip]{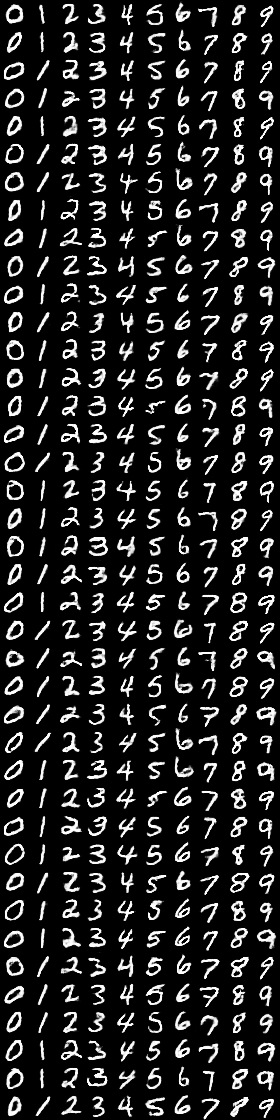}
\put(19,-11){epoch 15}
\end{overpic}
\begin{overpic}[width=0.225\textwidth,trim={0 952px 140px 0},clip]{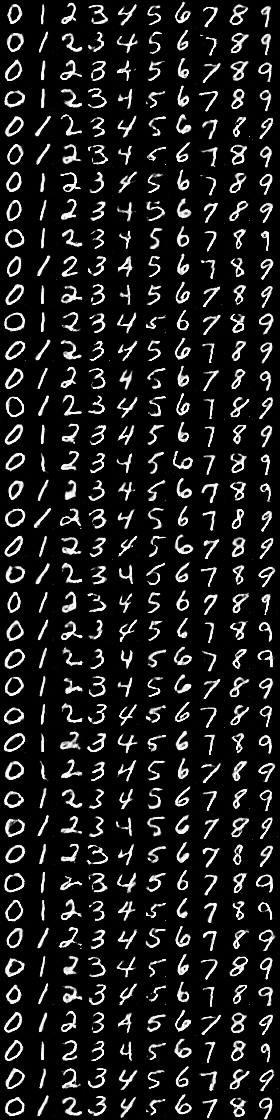}
\put(19,-11){epoch 20}
\end{overpic}\\
\vspace{13px}
\begin{overpic}[width=0.225\textwidth,trim={0 952px 140px 0},clip]{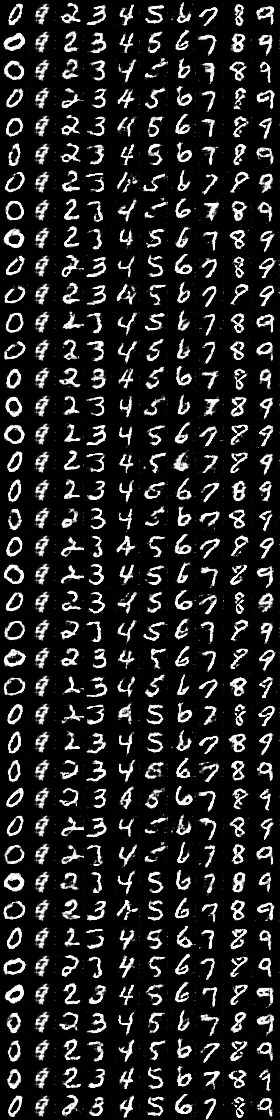}
\put(23,-11){epoch 5}
\put(-15,26){\rotatebox{90}{Rational}}
\end{overpic}
\begin{overpic}[width=0.225\textwidth,trim={0 952px 140px 0},clip]{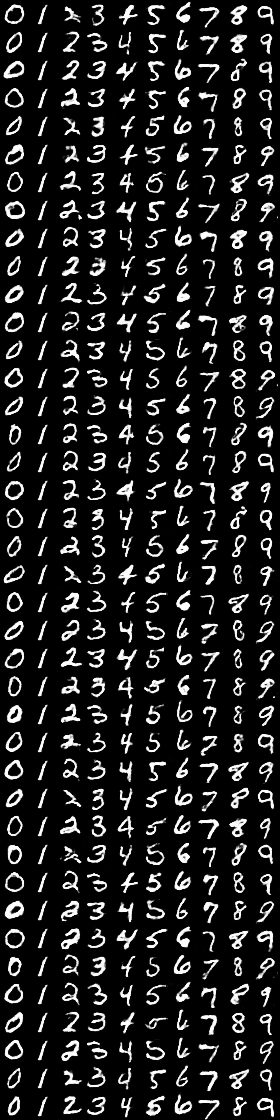}
\put(19,-11){epoch 10}
\end{overpic}
\begin{overpic}[width=0.225\textwidth,trim={0 952px 140px 0},clip]{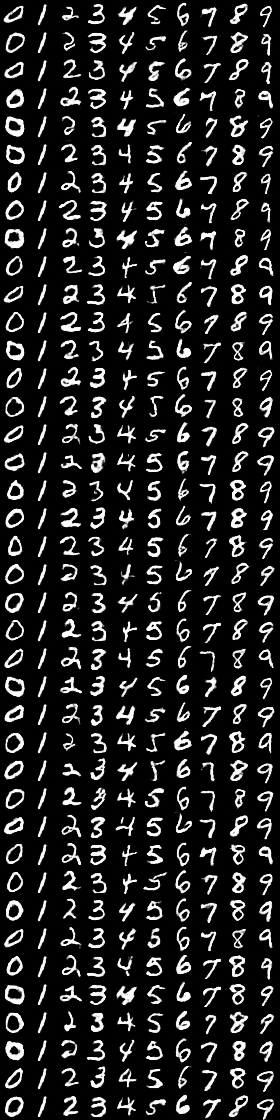}
\put(19,-11){epoch 15}
\end{overpic}
\begin{overpic}[width=0.225\textwidth,trim={0 952px 140px 0},clip]{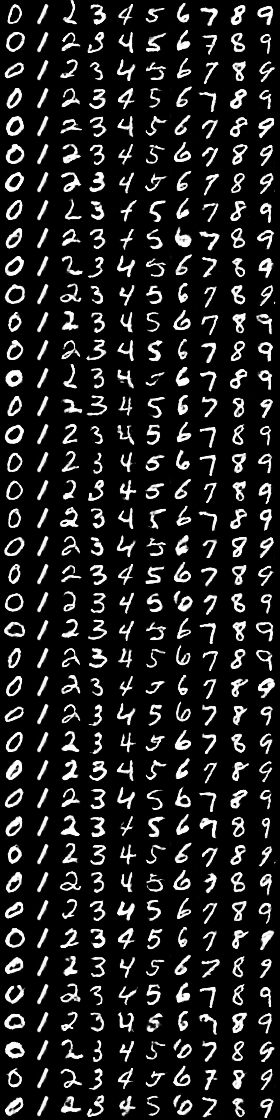}
\put(19,-11){epoch 20}
\end{overpic}
\end{minipage}%
\begin{minipage}{0.18\textwidth}
\centering
\begin{overpic}[width=\textwidth,trim={0 756px 140px 0},clip]{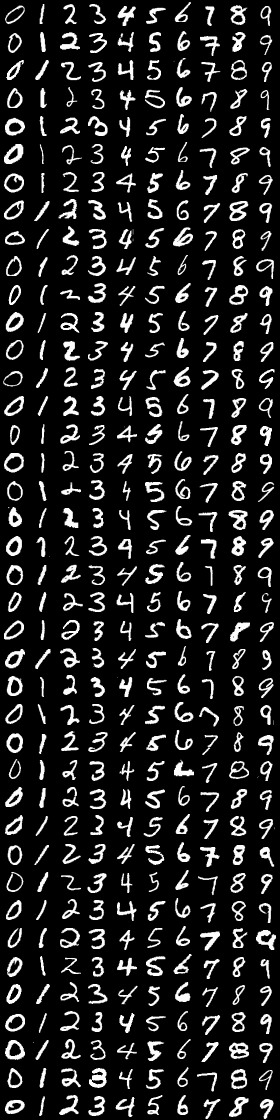}
\put(0.5,-5){MNIST images}
\end{overpic}
\end{minipage}
\vspace{0.5cm}
\caption{Digits generated by a ReLU (top) and rational (bottom) auxiliary classifier generative adversarial network. The right panel contains samples from the first five classes of the MNIST dataset for comparison.}
\label{fig_MNIST_gan}
\end{figure}

We show images of digits from the first five classes generated by a ReLU and rational GANs at different epochs of the training in \cref{fig_MNIST_gan} (the samples are generated randomly and are not manually selected). We observe that a rational network can generate realistic images with a broader range of features than the ReLU network, as illustrated by the presence of bold numbers at the epoch 20 in the bottom panel of \cref{fig_MNIST_gan}. 

\begin{figure}[htbp]
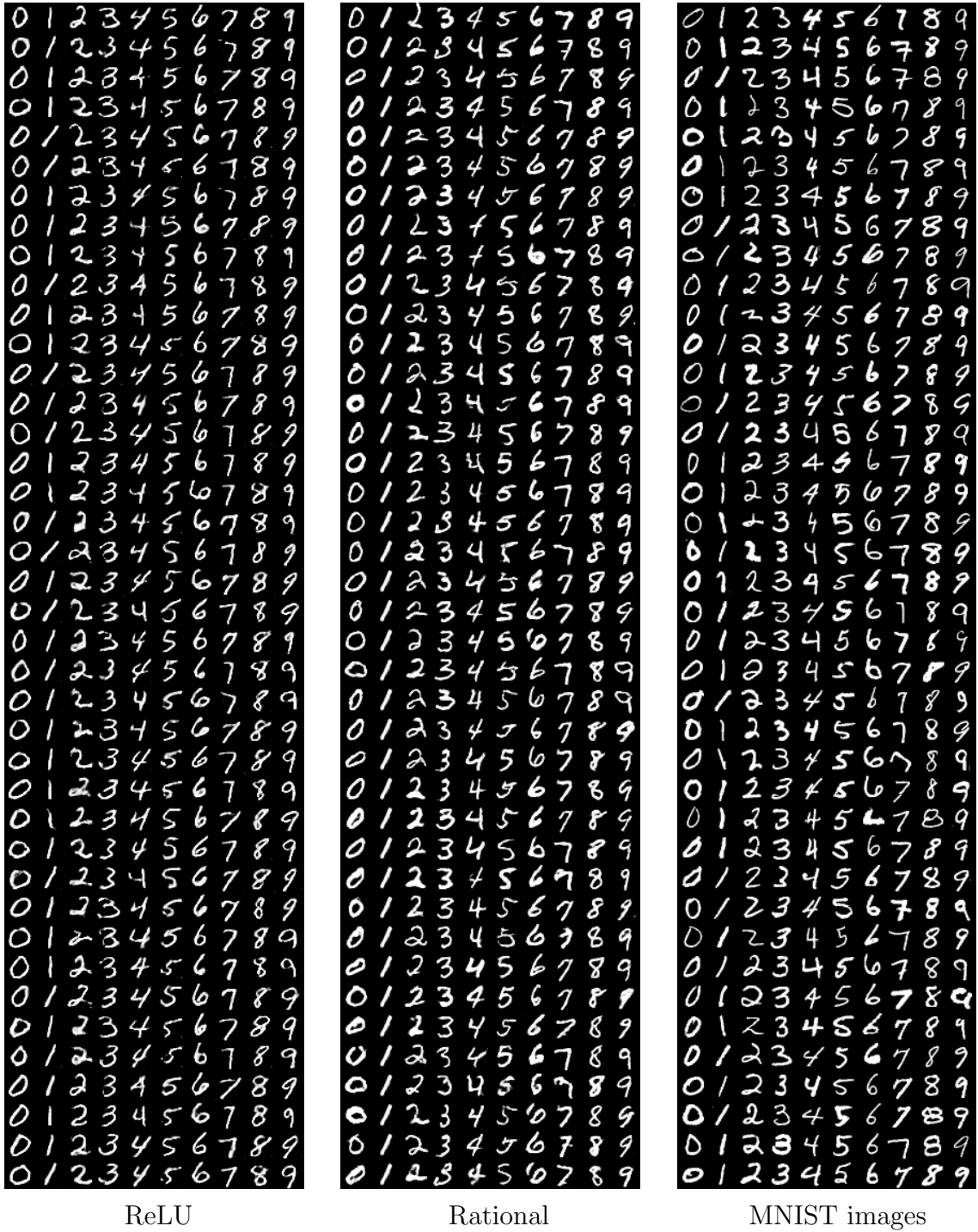

\centering
\begin{overpic}[width=0.3\textwidth,trim={0 0 0 0},clip]{Figure/Rational/gan/relu/plot_epoch_020_generated.png}
\put(10,-3){ReLU}
\end{overpic}
\hspace{0.3cm}
\begin{overpic}[width=0.3\textwidth,trim={0 0 0 0},clip]{Figure/Rational/gan/rat/plot_epoch_020_generated.png}
\put(9,-3){Rational}
\end{overpic}
\hspace{0.3cm}
\begin{overpic}[width=0.3\textwidth,trim={0 0 0 0},clip]{Figure/Rational/gan/mnist_images.png}
\put(6,-3){MNIST images}
\end{overpic}
\vspace{0.5cm}
\caption{Forty images generated by a ReLU network and a rational network after 20 epochs, together with real images from the MNIST dataset.}
\label{fig_MNIST_gan_supp}
\end{figure}

We report in~\cref{fig_MNIST_gan_supp} samples of the $10$ classes present in the MNIST dataset (right) and images generated at the $20$th epoch by the GAN with ReLU/Leaky ReLU units (left) and rational activation functions (middle). We observe that the digits one generated by the rational network are identical, suggesting that the rational GAN suffers from mode collapse. It should be noted that generative adversarial networks are notoriously tricky to train~\cite{goodfellow2016deep}. The hyper-parameters of the reference model are intensively tuned for a piecewise linear activation function (as shown by the use of Leaky ReLU in the discriminator network). Moreover, many stabilization methods have been proposed to resolve the mode collapse and non-convergence issues in training, such as Wasserstein GAN~\cite{arjovsky2017wasserstein}, Unrolled Generative Adversarial Networks~\cite{metz2016unrolled}, and batch normalization~\cite{ioffe2015batch}. These techniques could be explored and combined with rational networks to address the mode collapse issue observed in this experiment.

\dobib
 % Rational networks

\setcounter{chapter}{4}
\renewcommand{\thefootnote}{\fnsymbol{footnote}}
\chapter[Data-driven discovery of Green's functions with deep learning]{Data-driven discovery of Green's functions with deep learning\footnotemark} \label{chap_data_green}

\footnotetext{This chapter is based on a paper with Christopher Earls and Alex Townsend~\cite{boulle2021data}, published in Scientific Reports. Earls and Townsend had an advisory role; I designed the deep learning method, performed the numerical experiments, and was the lead author in writing the paper.}

\renewcommand*{\thefootnote}{\arabic{footnote}}
\setcounter{footnote}{0}

Deep learning (DL) holds promise as a scientific tool for discovering elusive patterns within the natural and technological world~\cite{goodfellow2016deep,lecun2015deep}. These patterns hint at undiscovered partial differential equations (PDEs) that describe governing phenomena within biology and physics. From sparse and noisy laboratory observations, we aim to learn mechanistic laws of nature~\cite{brunton2020machine,karniadakis2021physics}. Recently, scientific computing and machine learning have successfully converged on PDE discovery~\cite{Brunton,Rudy,schaeffer2017learning,zhang2020data}, PDE learning~\cite{feliu2020meta,gin2020deepgreen,li2020neural,lu2021learning,raissi2018deep,raissi2020hidden}, and symbolic regression~\cite{schmidt2009distilling,Udrescu2020} as promising means for applying machine learning to scientific investigations. These methods attempt to discover the coefficients of a PDE model or learn the operator that maps excitations to system responses. The recent DL techniques addressing the latter problem are based on approximating the solution operator associated with a PDE by a neural network (NN)~\cite{feliu2020meta,gin2020deepgreen,li2020neural,lu2021learning,raissi2018deep}. While excellent for solving PDEs, we consider them as ``black box'' and focus here on a data-driven strategy that improves human understanding of the governing PDE model.

We then offer a radically different, alternative approach that is backed by theory~\cite{boulle2021learning} and infuse an interpretation in the model by learning well-understood mathematical objects that imply underlying physical laws. We devise a DL method, employed for learning the Green's functions~\cite{stakgold2011green} associated with unknown governing linear PDEs, and train the neural networks by collecting physical system responses from random excitation functions drawn from a Gaussian process (GP). The empirically derived Green's functions relate the system's response (or PDE solution) to a forcing term, and can then be used as a fast reduced-order PDE solver. The existing graph kernel network~\cite{li2020neural} and DeepGreen~\cite{gin2020deepgreen} techniques also aim to learn solution operators of PDEs based on Green's functions. While they show competitive performance in predicting the solution of the PDE for new forcing functions, the errors between the exact and learned Green's functions are relatively large, which makes the extraction of qualitative and quantitative features of the physical system challenging.

Our secondary objective is to study the discovered Green's functions for clues regarding the physical properties of the observed systems. Our approach relies on the rational neural networks introduced in the previous chapter, which have higher approximation power than standard networks and carry human-understandable features of the PDE, such as shock and singularity locations, as we shall see later.

In this chapter, we use techniques from deep learning to discover the Green's function of linear differential equations $\L u = f$ from input-output pairs $(f,u)$, as opposed to directly learning $\L$, or model parameters.
In this sense, our approach is agnostic to the forward PDE model, but nonetheless offers insights into its physical properties. There are several advantages to learning the Green's function. First, once the Green's function is learned by a neural network, it is possible to compute the solution, $u$, for a new forcing term, $f$, by evaluating an integral (see~\cref{eq_Green_integral}); which is more efficient than training a new NN. Second, the Green's function associated with $\mathcal{L}$ contains information about the operator, $\mathcal{L}$, and the type of boundary constraints that are imposed; which helps uncover mechanistic understanding from experimental data. Finally, as discussed in \cref{chapt_PDE_learning}, it is easier to train NNs to approximate Green's functions, which are square-integrable functions under sufficient regularity conditions~\cite{dong2009green,gruter1982green,stakgold2011green}, than trying to approximate the action of the linear differential operator, $\L$, which is not bounded~\cite{Kreyszig}. Also, any prior mathematical and physical knowledge of the operator, $\L$, can be exploited in the design of the NN architecture, which could enforce a particular structure such as symmetry of the Green's function.

\section{Learning Green's functions}

We consider linear differential operators, $\L$, defined on a bounded domain $\Omega\subset\R^d$, where $d\in\{1,2,3\}$ denotes the spatial dimension. The aim of our method is to discover properties of the operator, $\L$, using $N$ input-output pairs $\{(f_j,u_j)\}_{j=1}^N$, consisting of forcing functions, $f_j:\Omega\to\R$, and system responses, $u_j:\Omega\to\R$, which are solutions to the following equation:
\begin{equation} \label{eq_problem}
\L u_j = f_j, \quad \mathcal{D}(u_j,\Omega) = g,
\end{equation}
where $\mathcal{D}$ is a linear operator acting on the solutions, $u$, and the domain, $\Omega$; with $g$ being the constraint. We assume that the forcing terms have sufficient regularity, and that the operator, $\D$, is a constraint so that \cref{eq_problem} has a unique solution~\cite{stakgold2011green}. An example of constraint is the imposition of homogeneous Dirichlet boundary conditions on the solutions: $\mathcal{D}(u_j,\Omega) \coloneqq u_{j}|_{\partial\Omega}=0$. Note that boundary conditions, integral conditions, jump conditions, or non-standard constraints, are all possible (see \cref{sec_linear_const}).

\subsection{Definitions}
A Green's function~\cite{arfken2011mathematical,evans10,myint2007linear,stakgold2011green} of the operator, $\L$, is defined as the solution to the following equation:
\[\L G(x,y) = \delta(y-x),\quad x,y\in\Omega,\]
where $\L$ is acting on the function $x\mapsto G(x,y)$ for fixed $y\in\Omega$, and $\delta(\cdot)$ denotes the Dirac delta function. The Green's function is well-defined and unique under mild conditions on $\L$, and suitable solution constraints imposed via an operator, $\mathcal{D}$ (see~\cref{eq_problem})~\cite{stakgold2011green}. Moreover, if $(f,u)$ is an input-output pair, satisfying \cref{eq_problem} with $g=0$, then
\[u(x) = \int_{\Omega}G(x,y)f(y)\d y, \quad x\in\Omega.\]
Therefore, the Green's function associated with $\mathcal{L}$ can be thought of as the right inverse of $\mathcal{L}$. 

Let $u_{\text{hom}}$ be the homogeneous solution to~\eqref{eq_problem}, so that
\[\L u_{\text{hom}} = 0, \quad \mathcal{D}(u_{\text{hom}},\Omega) = g.\]
Using superposition, we can construct solutions, $u_j$, to \cref{eq_problem} as $u_j = \tilde{u}_j+u_{\text{hom}}$, where $\tilde{u}_j$ satisfies
\[\L \tilde{u}_j = f_j, \quad \mathcal{D}(\tilde{u}_j,\Omega) = 0.\]
Then, the relation between the system's response, $u_j$, and the forcing term, $f_j$, can be expressed via the Green's function as
\[u_j(x) = \int_{\Omega}G(x,y)f_j(y)\d y + u_{\text{hom}}(x),\quad x\in\Omega.\]
In this chapter, we focus on learning Green's functions and homogeneous solutions from a fixed boundary constraint $g$ but one could also approximate a second Green's function associated with $u_{\text{hom}}$ from multiple boundary constraints. Therefore, we train two NNs: $\N_G:\Omega\times\Omega\to\R\cup\{\pm\infty\}$ and $\N_{\text{hom}}:\Omega\to\R$, to learn the Green's function, and also the homogeneous solution associated with $\L$ and the constraint operator $\mathcal{D}$. Note that this procedure allows us to discover boundary conditions, or constraints, directly from the input-output data without imposing it in the loss function (which often results in training instabilities~\cite{wight2020solving}).

\subsection{Theoretical justification}
Our approach for learning Green's functions associated with linear differential operators has a theoretically rigorous underpinning. Indeed, we showed in \cref{chapt_PDE_learning} that uniformly elliptic operators in three dimensions have an intrinsic \emph{learning rate}, which characterizes the number of training pairs needed to construct an $\epsilon$-approximation in the $L^2$-norm of the Green's function, $G$, with high probability, for $0<\epsilon<1$. The number of training pairs depends on the quality of the covariance kernel used to generate the random forcing terms, $\{f_j\}_{j=1}^N$. Our choice of covariance kernel (\cref{sec_generation_data}) is motivated by the GP quality measure (cf.~\cref{sec_quality_kernel}), to ensure that our set of training forcing terms is sufficiently diverse to capture the action of the solution operator, $f\mapsto u(x) = \int_{\Omega}G(x,y)f(y)\d y$, on a diverse set of functions.

Similarly, the choice of rational NNs to approximate the Green's function, and the homogeneous solution, is justified by the higher approximation power of these networks over ReLU as observed in~\cref{chapt_rational}. Other adaptive activation functions have been proposed for learning or solving PDEs with NNs~\cite{jagtap2020adaptive}, but they are only motivated by empirical observations. Both theory and experiments support rational NNs for regression problems. The number of trainable parameters, consisting of weight matrices, bias vectors, and rational coefficients, needed by a rational NN to approximate smooth functions within $0<\epsilon<1$, can be completely characterized~\cite{boulle2020rational}. This motivates our choice of NN architecture for learning Green's functions.

\section{Deep learning method}

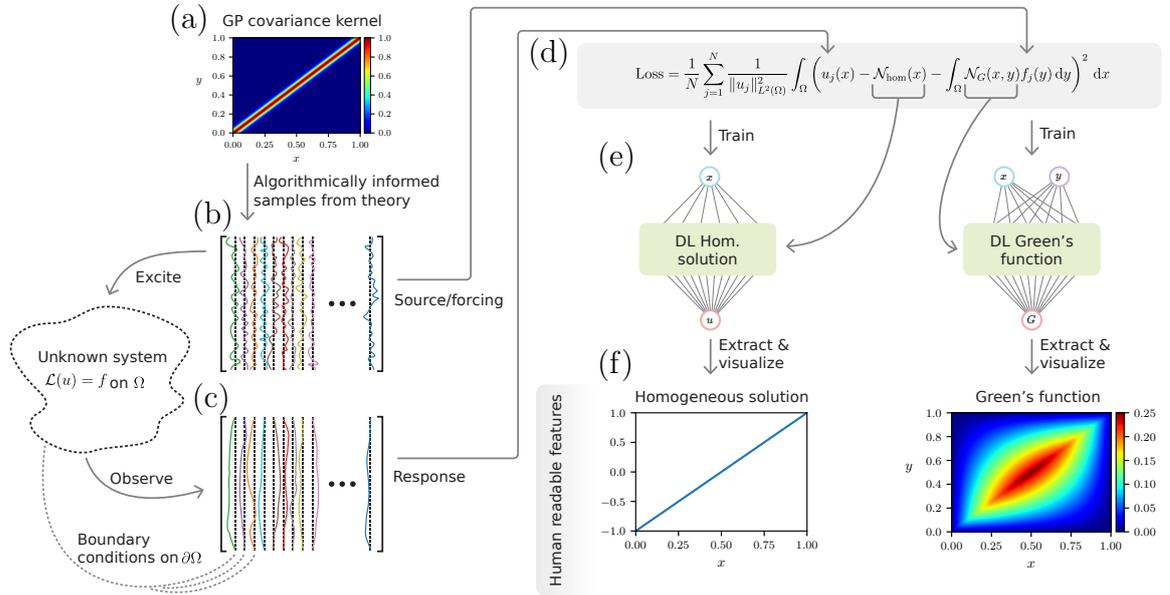
\begin{figure}[htbp]
\centering
\vspace{0.5cm}
\begin{overpic}[width=\textwidth]{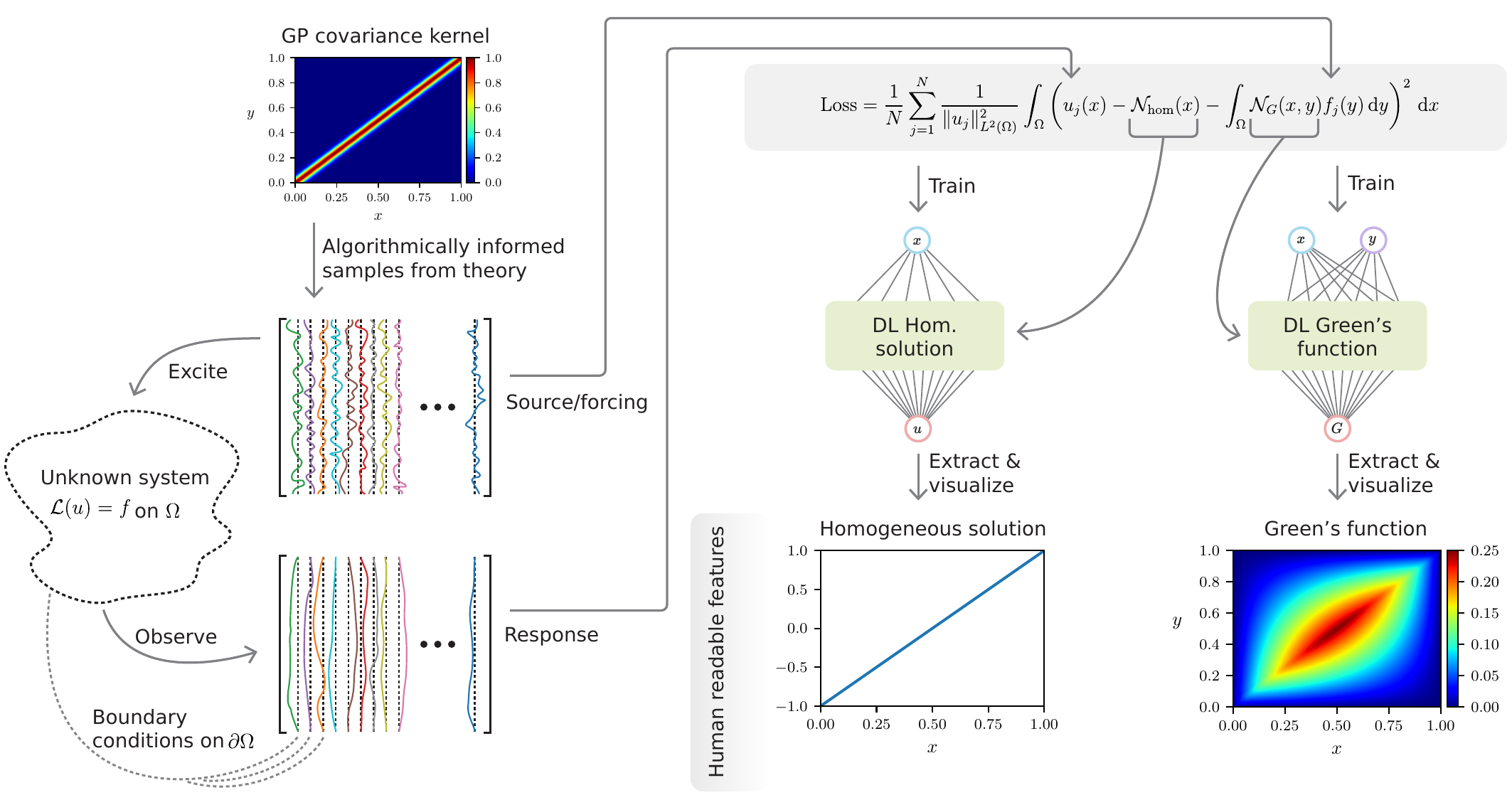}
\put(14,50){(a)}
\put(16,33){(b)}
\put(16,17){(c)}
\put(45,47){(d)}
\put(51,38){(e)}
\put(51,20){(f)}
\end{overpic}
\caption{Schematic of our DL method for learning Green's functions from input-output pairs.  (a) The covariance kernel of the Gaussian process, which is used to generate excitations. (b) The system's response to each excitation is computed and recorded (c). (d) A loss function is minimized to train rational NNs (e). (f) The learned Green's function and homogeneous solution are visualized by sampling the NNs.}
\label{fig_idea}
\end{figure}

Our DL approach (see~\cref{fig_idea}) begins with excitations (or forcing terms), $\{f_j\}_{j=1}^N$, sampled from a Gaussian process having a carefully designed covariance kernel, and corresponding system responses, $\{u_j\}_{j=1}^N$  (see~\cref{chap_svd}). It is postulated that there is an unknown linearized governing PDE so that $\L u_j = f_j$. The selection of random forcing terms is theoretically justified by~\cref{chapt_PDE_learning} and enables us to learn the dominant eigenmodes of the solution operator, using only a small number, $N$, of training pairs. The Green's function, $G$, and homogeneous solution, $u_{\text{hom}}$, which encodes the boundary conditions associated with the PDE, satisfy
\begin{equation} \label{eq_Green_integral}
u_j(x) = \int_{\Omega}G(x,y)f_j(y)\d y + u_{\text{hom}}(x),\quad x\in\Omega,
\end{equation}
and are approximated by two rational neural networks: $\N_G$ and $\N_{\text{hom}}$.

The parameters of the NNs representing the Green's function and homogeneous solution are simultaneously learned through minimization of the loss function displayed in \cref{fig_idea}(d). We discretize the integrals in the loss function at the specified measurement locations $\{x_i\}_{i=1}^{N_u}$, within the domain, $\Omega$, and forcing term sample points, $\{y_i\}_{i=1}^{N_f}$, respectively, using a quadrature rule. 

In this section, we detail the deep learning method used to learn Green's functions. Our DL technique is data-driven and requires minimal by-hand parameter tuning. In fact, all the numerical examples described in this chapter are performed using a single rational NN architecture, initialization procedure, and optimization algorithm\footnote{All data and codes used in this chapter are publicly available on the GitHub and Zenodo repositories at \url{https://github.com/NBoulle/greenlearning/}~\cite{boulleZenodo} to reproduce the numerical experiments and figures. A software package, including additional examples and documentation, is also available at  \url{https://greenlearning.readthedocs.io/}.}.

\subsection{Generating the training data} \label{sec_generation_data}

We create a training dataset, consisting of input-output functions, $\{(f_j\,u_j)\}$ for $1\leq j \leq N$, in three steps: (1) Generating the forcing terms by sampling random functions from a Gaussian process, (2) Solving \cref{eq_problem} for the generated forcing terms, and (3) Sampling the forcing terms, $f_j$, at the points $\{y_1,\ldots,y_{N_f}\}\subset\Omega$ and the system's responses, $u_j$, at $\{x_1,\ldots,x_{N_u}\}\subset\Omega$. Here, $N_f$ and $N_u$ are the forcing and solution discretization sizes, respectively. We recommend that all the forcing terms are sampled on the same grid and similarly for the system's responses. This minimizes the number of evaluations of $\N_G$ during the training phase and reduces the computational and memory costs of training. 

The spatial locations of points $\{y_i\}$ and the forcing discretization size, $N_f$, are chosen arbitrarily to train the NNs as the forcing terms are assumed to be known over $\Omega$. In practice, the number, $N_u$, and location of the measurement points, $\{x_i\}$, are imposed by the nature of the experiment, or simulation, performed to measure the system's response. When $\Omega$ is an interval, we always select $N_f=200$, $N_u=100$, and equally-spaced sampled points for the forcing and response functions. 

Unless otherwise stated, the training data comprises $N=100$ forcing and solution pairs, where the forcing terms are drawn at random from a Gaussian process, $\mathcal{GP}(0,K_{\text{SE}})$, where $K_{\text{SE}}$ is the squared-exponential covariance kernel~\cite{rasmussen2006gaussian} defined as
\begin{equation} \label{eq_kernel}
K_{\text{SE}}(x,y) = \exp\left(-\frac{|x-y|^2}{2\ell^2}\right),\quad x,y\in\Omega.
\end{equation}
As discussed in \cref{sec_cov_kernel}, the parameter $\ell>0$ in \cref{eq_kernel} is called the length-scale parameter, and characterizes the correlation between the values of $f\sim \mathcal{GP}(0,K_{\text{SE}})$ at $x$ and $y$ for $x,y\in \Omega$. A small parameter, $\ell$, yields highly oscillatory random functions, $f$, and determines the ability of the GP to generate a diverse set of training functions. This last property is crucial for capturing different modes within the operator, $\L$, and for learning the associated Green's function accurately~\cite{boulle2021learning}. Other possible choices of covariance kernels include the periodic kernel~\cite{rasmussen2006gaussian}:
\[K_{\text{Per}}(x,y)=\exp\left(-\frac{2\sin^2(\pi|x-y|)}{\ell^2}\right),\quad x,y\in\Omega,\]
which is used to sample periodic random functions for problems with periodic boundary conditions (\cref{fig_mean}(b)). Another possibility is a kernel from the Mat\'ern family~\cite{rasmussen2006gaussian} or the Jacobi kernel introduced in \cref{sec_Jacobi_kernel}.

\begin{figure}[htbp]
\centering
\vspace{0.5cm}
\begin{overpic}[width=0.8\textwidth]{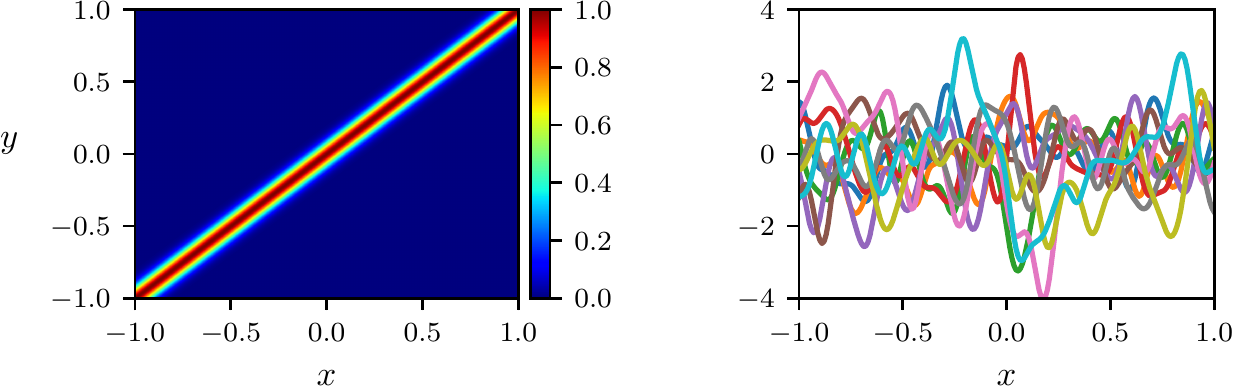}
\put(-1,30){(a)}
\put(54,30){(b)}
\end{overpic}
\caption{Random forcing terms. (a) Squared exponential covariance kernel $K_{\textup{SE}}$ on $[-1,1]^2$ with normalized length-scale $\lambda=0.03$ (b) 10 functions sampled from the Gaussian process $\mathcal{GP}(0,K_{\textup{SE}})$.}
\label{fig_covariance_kernel}
\end{figure}

When $\Omega$ is an interval $[a,b]$, we introduce a normalized length-scale parameter $\lambda = \ell/(b-a)$, so that the method described does not depend on the length of the interval. In addition, we choose $\lambda=0.03$, so that the length-scale, $\ell$, is larger than the forcing spatial discretization size, which allows us to adequately resolve the functions sampled from the GP with the discretization. More precisely, we make sure that $\ell\geq (b-a)/N_f$ so that $1/N_f \leq \lambda$. In \cref{fig_covariance_kernel}, we display the squared-exponential covariance kernel on the domain $\Omega=[-1,1]$, along with ten random functions sampled from $\mathcal{GP}(0,K_{\text{SE}})$.

When~\eqref{eq_problem} is a boundary-value problem, we generate training pairs by solving the PDE with a spectral method~\cite{trefethen2000spectral} using the Chebfun software system~\cite{driscoll2014chebfun}, written in MATLAB, and using a tolerance of $5\times 10^{-13}$. We also solve the homogeneous problem with zero-forcing, to compare the learned and exact homogeneous solutions. The exact homogeneous solution is not included in the training dataset. When the homogeneous solution is zero, the solutions, $\{u_j\}_{j=1}^{N}$, and forcing terms, $\{f_j\}_{j=1}^N$, are rescaled, so that $\max_{1\leq j\leq N} \|u_j\|_{L^\infty(\Omega)}=1$. By doing this, we facilitate the training of the NNs by avoiding disproportionately small-scale or large-scale data. In the presence of real data, with no known homogeneous solution, one could instead normalize the output of the NNs, $\N_G$ and $\N_{\text{hom}}$, to facilitate the training procedure.

\subsection{Rational neural networks} \label{sec_rational_net}

As introduced in \cref{chapt_rational}, rational NNs consist of NNs with adaptive rational activation functions $x\mapsto\sigma(x) = p(x)/q(x)$, where $p$ and $q$ are two polynomials, whose coefficients are trained at the same time as the other parameters of the networks, such as the weights and biases. These coefficients are shared between all the neurons in a given layer but generally differ between the network's layers. This type of network was proven to have better approximation power than standard Rectified Linear Unit (ReLU) networks~\cite{glorot2011deep,yarotsky2017error}, which means that they can approximate smooth functions more accurately with fewer layers and network parameters (see~\cref{sec_th_result}). It is also observed in \cref{sec_experiments} that rational NNs require fewer optimization steps in practice and therefore can be more efficient to train than other activation functions.

The NNs, $\N_G$ and $\N_{\text{hom}}$, which approximate the Green's function and homogeneous solution associated with \cref{eq_problem}, respectively, are chosen to be rational NNs with 4 hidden layers and 50 neurons in each layer. We choose the polynomials, $p$ and $q$, within the activation functions to be of degree 3 and 2, respectively, and initialize the coefficients of all the rational activation functions so that they are the best $(3,2)$ rational approximant to a ReLU (see \cref{sec_experiments} for details). The motivation is that the flexibility of the rational functions brings extra benefit in the training and accuracy over the ReLU activation function. We highlight that the increase in the number of trainable parameters, due to the adaptive rational activation functions, is only linear with respect to the number of layers and negligible compared to the total number of parameters in the network as:
\[\text{number of rational coefficients} = 7\times \text{number of hidden layers} = 28.\]
The weight matrices of the NNs are initialized using Glorot normal initializer~\cite{glorot2010understanding}, while the biases are initialized to zero.

Another advantage of rational NNs is the potential presence of poles, \emph{i.e.}, zeros of the polynomial $q$. While the initialization of the activation functions avoids training issues due to potential spurious poles, the poles can be exploited to learn physical features of the differential operator (see \cref{sec_singularity}). Therefore, the architecture of the NNs also supports the aim of a human-understandable approach for learning PDEs. In higher dimensions, such as $d = 2$ or $d =3$, the Green's function is not necessarily bounded along the diagonal, \emph{i.e.}, $\{(x,x),\, x\in\Omega\}$; thus making the poles of the rational NNs crucial.

Finally, we emphasize that the enhanced approximation properties of rational NNs make them ideal for learning Green's functions and, more generally, approximating functions within regression problems. These networks may also be of benefit to other approaches for solving and learning PDEs with DL techniques, such as DeepGreen~\cite{gin2020deepgreen}, Neural operator~\cite{li2020neural}, Fourier neural operator~\cite{li2020fourier}, DeepONet~\cite{lu2021learning}, and PINNs~\cite{raissi2019physics}.

\subsection{Loss function} \label{sec_loss}

The NNs, $\N_G$ and $\N_{\text{hom}}$, are trained by minimizing a mean square relative error (in the $L^2$-norm) regression loss, defined as:
\begin{equation} \label{eq_loss}
\text{Loss} = \frac{1}{N}\sum_{j=1}^N\frac{1}{\|u_j\|_{L^2(\Omega)}^2}\int_{\Omega} \left(u_j(x) - \mathcal{N}_{\text{hom}}(x) - \int_\Omega \mathcal{N}_G(x,y)f_j(y)\,\textup{d} y\right)^2\,\textup{d} x.
\end{equation}
Unless otherwise stated, the integrals in \cref{eq_loss} are discretized by a trapezoidal rule~\cite{suli2003introduction} using training data values that coincide with the forcing discretization grid, $\{y_i\}_{i=1}^{N_f}$, and measurement points, $\{x_i\}_{i=1}^{N_u}$. As an example, for $1\leq j\leq N$, the squared $L^2$-norm of $u_j$, on a one-dimensional domain $\Omega=[a,b]\subset\R$, is approximated as
\[\|u_j\|_{L^2(\Omega)}^2 = \int_{a}^b u_j(x)^2\d x\approx \sum_{i=2}^{N_u}\frac{u_j(x_{i-1})^2+u_j(x_{i})^2}{2}\Delta_{x_i},\]
where $\Delta_{x_i}=x_i-x_{i-1}$ is the length of the $i$th subinterval $[x_{i-1},x_i]$.

Later in \cref{sec_location}, we compare the results obtained by using trapezoidal integration, described above, and a Monte-Carlo integration~\cite{binder2012monte}:
\[\|u_j\|_{L^2(\Omega)}^2 \approx \frac{b-a}{N_u}\sum_{i=1}^{N_u}u_j(x_{i})^2,\]
which has a lower convergence rate to the integral with respect to the number of points, $N_u$. This integration technique is, however, particularly suited for approximating integrals in high dimensions, or with complex geometries~\cite{binder2012monte}. One could also use a mesh of the domain and compute the integrals with a quadrature rule on each cell.

It is also possible to incorporate some prior knowledge about the Green's function in the loss function, by adding a penalty term. If the differential operator is self-adjoint, then depending on the constraint operator $\mathcal{D}$, the associated Green's function is symmetric, \emph{i.e.}, $G(x,y) = G(y,x)$ for all $x,y\in\Omega$. In this case, one can train a symmetric NN $\N_G$ defined as 
\[
\N_G(x,y) = \mathcal{N}(x,y) + \mathcal{N}(y,x), \quad x,y\in\Omega.
\]
However, our numerical experiments reveal that the NNs can learn both boundary conditions and symmetry properties directly, from the training data, without additional constraints on the loss function or network architectures.

\subsection{Optimization algorithm} \label{sec_optimization}

The NNs are implemented with single-precision floating-point format within the TensorFlow DL library~\cite{abadi2016tensorflow}, and are trained\footnote{The numerical experiments are performed on a desktop computer with a Intel\textsuperscript{\tiny\textregistered} Xeon\textsuperscript{\tiny\textregistered} CPU E5-2667 v2 @ 3.30GHz and a NVIDIA\textsuperscript{\tiny\textregistered} Tesla\textsuperscript{\tiny\textregistered} K40m GPU.} using a two-step optimization procedure to minimize the loss function. First, we use Adam's algorithm~\cite{kingma2014adam} for the first 1000 optimization steps (or epochs), with default learning rate $0.001$ and parameters $\beta_1 = 0.9$, $\beta_2 = 0.999$. Then, we employ the limited memory BFGS, with bound constraints (L-BFGS-B) optimization algorithm~\cite{byrd1995limited,liu1989limited}, implemented in the SciPy library~\cite{virtanen2020scipy}, with a maximum of $5\times 10^4$ iterations. This training procedure is used by Lu \emph{et al.} to train physics-informed NNs (PINNs) and mitigate the risk of the optimizer getting stuck at poor local minima~\cite{lu2021deepxde}. The L-BFGS-B algorithm is also successful for PDE learning~\cite{raissi2018deep} and PDE solvers using DL techniques~\cite{lu2021deepxde,raissi2019physics}. Moreover, this optimization algorithm takes advantage of the smoothness of the loss function by using quasi-Newton approximations to second-order derivatives and often converges in fewer iterations than Adam's algorithm and other methods based on stochastic gradient descent~\cite{lu2021deepxde}. Within this setting, rational NNs are beneficial because the activation functions are smooth while maintaining an initialization close to ReLU.

\begin{figure}[htbp]
\centering
\begin{overpic}[width=0.5\textwidth, trim=0 0 0 0,clip]{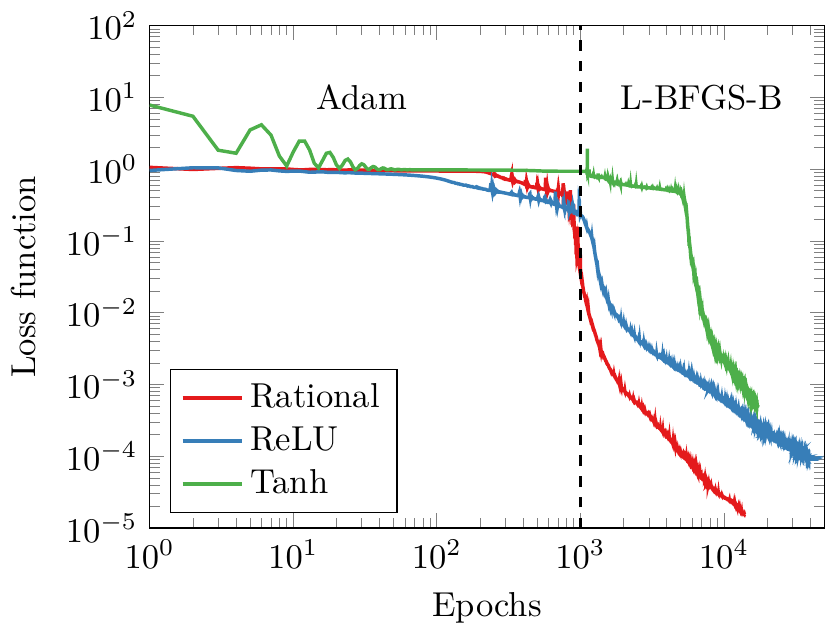}
\end{overpic}
\caption{Loss function during training. Loss function magnitudes for the ReLU, tanh, and rational NNs with respect to the number of epochs. The networks are trained to learn the Green's function of the Helmholtz operator with homogeneous Dirichlet boundary conditions and frequency $K=15$. Adam's optimizer is used until 1000 epochs (before the dashed line) and L-BFGS-B is employed thereafter.}
\label{fig_loss}
\end{figure}

In \cref{fig_loss}, we display the value of the loss function during the training of the NNs with different activation functions: rational, ReLU, and hyperbolic tangent (tanh). In this example, we aim to learn the Green's function of a high-frequency Helmholtz operator with homogeneous Dirichlet boundary conditions on $\Omega = [0,1]$:
\begin{equation} \label{eq_helmholtz}
\L u = \frac{d^2u}{dx^2}+K^2 u,\quad u(0)=u(1)=0,
\end{equation}
where $K=15$ denotes the Helmholtz frequency. Note that the operator defined by \cref{eq_helmholtz} is indefinite but invertible. We first remark in \cref{fig_loss} that the rational NN is easier to train than the other NNs, as it minimizes the loss function to $10^{-5}$ with $\approx 15000$ epochs, while a ReLU NN requires three times as many epochs to reach $10^{-4}$. We also see that the loss function for the ReLU and rational NN becomes more oscillatory~\cite{bengio2012practical} and harder to minimize before epoch 1000, while it converges much faster after switching to L-BFGS-B. In theory, one could introduce a variable learning rate that improves the behavior of Adam's optimizer~\cite{george2006adaptive,smith2017cyclical}. However, that introduces an additional parameter, which is not desirable in the context of PDE learning. We aim to design an adaptive and easy-to-use method that does not require extensive hyperparameter tuning. We also observe that the tanh NN has a similar convergence rate to the rational NN due to the smoothness of the activation function, but this network exhibits instability during training, as indicated by the high value of the loss function when the optimization terminates. Rational NNs do not suffer from this issue, thanks to the initialization close to a ReLU NN, as can be observed in \cref{fig_loss}, when focusing on the value of the loss function corresponding to the early optimization steps.

\subsection{Measuring the results} \label{sec_measure_res}

Once the NNs have been trained, we visualize the Green's functions by sampling the networks on a fine $1000\times 1000$ grid of $\Omega\times\Omega$. In the case where the exact Green's function $G_{\text{exact}}$ is known, we measure the accuracy of the trained NN, $\N_G$, using a relative error in the $L^2$-norm:
\begin{equation} \label{eq_rel_error}
\text{Relative Error} = 100\times \|G_{\text{exact}}-\N_G\|_{L^2(\Omega)}/\|G_{\text{exact}}\|_{L^2(\Omega)}.
\end{equation}
Here, we multiplied by 100 to obtain the relative error as a percentage (\%). This illustrates an additional advantage of using a Green's function formulation: we can create test case problems with known Green's functions and evaluate the method using relative error and offer performance guarantees on benchmark problems. The standard approaches in the literature often use best-case and worst-case examples as testing procedures and therefore do not guarantee that the solution operator is accurately learned. The ``worst-case'' examples can be misleading if they consist of functions with similar behavior to the forcing terms already included in the training dataset. Furthermore, since the space of possible forcing terms is of infinite dimension, it is not possible to evaluate the trained NNs on all these functions to obtain a true worst-case example.

\section{Robustness of the method} \label{sec_robustness}

We test the robustness of our DL method for learning Green's functions and homogeneous solutions of differential equations, with respect to the number of training pairs, the discretization of the solutions and forcing terms, and the noise perturbation of the training solutions, $\{u_j\}_{j=1}^N$. For consistency, we perform numerical experiments where we learn the Green's function of the Helmholtz operator with parameter $K=15$ and homogeneous Dirichlet boundary conditions (see \cref{eq_helmholtz}). The performance is measured using the relative error in the $L^2$-norm defined in \cref{eq_rel_error} between the trained network, $\N_G$, and the exact Green's function, $G_{\text{exact}}$, whose analytic expression is given by
\[G_{\text{exact}}(x,y) =
\begin{cases}
\frac{\sin(15x)\sin(15(y-1))}{15\sin(15)}, & \text{if } x\leq y,\\
\frac{\sin(15y)\sin(15(x-1))}{15\sin(15)}, & \text{if } x> y,\\
\end{cases}
\]
where $x,y\in[0,1]$.

\subsection{Influence of the activation function on the accuracy} \label{sec_act_accuracy}

We first compare the performances of different activation functions for learning the Green's functions of the Helmholtz operator by training the NNs, $\N_G$ and $\N_{\text{hom}}$, with rational, ReLU, and tanh activation functions. The numerical experiments are repeated ten times to study the statistical effect of the random initialization of the network weights and the stochastic nature of Adam's optimizer. The rational NN achieves a mean relative error of $1.2\%$ (with a standard deviation of $0.2\%$), while the ReLU NN reaches an average error of $3.3\%$ (with a standard deviation of $0.2\%$), which is about three times larger. Note that the ten times difference in the loss function between ReLU and Rational NNs, displayed in \cref{fig_loss}, is consistent with the factor of three in the relative error since the loss is a mean squared error and $\sqrt{10}\approx 3$. This indicates that the rational neural networks are not overfitting the training dataset. One of the numerical experiments with a tanh NN terminated early due to the training instabilities mentioned in \cref{sec_loss}, achieving a relative error of $99\%$. We excluded this problematic run when comparing the ReLU and rational NN's accuracy, limiting ourselves only to cases where the training was successful. The ReLU and rational NNs did not suffer from such issues and were always successful. The averaged relative error of the tanh NN, over the nine remaining experiments, is equal to $3.9\%$ (with a standard deviation of $1.4\%$), which is slightly worse than the ReLU NN, with higher volatility of the results.

\begin{figure}[htbp]
\centering
\vspace{0.5cm}
\begin{overpic}[width=0.8\textwidth]{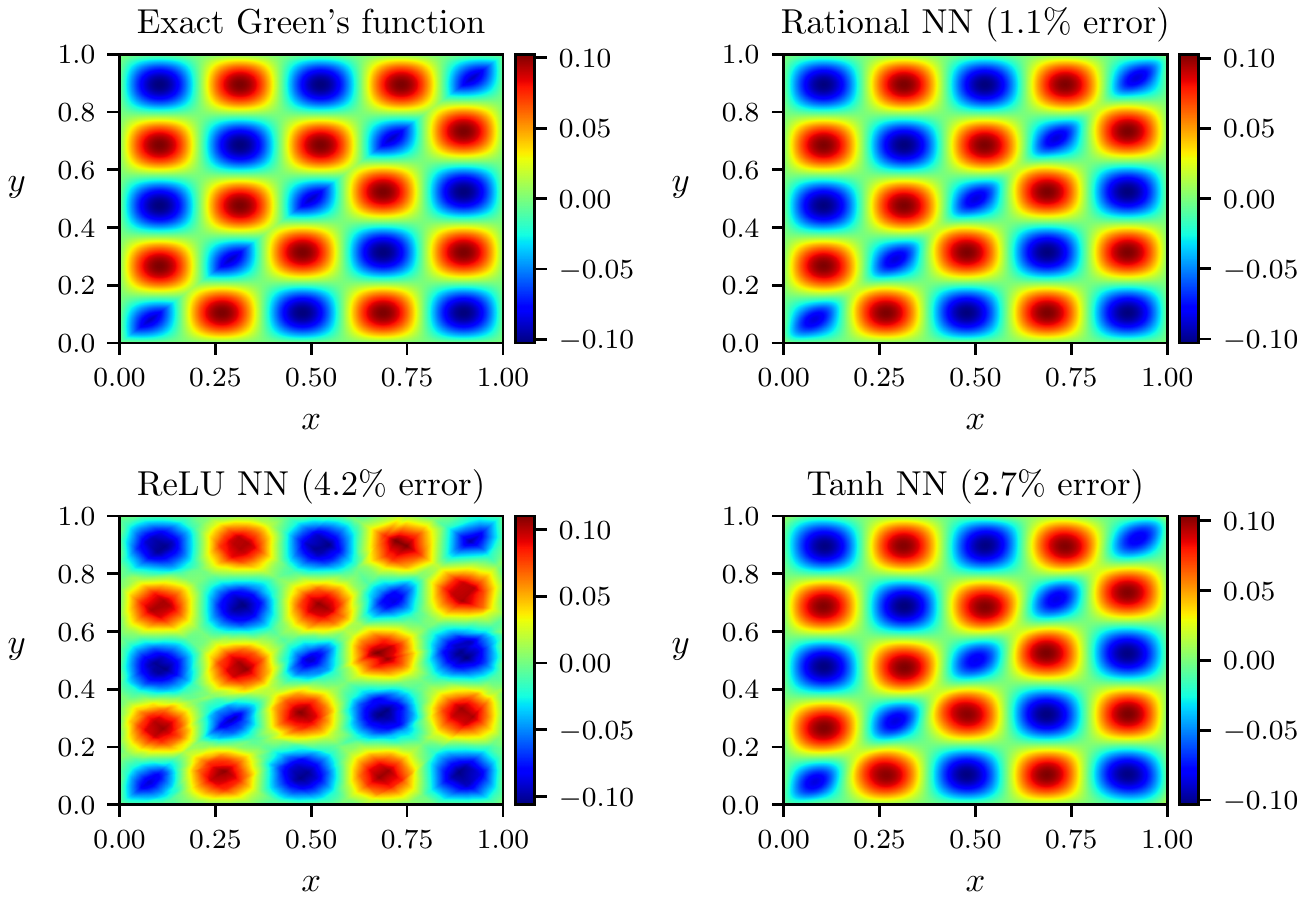}
\end{overpic}
\caption{Comparison of activation functions. Exact and learned Green's functions of the Helmholtz operator by a rational, ReLU, and tanh NN. The relative error in the $L^2$ norm is reported in the titles of the panels.}
\label{fig_activation_functions}
\end{figure}

The exact and learned Green's functions with rational, ReLU, and tanh NNs are displayed in \cref{fig_activation_functions}. We see that the rational and tanh NNs are smooth approximations of the exact Green's function, while visual artifacts are present for the ReLU NN as it is piecewise linear, despite its good accuracy.

\subsection{Number of training pairs and spatial measurements} 

\begin{figure}[htbp]
\centering
\vspace{0.5cm}
\begin{overpic}[width=\textwidth, trim=0 0 0 0,clip]{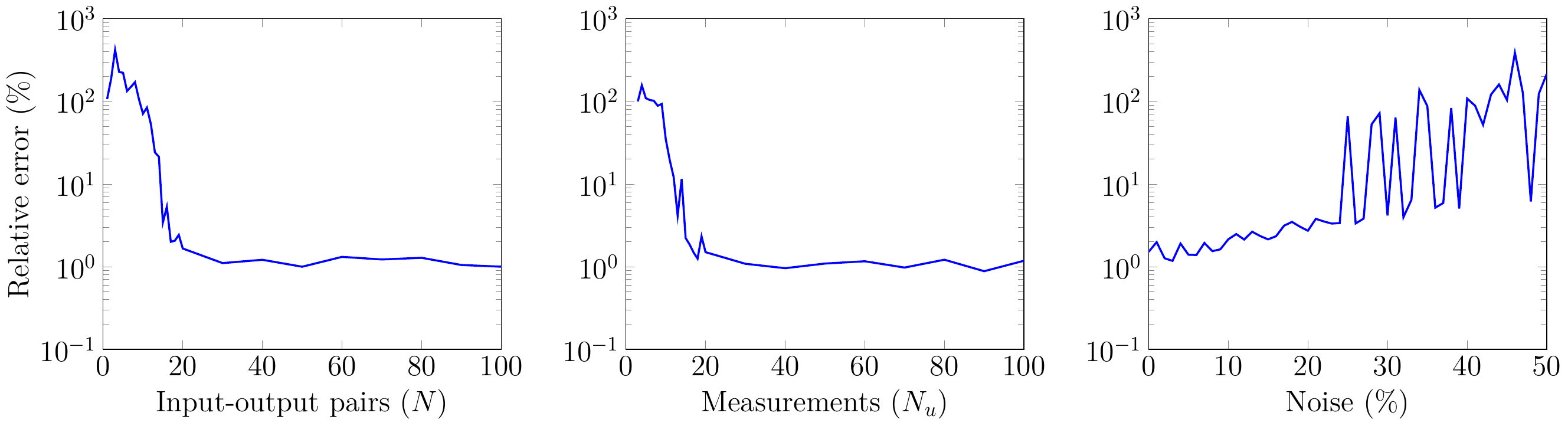}
\put(0,26){(a)}
\put(33.2,26){(b)}
\put(67,26){(c)}
\end{overpic}
\caption{Robustness of the method. Relative error of the learned Green's function of the Helmholtz operator with respect to the number of input-output pairs (a), spatial measurements (b), and level of Gaussian noise perturbation (c).}
\label{fig_robustness}
\end{figure}

This section describes our method's accuracy as we change the number of training pairs and the size of the spatial discretization. First, we fix the number of spatial measurements to be $N_u=100$, and then vary the number of input-output pairs, $\{(f_j,u_j)\}_{j=1}^N$, of the training dataset for the Helmholtz operator with Dirichlet boundary conditions (see \cref{eq_helmholtz}). As we increase $N$ from $1$ to $100$, we report the relative error of the Green's function learned by a rational NN in \cref{fig_robustness}(a). Next, in \cref{fig_robustness}(b), we display the relative error on the learned Green's function as we increase $N_u$ from $3$ to $100$, with $N=100$ input-output pairs. Note that we only perform the numerical experiments once since we obtained a low variation of the relative errors in \cref{sec_act_accuracy} when the networks, $\N_G$ and $\N_{\text{hom}}$, have rational activation functions. We observe similar behavior in \cref{fig_robustness}(a) and (b), where the relative error first rapidly (exponentially) decreases as we increase the number of functions in our dataset or spatial measurements of the solutions to the Helmholtz equations with random forcing terms. One important thing to notice is our method's ability to learn the Green's function of a high-frequency Helmholtz operator, with only $1\%$ relative error, using very few training pairs. The learning rate of our deep learning technique for this operator appears to be poly-logarithmic, \emph{i.e.},~the number of input-output pairs required to learn the Green's function within accuracy $0<\epsilon<1$ behaves like $\mathcal{O}(\text{polylog}(1/\epsilon))$, as predicted by the remark in \cref{sec_stable_H_reconst}.

The performance reaches a plateau at $N\approx 20$ and $N_u\approx 20$, respectively, and ceases to improve. However, the stagnation of the relative error for more numerous training data and spatial measurements is expected and can be explained by our choice of covariance kernel length-scale, which restricts the GP's ability to generate a wide variety of forcing terms. Following \cref{sec_generation_data}, we chose a normalized length-scale parameter $\lambda=0.03$, which yields approximately $20$ eigenvalues greater than $10^{-2}$. This issue can be resolved by decreasing the length-scale parameter and concomitantly increasing the forcing discretization size or choosing another covariance kernel with a less pronounced eigenvalue decay rate (see~\cref{sec_cov_kernel}). In summary, the number of spatial measurements should be larger than $1/\lambda$ to resolve the forcing terms and the number of input-output pairs should correspond to the number of covariance kernel eigenvalues greater than $10^{-2}$.

\subsection{Noise perturbation} 

The impact of noise in the training dataset on the accuracy of the learned Green's function is gauged experimentally by perturbing the system's response measurements with Gaussian noise as
\begin{equation} \label{eq_noise}
u_j^{\text{noise}}(x_i) = u_j(x_i)(1+\delta c_{i,j}),
\end{equation}
where the coefficients $c_{i,j}$ are i.i.d.~Gaussian random variables for $1\leq i\leq N_u$ and $1\leq j\leq N$, and $\delta$ denotes the noise level (in percent). We then vary the level of Gaussian noise perturbation from $0\%$ to $50\%$, train the NNs, $N_G$ and $N_{\text{hom}}$, for each choice of the noise level, and report the relative error in \cref{fig_robustness}(c). We first observe a low impact of the noise level on the accuracy of the learned Green's function, as a perturbation of the system's responses measurements with $20\%$ noise only increases the relative error by a factor of $2$ from $1.5\%$ (no noise) to $2.7\%$. When the level of noise exceeds $25\%$, we notice large variations of the relative errors and associated higher volatility in results, characterized by a large standard deviation in error associated with repeated numerical experiments. We consider our DL approach relatively robust to noise in the training dataset.

\subsection{Location of the measurements} \label{sec_location}

As described in \cref{sec_generation_data}, by default, we use a uniform grid for spatial measurements of the training dataset, and thus we discretize the integrals in the loss function (cf.~\cref{eq_loss}) using a trapezoidal rule. We conducted additional numerical experiments on the Helmholtz example to study the influence of the measurements' location and quadrature rule on the relative error of the learned Green's function. We report the relative errors between the learned and exact Green's functions in \cref{tab_quad}, using a Monte-Carlo or a trapezoidal rule to approximate the integrals and uniform or random spatial measurements. In the latter case, the measurement points $\{x_i\}_{i=1}^{N_u}$ are independently and identically sampled from a uniform distribution, $\mathcal{U}(0,1)$, where $\Omega=[0,1]$ is the domain. We find that the respective relative errors vary between $0.96\%$ and $1.3\%$. Therefore, we do not observe statistically significant differences in the relative error computed by rational NNs. These results support the claim that our method is relatively robust to the type of spatial measurements in the training dataset.

\begin{table}[htbp]
\centering
\caption{Choice of quadrature rules. Relative error of the Green's function of the Helmholtz operator with frequency $K=15$ learned by a rational NN with respect to the type of spatial measurements and quadrature rule (Monte-Carlo or trapezoidal rule) used.}
\label{tab_quad}
\begin{tabular}{c|cc}
\hline
Spatial measurements & Monte-Carlo & Trapezoidal rule \\
\hline
Random & $1.1\%$ & $1.3\%$\\
Uniform & $1.3\%$ & $0.96\%$ \\
\hline
\end{tabular}
\end{table}

\subsection{Missing measurements data} \label{sec_missing_data}

\begin{figure}[htbp]
\centering
\vspace{0.5cm}
\begin{overpic}[width=0.8\textwidth]{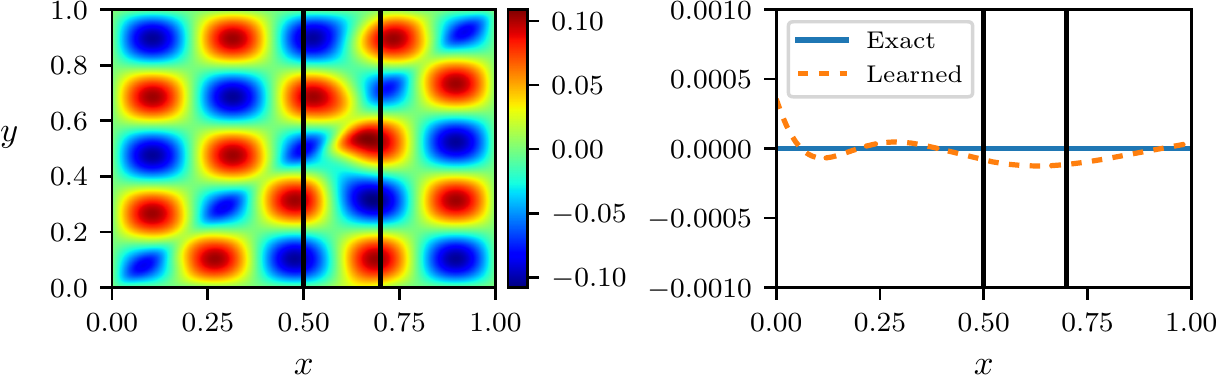}
\put(-1,30){(a)}
\put(50.3,30){(b)}
\end{overpic}
\caption{Gap in measurements. (a) Green's function of the Helmholtz operator and its homogeneous solution (b) learned by a rational NN with no measurement points for $x\in[0.5,0.7]$. The space between the vertical black lines indicates where there is a lack of spatial measurements.}
\label{fig_gap}
\end{figure}

Since experimental data may be partially corrupted or unavailable at some spatial locations, we assess our method's accuracy with respect to missing measurement data in the training dataset. We consider the high-frequency Helmholtz operator, defined on the domain $\Omega = [0,1]$ by \cref{eq_helmholtz}, with homogeneous Dirichlet boundary conditions. We introduce a gap in the spatial measurements located at $x\in[0.5,0.7]$ by sampling the system's responses, $\{u_j\}_{j=1}^N$, uniformly on the domain $[0,0.5]\cup [0.7,1]$. Note that the forcing terms, $\{f_j\}_{j=1}^N$, are still sampled uniformly on the whole domain since they are assumed to be known. The Green's function and homogeneous solution learned by the rational NNs are displayed in \cref{fig_gap}(a) and (b), respectively. Surprisingly, we find that the NN, $\N_G$, can capture the high-frequency pattern of the Green's function and achieves a relative error of $8.2\%$, despite the large gap within the measurement data for $x\in[0.5,0.7]$. Another interesting outcome of this numerical experiment is that the lack of spatial measurements in a specific interval does not influence the accuracy of our method outside this location, \emph{i.e.}, for $x\in [0,0.5]\cup [0.7,1]$ and $y\in [0,1]$ in this example. This phenomenon might be explained by the existence of non-local effects when expressing the solution operator associated with the PDE as an integral operator.

\section{Human-understandable features} \label{sec_features}

\begin{figure}[htbp]
\centering
\vspace{0.5cm}
\begin{overpic}[width=\textwidth, trim=0 0 0 0,clip]{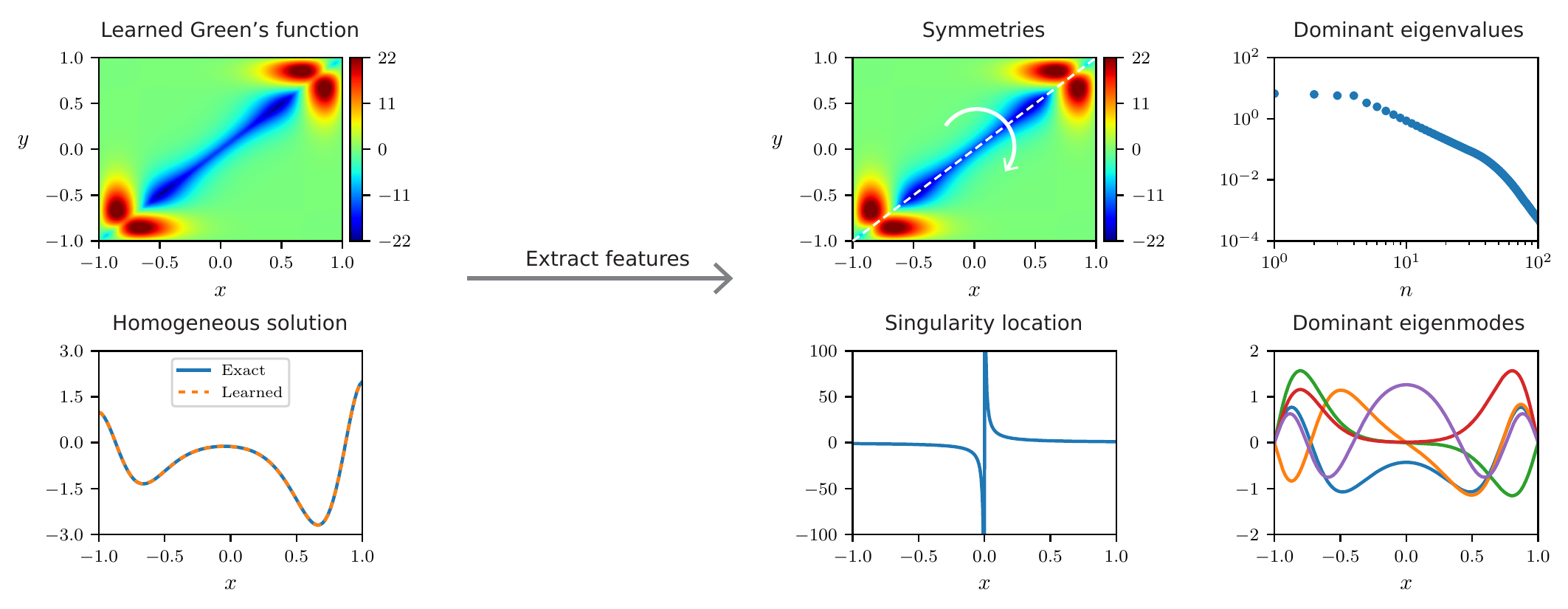}
\put(0.5,37){(a)}
\put(0.5,18){(b)}
\put(48,37){(c)}
\put(48,18){(d)}
\put(75,37){(e)}
\put(75,18){(f)}
\end{overpic}
\caption{Feature extraction from learned Green's functions. The NNs for the learned Green's function (a) and homogeneous solution (b) enable the extraction of qualitative and quantitative features associated with the differential operator. For example, the  symmetries in the Green's function reveal PDE invariances (c), poles of rational NNs identify singularity type and location (d), the dominant eigenvalues (e) and eigenmodes (f) of the learned Green's function are related to the eigenvalues and eigenmodes of the differential operator.}
\label{fig_schema}
\end{figure}
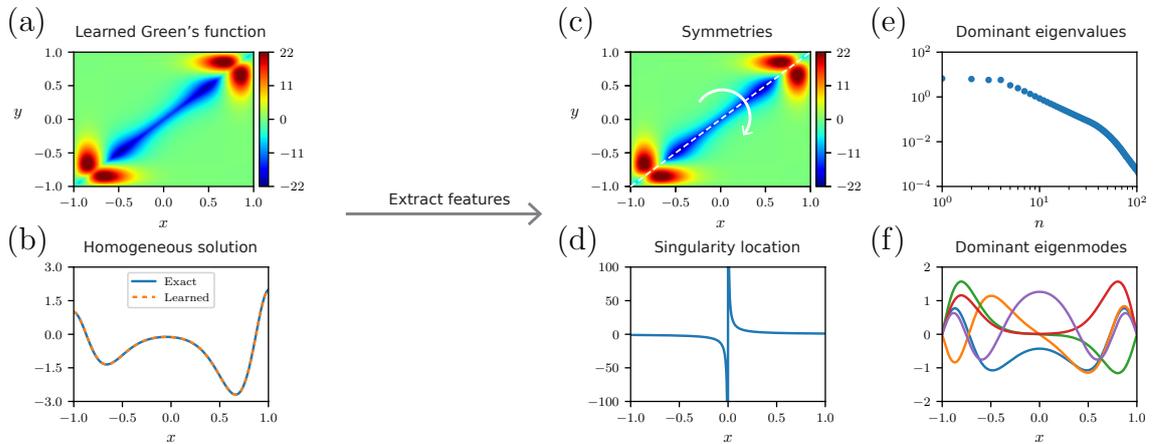

The trained NNs contain both the desired Green's function and homogeneous solution, which we evaluate and visualize to glean novel insights concerning the underlying governing PDE (\cref{fig_schema}). In this way, we achieve one part of our human interpretation goal: finding a link between the properties of the Green's function and that of the underlying differential operator and solution constraints. 

As an example, if the Green's function is symmetric, \emph{i.e.}, $G(x,y)=G(y,x)$ for all $x,y\in\Omega$, then the operator $\L$ is self-adjoint. Another aspect of human interpretability is that the poles of the trained rational NN tend to cluster in a way that reveal the location and type of singularities in the homogeneous solution, discussed further in \cref{sec_singularity} below. Finally, there is a direct correspondence between the dominant eigenmodes and eigenvalues (as well as the singular vectors and singular values) of the learned Green's function and those of the differential operator. The correspondence gives insight into the important eigenmodes that govern the system's behavior (see~\cref{sec_eig,sec_svd} below). This section highlights that several features of the differential operators can be extracted from the learned Green's function, which supports our aim of uncovering mechanistic understanding from input-output pairs of forcing terms and solutions.

\subsection{Linear constraints and symmetries} \label{sec_linear_const}

We first remark that boundary constraints, such as the constraint operator, $\D$, of \cref{eq_problem}, can be recovered from the Green's function, $G$, of the differential operator, $\L$. Let $f\in C_c^\infty(\Omega)$ be any infinitely differentiable function with compact support on $\Omega$, and $u$ be the solution to \cref{eq_problem} with forcing term, $f$, such that
\[u(x) = \int_{\Omega} G(x,y)f(y)\d y + u_{\text{hom}}(x),\quad x \in\Omega.\]
Under sufficient regularity conditions, the linearity of the operator, $\D$, implies that $\D(G(\cdot,y),\Omega)=0$ for all $y\in\Omega$. For instance, if $\D$ is the Dirichlet operator: $\D(u,\Omega) = u_{|\partial\Omega}$, then the Green's function satisfies $G(x,y) = 0$ for all $x\in\partial\Omega$.

\begin{figure}[htbp]
\centering
\vspace{0.5cm}
\begin{overpic}[width=0.8\textwidth]{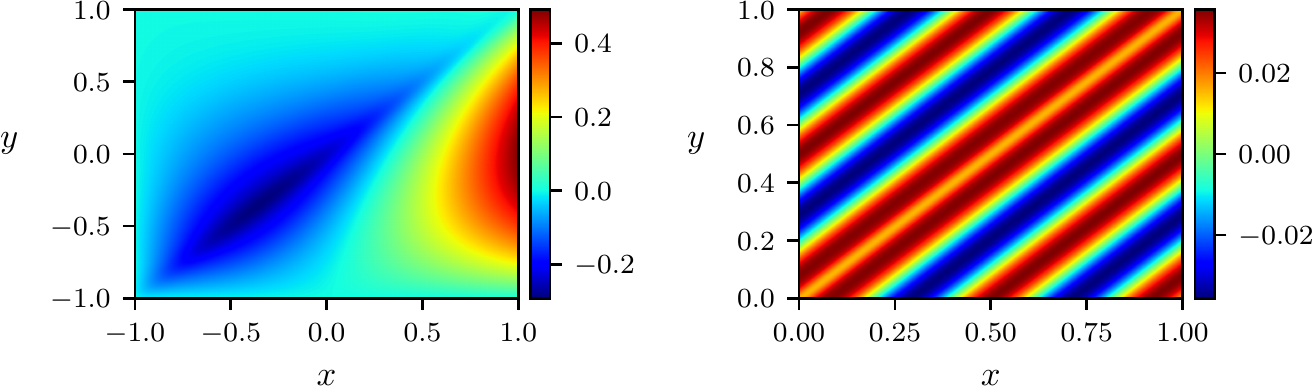}
\put(-2,30){(a)}
\put(50,30){(b)}
\end{overpic}
\caption{Extraction of linear constraints. (a) Learned Green's functions of a second-order differential operator with an integral constraint defined in \cref{eq_mean_condition}. (b) Green's function of the Helmholtz operator with periodic boundary conditions learned by a rational NN.}
\label{fig_mean}
\end{figure}

As an example, we display in \cref{fig_mean}(a) the learned Green's function of the following second-order differential operator on $\Omega=[-1,1]$ with an integral constraint on the solution:
\begin{equation} \label{eq_mean_condition}
\L u=\frac{d u^2}{dx^2}+x^2u, \quad u(-1)=1,\quad \int_{-1}^1 u(x)\d x=2.
\end{equation}
We observe that $G(-1,y)=0$ for all $y\in[-1,1]$ and one can verify that the relation $\int_{-1}^1 G(x,y)\d x=0$ holds for any $y\in[-1,1]$. In a second example, we learn the Green's function of the Helmholtz operator on $\Omega=[0,1]$ with frequency $K=15$ and periodic boundary conditions: $u(0)=u(1)$. One can see in \cref{fig_mean}(b) that the Green's function itself is periodic and that $G(0,y)=G(1,y)$ for all $y\in[0,1]$, as expected. The periodicity of the Green's function in the $y$-direction: $G(x,0)=G(x,1)$ for $x\in[0,1]$, is due to the fact that the Helmholtz operator is self-adjoint, which implies symmetry in the associated Green's function. Furthermore, any linear constraint $\mathcal{C}(u)=0$ such as linear conservation laws or symmetries~\cite{olver2000applications}, satisfied by all the solutions to \cref{eq_problem}, under forcing $f\in C_c^\infty(\Omega)$, is also satisfied by the Green's function, $G$, such that $\mathcal{C}(G(\cdot,y))=0$ for all $y\in\Omega$, and is therefore witnessed by the Green's function.

\subsection{Eigenvalue decomposition} \label{sec_eig}

Let $\L$ be a self-adjoint operator and consider the following eigenvalue problem:
\begin{equation} \label{eq_eig_L}
\L v = \lambda v,\quad \D(v,\Omega) = 0,
\end{equation}
where $v$ is an eigenfunction of the differential operator, $\L$, satisfying the homogeneous constraints with associated eigenvalue, $\lambda> 0$. The eigenfunction, $v$, can be expressed using the Green's function, $G$, of $\L$ as
\[v(x) = \lambda\int_\Omega G(x,y)v(y)\d y,\quad x\in\Omega,\]
which implies that $v$ is also an eigenfunction of the integral operator with kernel $G$, but with eigenvalue $1/\lambda$. Consider now the eigenvalue problem associated with the Green's function, itself:
\[\int_\Omega G(x,y) w(y)\d y = \mu w(x), \quad x\in\Omega,\]
where $\mu>0$. Then, we find that $(w,1/\mu)$ are solutions to the eigenvalue problem~\eqref{eq_eig_L}. Consequently, the differential operator, $\L$, and integral operator with kernel, $G$, share the same eigenfunctions, but possess reciprocal eigenvalues~\cite{stakgold2011green}. Thus, we can effectively compute the lowest eigenmodes of $\L$ from the learned Green's function.

\begin{figure}[htbp]
\centering
\vspace{0.5cm}
\begin{overpic}[width=0.8\textwidth]{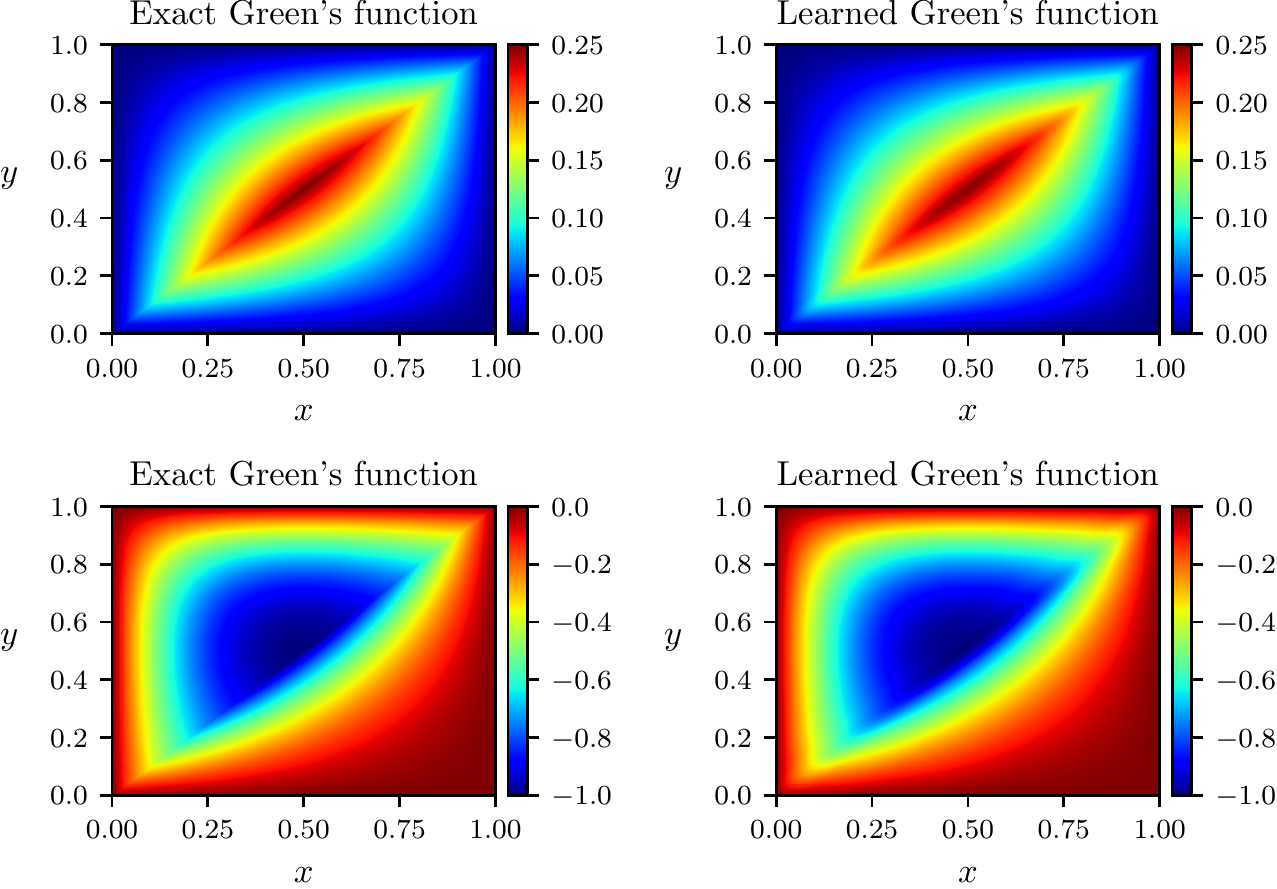}
\put(-2,67){(a)}
\put(-2,31){(b)}
\end{overpic}
\caption{Laplace and advection-diffusion operators. Exact and learned Green's functions of the Laplace (a) and advection-diffusion (b) operators.}
\label{fig_laplace_rational}
\end{figure}

We now evaluate our method's ability to accurately recover the eigenfunctions of the Green's function that are associated with the largest eigenvalues, in magnitude, from input-output pairs. We train a NN to learn the Green's function of the Laplace operator $\L u=-d^2 u/dx^2$ on $[0,1]$, with homogeneous Dirichlet boundary conditions, and numerically compute its eigenvalue decomposition. In \cref{fig_laplace_rational}(a), we display the learned and exact Green's function, whose expression is given for $x,y\in[0,1]$ by 
\[
G_{\text{exact}}(x,y) = 
\begin{cases}
x(1-y), & \text{if } x\leq y,\\
y(1-x), & \text{if } y < x.
\end{cases}
\]

\begin{figure}[htbp]
\centering
\vspace{0.5cm}
\begin{overpic}[height=0.56\textwidth, trim=0 0 0 0,clip]{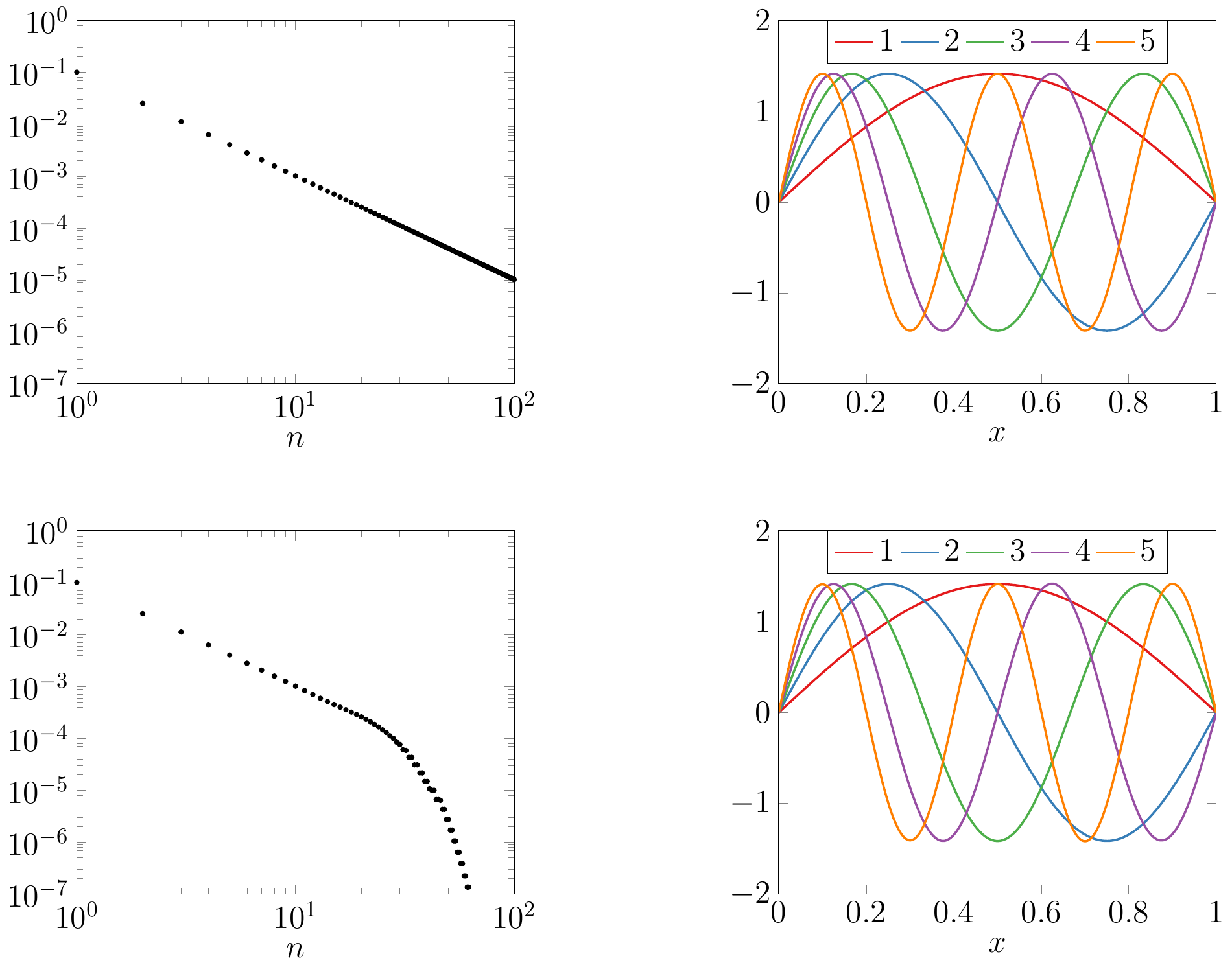}
\put(-5,77){(a)}
\put(-5,35.5){(b)}
\put(22.5,80){$D$}
\put(80,80){$V$}
\end{overpic}
\caption{Eigenvalue decomposition. The first 100 largest eigenvalues and first five eigenfunctions of the exact (a) and learned (b) Green's functions of the Laplace operator. The eigenvalues are represented in the left panels, while the right panels illustrate the first five eigenfunctions of the Green's function.}
\label{fig_laplace_eig}
\end{figure}
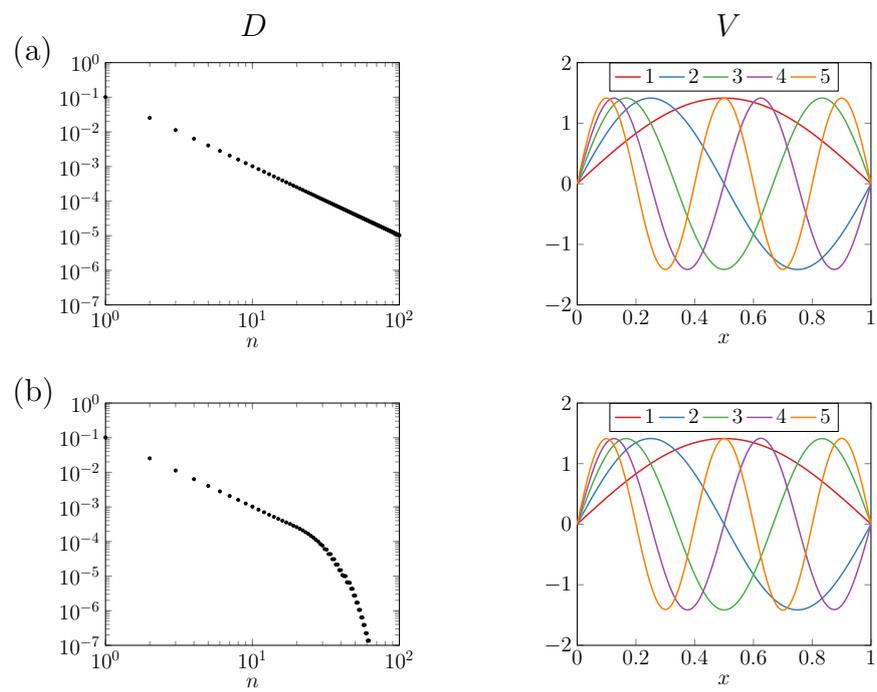

The one hundred largest eigenvalues in magnitude, along with the corresponding first five eigenfunctions, are visualized for the exact and learned Green's functions in \cref{fig_laplace_eig}. Note that the eigenvectors of the learned Green's functions are normalized and flipped to match the ones of the exact Green's function because eigenfunctions are unique up to a scalar multiple when the eigenvalues are all distinct. We find that we can recover the largest eigenvalues and eigenfunctions of the learned Green's function and that the first 20 largest eigenvalues remain accurate. Therefore, the approximation error between the learned and exact Green's functions mainly affects the smallest eigenvalues. This is an essential feature of our method since the dominant eigenmodes of the differential operator $\L$ are associated with the largest eigenvalues of the Green's functions, which can be learned accurately. The exponential decay of the smallest eigenvalues of the learned Green's function in the left panel of \cref{fig_laplace_eig}(b) is because the rational NN is a smooth approximation to the exact Green's function.

\subsection{Singular value decomposition} \label{sec_svd}

When the Green's function of the differentiation operator, $\L$, is square-integrable, its associated Hilbert--Schmidt integral operator admits a singular value decomposition (SVD) (see~\cref{sec_HS}). Then, there exist a positive sequence $\sigma_1\geq \sigma_2\geq \cdots> 0$, and two orthonormal bases, $\{\phi_j\}$ and $\{\psi_j\}$, of $L^2(\Omega)$ such that
\begin{equation} \label{eq_svd}
u(x) = \int_{\Omega}G(x,y)f(y)\d y + u_{\text{hom}}(x) = \sum_{\substack{j=1\\\sigma_j>0}}^\infty \sigma_j \langle\phi_j,f\rangle\psi_j(x) + u_{\text{hom}}(x),\quad x\in\Omega,
\end{equation}
where $u$ is the solution to \cref{eq_problem} with forcing term $f$, and $\langle\cdot,\cdot\rangle$ denotes the inner product in $L^2(\Omega)$. Therefore, the action of the solution operator $f\mapsto u$ can be approximated using the SVD of the integral operator. Similarly to \cref{sec_eig} with the eigenvalue decomposition, the dominant terms in the expansion of \cref{eq_svd} are associated with the largest singular values of the integral operator.

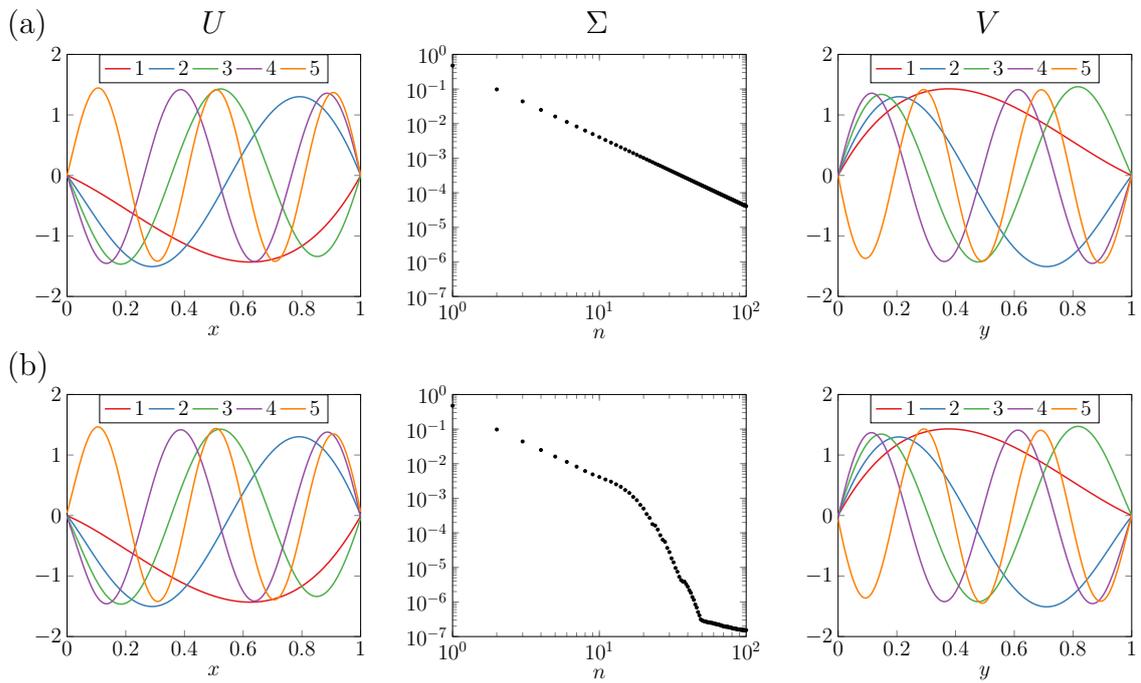
\begin{figure}[htbp]
\centering
\vspace{0.5cm}
\begin{overpic}[height=0.56\textwidth, trim=0 0 0 0,clip]{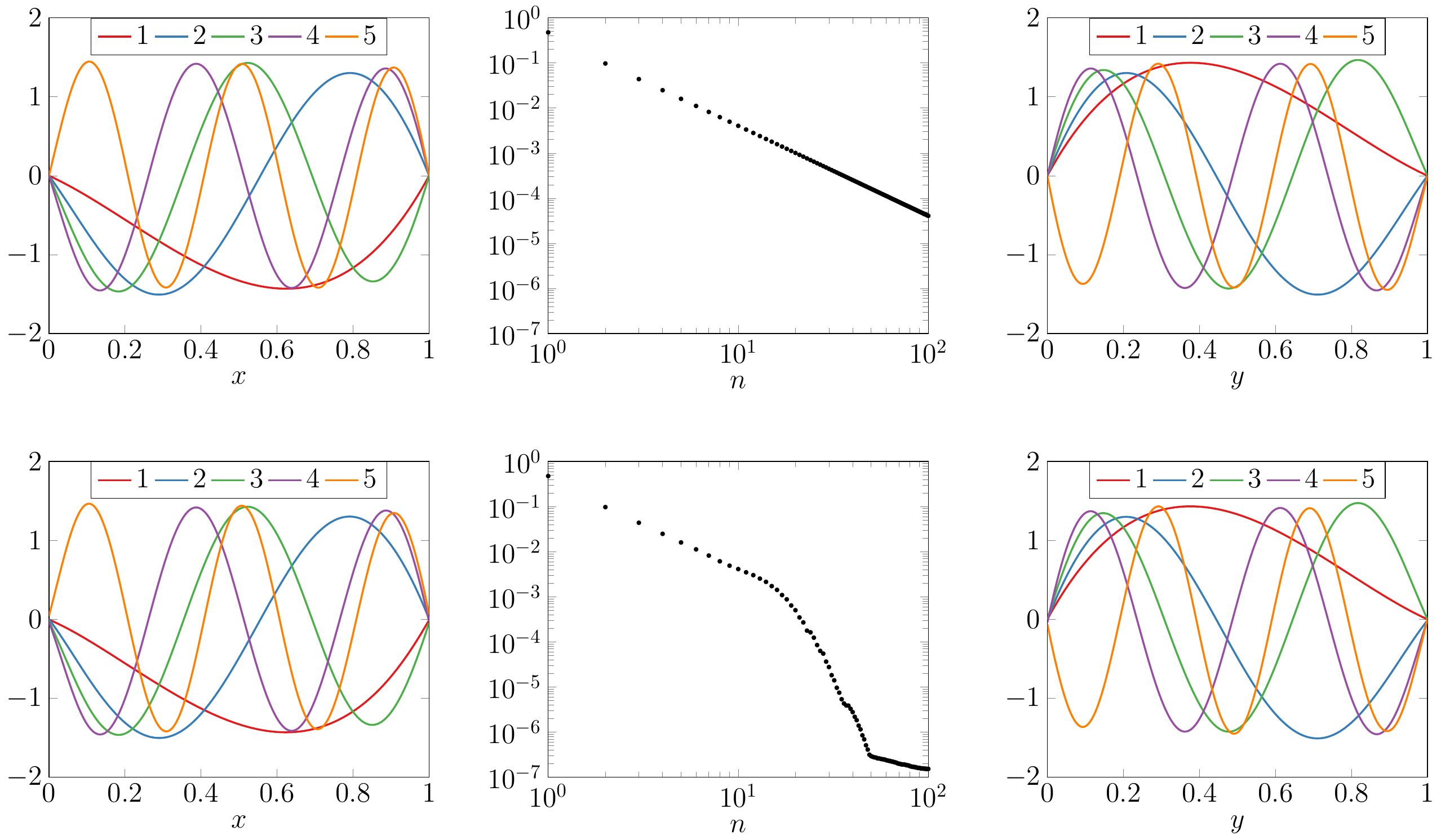}
\put(-2,59){(a)}
\put(-2,28){(b)}
\put(15.5,59){$U$}
\put(50,59){$\Sigma$}
\put(85,59){$V$}
\end{overpic}
\caption{Singular value decomposition. Singular value decomposition of the exact (a) and learned (b) Green's functions of the advection-diffusion operator defined by \cref{eq_adv_diff}. The left and right panels, respectively, show the first five left and right singular vectors, $\{\phi\}_{n=1}^5$ and $\{\psi\}_{n=1}^5$, of the exact and learned Green's functions. The singular values of the Green's functions are plotted in the middle panel.}
\label{fig_adv_diff_SVD}
\end{figure}

We now show that one can accurately recover the first singular values and singular vectors from the Green's function learned by a rational NN. We train a rational NN to learn the Green's function of an advection-diffusion operator $\L$ on $\Omega=[0,1]$ with Dirichlet boundary conditions, defined as
\begin{equation} \label{eq_adv_diff}
\L u=\frac{1}{4}\frac{d^2 u}{dx^2}+\frac{du}{dx}+u,\quad u(0)=1,\, u(1)=-2.
\end{equation}
The learned Green's function is illustrated in \cref{fig_laplace_rational}(a), next to the exact Green's function given by:
\[
G_{\text{exact}}(x,y) = 
\begin{cases}
4x(y-1)\exp(-2(x-y)), & \text{if } x\leq y,\\
(x-1)y, & \text{if } y < x,
\end{cases}
\]
for $x,y\in[0,1]$. In \cref{fig_adv_diff_SVD}, we display the first five left and right singular vectors and the singular values of the exact and learned Green's functions. We observe that the first fifteen singular values of the learned Green's functions are accurate. This leads us to conclude that our method enables the construction of a low-rank representation of the solution operator associated with the differential operator, $\L$, and allows us to compute and analyze its dominant modes.

\subsection{Schr\"odinger equation with double-well potential}

\begin{figure}[htbp]
\centering
\vspace{0.5cm}
\begin{overpic}[width=0.8\textwidth]{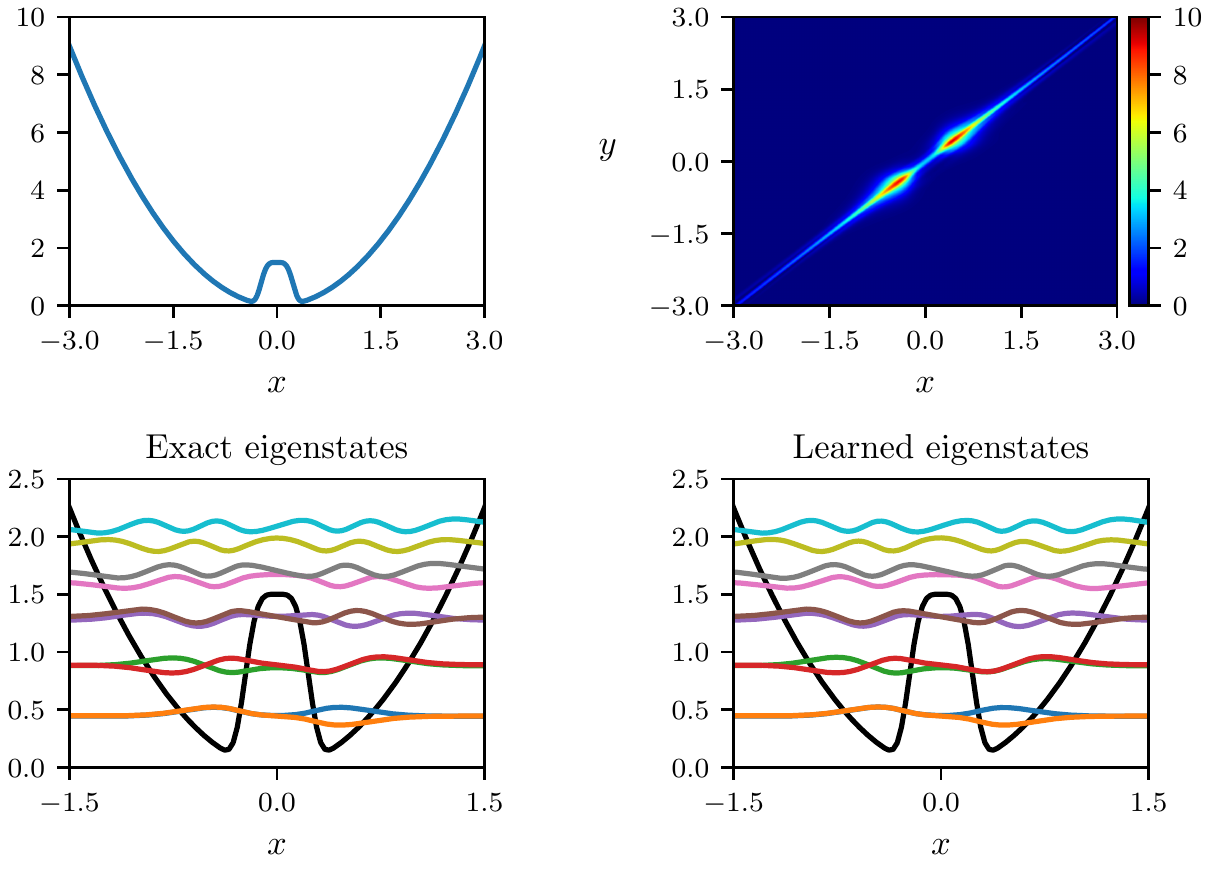}
\put(-4,71){(a)}
\put(48,71){(b)}
\put(-4,33){(c)}
\put(48,33){(d)}
\end{overpic}
\caption{Schr\"odinger equation. (a) Double well potential $V(x)=x^2+1.5\exp(-(4x)^4)$. (n) Learned Green's function of the Schr\"odinger equation with potential $V(x)$. (c) First ten exact eigenstates computed numerically from the Schr\"odinger operator and (d) eigenstates computed from the learned Green's function displayed in (b). The eigenfunctions are shifted by an amount corresponding to the eigenvalue. The double-well potential is shown as a black curve.}
\label{fig_schrodinger_rational}
\end{figure}
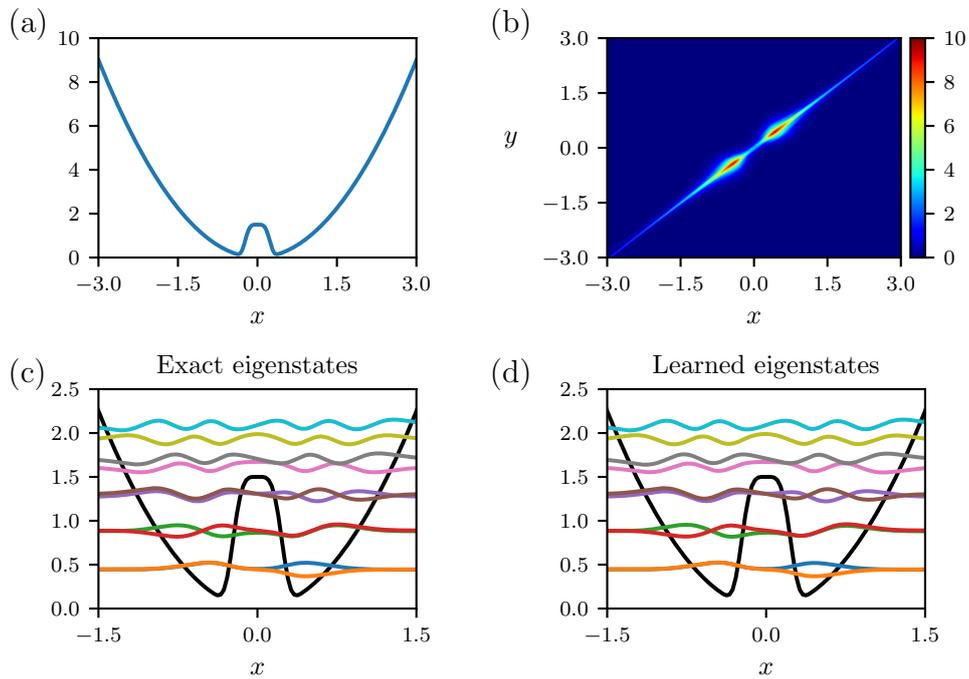

We highlight the ability of our DL method to learn physical features of an underlying system by considering the steady-state one-dimensional Schr\"odinger operator on $\Omega=[-3,3]$:
\[\L(u)=-h^2\frac{d^2u}{dx^2}+V(x)u,\quad u(-3)=u(3)=0,\]
with double-well potential $V(x) = x^2+1.5\exp(-(4x)^4)$ and $h=0.1$~\cite{trefethen2017exploring}. The potential $V(x)$ is illustrated in \cref{fig_schrodinger_rational}, along with the Green's function learned by the rational NN from pairs of forcing terms and the system's responses. First, the shape of the well potential can be visualized along the diagonal of the Green's function in \cref{fig_schrodinger_rational}(b). Next, in \cref{fig_schrodinger_rational}, we compute the first ten eigenstates of the Schr\"odinger operator in Chebfun~\cite{driscoll2014chebfun} and plot them using a similar representation as~\cite[Figure~6.9]{trefethen2017exploring}. Similarly to \cref{sec_eig}, we compute the eigenvalue decomposition of the Green's function learned by a rational NN and plot the eigenstates (shifted by the corresponding eigenvalues) in \cref{fig_schrodinger_rational}. Note that the eigenvalues of the operator and the Green's functions are reversed. We observe a perfect agreement between the first ten exact and learned eigenstates. These energy levels capture information about the states of atomic particles modeled by the Schr\"odinger equation.

\subsection{Singularity location and type} \label{sec_singularity}

The input-output function of a rational NN is a high-degree rational function, which means that it has poles (isolated points for which it is infinite). In rational function approximation theory, it is known that the poles of a near-optimal rational approximant tend to cluster near a function's singularities~\cite{trefethen2021exponential}. The clustering of the poles near the singularity is needed for the rational approximant to have excellent global approximation~\cite{stahl1992best,stahl1993best}. Moreover, the type of clustering (algebraic, exponential, beveled exponential) can reveal the type of singularity (square-root, blow-up, non-differentiable) at that location. This feature of rational approximants is used in other settings~\cite{beyene1999pole}.  

\begin{figure}[htbp]
\centering
\vspace{0.5cm}
\begin{overpic}[width=\textwidth]{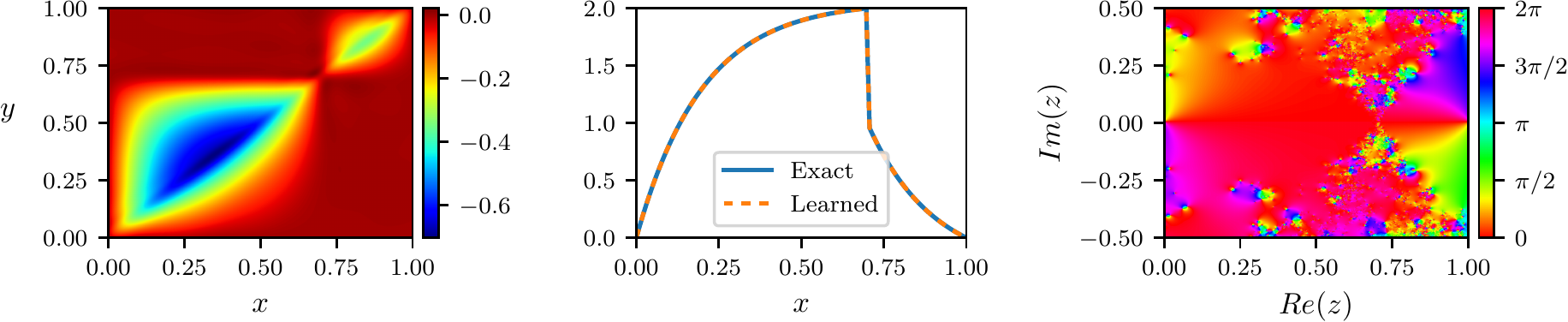}
\put(-1,20){(a)}
\put(33,20){(b)}
\put(65,20){(c)}
\end{overpic}
\caption{Singularity location. (a) Learned Green's function of a second-order differential operator with a jump condition at $x=0.7$. Homogeneous solution of the operator with jump condition (b) and argument of the rational NN representing the homogeneous solution in the complex plane (c).}
\label{fig_jump}
\end{figure}

We show that the rational NNs also cluster poles in a way that identifies its location and type. In~\cref{fig_jump}(c), we display the complex argument of the trained rational NN for the Green's function of a second-order differential operator with a jump condition, defined on $\Omega=[0,1]$ as
\[\L u = 0.2\frac{d^2 u}{dx^2}+\frac{d u}{dx},\quad u(0)=u(1)=0,\, u(0.7^{-})=2,\, u(0.7^{+})=1.\]
These diagrams are known as phase portraits and are useful for illustrating complex analysis~\cite{wegert2012visual}. A pole of the rational function can be identified as a point in the complex plane for which the full colormap goes around that point in a clockwise fashion. In  particular, in~\cref{fig_jump}(c), we see that the poles of the rational function cluster quite closely to the real-line (where $Im(z) = 0$) at $x = 0.7$. If the clustering is examined more closely, it may be possible to reveal that the singularity in the Green's function at $x = 0.7$ is due to a jump condition. 

\begin{figure}[htbp]
\centering
\vspace{0.5cm}
\begin{overpic}[width=0.8\textwidth]{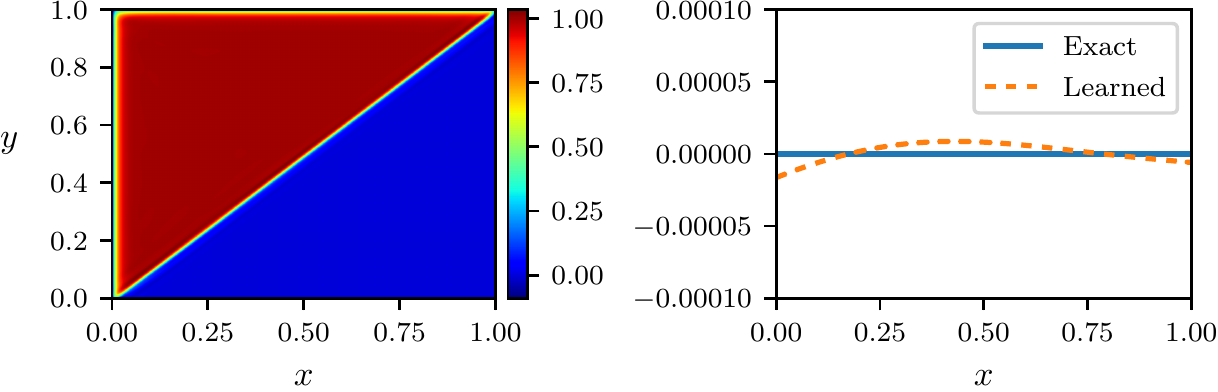}
\put(-1,33){(a)}
\put(50,33){(b)}
\end{overpic}
\caption{Boundary layer. Learned Green's function (a) and homogeneous solution (b) to a differential equation with a boundary layer around $x=0$.}
\label{fig_boundary_layer}
\end{figure}

Rational NNs are also important for resolving Green's function with boundary layers as the NN can resolve the boundary layer by clustering its poles in the complex plane. In~\cref{fig_boundary_layer}, we see a learned Green's function of a differential equation with a boundary layer at $x = 0$ with $\nu=10^{-2}$:
\[\L u=-\nu\frac{d^2 u}{dx^2}-\frac{d u}{dx}, \quad u(0)=u(1)=0, \quad \Omega=[0,1].\]
The analytical expression for the Green's function is given by the following equation:
\[
G_{\text{exact}}(x,y) = 
\begin{cases}
\frac{1}{e^{1/\nu}-1}\left(1-e^{-x\nu}\right)\left(e^{1/\nu}-e^{y/\nu}\right), & \text{if } x\leq y,\\
\frac{1}{e^{1/\nu}-1}\left(1-e^{(1-x)/\nu}\right)\left(1-e^{y/\nu}\right), & \text{if } y < x,
\end{cases}
\]
While the Green's function is not smooth, our rational NN still resolves it with relatively good accuracy, as shown by the sharp interface along the diagonal.

\section{Viscous shock and multiphysics examples}

In this section, we focus on two physical models and analyse the Green's functions discovered by our deep learning approach.

\subsection{Viscous shock}

\begin{figure}[htbp]
\centering
\vspace{0.5cm}
\begin{overpic}[width=\textwidth]{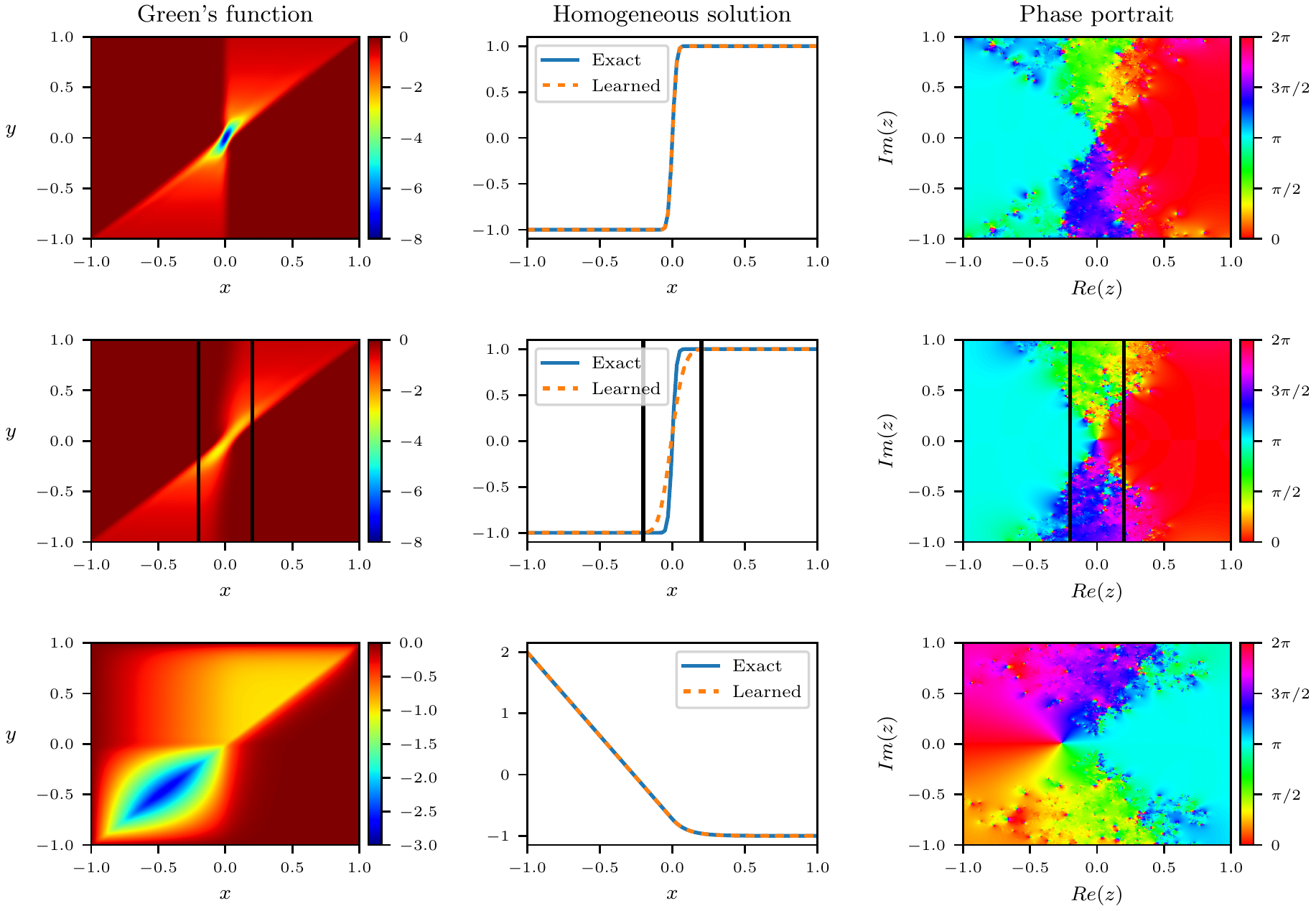}
\put(0,67){(a)}
\put(34,67){(b)}
\put(65,67){(c)}
\put(0,44){(d)}
\put(34,44){(e)}
\put(65,44){(f)}
\put(0,21){(g)}
\put(34,21){(h)}
\put(65,21){(i)}
\end{overpic}
\caption{Green's functions learned by rational neural networks. (a) Green's function of a differential operator with a viscous shock at $x=0$, learned by a rational NN. (b) Learned and exact (computed by a classical spectral method) homogeneous solution to the differential equation with zero forcing term. (c) Phase portrait of the homogeneous rational NN evaluated on the complex plane. (d)-(f) Similar to (a)-(c), but without any system's response measurements in $x\in[-0.2,0.2]$ (see vertical black lines) near the shock. (g) Learned Green's function and homogeneous solution (h) of an advection-diffusion operator with advection occurring for $x\geq 0$. (i) Phase portrait of the homogeneous NN on the complex plane.}
\label{fig_experiments}
\end{figure}
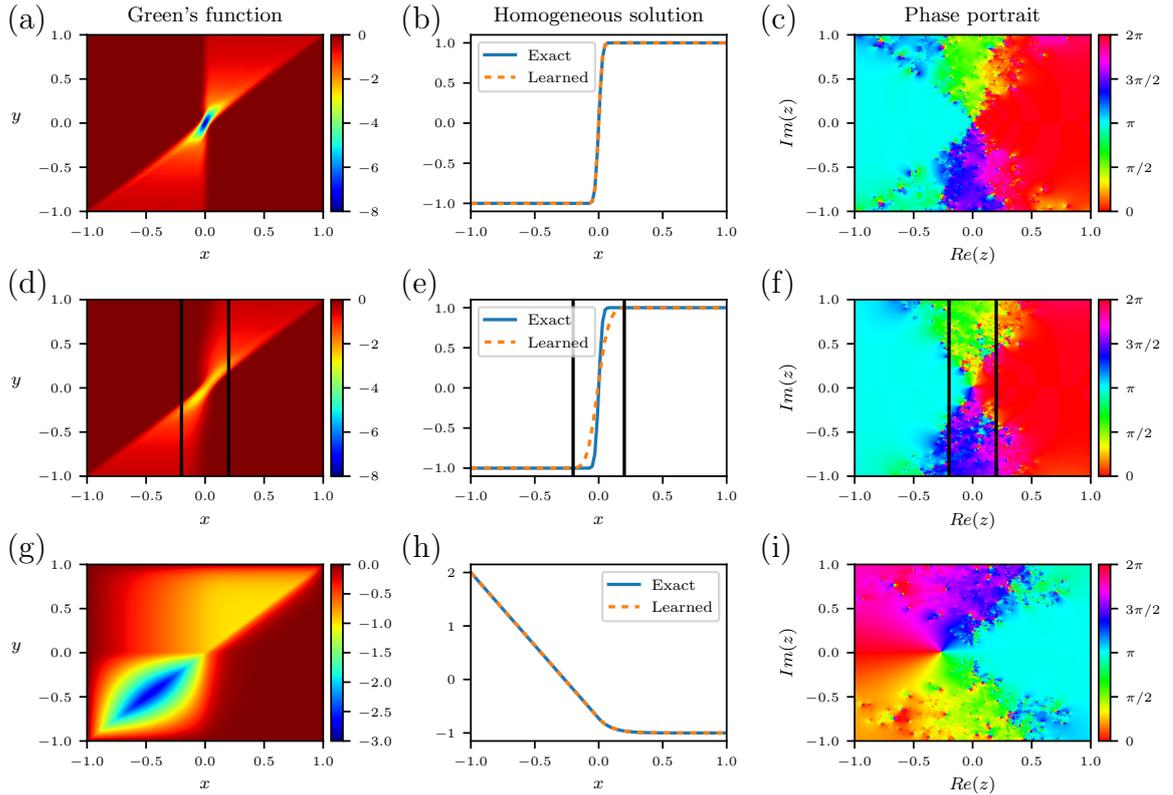

As a first example, we consider a second-order differential operator having suitable variable coefficients to model a viscous shock at $x=0$~\cite{lee1997fast}:
\[\L u = 10^{-3}\frac{d^2 u}{dx^2}+2x\frac{du}{dx}, \quad u(-1) = -1,\, u(1)=1.\]
The system's responses are obtained by solving the PDE, with Dirichlet boundary conditions, using a spectral numerical solver for each of the $N=100$ random forcing terms, sampled from a GP having a squared-exponential covariance kernel~\cite{boulle2021learning}. The learned Green's function is displayed in \cref{fig_experiments}(a) and satisfies the following symmetry relation: $G(x,y) = G(-x,-y)$, indicating the presence of a reflective symmetry group within the underlying PDE. Indeed, if $u$ is a solution to $\L u=f$ with homogeneous boundary conditions, then $u(-x)$ is a solution to $\L v=f(-x)$. We also observe in \cref{fig_experiments}(b) and (c) that the homogeneous solution is accurately captured and that the poles of the homogeneous rational NN cluster near the real axis around $x=0$: the location of the singularity induced by the shock (cf.~\cref{sec_singularity}).

Next, we reproduce the same viscous shock numerical experiment, except that this time we remove measurements of the system's response from the training dataset in the interval $[-0.2,0.2]$: adjacent to the shock front. By comparing \cref{fig_experiments}(d)-(f) and \cref{fig_experiments}(a)-(c), we find that the Green's function and homogeneous solution, learned by the rational NNs, may not be affected in the region outside of the interval with missing data. In some cases, the NNs can still accurately capture the main features of the Green's function and homogeneous solution in the region lacking measurements. The robustness of our method to noise perturbation and corrupted or missing data is of significant interest and promising for real applications with experimental data.

\subsection{Advection-diffusion operator}

We next apply our DL method to discover the Green's function and homogeneous solution of an advection-diffusion operator, where the advection is dominant only within the right half of the domain:
\[\L u = 0.1\frac{d^2 u}{dx^2}+\mathbb{I}_{(x\geq 0)}\frac{du}{dx},\quad u(-1)=2,\, u(1)=-1,\]
on $\Omega = [-1,1]$. Here, $\mathbb{I}_{(x\geq 0)}$ denotes the characteristic function on $x\geq 0$. The resulting equation is diffusive on the left half of the domain, while the advection is turned on for $x\geq 0$. The output of the Green's function NN is plotted in \cref{fig_experiments}(g), where we observe the disparate spatial behaviors of the dominant physical mechanisms. This can be recognized when observing the restriction of the Green's function to the subdomain $[-1,0]\times[-1,0]$, where the observed solution is reminiscent of the Green's function for the Laplacian; thus indicating that the PDE is diffusive on the left half of the domain. Similarly, the restriction of the learned Green's function to $[0,1]\times[0,1]$ is characteristic of advection. 

In \cref{fig_experiments}(h) and (i), we display the homogeneous solution NN, along with the phase of the rational NN, evaluated on the complex plane. The agreement between the exact and learned homogeneous solution illustrates the ability of the DL method to accurately capture the behavior of a system within ``multiphysics'' contexts. The choice of rational NNs is crucial here: to deepen our understanding of the system, as the poles of the homogeneous rational NN characterize the location and type of singularities in the homogeneous solution. Here the change in behavior of the differential operator from diffusion to advection is delineated by the location of the poles of the rational NN. 

\section{Two-dimensional operators and systems}

Our deep learning technique for learning Green's functions generalizes well in two dimensions and for systems of linear partial differential equations as we will see in this section.

\subsection{Differential operators in two dimensions} \label{sec_dimension_2}

We demonstrate the ability of our method to learn Green's functions associated with two-dimensional operators by repeating the numerical experiment of~\cite{li2020neural}, which consists of learning the Green's function of the Poisson operator on the unit disk $\Omega=D(0,1)$, with homogeneous Dirichlet boundary conditions:
\[\L u = \nabla^2 u,\quad u_{|\partial D(0,1)} = 0.\]
This experiment is a good benchmark for PDE learning techniques as the analytical expression of the Green's function in Cartesian coordinates can be expressed as~\cite{myint2007linear}:
\[G_{\text{exact}}(x,y,\tilde{x},\tilde{y}) = \frac{1}{4\pi}\ln\left(\frac{(x-\tilde{x})^2+(y-\tilde{y})^2}{(x\tilde{y}-\tilde{x}y)^2+(x\tilde{x}+y\tilde{y}-1)^2}\right),\]
where $(x,y),(\tilde{x},\tilde{y})\in D(0,1)$.

\begin{figure}[htbp]
\centering
\vspace{0.5cm}
\begin{overpic}[width=\textwidth]{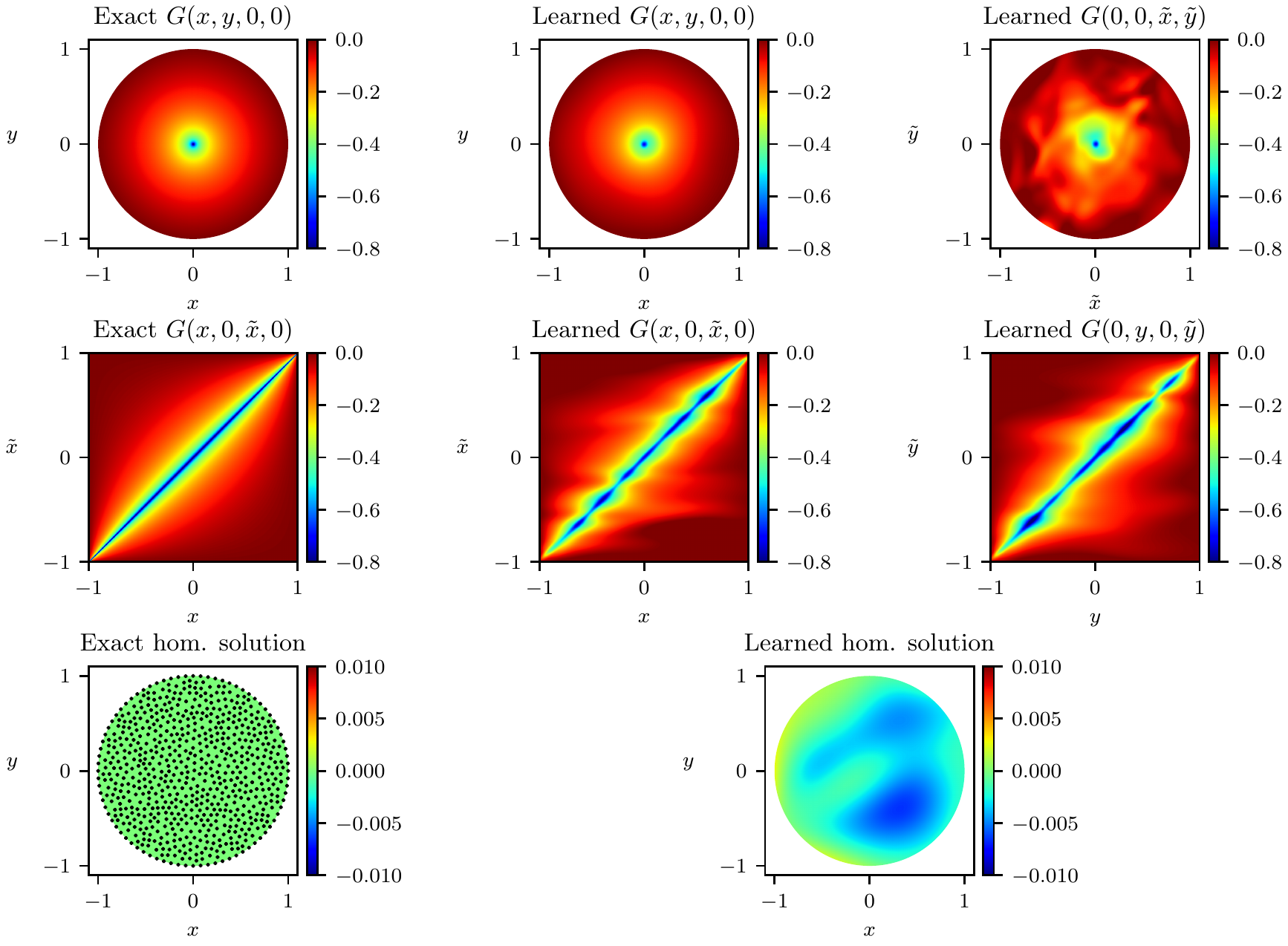}
\put(-1,71){(a)}
\put(33.5,71){(d)}
\put(68.5,71){(e)}
\put(-1,47){(b)}
\put(33.5,47){(f)}
\put(68.5,47){(g)}
\put(-1,23){(c)}
\put(51.5,23){(h)}
\end{overpic}
\caption{Poisson equation on the disk. Exact (a)-(b) and learned (d)-(f) Green's function of the Poisson operator on the unit disk, evaluated at two-dimensional slices. The colorbar is scaled to remove the singularity of the Green's function at $(x,y)=(\tilde{x},\tilde{y})$. (c) Exact homogeneous solution with sample points for the training functions and (h) homogeneous solution learned by the rational NN.}
\label{poisson_disk}
\end{figure}
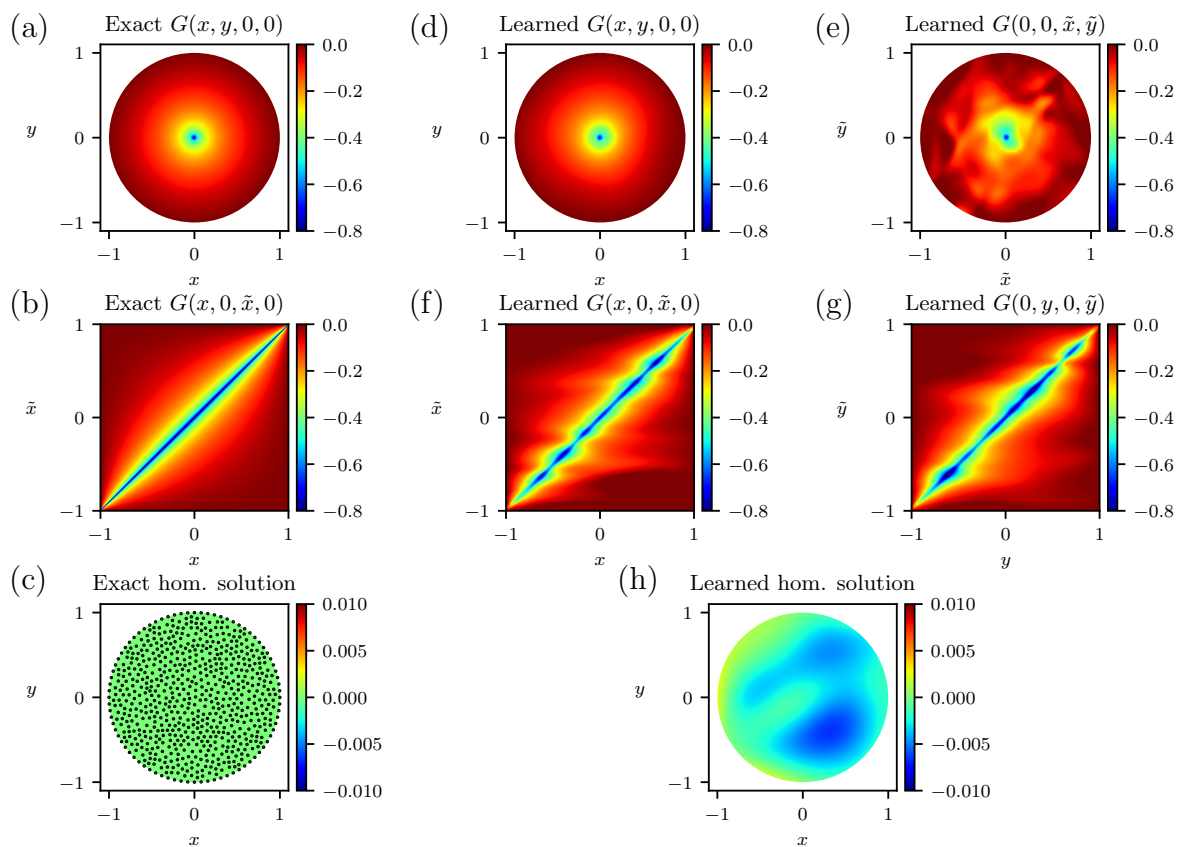

The training dataset for this numerical example is created as follows. First, we generate $N=100$ random forcing terms using the command \texttt{randnfundisk} of the Chebfun software~\cite{driscoll2014chebfun,filip2019smooth,wilber2017computing} with a frequency parameter of $\lambda = 0.2$, and then solve the Poisson equation, with corresponding right-hand sides, using a spectral method. Then, the forcing terms and system responses (\emph{i.e.}~solutions) are sampled at the $N_u=N_f=673$ nodes of a disk mesh, generated using the Gmsh software~\cite{geuzaine2009gmsh}. Moreover, the mesh structure ensures that the repartition of the sample points is approximately uniform in the disk (\cref{poisson_disk}(c)) and that the boundary is accurately captured.

The Green's function and homogeneous rational NNs have four hidden layers and width of $50$ neurons, with 4 and 2 input nodes, respectively, as the Green's function is defined on $\Omega\times\Omega$. The two-dimensional integrals of the loss function~\eqref{eq_loss} are discretized using uniform quadrature weights: $w_i = \pi/N_f$ for $1\leq i\leq N_f$. In \cref{poisson_disk}(d)-(g), we visualize four two-dimensional slices of the learned Green's function together with two slices of the exact Green's function in panels (a) and (b). Because of the symmetry in the Green's function, due to the self-adjointness of $\L$ and the boundary constraints, the exact Green's function satisfies $G(x,y,0,0) = G(0,0,x,y)$ for $(x,y)\in D(0,1)$. Therefore, we compare \cref{poisson_disk}(a) to \cref{poisson_disk}(d)-(e), and similarly for \cref{poisson_disk}(b) and \cref{poisson_disk}(f)-(g). We observe that the Green's function is accurately learned by the rational NN, which preserves low approximation errors near the singularity at $(x,y)=(\tilde{x},\tilde{y})$, contrary to the neural operator technique~\cite{li2020neural}. The visual artifacts present in \cref{poisson_disk}(e)-(g) are likely due to the low spatial discretization of the training data. One could increase the number of spatial measurements or use a high-order quadrature rule.

\subsection{System of differential equations} \label{sec_system}

The method extends also naturally to systems of differential equations. Let $f = \begin{bmatrix}
f^1 & \cdots & f^{n_f}
\end{bmatrix}^\top:\Omega\to\R^{n_f}$ be a vector of $n_f$ forcing terms and 
$u = \begin{bmatrix}
u^1 & \cdots & u^{n_u}
\end{bmatrix}^\top:\Omega\to\R^{n_u}$ be a vector of $n_u$ system responses such that
\begin{equation} \label{eq_system}
\L
\begin{bmatrix}
u^1\\
\vdots\\
u^{n_u}
\end{bmatrix}
=\begin{bmatrix}
f^1 \\ 
\vdots \\
f^{n_f}
\end{bmatrix},\quad 
\mathcal{D}\left(
\begin{bmatrix}
u^1\\
\vdots\\
u^{n_u}
\end{bmatrix},\Omega\right)=
\begin{bmatrix}
g^1\\
\vdots\\
g^{n_u}
\end{bmatrix}.
\end{equation}
The solution to \cref{eq_system} with $f=0$ is called the homogeneous solution and denoted by $u_{\text{hom}}=\begin{bmatrix}
u_{\text{hom}}^1 & \cdots & u_{\text{hom}}^{n_u}
\end{bmatrix}^\top$. Similarly to the scalar case, we can express the relation between the system's response and the forcing term using Green's functions and an integral formulation as
\begin{equation} \label{eq_system_loss}
u^i(x)=
\sum_{j=1}^{n_f}
\int_\Omega G_{i,j}(x,y)f^j(y)\d y+u_{\text{hom}}^{i}(x),\quad x\in\Omega,
\end{equation}
for $1\leq i\leq n_u$. Here, $G_{i,j}:\Omega\times\Omega\to \R\cup\{\pm\infty\}$ is a component of the \emph{Green's matrix} for $1\leq i\leq n_u$ and $1\leq j\leq n_f$, which consists of a $n_u\times n_f$ matrix of Green's functions:
\[G(x,y) = \begin{bmatrix}
G_{1,1}(x,y) & \cdots & G_{1,n_f}(x,y)\\
\vdots & \ddots & \vdots \\
G_{n_u,1}(x,y) & \cdots & G_{n_u,n_f}(x,y)
\end{bmatrix},\quad x,y\in\Omega.\]
Following \cref{eq_system_loss}, we remark that the differential equations decouple, and therefore we can learn each row of the Green's function matrix independently. That is, for each row $1\leq i\leq n_u$, we train $n_f$ NNs to approximate the components $G_{i,1},\ldots,G_{i,n_f}$, and one NN to approximate the $i$th component of the homogeneous solution, $u_{\text{hom}}^i$.

\begin{figure}[htbp]
\centering
\vspace{0.5cm}
\begin{overpic}[width=\textwidth]{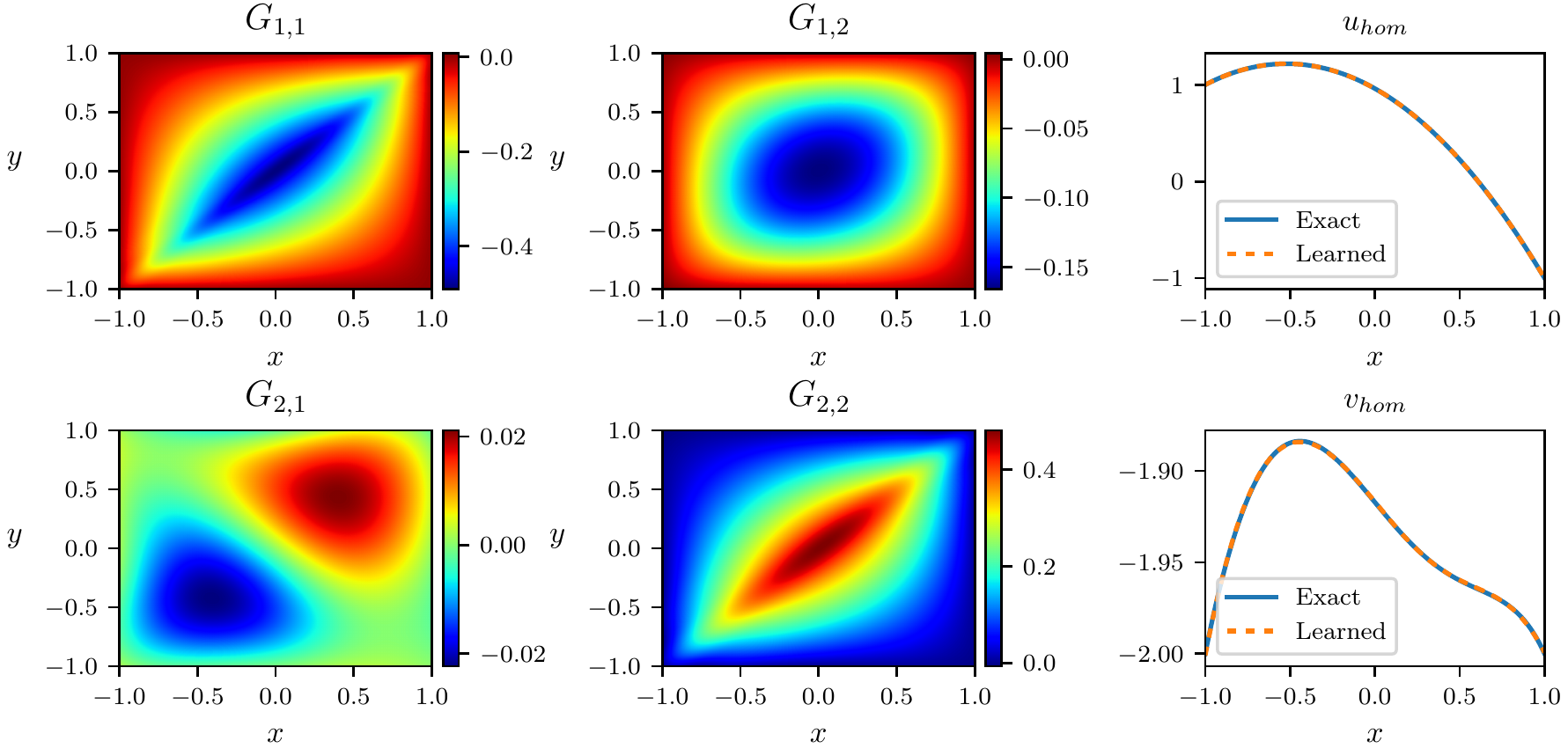}
\put(-1,45){(a)}
\put(71,45){(b)}
\end{overpic}
\caption{Green's matrix of system of ODEs. (a) Matrix of Green's function learned from the system of ordinary differential equations~\eqref{eq_ODE}. (b) Homogeneous solutions associated with the system of ODEs.}
\label{fig_ODE}
\end{figure}

As an example, we consider the following system of ordinary differential equations (ODEs) on $\Omega = [-1,1]$:
\begin{subequations} \label{eq_ODE}
\begin{align}
\frac{d^2 u}{d x^2} - v &= f^1,\\
\frac{-d^2 v}{d x^2} +xu &= f^2,
\end{align}
\end{subequations}
with boundary conditions: $u(-1) = 1$, $u(1) = -1$, $v(-1)=v(1) =-2$. In \cref{fig_ODE}, we display the different components of the Green's matrix and the exact solution (computed by a spectral method), along with the learned homogeneous solutions. We find that the Green's function matrix provides insight on the coupling between the two system variables, $u$ and $v$, as shown by the diagonal components $G_{1,2}$ and $G_{2,1}$ of the Green's matrix in \cref{fig_ODE}(a). Similarly, the components $G_{1,1}$ and $G_{2,2}$ are characteristic of diffusion operators, which appear in \cref{eq_ODE}. In this case, the Green's matrix can be understood as a $2\times 2$ block inverse~\cite{lu2002inverses} of the linear operator, $\L$.

\section{Nonlinear and vector-valued equations}

We can also discover Green's functions from forcing terms and concomitant solutions to nonlinear differential equations possessing semi-dominant linearity as well as Green's functions associated with vector-valued equations.

\subsection{Linearized models of nonlinear operators} \label{sec_nonlinear}

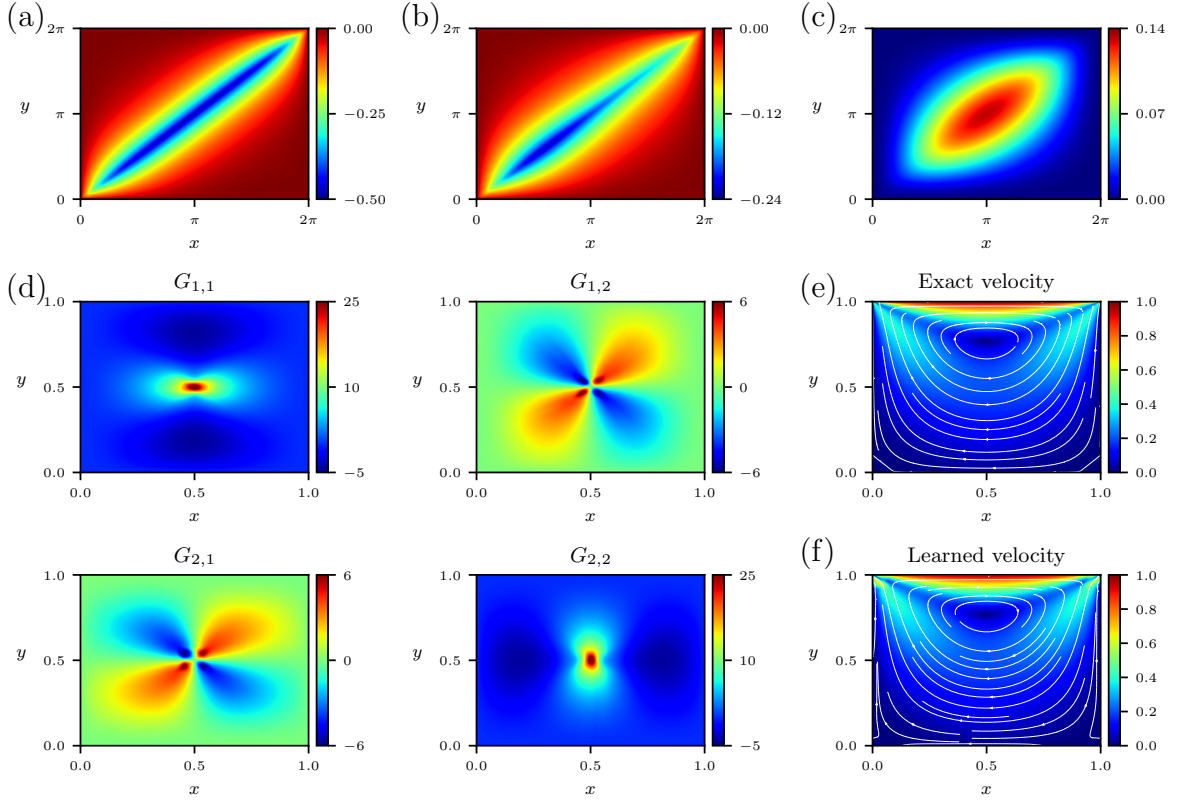
\begin{figure}[htbp]
\centering
\vspace{0.5cm}
\begin{overpic}[width=\textwidth]{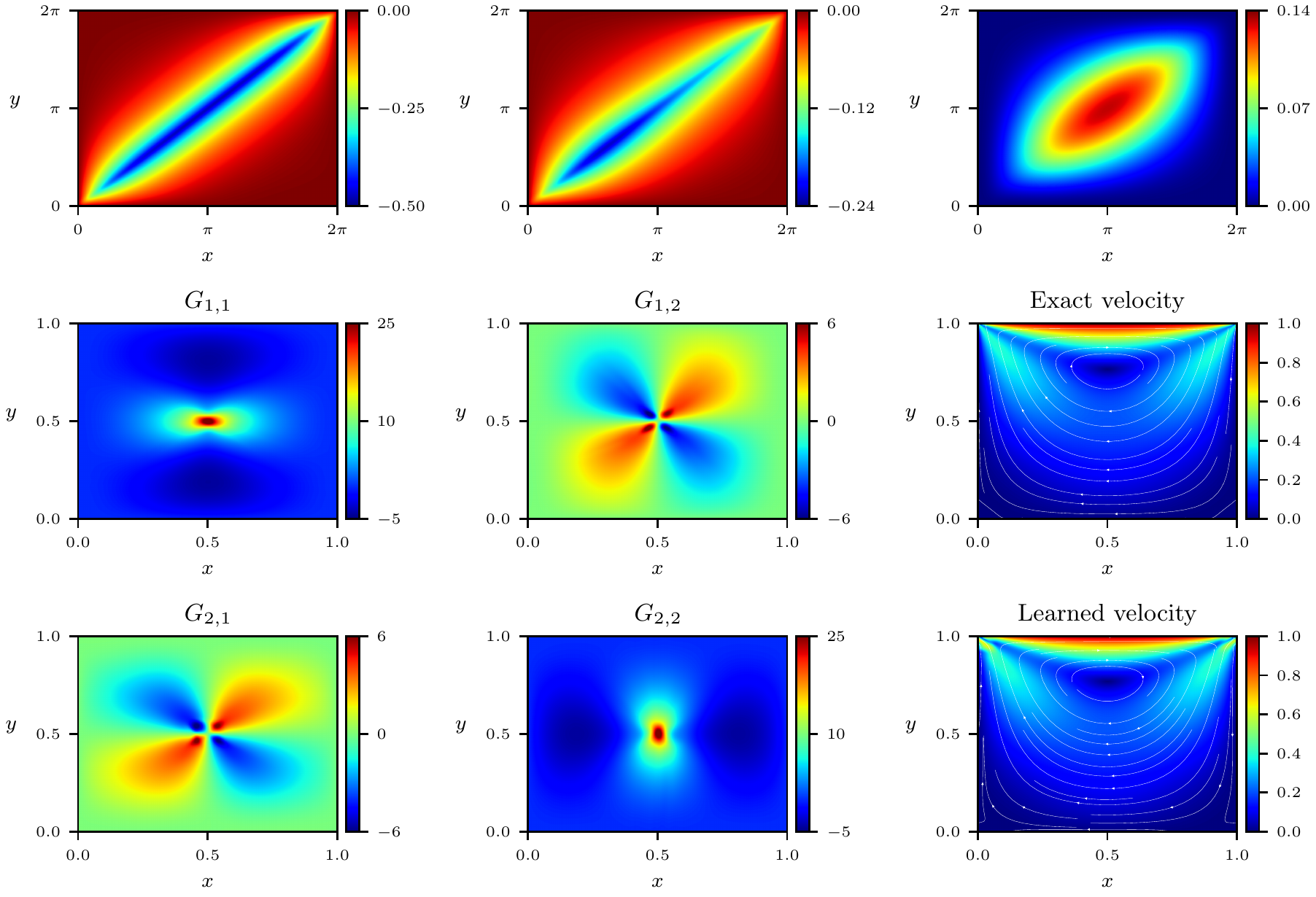}
\put(-0.5,67.5){(a)}
\put(33.5,67.5){(b)}
\put(68,67.5){(c)}
\put(-0.5,44){(d)}
\put(68,44){(e)}
\put(68,21){(f)}
\end{overpic}
\caption{Linearized models and Stokes flow. (a)-(c) Green's functions of three differential operators: Helmholtz, Sturm--Liouville, and biharmonic, with cubic nonlinearity. (d) Matrix of Green's functions of a two-dimensional Stokes flow in a lid-driven cavity, evaluated at the two-dimensional slice $(x,y,0.5,0.5)$. Velocity magnitude and streamlines of the exact (e) and learned (f) homogeneous solution to the Stokes equations with zero applied body force.}
\label{fig_stokes}
\end{figure}

We demonstrate that our DL method can be used to linearize and extract Green's functions from nonlinear boundary value problems of the form
\[\L u+\epsilon\N(u) = f, \quad \D(u,\Omega) = g,\]
where $\L$ denotes a linear operator, $\N$ is a nonlinear operator, and $\epsilon < 1$ is a small parameter controlling the nonlinearity.
We demonstrate this ability on the three nonlinear boundary value problems, dominated by the linearity, used in \cite{gin2020deepgreen}. In \cref{fig_stokes}(a)-(c), we visualize the Green's function NNs of three operators with cubic nonlinearity considered in~\cite{gin2020deepgreen}.

First, \cref{fig_stokes}(a) illustrates the learned Green's function of a cubic Helmholtz system on $\Omega=[0,2\pi]$ with homogeneous Dirichlet boundary conditions:
\[\frac{d^2 u}{dx^2}+\alpha u+\epsilon u^3 = f(x),\]
where $\alpha=-1$ and $\epsilon = 0.4$. Next, in \cref{fig_stokes}(b), we consider a nonlinear Sturm--Liouville operator of the form:
\[[-p(x)u']'+q(x)(u+\epsilon u^3) = f(x),\quad u(0) = u(2\pi) = 0,\]
with $p(x) = 0.4\sin(x)-3$, $q(x)=0.6\sin(x)-2$, and $\epsilon=0.4$. The notation $u'$ denotes the derivative with respect to $x$, $du/dx$. Finally, the example represented in \cref{fig_stokes}(c) is the learned Green's function of a nonlinear biharmonic operator:
\[[-p(x)u'']''+q(u+\epsilon u^3) = f(x),\quad u(0) = u(2\pi) = 0,\]
where $p=-4$, $q=2$, and $\epsilon = 0.4$.

The nonlinearity does not prevent our method from discovering a Green's function of an approximate linear model, from which one can understand features such as symmetry and boundary conditions. This property is crucial for tackling time-dependent problems, where the present technique may be extended and applied to uncover linear propagators.

\subsection{Lid-driven cavity problem} \label{sec_lid_driven}

\begin{figure}[htbp]
\centering
\vspace{0.5cm}
\begin{overpic}[width=0.8\textwidth]{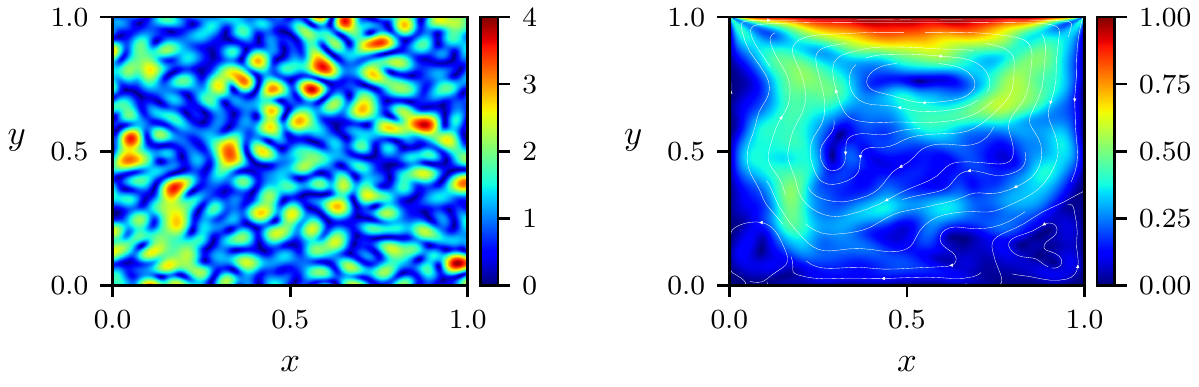}
\put(-1,31){(a)}
\put(50,31){(b)}
\end{overpic}
\caption{Training functions for Stokes flow. (a) Magnitude of a random applied body force used as a forcing term in the Stokes equations. (b) Velocity magnitude and streamlines of the system's response.}
\label{fig_stokes_training}
\end{figure}

Finally, we consider a classical benchmark in fluid dynamics consisting of Stokes flow in a two-dimensional lid-driven cavity problem~\cite{elman2014finite}. We aim to discover the matrix of Green's functions of the Stokes flow~\cite{blake1971note}, which is modelled by the following system of equations on the domain $\Omega = [0,1]^2$,
\begin{align*}
\mu\nabla^2\mathbf{u}-\nabla p&=\mathbf{f},\\
\nabla\cdot\mathbf{u}&=0.
\end{align*}
Here, $\mathbf{u}=(u_x,u_y)$ is the fluid velocity, $p$ is the pressure, $\mathbf{f}=(f_x,f_y)$ is an applied body force (\emph{i.e.}~a forcing term), and $\mu=1/100$ is the dynamic viscosity. The fluid velocity satisfies no-slip boundary conditions on the walls, except on the top wall where $\mathbf{u}=(1,0)$. We first generate one hundred forcing terms, $\mathbf{f}$, with two smooth random components, in the Chebfun software~\cite{driscoll2014chebfun,filip2019smooth} using the \texttt{randnfun2} command with wavelength parameter $\lambda = 0.1$. The Stokes equations are then discretized with Taylor--Hood finite elements~\cite{boffi2013mixed,taylor1973numerical} for the velocity and pressure on a mesh with $96\times 96$ square cells and subsequently solved using the Firedrake finite element library~\cite{rathgeber2016firedrake}. We illustrate in \cref{fig_stokes_training} an example of applied body force and velocity solution obtained by solving the system of PDEs. We then create the training dataset for the NNs by sampling the applied body forces and corresponding velocity solutions, $\mathbf{u}$, on a regular $25\times 25$ grid. 

In this context, the relation between the system's responses and the forcing terms can be expressed using a Green's matrix, which consists of a two-by-two matrix of Green's functions and whose components reveal features of the underlying system such as symmetry and coupling (\cref{fig_stokes}(d) and \cref{sec_system}). The four Green's functions and two homogeneous NNs have the same architecture as the one described in \cref{sec_rational_net}, except that they have respectively four and two input nodes (instead of two and one) due to the current spatial dimension. \cref{fig_stokes}(e) and (f) illustrate that the homogeneous solution to the Stokes equation is accurately captured by the homogeneous rational NN, despite the corner singularities and coarse measurement grid. The four components of the Green's matrix for the Stokes flow are evaluated on the two-dimensional slice $(x,y,0.5,0.5)$, for $x,y\in[0,1]$, and displayed in \cref{fig_stokes}(d). This figure allows us to visualize the system's response to a point force, $\mathbf{f}=(f_x,f_y)$, located at $(0.5,0.5)$, with the system's response being denoted as $\mathbf{u}=(u_x,u_y)$, where
\begin{align*}
u_x(x,y) &= G_{1,1}(x,y,0.5,0.5)f_x+G_{1,2}(x,y,0.5,0.5)f_y,\\
u_y(x,y) &= G_{2,1}(x,y,0.5,0.5)f_x+G_{2,2}(x,y,0.5,0.5)f_y,
\end{align*}
for $x,y\in[0,1]$. The visualization of the $G_{2,2}$ component in \cref{fig_stokes}(d), corresponding to the system's response to a unitary vertical point force $\mathbf{f}=(0,1)$ is reminiscent of~\cite[Figure~1]{ekiel2018stokes}. 

\begin{figure}[htbp]
\centering
\vspace{0.5cm}
\begin{overpic}[width=0.7\textwidth]{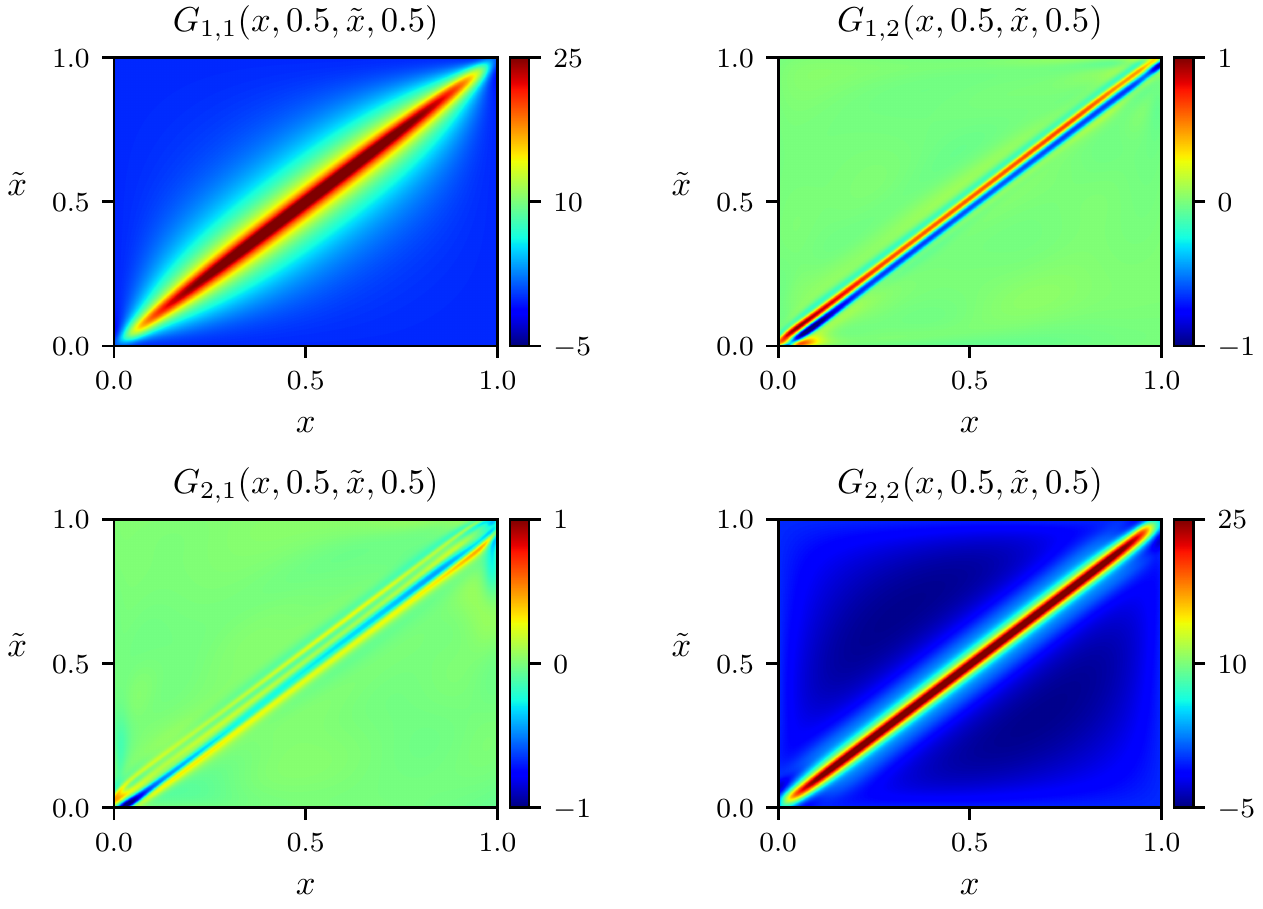}
\end{overpic}
\caption{2nd Green's matrix slice of Stokes flow. The four components of the Green's matrix learned by a rational neural network evaluated at the two-dimensional slice $(x,0.5,\tilde{x},0.5)$.}
\label{fig_slice_stokes_2}
\end{figure}

\begin{figure}[htbp]
\centering
\vspace{0.5cm}
\begin{overpic}[width=0.7\textwidth]{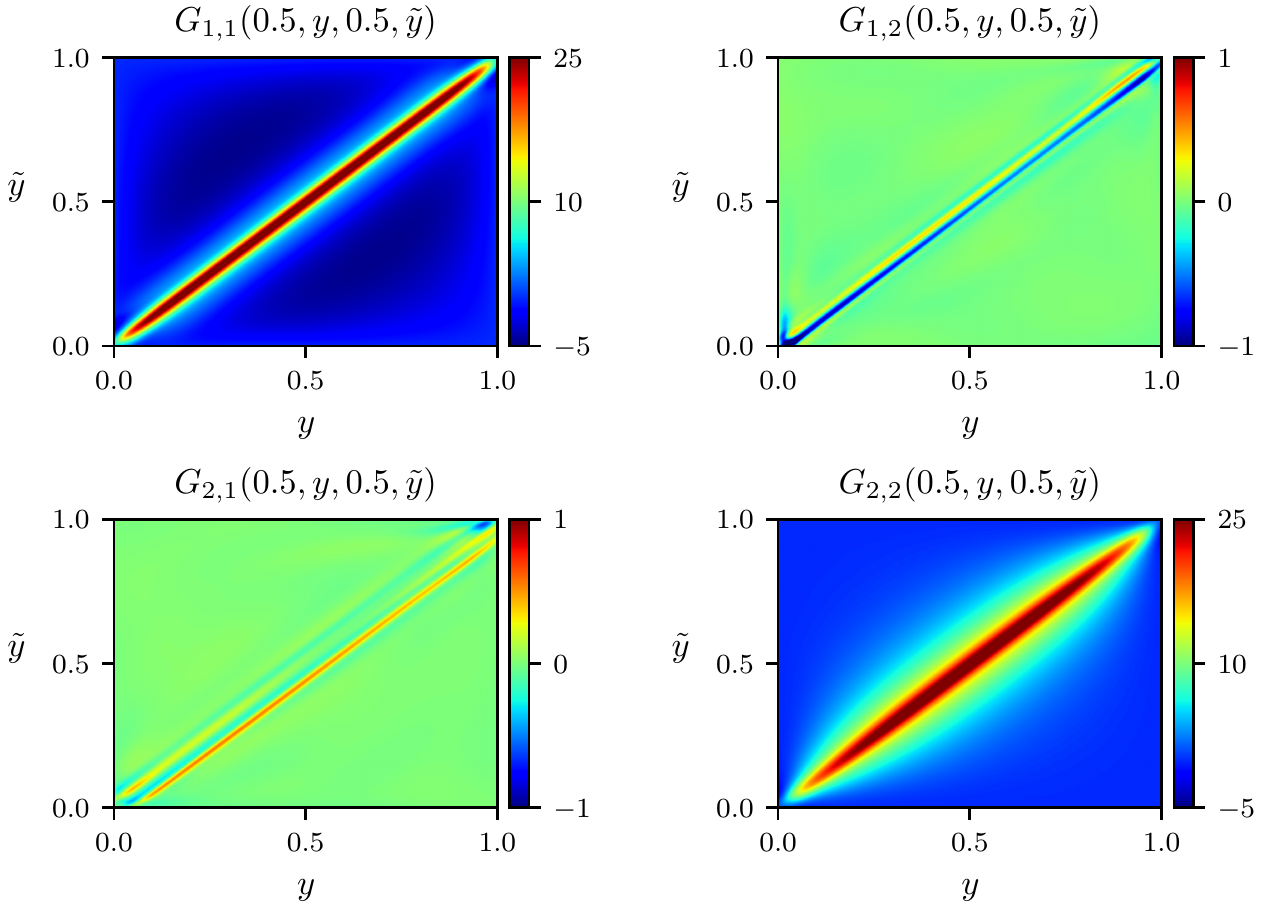}
\end{overpic}
\caption{3rd Green's matrix slice of Stokes flow. The four components of the Green's matrix learned by a rational neural network evaluated at the two-dimensional slice $(0.5,y, 0.5, \tilde{y})$.}
\label{fig_slice_stokes_3}
\end{figure}

\begin{figure}[htbp]
\centering
\vspace{0.5cm}
\begin{overpic}[width=0.7\textwidth]{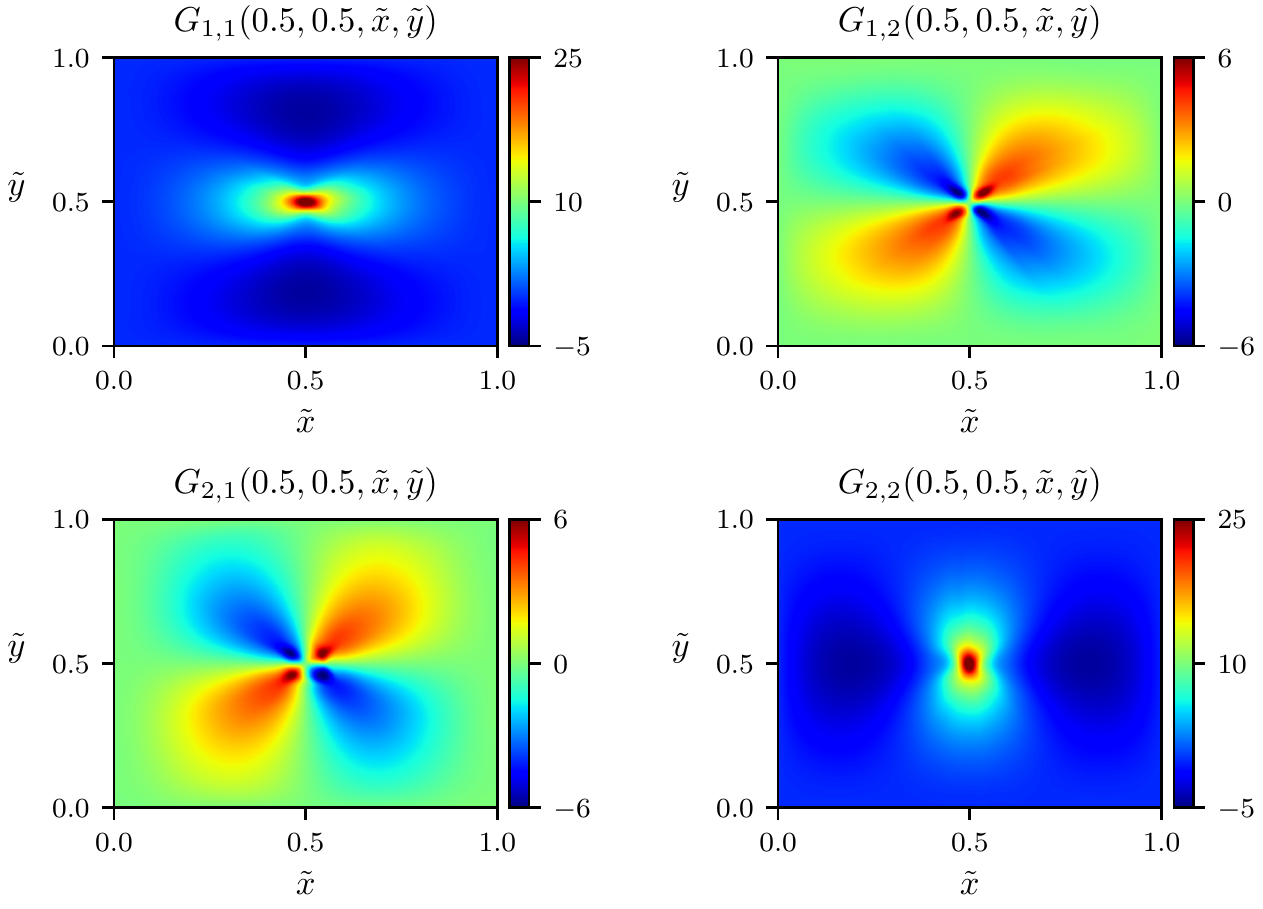}
\end{overpic}
\caption{4th Green's matrix slice of Stokes flow. The four components of the Green's matrix learned by a rational neural network evaluated at the two-dimensional slice $(0.5,0.5,\tilde{x},\tilde{y})$.}
\label{fig_slice_stokes_4}
\end{figure}

Finally, we evaluate the components of the Green's matrix at three other two-dimensional slices: $(x,0.5,\tilde{x},0.5)$, $(0.5,y,0.5,\tilde{y})$, $(0.5,0.5,\tilde{x},\tilde{y})$ and display them respectively in \cref{fig_slice_stokes_2,fig_slice_stokes_3,fig_slice_stokes_4}. These figures illustrate the different symmetries of the Green's matrix, which are captured by the rational NNs. As an example, we see in \cref{fig_slice_stokes_2,fig_slice_stokes_3} that $G_{1,1}(x,0.5,\tilde{x},0.5) = G_{2,2}(0.5,x,0.5,\tilde{x})$ and $G_{2,2}(x,0.5,\tilde{x},0.5) = G_{1,1}(0.5,x,0.5,\tilde{x})$, for $x,\tilde{x}\in [0,1]$. Similarly, we find in \cref{fig_slice_stokes_4} that $G_{1,1}(0.5,0.5,\tilde{x},\tilde{y})=G_{1,1}(0.5,0.5,\tilde{y},\tilde{x})$ and $G_{1,2}(0.5,0.5,\tilde{x},\tilde{y})=G_{2,1}(0.5,0.5,\tilde{x},\tilde{y})$, for $\tilde{x},\tilde{y}\in [0,1]$. The $G_{1,2}$ and $G_{2,1}$ components of the Green's matrix in \cref{fig_slice_stokes_2} highlight a singularity along the diagonal $(x,0.5,x,0.5)$ for $x\in [0,1]$. However, this singularity does not prevent the rational NNs from accurately learning the different components of the Green's matrix displayed in \cref{fig_stokes}(d) and \cref{fig_slice_stokes_2,fig_slice_stokes_3,fig_slice_stokes_4}.

\section{Time-dependent equations} \label{sec_time_dep}

In this section, we show that one can use a time-stepping scheme to discretize a time-dependent PDE and learn the Green's function associated with the time-propagator operator $\tau:u_n \to u_{n+1}$, where $u_n$ is the solution of the PDE at time $t = n\Delta t$ for a fixed time step $\Delta t$. As an example, we consider the time-dependent Schr\"odinger equation with a harmonic trap potential $V(x) = x^2$ given by
\begin{equation} \label{eq_schrodinger_time}
\mathrm{i}\frac{\partial \psi(x,t)}{\partial t} = -\frac{1}{2}\frac{\partial^2\psi(x,t)}{\partial x^2} + x^2\psi(x,t),\quad x\in[-3,3],
\end{equation}
with homogeneous Dirichlet boundary conditions. We use a Crank--Nicolson time-stepping scheme with time step $\Delta t = 2\times10^{-2}$ to discretize \cref{eq_schrodinger_time} in time and obtain
\[\mathrm{i}\frac{\psi_{n+1}-\psi_{n}}{\Delta t} = \frac{1}{2}\left[-\frac{1}{2}\frac{d^2\psi_{n+1}}{d x^2} + x^2\psi_{n+1}-\frac{1}{2}\frac{d^2\psi_{n}}{d x^2} + x^2\psi_{n}\right].\]
Our training dataset consists of one hundred random initial forcing functions $\psi_n$ at time $t$ and associated response $\psi_{n+1}$ at time $t+\Delta t$. The functions $\psi_n$ have real and imaginary parts sampled from a Gaussian process with periodic kernel and length-scale parameter $\lambda=0.5$ (see \cref{sec_generation_data}), and multiplied by the Gaussian damping function $g(x) = e^{-x^6/20}$ to ensure that the functions decay to zero before reaching the domain boundaries. We then train a rational neural network to learn the Green's function $G$ associated with the time-propagator operator such that
\[\tau(\psi_n)(x) = \int_{-3}^3 G(x,y)\psi_n(y) \d y= \psi_{n+1}(x),\quad x\in [-3,3].\]
Note that since $\psi$ takes complex values, we in fact split \cref{eq_schrodinger_time} into a system of equations for the real and imaginary parts of $\psi$, and learn the Green's matrix associated with the system (see \cref{sec_system}). 

\begin{figure}[htbp]
\centering
\vspace{0.5cm}
\begin{overpic}[width=0.95\textwidth]{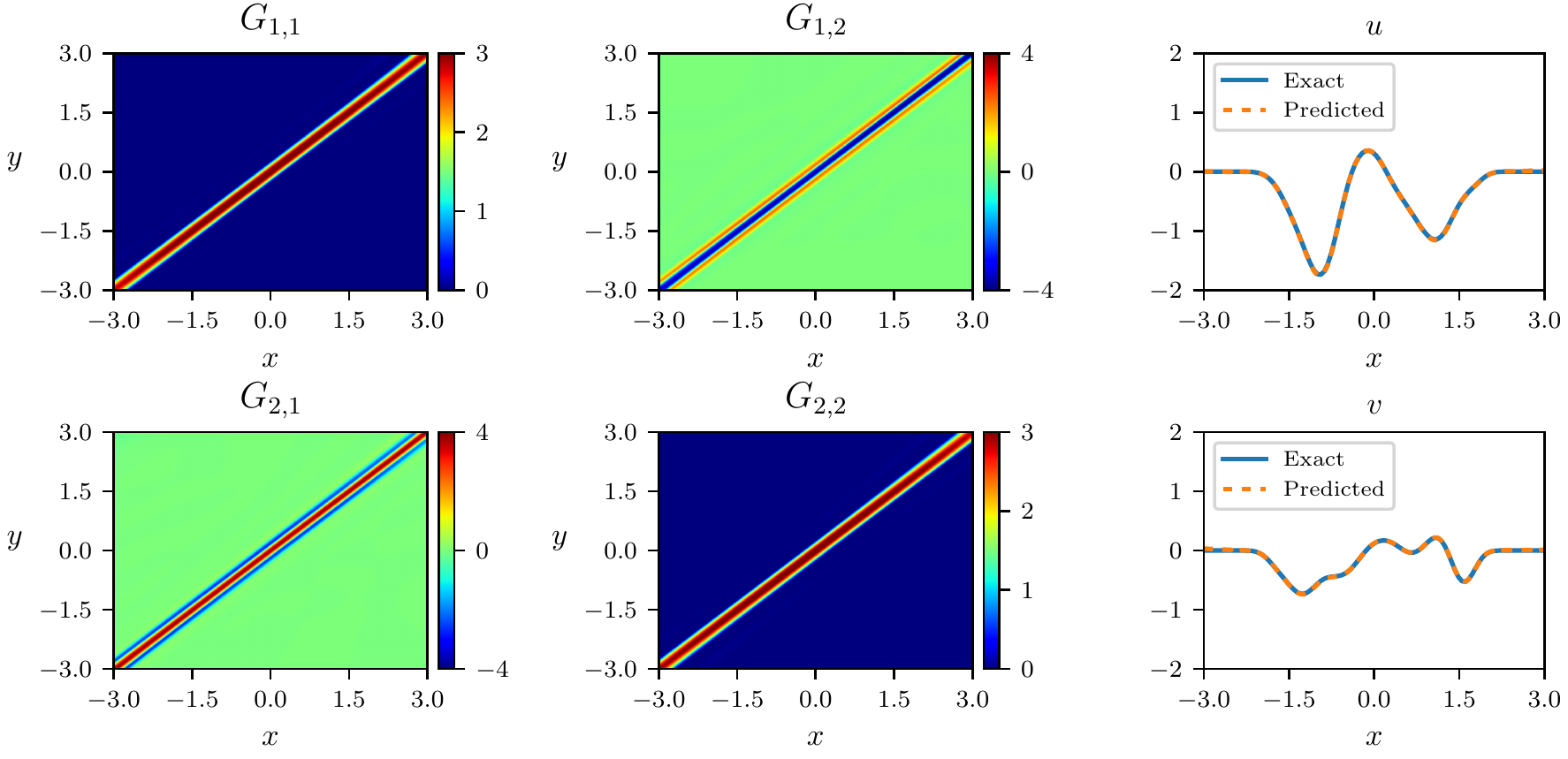}
\put(-1,45){(a)}
\put(69,45){(b)}
\end{overpic}
\caption{Green's matrix of the time-dependent Schr\"odinger equation. (a) The four components of the Green's matrix for the time propagator operator of the time-dependent Schr\"odinger equation discretized using a time-stepping scheme. (b) Real and imaginary components of the worst case prediction of the solution at the next time step.}
\label{fig_schrodinger_time}
\end{figure}

We report the Green's matrix of the time-propagator operator for the Schr\"odinger equation in \cref{fig_schrodinger_time}(a) and observe that the four components are dominated by the diagonal, which is expected for a small time-step. Additionally, we evaluate the accuracy of the learned Green's functions by generating a testing dataset with one hundred initial functions $\psi_{n}$, sampled from the same distribution, and associated solution $\psi_{n+1}$ at time $t+\Delta t$. We then compute the average (over the one hundred test cases) relative error in the $L^2$ norm between the exact solution $\psi_{n+1}$ and the one predicted using the learned Green's functions, $\psi_{n+1}^{\text{pred}}$, as
\[\textup{relative error} = \|\psi_{n+1}-\psi_{n+1}^{\text{pred}}\|_{L^2([-3,3])} / \|\psi_{n+1}\|_{L^2([-3,3])},\]
where $\psi_{n+1}^{\text{pred}}$ is defined as
\[\psi_{n+1}^{\text{pred}}(x) = \int_{-3}^3 G(x,y) \psi_n(y)\d y,\quad x\in[-3,3].\]
Finally, we obtain an average relative error of $1.3\%$ with standard deviation $0.2\%$ across the 100 test cases, confirming the good accuracy of our method. We display the worst-case prediction of the solution $\psi_{n+1}$ in \cref{fig_schrodinger_time}(b).

\dobib
 % Green learning

\chapter*{Conclusions}
\addcontentsline{toc}{chapter}{Conclusions}

This thesis derived theoretical results and a practical deep learning algorithm for approximating Green's functions associated with linear partial differential equations (PDEs) from pairs of forcing terms and solutions to a PDE.

By generalizing the randomized singular value decomposition (SVD) to Hilbert--Schmidt (HS) operators in \cref{chapt_PDE_learning}, we showed that one can rigorously learn the Green's function associated with an elliptic PDE in three dimensions. We derived a learning rate associated with elliptic partial differential operators in three dimensions and bounded the number of input-output training pairs required to recover a Green's function approximately with high probability. The random forcing functions are sampled from a Gaussian process (GP) with mean zero and are characterized by the associated covariance kernel. One practical outcome of this work is a measure for the quality of covariance kernels, which may be used to design efficient GP kernels for PDE learning tasks.

We then explored the practical extensions of the randomized SVD to Gaussian random vectors with correlated entries (\emph{i.e.},~nonstandard covariance matrices) and HS operators in \cref{chap_svd}. This chapter motivates new computational and algorithmic approaches for constructing the covariance kernel based on prior information to compute a low-rank approximation of matrices and impose properties on the learned matrix and random functions from the GP. We performed numerical experiments to demonstrate that covariance matrices with prior knowledge can outperform the standard identity matrix used in the literature and lead to near-optimal approximation errors. In addition, we proposed a covariance kernel based on weighted Jacobi polynomials, which allows the control of the smoothness of the random functions generated and may find practical applications in PDE learning~\cite{boulle2021data,boulle2020rational} as it imposes prior knowledge of Dirichlet boundary conditions. The algorithm presented in this chapter is limited to matrices and HS operators and does not extend to unbounded operators such as differential operators. Additionally, the theoretical bounds only offer probabilistic guarantees for Gaussian inputs, while sub-Gaussian distributions~\cite{kahane1960proprietes} of the inputs would be closer to realistic application settings.

Motivated by the theoretical results obtained in \cref{chapt_PDE_learning,chap_svd}, we wanted to design an efficient deep learning architecture for learning Green's functions. In \cref{chapt_rational}, we investigated rational neural networks, which are neural networks with smooth trainable activation functions based on rational functions. We proved theoretical statements quantifying the advantages of rational neural networks over ReLU networks. In particular, we remarked that a composition of low-degree rational functions has a good approximation power but a relatively small number of trainable parameters. Therefore, we showed that rational neural networks require fewer nodes and exponentially smaller depth than ReLU networks to approximate smooth functions to within a certain accuracy. This improved approximation power has practical consequences for large neural networks, given that a deep neural network is computationally expensive to train due to expensive gradient evaluations and slower convergence. The experiments conducted in the chapter demonstrate the potential applications of these rational networks for solving PDEs and generative adversarial networks. The practical implementation of rational networks is straightforward in the TensorFlow framework and consists of replacing the activation functions by trainable rational functions. The main benefits of rational NNs are their fast approximation power, the trainability of the activation parameters, and the smoothness of the activation function outside poles.

Our primary objective in \cref{chap_data_green} was to uncover mechanistic understanding from input-output data using a human-understandable representation of an underlying hidden differential operator. This representation took the form of a rational NN for the Green's function. We extensively described all the physical features of the operator that can be extracted and discovered from the learned Green's function and homogeneous solutions, such as linear conservation laws, symmetries, shock front and singularity locations, boundary conditions, and dominant modes. Our deep learning method for learning Green's functions and extracting human-understandable properties of partial differential equations benefits from the adaptivity of rational neural networks and its support for qualitative feature detection and interpretation. We successfully tested our approach with noisy and sparse measurements as training data in one and two dimensions. The design of our network architecture, and the covariance kernel used to generate the system forcings, are guided by rigorous theoretical statements, obtained in \cref{chapt_PDE_learning,chap_svd,chapt_rational}, that offer performance guarantees. This shows that our proposed deep learning method may be used to discover new mechanistic understanding with machine learning.

The deep learning method naturally extends to the case of three spatial dimensions but these systems are more challenging due to the GPU memory demands required to represent the six-dimensional inputs used to train the neural network representing the Green's function. However, alternative optimization algorithms than the one we used, such as mini-batch optimization~\cite{kingma2014adam,li2014efficient}, may be employed to alleviate the computational expense of the training procedure. While our method is demonstrated on linear differential operators, it can be extended to nonlinear, time-dependent problems that can be linearized using an implicit-explicit time-stepping scheme~\cite{ascher1997implicit,pareschi2005implicit} or an iterative method~\cite{kelley1995iterative}. This process should allow us to learn Green's functions of linear time propagators and understand physical behavior in time-dependent problems from input-output data such as the time-dependent Schr\"odinger equation. The numerical experiments conducted in \cref{chap_data_green} highlight that our approach can generalized to discover Green's functions of some linearization of a nonlinear differential operator.

There are many future research directions exploring the potential applications of rational networks beyond PDE learning, in fields such as image classification, time series forecasting, and generative adversarial networks. These applications already employ nonstandard activation functions to overcome various drawbacks of ReLU. Another exciting and promising field is the numerical solution and data-driven discovery of partial differential equations with deep learning. We believe that popular techniques such as physics-informed neural networks~\cite{raissi2019physics} could benefit from rational NNs to improve the robustness and performance of PDE solvers, both from a theoretical and practical viewpoint.

Finally, while the ideas present in \cref{chapt_PDE_learning} have been recently applied to derive a learning rate for Green's functions associated with parabolic PDEs~\cite{boulle2022parabolic}, obtaining theoretical results for more general classes of PDEs, such as hyperbolic, fractional, or stochastic PDEs, remain highly challenging. Such studies are essential to understand which mathematical models can be learned from data, obtain performance guarantees, and, more generally, deepen our knowledge of PDE learning techniques.

\dobib

%now enable appendix numbering format and include any appendices
%\appendix
%\subfile{appendix1}

%next line adds the Bibliography to the contents page
%\subfile{biblio}
\bibliographystyle{siam}  %use the plain bibliography style
\bibliography{references}  %use a bibtex bibliography file refs.bib

\end{document}